\documentclass[12pt]{article}
\usepackage[utf8]{inputenc}

\usepackage{amsmath,amssymb,amsthm,mathtools,tikz-cd, bbm}
\usepackage[margin=1in]{geometry}
\usepackage{hyperref}
\hypersetup{
    colorlinks=true,
    linkcolor=blue,
    filecolor=magenta,      
    urlcolor=cyan,
}
\newtheorem{theorem}{Theorem}[section]
\newtheorem{lemma}[theorem]{Lemma}
\newtheorem{corollary}[theorem]{Corollary}
\newtheorem{sublemma}{Sublemma}[theorem]
\newtheorem{proposition}[theorem]{Proposition}
\newtheorem{claim}{Claim}
\newtheorem{subclaim}{Subclaim}[claim]
\usepackage{etoolbox}
\AtEndEnvironment{theorem}{\setcounter{claim}{0}}
\AtEndEnvironment{lemma}{\setcounter{claim}{0}}
\AtEndEnvironment{proposition}{\setcounter{claim}{0}}

\theoremstyle{remark}
\newtheorem{remark}[theorem]{Remark}

\theoremstyle{definition}
\newtheorem{definition}[theorem]{Definition}
\newtheorem{example}[theorem]{Example}

\renewenvironment{proof}[1][Proof.]{\begin{trivlist}
\item[\hskip \labelsep {\textit{#1}}]}{\end{trivlist}}

\newcommand{\<}{\langle}
\newcommand{\hd}{\text{HOD}}
\newcommand{\conc}{{}^\frown}

\newcommand{\defeq}{=_\text{def}}
\newcommand{\mfc}{\mathfrak{C}}
\newcommand{\adp}{\mathsf{AD^+}}
\newcommand{\ad}{\mathsf{AD}}

\newcommand{\R}{\mathbb{R}}
\newcommand{\adr}{\ad_\R}

\newcommand{\hpc}{\textsf{HPC}}
\newcommand{\C}{\mathbb{C}}
\newcommand{\lh}{\text{lh}}
\newcommand{\ult}{\text{Ult}}
\newcommand{\crit}{\text{crit}}
\newcommand{\dom}{\text{dom}}
\newcommand{\ran}{\text{ran}}
\newcommand{\rs}{\text{rs}}

\newcommand{\pred}{\text{-pred}}
\newcommand{\restrict}{\upharpoonright}
\newcommand{\is}{\trianglelefteq}
\newcommand{\isneq}{\triangleleft}
\newcommand{\tree}[1]{\mathcal{#1}}
\newcommand{\mtree}[1]{{\mathbb{#1}}}

\newcommand{\id}{\text{id}}

\newcommand{\nul}{\lambda}
\newcommand{\nuz}{\hat{\lambda}}

\newcommand{\ihat}{\hat\imath}

\def\ihat{\hat{\imath}}
\def\ohat{\hat{o}}

\def\OR{\textrm{OR}}

\def\hc{\textrm{HC}}

\def\nul{\lambda}
                                               
\def\itS{\mathcal{S}}
\def\itT{\mathcal{T}}
\def\itU{\mathcal{U}}
\def\itV{\mathcal{V}}
\def\itW{\mathcal{W}}

\newcommand{\inup}{\mathbin{\rotatebox[origin=c]{90}{$\in$}}}
\makeatletter
\tikzset{
  edge node/.code={%
      \expandafter\def\expandafter\tikz@tonodes\expandafter{\tikz@tonodes #1}}}
\makeatother
\tikzset{
  is/.style={
    draw=none,
    edge node={node [sloped, allow upside down, auto=false]{$\is$}}},
  Is/.style={
    draw=none,
    every to/.append style={
      edge node={node [sloped, allow upside down, auto=false]{$\is$}}}
  }
}
\tikzset{
  eq/.style={
    draw=none,
    edge node={node [sloped, allow upside down, auto=false]{$=$}}},
  Eq/.style={
    draw=none,
    every to/.append style={
      edge node={node [sloped, allow upside down, auto=false]{$=$}}}
  }
}

\title{Full normalization for mouse pairs}
\author{Benjamin Siskind and John Steel}

\begin{document}

\maketitle
\setcounter{section}{-1}
\section{Introduction}
A \textit{mouse pair} is a premouse $M$ together with an iteration strategy $\Sigma$ 
for $M$ with certain condensation
properties. The notion is isolated in \cite{nitcis}. That book proves 
 a comparison theorem for mouse pairs, and 
shows that many of the basic results of inner model theory can 
be stated in their proper general
form  by considering mouse pairs instead of just mice. For example, we have
 the full Dodd-Jensen property for mouse pairs,
and thus a wellfounded mouse pair order, whereas these both fail if we 
consider iterable premice in isolation. 
These and other results seem to indicate that mouse pairs, and not just their 
mouse components, are the fundamental objects of study in inner model theory.

One important technical device employed in \cite{nitcis} is
embedding normalization. Given a stack $s$ of normal trees on $M$
 with last model $P$,
 there is a natural attempt to build a minimal normal tree $W(s)$ such that
 $P$
 embeds into its last model. $W(s)$ is called the
 {\em embedding normalization} of $s$.  An iteration strategy $\Sigma$ for
 $M$ {\em normalizes well} iff whenever $s$ is by $\Sigma$, then $W(s)$ exists and is
 by $\Sigma$. It is one of the defining properties of mouse pairs that their
strategy components normalize well. \cite{nitcis} shows
that full background extender constructions  done in an appropriate universe yield 
mouse pairs.\footnote{ Embedding normalization
was first studied systematically by Schlutzenberg and Steel, 
independently and then partly jointly. The good behavior of nice 
strategies on infinite stacks is due to Schlutzenberg (see \cite{farmer}).
Some of this work was later re-cast by Jensen and extended to his 
$\Sigma^*$-elementary iterations in \cite{jensen}.}

Our main result here is that, assuming $\adp$, if $(M,\Sigma)$ is a mouse pair,
and $s$ is a stack of normal trees on  $M$ by $\Sigma$ with last model $P$,
then in fact there is a
normal tree $X(s)$ by $\Sigma$ whose last model is equal to $P$.\footnote{So in what sense
then was $W(s)$ minimal? The answer has to do with how the extenders used in $s$
get associated to extenders used in $W(s)$, which is more direct than the way
they are associated to extenders used in $X(s)$. $W(s)$ is minimal, granted that
we demand the more direct connection. See \ref{factor lemma} below.} We call
$X(s)$ the {\em full normalization} of $s$, and say that $\Sigma$
fully normalizes well. Special cases of 
this theorem were proved in \S 6.1 of \cite{nitcis}.\footnote{What we actually prove
here is Theorem \ref{fullnormalizationeheorem}, which does not cover certain anomalous
stacks $s$. The complete proof involves more bookkeeping, but no further ideas.}
In contrast to embedding normalization, the proof goes beyond iteration tree 
combinatorics. The sort of
phalanx comparison typical of condensation proofs for ordinary mice
comes into play, even in the case $s$ consists of two normal trees, each
using only one extender.

Our results on full normalization figure heavily in the construction of
optimal Suslin representations for mouse pairs given in \cite{mouse.suslin}.
As one would expect, such Suslin representations are useful. For example,
\cite{mouse.suslin} uses them to characterize the {\em Solovay sequence}
$\langle \theta_\alpha \mid \alpha < \Omega \rangle$
in terms of the cutpoint Woodin cardinals of $\hd$, 
assuming $\adr$ and a
natural mouse capturing hypothesis.\footnote{The mouse capturing hypothesis is
$\hpc$, or {\em HOD Pair Capturing}. It simply asserts that the iteration
strategies of mouse pairs are Wadge cofinal in the Suslin-co-Suslin sets.
\cite{mouse.suslin} shows that assuming $\adr$ + $\hpc$, 
the following are equivalent: (i) $\delta$ is a cutpoint Woodin cardinal of
$\hd$, and (ii) $\delta = \theta_0$, or $\delta = \theta_{\alpha+1}$ for
some $\alpha$.} The papers \cite{mouse.suslin} and 
\cite{jacksonsargsyansteel} show that under the same hypotheses, the Suslin
cardinals are precisely the cardinalities of cutpoints in $\hd$.\footnote{$\kappa$
is a cutpoint of $\hd$ iff there is no extender $E$ on the $\hd$-sequence
such that $\crit(E) < \kappa \leq \lh(E)$.} For further applications
of our results on full normalization, see \cite{mouse.suslin} and 
\cite{jacksonsargsyansteel}.

  The idea of our proof 
  that there is a
normal tree $X(s)$ by $\Sigma$ whose last model is equal to $P$ is roughly as follows.
 It is not too hard to
define $X(s)$, by extending the definition of $X(s)$ in the special cases 
covered by \cite{nitcis}
to arbitrary stacks $s$. 
 The main problem is to
show that $X(s)$ is by $\Sigma$. Suppose then that
$\tree{S} = X(s)\restrict \lambda$ is by
$\Sigma$, and $X(s)$ picks $b = [0,\lambda)_X$, while $\Sigma(\tree{S}) = c$.
We must show that $b=c$. For this, we compare the phalanxes of the trees 
$\tree{S}\conc b$ and 
$\tree{S}\conc c$, using $\Sigma$ to iterate the  phalanx
$\Phi(\tree{S}\conc c)$. The strategy for iterating the phalanx
$\Phi(\tree{S} \conc b)$ comes from pulling back the strategy
for $\Phi(W(s))$ induced by $\Sigma$,
under
a natural embedding $\Psi \colon X(s) \to W(s)$ that comes out of the definition
of $X(s)$. Here we face one of our main new problems: unless
$X(s) = W(s)$, $\Psi$ is not actually a tree embedding, but something weaker.
A fair amount of our work is devoted to isolating the properties of
$\Psi$ in the notion of a {\em weak tree embedding}, and showing
that if $\Gamma \colon \tree{T} \to \tree{U}$ is a weak tree embedding,
then we can use $\Gamma$ to pull back strategies for $\Phi(\tree{U})$ to
strategies for $\Phi(\tree{T})$.

There is a second issue, one that also comes up in the
phalanx comparisons of \cite{nitcis}.
In order to show that $b=c$,
we need to use the full Dodd-Jensen
property, and so we must compare the last models of our phalanxes as mouse pairs.
This means we must use something like of the mouse pair comparison process
developed in \cite{nitcis}, comparing both phalanxes with the levels of a common 
background construction. To ensure that our comparison process doesn't terminate 
in a trivial way (by applying the same extender to a model common to both trees), 
at certain stages we lift a phalanx.

This much follows the comparison processes in
the proofs of solidity, universality, condensation,
in \cite{nitcis}, where the resulting systems are called 
{\em pseudo-iteration trees}.\footnote{ The argument closest to the one we
give here is the proof that {\sf UBH} holds in lbr hod mice, Theorem 7.3.2 of
\cite{nitcis}. See also \cite{trang}.} However, here 
our phalanxes are all of the
form $\Phi(\tree{T})$, for some iteration tree $\tree{T}$,  and this enables us to
lift them in a different way. Namely,
we can use one step of the embedding normalization process, lifting $\tree{T}$
to $W(\tree{T},F)$.
The resulting system is best viewed as a tree of normal iteration trees,
 something we shall
call a \textit{meta-iteration tree}, or \textit{meta-tree}. 

     The meta-tree notion evolved from the work of Jensen, Schlutzenberg, and Steel on
embedding normalization.
Its full, general form is due to Schlutzenberg. (See \cite{jensen} and \cite{farmer}.
 Those papers
use somewhat different terminology for meta-trees and their associated apparatus.)
Meta-trees provide a very convenient framework for thinking about certain aspects
of iteration tree combinatorics, and in particular, they help a lot here.
The general results about meta-trees that we shall need are
due to Schlutzenberg and Siskind. We shall state those results and
outline their proofs here, but the reader should see \cite{farmer} and
\cite{associativity} for an in-depth treatment. The general theory of meta-iteration trees
has other applications that are described in those papers.

In section 1 of the paper we review embedding normalization. The one new
result here is a factoring lemma for tree embeddings.
Section  2 is devoted to the general theory of meta-trees.

Section 3 proves a general comparison theorem for the ``tail normal
components" of nice meta-strategies, Theorem \ref{main comparison theorem}. 
This is where we show that the appropriate
phalanx comparisons terminate.
One can think of this as a strategy comparison
theorem for phalanxes of the form $\Phi(\tree{S})$.\footnote{Since not all phalanxes
 are of this
form, it is not a truly general strategy comparison theorem for phalanxes.}
We then use our tree-phalanx comparison theorem to characterize those meta-strategies
that are induced by an ordinary iteration strategy.\footnote{An ordinary
strategy $\Sigma$ induces a meta-strategy $\Sigma^*$ via:
$\langle \mtree{S}_\xi \mid \xi < \lambda \rangle
  \text{ is by }\Sigma^* \Leftrightarrow\text{$\forall \xi < \lambda$ (every tree occurring in $\mtree{S}_\xi$ is by }\Sigma)$. See \ref{meta-iterability}.} It turns
out every sufficiently nice meta-strategy is induced in this way. (That is Theorem
\ref{induced strategy theorem}.)  The moral one might draw is that meta-strategies
are not something fundamentally new, but rather a useful way of organizing
constructions and proofs to do with ordinary strategies.

The main step toward Theorem \ref{induced strategy theorem} is Lemma
\ref{induced strategy lemma}, which is a kind of uniqueness theorem
for ordinary iteration strategies $\Sigma$ whose induced meta-strategies
$\Sigma^*$ behave well.

Section 4 contains the definitions of $X(s)$, weak tree embeddings,
and the weak tree embedding from $X(s)$ to $W(s)$.
We then use the meta-strategy uniqueness results of section 3 to show
that if $(M,\Sigma)$ is a mouse pair,
then  $\Sigma$ condenses to itself
under weak tree embeddings. This is Theorem \ref{vshctheorem}, one of
the central results of the paper. It follows easily from this
{\em very strong hull condensation} property that $\Sigma$
fully normalizes well.\footnote{ The results in section 4 are essentially
due to the second author. 
He outlined proofs of them that made use of pseudo-iteration trees 
in \cite{localhod}. The conversion of those outlines to full proofs that we
present here, making use of meta-iteration trees, is due to the first author.
We shall attribute the results of sections 2 and 3 when we get to them.}

We assume that the reader is familiar with \cite{nitcis}, especially
Chapter 6 (on embedding normalization) and Chapter 9 (on phalanx comparison).
We shall review and re-do much of this work, however, so this prerequisite
is perhaps less burdensome than it may seem. Our unexplained general
inner-model-theoretic notation is laid out in Chapters 2-4 of \cite{nitcis}.
In the rest of the paper,
   by \lq\lq premouse" we mean Jensen-indexed pure extender or least branch pfs
premouse. ``pfs" stands for ``projectum-free spaces", a variant on the standard
fine structure that is described in Chapter 4 of \cite{nitcis}.

\section{Tree embeddings and embedding normalization}

      We review some material from \cite{nitcis}, while adding a 
few things to it. 

     We shall be using the projectum-free spaces fine structure of
\cite{nitcis}, for reasons explained there.\footnote{See the beginning of
Chapter 4 of \cite{nitcis}.} Our iteration trees on such
premice will be linear stacks of {\em plus trees}. Briefly:
for $F$ an extender of the sort that might appear on the sequence
of a pfs premouse\footnote{The extenders are Jensen-indexed, so that
$E$ is indexed in $M$ at $\lh(E) = \lambda(E)^{+,N}$, where $N=
\ult(M||\lh(E),E)$ and $\lambda(E) = i_E^{M||\lh(E)}(\crit(E))$.},
$E^+$ is the extender of $E$-then-$D$, where $D$ is the order zero
measure of $\ult(M||\lh(E),E)$ on $\lambda(E)$. $F$ has {\em plus type}
iff $F=E^+$ for some such $E$. If $F=E^+$, then $F^- = E$. If $F$ is not
of plus type, $F^- = F$. In both cases we let $\lh(F) = \lh(F^-)$,
and set $\nuz(F) = \lambda(F^-)$.   The {\em extended $M$-sequence}
consists of all $E$ and $E^+$ on the $M$-sequence. 

A {\em plus tree}
on $M$ is a quasi-normal iteration tree $\tree{T}$ on $M$ that is allowed
to use an extender from the extended sequence of $M_\alpha^{\tree{T}}$
to form $M_{\alpha+1}^{\tree{T}}$. $\tree{T}$ is {\em $\lambda$-separated}
if it always uses extenders of plus type, and {\em $\lambda$-tight} if
it never uses such extenders. Roughly, $\tree{T}$ is {\em quasi-normal}
iff the extenders it uses are $\nuz$-nondecreasing, and it always applies
its extenders to the longest possible initial segment of the earliest
possible model. $\tree{T}$ is {\em normal} iff it is quasi-normal and
the extenders it uses are (strictly) $\lh$-increasing.

\subsection{Tree embeddings} Tree embeddings in general were isolated
 by the second author and Schlutzenberg.  They arise naturally in the context
  of embedding normalization and will feature prominently in the rest
 of the paper. The precise definition must be tailored to the context of 
interest: we will have to look at tree embeddings between
 arbitrary \textit{plus trees}. This variation was worked out by the second author in \cite{nitcis}. 
  
We first establish some notation for applications of the Shift 
Lemma (that is, copy maps). As usual, for $E$ an extender on the 
extended $M$-sequence, the \textit{domain of $E$} is
  $\dom(E) = M||(\crit (E)^+)^{M|\lh(E)}$.
  
  \begin{definition} Let $\varphi \colon M \to N$ and $\pi \colon
   P \to Q$ be nearly elementary, $E$ an extender
on the extended $M$-sequence, and $F$ an extender on the extended $N$-sequence. 
We say \textit{the Shift Lemma applies to $(\varphi, \pi, E, F)$}
  iff \begin{enumerate}
      \item[i.] $\dom(E) = P||({\crit (E)^+})^{P}$,\footnote{We allow $P=\dom(E)$ here.}
      \item[ii.] $\varphi \restrict (\dom(E) \cup
 \{\dom(E) \}) = \pi \restrict (\dom(E) \cup \{\dom(E) \})$, and
      \item[iii.] either \begin{enumerate}
          \item[(a)] $F=\varphi(E)$ or
          \item[(b)] $E$ is not of plus type and $F=\varphi(E)^+$.
      \end{enumerate}
  \end{enumerate}
In this situation, the Shift Lemma (cf. \cite{nitcis}, Corollary 2.5.20)
 gives that there is a unique nearly elementary map $\sigma: Ult(P, E)\to Ult(Q,F)$ such that
\begin{enumerate}
\item $\sigma \circ  i^{P}_E= i^{Q}_{F} \circ\pi$ and
\item $\sigma \restrict \varepsilon(E) = \varphi \restrict \varepsilon(E)$.
\end{enumerate}
We call this $\sigma$ the \textit{copy map associated to $(\varphi, \pi, E, F)$}.\footnote{If hypothesis (iii)(a) obtains, then the existence of $\sigma$ follows literally by the Shift Lemma. If hypothesis (iii)(b) obtains, then $\sigma=i^{\ult(M, F^-)}_\mu\circ \sigma^-$ where $\sigma^-$ is the ordinary copy map associated to $(\varphi, \pi, E, F^-)$ and $\mu$ is the Mitchell-order zero normal measure on $\lambda(F^-)$. Since $\crit(i^{\ult(M, F^-)}_\mu)=\lambda(F^-)$ and $\sigma^-$ satisfies (2) and the appropriate version of (1), this $\sigma$ does actually satisfy (1) and (2). Uniqueness of this map is guaranteed as usual, since $\ult(P, E)$ is the hull of points in the range of $i_E^P$ together with $\varepsilon(E)=\lambda(E)$. Also note that in this case, if $\sigma^-$ is elementary (for example when $\langle \varphi, \pi\rangle: (P, E)\to^* (Q, F^-)$), then so is $\sigma$.} 
\end{definition}

It will be convenient to use all this terminology even when $\dom(E)\is P$ but 
clause (i), above, fails. In this case, letting $\bar{P}\is P$ the least
 initial segment $R$ of $P$ such that $\dom(E)\is R$ and $\rho(R)\leq \crit(E)$, 
we have that $\pi(\bar{P})$ is the least initial segment $S$ of $Q$ such 
that $\dom(\varphi(E))\is S$ and $\rho(S)\leq \crit(\varphi(E))$ by near 
elementarity of $\pi$. We'll say that the \textit{Shift Lemma applies 
to }$(\varphi, \pi, E, F)$ when it applies to $(\varphi, \pi\restrict\bar{P}, E, F)$ 
and let the copy map associated to $(\varphi, \pi, E, F)$ be that 
associated to $(\varphi, \pi\restrict\bar{P}, E, F)$.

Our definition of a tree embedding between two plus trees will 
actually be an ostensible weakening of that in \cite{nitcis}. 
The advantage for this definition is its relative ease of verification.
 We shall verify it is equivalent to the definition from \cite{nitcis} shortly.
 \begin{definition}\label{treembdef} Let $\itS$ and $\itT$
 be plus trees
 on a premouse $M$.
 A \textit{tree embedding}
  $\Phi:\tree{S}\to \tree{T}$ is a system
 $\langle v,u, \{s_\xi\}_{\xi<\lh\tree{S}}, \{t_\zeta\}_{\zeta+1<\lh(\tree{S})}\rangle$ such that
\begin{enumerate}
    \item  $v:\lh(\tree{S})\to \lh(\tree{T})$ is tree-order preserving,
     $u:\{\eta\,|\,\eta+1<\lh (\tree{S})\}\to \lh (\tree{T})$, 
     $v(\xi)=\sup\{u(\eta)+1\,|\, \eta<\xi\}$, and for all $\xi+1<\lh(\tree{S})$, $v(\xi)\leq_\tree{T}u(\xi)$;
    \item For all $\xi$ and $\eta\leq_\tree{S}\xi$,
    \begin{enumerate}
        \item $s_\xi: M^\tree{S}_\xi\to M^\tree{T}_{v(\xi)}$ is nearly elementary and 
        $s_0 = id_{M^\tree{S}_0}$,
        \item $\hat\imath^\tree{T}_{v(\eta),v(\xi)}\circ s_\eta=
        s_\xi\circ \hat\imath^\tree{S}_{\eta,\xi}$,\footnote{We are
         using $\hat i^\tree{T}$ to denote the possibly partial branch embeddings 
         of $\tree{T}$, as in \cite{nitcis}.} and
         \item if $\xi+1<\lh(\tree{S})$, then $t_\xi= \hat\imath^\tree{T}_{v(\xi),u(\xi)}\circ s_\xi$ with $M_\xi^\tree{S}|\lh(E_\xi^\tree{S})\is \dom(t_\xi)$;
         \end{enumerate}
        \item for all $\xi+1 < \lh(\itS)$, letting $\eta=\tree{S}\pred(\xi+1)$ and
     $\eta^*=\tree{T}\pred(u(\xi)+1)$, 
    \begin{enumerate}
    \item either $E^\tree{
    T}_{u(\xi)}=t_\xi^\Phi(E_\xi^\tree{S})$ or else $E^\tree{S}_\xi$ is not of plus type and $E^\tree{
    T}_{u(\xi)}=t_\xi^\Phi(E_\xi^\tree{S})^+$,
    \item $\eta^*\in[v(\eta),u(\eta)]_\tree{T}$,
        \item $s_{\xi+1}^\Phi\restrict \varepsilon(E_\xi^\tree{S})=t^\Phi_\xi\restrict\varepsilon(E_\xi^\tree{S})$.
        \end{enumerate}
    \end{enumerate}
    \end{definition}

Clause (1) implies that $v(0)=0$ but also that $v$ is continuous at limit
 ordinals $\lambda$ so that clause (2)(b) then gives that $s_\lambda : M_\lambda^{\tree{T}} \to
M_{v(\lambda)}^{\tree{U}}$ is the direct limit of the $s_\eta$ for
$\eta <_T \lambda$ sufficiently large. 

Clause (3) together with the commutativity clause (2)(b) guarantee that the 
all of the $s$-maps are actually copy maps and, in fact, elementary,
 not just nearly elementary. To see this, we just need to see 
that we've maintained enough agreement so that the Shift Lemma
 applies at successors. Then we can repeat the proof of the Copy
 Lemma of \cite{nitcis} (Lemma 4.5.17) to get that all the
 $s$-maps are elementary, by induction.

\begin{proposition}
Let $\Phi:\tree{S}\to \tree{T}$ be a tree embedding. Then 
for all $\xi<\lh(\tree{S})$, \begin{enumerate}
    \item for all $\eta<\xi$, $s_\xi\restrict\varepsilon(E_\eta^\tree{S})=s_\eta\restrict\varepsilon(E_\eta^\tree{S})$;
\item if $\xi+1 < \lh(\itS)$, letting $\eta=\tree{S}\pred(\xi+1)$ and
     $\eta^*=\tree{T}\pred(u(\xi)+1)$,
\begin{enumerate}
   \item $\eta^*$ is the least $\zeta\in [v(\eta), u(\eta)]_\tree{T}$ 
such that $\zeta=u(\eta)$ or else $\crit(\hat\imath^\tree{T}_{\zeta, u(\eta)})>\crit(E_{u(\xi)}^\tree{T})$,
    \item the Shift Lemma applies to $(t_\xi, \hat\imath^\tree{T}_{v(\eta),
 \eta^*}\circ s_\eta, E_\xi^\tree{S}, E_{u(\xi)}^\tree{T})$ and $s_{\xi+1}$ is the associated copy map; and
\end{enumerate}
\item $s_\xi$ is elementary.
\end{enumerate}
\end{proposition}

\begin{proof}[Proof sketch.]
(1) and (2) can be verified by induction, simultaneously. (3) can 
be proved by a separate induction using that (1) and (2) hold 
for all $\xi<\lh(\tree{S})$, using the analysis from \cite{nitcis} 
about when extenders of an iteration tree very close to the models 
to which they are applied. As mentioned above, this is essentially 
the proof of the Copy Lemma of \cite{nitcis}. We omit further detail. 
\end{proof}

This proposition implies that our current definition is equivalent to that of \cite{nitcis}.
 
  \begin{definition}
  If $\tree{T}$ is a plus tree, and $\alpha < \lh(\tree{T})$, then
\begin{align*}
\nul_\alpha^\tree{T} & = \sup \{ \lambda(E_\beta^\tree{T})\,|\,\beta < \alpha \}\\
      &  = \sup \{ \lambda(E_\beta^\tree{T}) \,|\, \beta <_T \alpha\} \end{align*}
is the sup of the Jensen generators of $M_\alpha^\tree{T}$.
\end{definition}

The agreement of the $s$ and $t$ maps in a tree embedding is given by

\begin{lemma}\label{pshullagree} Let $\langle v, u, \langle s_\beta \mid \beta < \lh \itT\rangle,
 \langle t_\beta \mid \beta+1 < \lh \itT\rangle\rangle$ be a tree embedding of $\itT$ into $\itU$; then
\begin{itemize}
\item[(a)] if $\alpha +1 < \lh(\itT)$, then $t_\alpha$ agrees with $s_\alpha$ on
$\nul_\alpha^{\itT}$,
\item[(b)] if $\beta < \alpha < \lh(\itT)$, then $s_\alpha$ agrees with $t_\beta$ on
        $\lh(E_\beta^{\itT}) +1$, and
\item[(b)] if $\beta < \alpha < \lh(\itT)$, then $s_\alpha$ agrees with $s_\beta$
on $\lambda_\beta^{\itT}$
\end{itemize}
\end{lemma}

We also have a formula for the point of application
 $\eta^* = T\pred(u(\xi)+1)$ 
occurring in clause (4) of 
Definition \ref{treembdef}, namely
\[
T\pred(u(\xi)+1) = \text{least $\gamma \in [v(\eta), u(\eta)]_T$ 
such that $\crit \ihat^\tree{T}_{\gamma,u(\eta)} > 
\ihat^\itT_{v(\eta),\gamma} \circ s_\eta(\mu)$,}
\]
where
\[
\eta = \tree{S}\pred(\xi +1) \mbox{  and  } \mu = \crit(E_\xi^{\itS}).
\]
These facts are proved in \cite{nitcis}. The proof uses the following
elementary fact about iteration trees.\footnote{See Proposition 8.2.1 of \cite{nitcis}.
In reading it, recall that the convention of \cite{nitcis} is that when $P$ is
    active, then $P||o(P) \isneq P$. }

\begin{proposition}
\label{unravelonbranch}
Let $\itS$ be a normal tree, let $\delta \le_S \eta$, and
suppose that $P \unlhd M_\eta^{\itS}$, but $P \not\is M_\sigma^{\itS}$ whenever
$\sigma <_S \delta$. Suppose also that $P \in \ran(\ihat^{\itS}_{\delta,\eta})$. Let
\begin{align*}
\alpha & = \text{ least $\gamma$ such that $P \unlhd M_\gamma^{\itS}$}\\
   &= \text{ least $\gamma$ such that $o(P) < \lh(E_\gamma^{\itS})$ or
    $\gamma = \eta$,}
\end{align*}
and
\[
\beta = \text{least $\gamma \in [0,\eta]_S$ such that }
                        o(P) < \crit(\ihat^\itS_{\gamma,\eta}) \mbox{ or } \gamma = \eta.
\] Then $\beta \in [\delta,\eta]_S$, and
\begin{itemize}
\item[(a)] either  $\beta = \alpha$, or $\beta = \alpha +1$, and
$\lambda(E_\alpha^{\itS}) \le o(P) < \lh(E_\alpha^{\itS})$;
\item[(b)] if $P = \dom(E_\xi^{\itS})$, then $S\pred(\xi+1) = \alpha = \beta$.
\end{itemize}
 (We allow $\delta = \eta$, with the understanding $\ihat_{\delta,\delta}$ is the identity.)
\end{proposition}

This proposition or its (very short) proof show up in many arguments
to do with extending tree embeddings. Commonly, one has
$\eta = u(\xi)$, $\delta = v(\xi)$, and $P = t_\xi(\bar{P})$.

\begin{definition}\label{extendedtreeemb}
Let $\tree{S}$ be of successor length $\gamma+1$ and $\Phi:\tree{S}\to \tree{T}$ a 
tree embedding where $\tree{T}$ has successor length $\delta+1$. 
If $v(\gamma)\leq_\tree{T} \delta$, we let $u(\gamma)=\delta$ and 
$t_\gamma=\hat\imath^\tree{T}_{v(\gamma),\delta}\circ s_\gamma$ and call 
the resulting system an \textit{extended tree embedding}. An extended tree embedding 
$\Phi$ is \textit{non-dropping} if $(v(\gamma),\delta]_\tree{T}$ doesn't drop.
\end{definition}

Note that if $\Phi:\tree{S}\to \tree{T}$ is a tree embedding and $\tree{S}$ has
 successor length $\delta+1$, then we can always view $\Phi$ as a non-dropping 
 extended tree embedding $\Phi:\tree{S}\to \tree{T}\restrict v(\delta)+1$. 
 
 \begin{remark} On the
  other hand, it is not always possible to extend tree embeddings defined on trees
  of limit length, even if our trees $\itS$ and $\itT$ are by the same nice
  iteration strategy. 
  
   For example, assume $M_1^\#$ exists and let $\Lambda$ be the iteration strategy
    for $M_1$. Toward a contradiction, suppose for every $\tree{S},\tree{T}$ of limit
     lengths by $\Lambda$ such that there is $\Phi: \tree{S}\to \tree{T}$ with 
     $\ran(v^\Phi)$ cofinal in $\lh \tree{T}$, $v^\Phi ``\Lambda(\tree{S})\subseteq 
     \Lambda(\tree{T})$. 
    Now let $\tree{T}\in M_1$ be a tree by $\Lambda$ of height ${\delta^+}^{M_1}$ 
which has no branch in $M_1$, where $\delta$ is the Woodin cardinal of $M_1$ (that
 such a tree exists is due to Woodin, see Lemma 1.1 of \cite{schindler steel}).
  Let $g$ be $Col(\omega, \delta)$ generic over $M_1$ and $h$ generic for the Namba 
  forcing over $M_1[g]$. The restriction of $\Lambda$ to countable trees which are
   in $M_1[g][h]$ is in $M_1[g][h]$ since Namba forcing adds no reals and $M_1$ contains
    the restriction of $\Lambda$ to trees of length $\delta$ which are in $M_1$. Now, 
    in $M_1[g][h]$, $\lh \tree{T}$ has countable cofinality. We can take a Skolem hull 
    to get $\tree{S}$ countable in $M_1[g][h]$ and $\Phi: \tree{S}\to \tree{T}$ with
     $\ran(v^\Phi)$ cofinal in $\lh \tree{T}$. Since $\tree{S}$ is countable and by
      $\Lambda$, $\Lambda(\tree{S})\in M_1[g][h]$. So, by assumption, 
      $v^\Phi`` \Lambda(\tree{S})\subseteq \Lambda(\tree{T})$; so $\Lambda(\tree{T})$ is 
      just the downwards closure of $v^\Phi``\Lambda(\tree{S})$ (in $M_1[g][h]$). This
       identification of $\Lambda(\tree{T})$ was independent of our choice 
       of $g,h, \tree{S}$, so we get $\Lambda(\tree{T})\in M_1$, a contradiction.
       \end{remark}

The tree embeddings that appear naturally 
in embedding normalization can be viewed as extended tree embeddings.
We deal almost exclusively
 with extended tree embeddings in the rest of the paper.

\begin{definition}\label{treeembagreedef}
For tree embeddings or extended tree embeddings 
$\Phi:\tree{S}\to \tree{T},\Gamma:\tree{U}\to \itV$, we put
 \[
 \Phi\restrict\xi+1\approx \Gamma \restrict\xi+1
 \] iff
  $\tree{S}\restrict\xi+1=\tree{U}\restrict\xi+1$,
   $v^\Phi\restrict\xi+1 = v^{\Gamma}\restrict\xi+1$, and
    $\tree{T}\restrict v^\Phi(\xi)+1=\tree{V}\restrict v^{\Gamma}(\xi)+1$.
    \end{definition}
If  $\Phi\restrict\xi+1\approx \Gamma \restrict\xi+1$, then
$s_\eta^\Phi= s_\eta^{\Gamma}$ for
      $\eta\leq \xi$, $u^\Phi\restrict\xi=u^{\Gamma}\restrict\xi$, and $t_\eta^\Phi=t_\eta^\Gamma$ for $\eta<\xi$.
        It does \textit{not} imply that $u^\Phi(\xi)= u^{\Gamma}(\xi)$ 
        (even when $\Phi, \Phi^*$ are extended tree embeddings). Intuitively, $u^\Phi(\xi)$
    and $t^\Phi_\xi$ are telling us how to inflate $E_\xi^\tree{S}$, while
        $\Phi \restrict
        \xi +1$ is the part of $\Phi$ that acts on $\tree{S} \restrict \xi+1$, which
does not include $E_\xi^\tree{S}$.
        
We will sometimes write $u_\alpha$ for $u(\alpha)$ 
and $v_\alpha$ for $v(\alpha)$.\footnote{This is a notational convention introduced by Schlutzenberg which cleans up some type-setting issues, as the reader may notice.}


\subsection{Direct limits under tree embeddings}

If $\Phi \colon \itS \to \itT$ and $\Psi \colon \itT \to \itU$ are
(extended) tree embeddings, then $\Psi \circ \Phi \colon \itS \to \itU$
is the tree embedding obtained by component-wise composition, in the obvious way.

\begin{definition}\label{directedsystemoftrees} Let $M$ be a premouse.
A \textit{directed system of plus trees} on $M$ is a system 
$\mathcal{D}=\langle\{\tree{T}_a\}_{a\in A},\{\Psi_{a,b}\}_{a\preceq b}\rangle$, 
where $\preceq$ is a directed partial order on some set $A$ and 
\begin{enumerate}
    \item[(a)] for any $a\in A$, $\tree{T}_a$ is a plus tree on $M$ of successor length,
    \item[(b)] for any $a,b\in A$ with 
    $a\prec b$, $\Psi_{a,b}: \tree{T}_a\to \tree{T}_b$ is an extended tree embedding,
    \item[(c)] for any $a,b,c\in A$ such that $a\preceq b\preceq c$, 
    $\Psi_{a,c}= \Psi_{b,c}\circ \Psi_{a,b}$.
\end{enumerate} 
\end{definition}

It follows from (c) that 
$\Psi_{a,a}$ is the identity extended tree embedding on $\tree{T}_a$.

Let $\mathcal{D}=\langle\{\tree{T}_a\}_{a\in A},\{\Psi_{a,b}\}_{a\preceq b}\rangle$
 be a directed system of plus trees on $M$, where $\preceq$ is a directed 
 partial order on $A$. Assuming the $t$-maps of our tree embeddings behave
 properly,
we shall define the direct limit of 
  $\mathcal{D}$, which we denote  $\lim\mathcal{D}$, as an algebraic structure. We then show
  that if all models in $\lim\mathcal{D}$ are wellfounded,
  then it is (isomorphic to) a plus tree on $M$.
In components, we shall have
\[
 \lim\mathcal{D} = \langle D, \leq, \leq^*, \{M_x\}_{x\in D}, \{E_x\}_{x\in D}, 
\{\Gamma_a\}_{a\in A}\rangle.
\]
Let us define these components.

Let
\[\Psi_{a, b}= \langle v_{a,b},u_{a,b},
 \{s^{a,b}_\gamma\}_{\gamma<\lh(\tree{T}_a)},
  \{t^{a,b}_\gamma\}_{\gamma<\lh( \tree{T}_a)}\rangle.
\]
A \textit{$u$-thread} is a partial function $x:A\rightharpoonup Ord$ such that
whenever $a \in \dom(x)$ and $a \prec b$, then $b \in \dom(x)$ and
$u_{a,b}(x(a)) = x(b)$. Since $\dom(x)$ is upward closed,
it must be $\prec$-cofinal. If $x$ and $y$ are $u$-threads, then
$x \sim y$ iff $\exists a (x(a) = y(a))$ iff $\forall a \in \dom(x) \cap \dom(y)
(x(a) = y(a))$. Every $u$-thread is equivalent to a unique {\em maximal} $u$-thread.
For any 
  $a\in A$ and $\gamma<\lh(\tree{T}_a)$, there is exactly one maximal $u$-thread $x$ such 
  that $x(a)=\gamma$, and 
  we write $x=[a,\gamma]_\mathcal{D}$ for it. 

$D$ is the set of all maximal $u$-threads. $\leq$ and $\leq^*$ will be certain partial 
orders with field $D$. Going forward, we write $\forall^*a \varphi(a)$ to abbreviate
 $\exists b \forall a\succeq b \,\varphi(a)$.

For $u$-threads $x,y$ we put 
\begin{align*}x\leq y &\Leftrightarrow 
\forall^*a \, x(a)\leq y(a)\\
x\leq^* y &\Leftrightarrow \forall^*a\, x(a)\leq_{\tree{T}_a} y(a).
\end{align*}
Since $\leq_{\tree{T}_a}$ is a refinement of the 
 order $\leq$ on ordinals, we get that $\leq^*$ is a refinement of $\leq$.
 $u$-maps preserve $\leq$ on ordinals, so $x \leq y$ iff for some
 (all) $a \in \dom(x) \cap \dom(y)$, $x(a) \leq y(a)$. But $u$-maps don't
 preserve tree-order everywhere, so the $\forall^*a$ quantifier is needed
 in the definition of $x \leq^* y$.

It's easy to see that $\leq$ is a linear order on $D$, but it could fail to be
 a wellorder. If it is a wellorder, we identify it with its order-type $\delta$. 
 In any case, we will think of $\langle D,\leq \rangle$ as the length of
  the direct limit. In the case that the direct limit produces an iteration tree,
   $\delta$ will really be its length.

We now define $\Gamma_a=\langle u_a, v_a, \{s^a_\gamma\}_{\gamma<\lh(\tree{T}_a)}, \{t^a_\gamma\}_{\gamma<\lh(\tree{T}_a)}\rangle$ along with $M_x$ and $E_x$ for $x$ such that $a\in \dom(x)$. We will actually only define the $u$-map and $t$-maps of the $\Gamma_a$; these determined the whole tree embedding in the case that the direct limit produces an iteration tree. So fix $a$ and $\gamma<\lh(\tree{T}_a)$. Let $x=[a,\gamma]_\mathcal{D}$. We set $u_a(\gamma)=x$. We'll actually leave $M_x$, $E_x$ and $t^a_\gamma$ undefined unless
the $t$-maps along $x$ are eventually total. So suppose we're in this case, i.e.
$\forall^*b \forall c\geq b\, (t^{b,c}_{x(b)}$ is total).
We define
\[
M_x = \lim \langle M^{\tree{T}_b}_{x(b)}, t^{b,c}_{x(b)}\,|\,
 b\text{ such that for all }c\succeq b, \,t^{b,c}_{x(b)}\text{ is total}\rangle.
 \]
  For any $b$ such that for all $c\succeq b$, $t^{b,c}_{x(b)}$ is total,
   we let $t^b_{x(b)}$ be the direct limit map and we put 
   $t^a_\gamma = t^b_{x(b)}\circ t^{a,b}_\gamma$ for any such $b$
    (this is independent of the choice of $b$). 
    
    Our assignment of the $E_x$ requires slightly more care. 
 Notice that 
    for any $u$-thread $x$,
    \[\forall^*b\forall c\succeq b E_{x(c)}^{\tree{T}_c}=t^{b,c}_{x(b)}(E_{x(b)}^{\tree{T}_b}).\]
    This is because the only way to have $E_{x(c)}^{\tree{T}_c}\neq t^{b,c}_{x(b)}(E_{x(b)}^{\tree{T}_b})$ is for $E_{x(b)}^{\tree{T}_b}$ to be \textit{not} of plus type and $E_{x(c)}^{\tree{T}_c}= t^{b,c}_{x(b)}(E_{x(b)}^{\tree{T}_b})^+$. But then since $E_{x(c)}^{\tree{T}_c}$ is now of plus type, for any $d\succeq c$, $E_{x(d)}^{\tree{T}_d}= t^{c,d}_{x(c)}(E_{x(c)}^{\tree{T}_c})$. So we may let let 
    \[
    E_x = t^b_{x(b)}(E^{\tree{T}_b}_{x(b)}),
    \] 
    for any $b$ such that $t^{b,c}_{x(b)}$ is total and $E_{x(c)}^{\tree{T}_c}=t^{b,c}_{x(b)}(E_{x(b)}^{\tree{T}_b})$ for all  $c\succeq b$. Again, this is independent of the choice of $b$.
This finishes the definition of $\lim \mathcal{D}$. 

We say that $\lim \mathcal{D}$ is \textit{wellfounded} iff     \begin{enumerate}
    \item for all $x\in D$, the model $M_x$ is defined and wellfounded,
\item $\leq$ is wellfounded,
  \item $\tree{U}=\langle M_x,E_x,\leq^*\rangle$ is a plus tree 
    (i.e. with models $M_x$, exit extenders $E_x$, and tree-order $\leq^*$).
\end{enumerate}
    
If $\lim \mathcal{D}$ is wellfounded, one can show that letting $v_a(x)=\sup\{u_a(y)+1\mid y<x\}$, we can define $s^a_\gamma$ to be the required copy maps so that   $\Gamma_a=\langle u_a, v_a, \{s^a_\gamma\}_{\gamma<\lh(\tree{T}_a)}, \{t^a_\gamma\}_{\gamma<\lh(\tree{T}_a)}\rangle$ is an extended tree embedding from $\tree{T}_a$ into $\tree{U}$ and
 $\Gamma_b\circ \Psi_{a,b}= \Gamma_a$ for every $a\preceq b$. Part of this is the analysis of successors in the $\leq^*$-order, below.

Perhaps surprisingly, we can drop conditions (2) and (3)
 in the definition of the wellfoundedness of the direct limit.

\begin{proposition}\label{direct limit characterization}
Let $\mathcal{D}$ be a directed system of plus trees. Then $\lim \mathcal{D}$ 
is wellfounded iff for every $u$-thread $x$, the models $M_x$ are defined
 and wellfounded.
\end{proposition}

Before we give a proof, we need the following observations about
 iterated applications of the Shift Lemma.

 \begin{lemma}\label{shift composition}
 Let $\pi_0: M_0\to M_1,$ $\pi_1: M_1\to M_2$ and
  $\sigma_0: N_0\to N_1$, $\sigma_1: N_1\to N_2$ be nearly elementary
   and let $E$ be on the extended $M_0$-sequence, $F$ on the extended $M_1$-sequence, and $G$ on the extended $M_2$-sequence.
 
 Suppose that the Shift Lemma applies to $(\pi_0, \sigma_0, E, F)$ and to $(\pi_1, \sigma_1, F, G)$.
   Let $\tau_0$ be the copy map associated
    to $(\pi_0, \sigma_0, E, F)$ and $\tau_1$ the copy map
associated to $(\pi_1, \sigma_1, F, G)$. 
 
 Then the Shift Lemma applies to $(\pi_1\circ \pi_0, \sigma_1\circ \sigma_0, E, G)$ and $\tau_1\circ\tau_0$ is the associated copy map.
\end{lemma}

Next we record how the copying interacts with direct limits.
 This is implicit in \cite{nitcis}.

\begin{lemma}\label{shift direct limits}
Let $\preceq$ be a directed partial order on a set $A$. Suppose we have 
directed systems of premice $\mathcal{M}=\langle \{M_a\}_{a\in A}, \{\pi_{a,b}\}_{a\preceq b}\rangle$ 
and  $\mathcal{N}=\langle \{N_a\}_{a\in A}, \{\sigma_{a,b}\}_{a\preceq b}\rangle$ and extenders
 $\{E_a\}_{a\in A}$ such that
  \begin{enumerate}
    \item[(a)]  for all $a\in A$, $E_a$ is on the extended $M_a$-sequence,
    \item[(b)] for all $a,b\in A$ such that $a\preceq b$, $\pi_{a,b}$ and $\sigma_{a,b}$ are nearly elementary, and
    \item[(c)] for all $a,b\in A$ such that $a\preceq b$, the 
    Shift Lemma applies to $(\pi_{a,b}, \sigma_{a,b}, E_a, E_b)$.
    \end{enumerate}
For $a,b\in A$ such that $a\preceq b$, let $\tau_{a,b}$ be 
copy map associated to $(\pi_{a,b}, \sigma_{a,b}, E_a, E_b)$. 
Let $M = \lim\mathcal{M}$, $N=\lim \mathcal{N}$,  $\pi_a:M_a\to M$ and
 $\sigma_a:N_a\to N$ be the direct limit maps, and $E$ the eventual common value of $\pi_a(E_a)$.\footnote{As in the discussion of the direct limit of a system of plus trees, we must have $\pi_{a,b}(E_a)=E_b$ on a $\preceq$-cofinal subset of $A$, so this makes sense.}

Let $\mathcal{P} = \langle \{\ult(N_a, E_a)\}_{a\in A}, \{\tau_{a,b}\}_{a\preceq b}\rangle$,
  $P = \lim\mathcal{P}$, $\tau_a:Ult(N_a, E_a)\to P$ the direct limit maps,
   and $j: N\to P$ the unique map such that for every $a\in A$, the following
   diagram commutes
\begin{center}
    \begin{tikzcd}
    N_a \arrow[r, "\sigma_{a}"] \arrow[d, "E_a"'] 
    & N \arrow[d,"j"]\\
    Ult(N_a, E_a) \arrow[r, "\tau_{a}"'] 
    & P.
    \end{tikzcd}
    \end{center}
Then $P= Ult(N,E)$, $j=i^N_E$, and for all $a\in A$, $\tau_a$ is the copy map
associated to $(\pi_a, \sigma_a, E_a, E)$.

Moreover, if $\tau_{a,b}$ is elementary for every $a,b\in A$ such that $a\preceq b$, then for every $a\in A$, $\tau_a$ is elementary.\footnote{For example, if $\langle \pi_{a,b}, \sigma_{a,b}\rangle:(N_a, E_a)\to^*(N_b, E_b)$ for every $a,b\in A$, then $\tau_a$ is elementary for every $a\in A$.}
\end{lemma}
 
 The following diagram illustrates the situation along chains of $\preceq$.

 \begin{center}
    \begin{tikzcd} 
    M_a \arrow[r, "\pi_{a,b}"] 
    & M_b \arrow[r] \arrow[rr, bend left, "\pi_{b}"]& {} \arrow[r, "\cdots" ,phantom]&  M\\
    E_a \arrow[r, "\mapsto", phantom]
 \arrow[u, "\scriptsize \inup" ,phantom]& E_b\arrow[r, "\mapsto", phantom]
\arrow[rr,mapsto, bend left]& {} \arrow[r, "\cdots" ,phantom]& E\arrow[u, "\scriptsize\inup" ,phantom]\\
    N_a \arrow[r, "\sigma_{a,b}"] \arrow[d, "E_a"]
    & N_b \arrow[d,"E_b"] \arrow[r]\arrow[rr, bend left, "\sigma_{b}"] &
     {} \arrow[r, "\cdots" ,phantom]& N \arrow[d,"j"]\\
    Ult(N_a, E_a) \arrow[r, "\tau_{a,b}"'] 
    & Ult(N_b, E_b) \arrow[r]  \arrow[rr, bend left,
     "\tau_{b}"]& {} \arrow[r, "\cdots" ,phantom]& P
    \end{tikzcd}
    \end{center}
    
    The main relevant fact here is that  $\pi_a$ and $\sigma_a$
    agree on $\dom(E_a) \cup \lbrace \dom(E_a) \rbrace$.
    This follows at once from the fact that $\pi_{b,c}$ agrees with
     $\sigma_{b,c}$
    on $\dom(E_b) \cup \lbrace \dom(E_b) \rbrace$, for all
    $b \prec c$. We leave the calculations that show everything else fits together
    properly to the reader.

\begin{proof}[Proof of Proposition \ref{direct limit characterization}.]
Again, this proposition amounts to saying (1) implies (2) and (3) in the above
 definition of when the direct limit is wellfounded. We first show (1) implies (2).

\begin{claim}\label{direct limit claim 1}
Let $x\in D$. Suppose $M_x$ is defined and that $\leq$ is illfounded
 below $x$. Then $M_{x}$ is illfounded.
\end{claim}
\begin{proof}
We define an order preserving embedding $f$ from $\leq\restrict x$ 
into the ordinals of $M_x$.

Let $x=[a,\alpha]$ and fix $y=[b,\beta]<x$. Without loss of generality,
 we may assume $b\preceq a$ and $\beta+1<\lh(\tree{T}_a)$ (this is just because 
 we can move to $c\geq a,b$ where we have $[c, u_{a,c}(\alpha)]=[a,\alpha]$ and
  $[c,u_{b,c}(\beta)]=[b,\beta]$, so that $u_{b,c}(\beta)+1<u_{a,c}(\alpha)+1\leq
  \lh (\tree{T}_c)$, since $[b,\beta]<[a,\alpha]$). We let \begin{equation*}
    f(y)  = t^{b}_{u_{a,b}(\alpha)}(\lh( E^{\tree{T}_a}_\beta)).
\end{equation*}

Clearly $f$ maps $y$ to an ordinal of $M_{x}$. It's easy to check
 that it is (strictly) order preserving, so $M_{x}$ is illfounded.
\hfill{$\qed$ Claim \ref{direct limit claim 1}}
\end{proof}

So now suppose the $(D,\leq)$ is a well-order and that the models $M_{x}$ exist and
 are wellfounded (i.e. (1) and (2)). We show (3) by induction on $(D,\leq)$. 

More specifically, for $x\in D$, we let $\mathcal{D}^{\leq x} = 
\langle \tree{T}_a\restrict (x(a)+1), \Psi_{a,b}\restrict(\tree{T}_a\restrict 
(x(a)+1))\,|\, a\preceq b \wedge a,b\in \dom (x)\rangle$. It's easy to see that
 $\lim \mathcal{D}^{\leq x} = (\lim \mathcal{D})\restrict x+1$, where for the 
 $\Gamma$ systems we mean that for any $a\in \dom (x)$, $(u_a)^{\mathcal{D}^{\leq x}}=
  u_a\restrict x(a)+1$ and the $t$-maps are the same.

\begin{claim}\label{direct limit induction}
For all u-threads $x\in D$,
\begin{enumerate}
    \item[(i)] $\lim \mathcal D^{\leq x}$ is wellfounded.

\item[(ii)] for all $a\in \dom (x)$, $\Gamma_a\restrict (\tree{T}_a\restrict x(a)+1)$
 is an extended tree embedding from $\tree{T}_a\restrict x(a)+1$ into
  $\lim \mathcal{D}\restrict x+1$ and 

\item[(iii)] for all $a\in \dom (x)$, all $b\succeq a$, and
\[(\Gamma_b\circ\Phi_{a,b})\restrict x(a)+1\approx \Gamma_a\restrict x(a)+1.\]

\end{enumerate}
\end{claim}
We first show that the $u$ maps of our system preserve tree-predecessors 
of successor $u$-threads on a tail, in the following sense.

\begin{claim}\label{direct limit claim 2}
For any $x\in D$ which has a $\leq$-successor $z$ in $D$, there is $y\in D$ such
 that \[y=\leq^*\pred (z).\] Moreover,
\[y=\leq^*\pred (z)\, \Leftrightarrow\forall^* a\, y(a)=\tree{T}_a\pred (x(a)+1).\]
\end{claim}

\begin{proof}
Fix $x$ and $z$ it's $\leq$-successor. First note that for any
 $a\in \dom (z)\cap \dom (x)$, $z(a)=x(a)+1$, since $z(a)\leq x(a)+1$ 
 (as $[a, x(a)+1]$ is a $u$-thread $>x$) but $z(a)\neq x(a)$ (since $z\neq x$).
  We'll show that there is a $u$-thread $y$ such that 
  $\forall ^* a y(a)= \tree{T}_a\pred(x(a)+1))$. 
Suppose that there is no immediate $\leq^*$-predecessor of $z$. We'll show that
 $\leq$ is illfounded, a contradiction.

We define sequences $\langle a_n\,|\,n\in \omega\rangle$,
 $\langle \beta_n\,|\,n\in \omega\rangle$ such that $a_n\prec a_{n+1}$ 
 and $\beta_n=\tree{T}_{a_n}\pred (x(a_n)+1)$ but $\beta_{n+1}<u_{a_n,a_{n+1}}(\beta_n)$. 
 Then, taking $y_n=[a_n, \beta_n]$ gives a witness to the illfoundedness $\leq$.

We start with any $a_0\in \dom (z)$ and take $\beta_0 = \tree{T}_{a_0}\pred (x(a_0)+1)$,
 as we must.

Given $a_n$ and $\beta_n=\tree{T}_{a_n}\pred(x(a_n)+1)$, let 
\[a_{n+1}>\xi_n \text{ least such that } u_{a_n,a_{n+1}}
(\beta_n)\neq \tree{T}_{a_{n+1}}\pred(x(a_{n+1})+1),\]
We have that such an $a_{n+1}$ exists, since otherwise $y_n=[a_n, b_n]$ is the 
immediate predecessor of $z$ in $\leq^*$. Now let 
$\beta_{n+1}= \tree{T}_{a_{n+1}}\pred (x(a_{n+1})+1)$. Since $\Psi_{a_n, a_{n+1}}$ is a
 tree embedding, we must have that $\beta_{n+1}\in 
 [v_{a_n, a_{n+1}}(\beta_n),u_{a_n,a_{n+1}}(\beta_n)]_{\tree{T}_{{a_{n+1}}}}$. 
 So since $\beta_{n+1}\neq u_{a_n,a_{n+1}}(\beta_n)$, we have 
 $\beta_{n+1}<u_{a_n, a_{n+1}}(\beta_n)$, as desired.
\hfill{$\qed$ Claim \ref{direct limit claim 2}}
\end{proof}

Let $z$ be a $u$-thread of successor rank, say $z$ is the $\leq$-successor
 of $x$ (i.e. the rank of $z$ is the rank of $x$ plus one). The observation made at 
 the start of the previous proof shows $z(a)=x(a)+1$ for all most all $a$. Fix such an
  $a$, so for all $b\succeq a$ $z(b)=x(b)+1$. So we have for all $b\succeq a$,
\begin{align*} z(b) &= x(b)+1\\ &=u_{a,b}(x(a))+1\\ &= v_{a,b}(x(a)+1)\\&= v_{a,b}(z(a)).
\end{align*} This shows that all successor $u$-threads are actually $v$-threads (defined 
in the obvious way) when $\leq$ is wellfounded. Even when $\leq$ is wellfounded, there 
may be $u$-threads which are not $v$-threads, so it was important to use $u$-threads in
 defining the direct limit.

Going forward, if 
$x$ is $u$-thread which is not the $\leq$-largest $u$-thread, we'll let $x+1$ be
 the $\leq$-successor of $x$.
\begin{proof}[Proof of Claim \ref{direct limit induction}.]
\setcounter{claim}{2}

We proceed by induction. We already know all the $M_y$ are defined and wellfounded 
and $\leq$ is wellfounded, so to show (i) we just need to see that $\langle M_y,E_y,
 \leq^*\restrict x\rangle$ is a plus iteration tree. 

In the base case, where $x$ is the minimum $u$-thread, this is unique tree of 
length one on the base model and (ii) and (iii) hold trivially. For the successor case, 
suppose we have (i)-(iii) for all $z\leq x$ and suppose that $x$ is not the last
 $u$-thread. So $x$ has a $\leq$-successor, $x+1$ and, appealing to Claim 
 \ref{direct limit claim 2}, we can take $y=\leq^*\pred(x+1)$. To show  (i) 
 we need to see that we're applying $E_x$ following the rules of quasi-normality, i.e.

\begin{subclaim}\label{direct limit subclaim1} $y$ is $\leq$-least such that
 $\crit (E_{x})<\lambda(E_{y})$ and for $P\is M_{y}$ least such that $M_{y}|\lh (E_{y})\is P$
  and $\rho(P)\leq \crit (E_{x})$,  \[M_{x+1}=\ult(P,E_{y}).\]
\end{subclaim}

\begin{proof}
By Claim \ref{direct limit claim 2}, we may take $a$ such that for all 
$b\succeq a$, $y(b)=\tree{T}_b \pred (x(b)+1)$.

For the appropriate choice of maps and models, we are now exactly in the 
situation of Lemma \ref{shift direct limits}. The rest of the Subclaim follows from the
 normality of each of the trees $\tree{T}_b$. We leave the details to the reader.
\hfill{$\qed$ Claim \ref{direct limit subclaim1}}
\end{proof}
$(ii)$ and $(iii)$ at $x+1$ easily follow.

Suppose now $x$ has limit rank and $(i)$-$(iii)$ hold for all $z<x$. We first need
 to see
\begin{subclaim}\label{direct limit subclaim2}
For $b=[0,x)_{\leq^*}\defeq \{y\,|\, y<^*x\}$, $b$ is $<$-cofinal in $x$, there are 
finitely many drops along $b$, and $M_{x}$ is the direct limit along $b$.
\end{subclaim}

\begin{proof}
To see that $b$ is cofinal, let $y<x$. Since $x$ has limit rank, $y+1<x$. Let $a$ be 
sufficiently large such that $y+1=[a,y(a)+1]$. Let $\gamma+1$ least such that 
$y(a)+1\leq\gamma+1\leq_{\tree{T}_a}x(a)$. We have that for all 
$b\succeq a$, \[y(b)<u_{a,b}(\gamma)+1=v_{a,b}(\gamma+1)\leq_{\tree{T}_b}v_{a,b}
(x(a))\leq_{\tree{T}_b} x(b),\] using here that $\Phi_{a,b}$ is a tree embeddings 
(and the $v$-maps of tree embeddings are tree-order preserving). So letting
 $z=[a,\gamma]+1$, we have that $y<z\leq^* x$. Since $x$ is not a successor, we actually
  have $y<z<^*x$, as desired.

Since the model $M_x$ is defined, there is an $a$ such that for all $b\succeq a$,
 $t^{a,b}_{x(a)}$ is total. Suppose first that there is some successor 
 $\eta<_{\tree{T}_a} x(a)$ such that $(\eta,x(a)]_{\tree{T}_a}$ doesn't drop. Then 
 for all $b\succeq a$, we have that $[v_{a,b}(\eta), x(b))_{\tree{T}_{b}}$ doesn't drop.
  There is some $u$-thread $z$ such that $z=[b,v_{a,b}(\eta)]$ for all sufficiently 
  large $\zeta'$. Now, any drops from $z$ to $x$ in the direct limit corresponds to a 
  drop in $[z(b), x(b))_{\tree{T}_b}$ for all sufficiently large $b$, so there are no 
  such drops.

In the remaining case, $x(a)$ is a successor ordinal and a drop in $\tree{T}_a$.
 Let $\beta=\tree{T}_a\pred (x(a))$. Letting $z$ the $u$-thread such that 
 $z(b)=v_{a,b}(x(a))$ for all sufficiently large $b$, we have $z<^* x$ and there can
  be no drops between $z$ and $x$ in the direct limit tree, just as before (since for
   all sufficiently large $b$, there are no drops in $(z(b), x(b)]_{\tree{T}_{b}}$,
    as $t^{a,b}_{x(a)}$ is total). By induction, this means there are only finitely 
    many drops.

Using our induction hypotheses $(ii)$ and $(iii)$ for $z<x$, it is straightforward 
to check that $M_x$ is the direct limit along $b$, so we leave it to the reader.
\hfill{$\qed$ Subclaim \ref{direct limit subclaim2}}

This gives us (i) and it is now easy to verify $(ii)$ and $(iii)$.
\end{proof}\hfill{$\qed$ Claim \ref{direct limit induction}}
\end{proof}
\hfill{$\qed$ Proposition \ref{direct limit characterization}}
\end{proof}

Notice that in the proof of Proposition \ref{direct limit characterization},
 we also verified (or rather left it to the reader to verify) that when the
  direct limit $\lim\mathcal{D}$ is wellfounded, then $\Gamma_a =
   \langle v_a, u_a, \{s^a_\xi\},\{t^a_\xi\}\rangle$ is an extended tree embedding from
    $\tree{T}_\xi $ into $\tree{U}$, the direct limit tree, and 
    $\Gamma_b\circ \Psi_{a,b}= \Gamma_a$ when $a\preceq b$.

 According to the following proposition,
 the direct limit we just defined is indeed the direct limit in the category of plus
 iteration trees of successor lengths and extended tree embeddings. 

\begin{proposition}\label{direct limit prop}
Let $\mathcal{D}= \langle \tree{T}_a, \Psi_{a,b}, \preceq\rangle$ be a directed 
system of trees, where $\preceq$ has field $A$.
Suppose there is a normal tree $\tree{S}$, and for each $a\in A$ an extended
 tree embedding $\Pi_a:\tree{T}_a\to \tree{S}$  such that
  whenever $a\preceq b$, $\Pi_b=\Psi_{a,b}\circ \Pi_a$; then the direct limit 
  $\lim\mathcal{D}$ is wellfounded, and there is a unique tree embedding 
  $\Pi:\lim \mathcal{D}\to \tree{S}$ such that $\Pi_a= \Pi\circ \Gamma_a$ for
   all $a\in A$.
\end{proposition}

We omit the straightforward proof.

\begin{remark}\label{direct limit remark} We can define the direct limit of a commuting 
system of trees under ordinary tree embeddings in the obvious way and verify the versions
of Propositions \ref{direct limit characterization} and \ref{direct limit prop}. It is
 easy to see that the direct limit of a system of trees under extended tree embeddings is
  either the same as the corresponding direct limit of under ordinary tree embeddings, or
   else is  some tree $\tree{U}$ with length $\gamma+1$ for a limit ordinal $\gamma$, and
    the corresponding direct limit under ordinary tree embeddings is just
     $\tree{U}\restrict\gamma$.
\end{remark}

\subsection{Embedding normalization and quasi-normalization}

   We begin with one-step embedding normalization.
Let $\tree{S}$ and $\tree{T}$ be normal trees of successor length on some
common base premouse $M$.  Let $F$ be on the sequence of last model of $\tree{T}$. Put
 \[
 \alpha=\alpha(F,\tree{T})<\lh (\tree{T})= \\\text{least $\gamma$ such that $F$ is
  on the sequence of $M^\tree{T}_\gamma$}
  \]
 and 
  \[
  \beta=\beta(F,\tree{T}) = \text{least $\gamma$ such that $\gamma=\alpha$ or 
  $\lambda(E^\tree{T}_\beta)>\crit (F)$.}
  \]
Suppose that $\tree{S}\restrict\beta+1=\tree{T}\restrict\beta+1$ and 
$\dom F\leq \lambda( E^\tree{S}_\beta)$, 
if $\beta+1<\lh (\tree{S})$.
In this case, we define a tree
\[
\tree{W}=W(\tree{S},\tree{T},F),
\]
a partial extended tree embedding
\[
\Phi^{W(\tree{S},\tree{T},F)} = \Phi \colon \tree{S} \to \tree{W},
\]
and a nearly elementary map
\[
\sigma^{W(\tree{S},\tree{T},F)} = \sigma\colon\ult(P, F) \to M_\infty^{\tree{W}},
\]
 where $P$ is the largest initial segment of the last model of $\tree{S}$ to which
we can apply $F$, and $M_\infty^{\tree{W}}$ is the last model of $\itW$.
 In general, we may reach illfounded models in forming $\tree{W}$, and we stop if so. 
 We say that $\tree{W}$ is wellfounded if we never reach illfounded models.
  If $\tree{S}$ and $\tree{T}$ are by a strategy $\Sigma$ which has
strong hull condensation, then $\tree{W}$ will be wellfounded.

We let $\tree{W}\restrict \alpha+1 =\tree{T}\restrict \alpha+1$ and
 $E^\tree{W}_\alpha=F$.
For the rest of $\tree{W}$, we consider cases. Let $Q$ be the initial segment
of $M_\beta^\itS$ to which $F$ is must be applied in a normal tree.

\paragraph{The dropping case.} $\beta+1 = \lh (\itS)$ and 
$Q\isneq M^\tree{S}_\beta$ or
$\beta + 1 < \lh (\itS)$, 
$Q\isneq M^\tree{S}_\beta|\lh (E^\tree{S}_\beta)$.  

In this case we have described all of $\tree{W}$ already:
 \[
 \tree{W}=\tree{T}\restrict\alpha+1 \conc \langle F\rangle,
 \]
and $\Phi \restriction \beta+1$ is just the identity
 on $\tree{S}\restrict\beta+1$. $\Phi$ is the associated extended tree embedding:
 $u(\beta)=\alpha+1$ and $t_\beta = \ihat_{v(\beta),u(\beta)}^\itW \circ s_\beta =
 i_F^Q \circ \mbox{id}$ = $i_F^Q$.

In this case $P=Q$, and  
$\ult(P,F) = M_{\alpha +1}^\itW = \ult(P,F)$ is the last model of $\tree{W}$, so
we set $\sigma=id$.

In the dropping case, $\Phi$ is total exactly when $\beta+1=\lh (\tree{S})$.

\paragraph{The non-dropping case.} Otherwise.

We define the $u$ and $v$ maps of $\Phi$ by
\begin{equation*}
    u(\xi) =
    \begin{cases*}
      \xi & if $\xi<\beta$, \\
      \alpha+1+(\xi-\beta) & if $\xi\geq \beta$,
    \end{cases*}
  \end{equation*}
and  
 \begin{equation*}
    v(\xi) =
    \begin{cases*}
      \xi & if $\xi\le \beta$, \\
      \alpha+1+(\xi-\beta) & if $\xi > \beta$.
    \end{cases*}
  \end{equation*}
  So $u = v$, except that $v(\beta) = \beta <_W u(\beta)=\alpha+1$.
  The remainder of $\itW$ and $\Phi$ are determined by our having set
  $E_\alpha^\itW = F$ and the rules for tree embeddings. 
  
  For example,
  if $\xi+1 < \lh (\itS)$, then $E_{u(\xi)}^\itW$ is defined to
  be $t_\xi(E_\xi^\itS)$, which then determines
  $M_{u(\xi)+}^\itW$ by normality.  Letting $\eta = S\pred(\xi+1)$, and 
  letting $\eta^*$
  be what normality dictates for $W\pred(u(\xi)+1)$, we must see that
  $\eta^*$ is properly located in $\itW$, namely that
  $v(\eta) \le_{W_\eta} \eta^* \le_{W_\eta} u(\eta)$.
  But it is easy to check that if $\eta \neq \beta$, then $v(\eta) = \eta^* = u(\eta)$,
  while if $\eta = \beta$, then $\eta^* = v(\eta) = \beta$ when $\crit(E_\eta)
   < \crit(F)$,
  and $\eta^* = u(\eta) = \alpha+1$ when $\crit(F) \le \crit(E_\eta^\itS)$.
  In either case, $t_\xi$ agrees with $\ihat_{v(\eta),\eta^*}^\itS \circ s_\eta$
  on $\dom(E_\xi^\itS$. Letting $Q$ be what $E_\xi^\itS$ applies to in $\itS$,
  and $Q^*$ what $E_{u(\xi)}^\itW$ applies to in $W$, and
  \[
  \pi = \ihat^\itW \circ s_\eta \restriction Q,
  \]
  we get that $\pi$ agrees with $t_\xi$ on $\dom(E_\xi^\itS)$, so the Shift 
  Lemma applies to $(t_\xi,\pi, E_\xi^\itS)$, and we can let
  \[
  s_{\xi+1} \colon \ult(Q,E_\xi^\itS) \to \ult(Q^*,E_{u(\xi)}^\itW)
  \]
  be the associated copy map.
  
  \begin{remark} The argument just sketched is a special case of a
  general lemma on extending tree embeddings given in \S8.2 of \cite{nitcis}.
  \end{remark}

In the non-dropping case, $\Phi$ is a total,
non-dropping extended tree embedding from $\itS$ to $\itW$.
We have $\lh (\tree{W})= \alpha+1+(\lh (\tree{S)}-\beta)$,
and $\ran (u)=[0,\beta)\cup [\alpha+1,\lh (\tree{W}))$.

Since $t_\beta=i^P_F$ and all the $t_\xi$ for $\xi\geq \beta$ agree
 with $t_\beta$ beyond $\dom (F)$, we have that $F$ is an initial factor
  of the extender of these $t_\xi$. We let $\sigma$ be the unique map 
  such that the last $t$ map factors as $\sigma\circ i^N_F$, where $N$ is
   the last model of $\tree{S}$ (using here that $F$ is total on the last
    model by our case hypothesis).

Next we review the one-step \textit{quasi-normalization}, 
which applies to plus trees. Let $\tree{T}$ a plus tree of successor
 length and $F$ an extender such that $F^-$ is on the sequence
 on the last model of $\tree{T}$.

Put
 \begin{enumerate}
     \item[] 
 $\alpha_0=\alpha_0(F,\tree{T})=$ least $\gamma$ such that 
 \begin{enumerate}\item[i.] $\alpha(\tree{T},F^-)\leq \alpha_0$ and 
\item[ii.] $\lh(F)<\hat\lambda(E_\gamma^\tree{T})$, or $E_\gamma^\tree{T}$ is of plus type, 
or $\gamma+1=\lh(\tree{T})$.
 \end{enumerate}
 
 and 
  \item[]
  $\beta=\beta(F,\tree{T}) =$ least $\gamma$ such that 
  $\crit (F)< \lh(E^\tree{T}_\beta)$.
  \end{enumerate}

It's easy to see that $\beta\leq\alpha_0$ and if $\beta<\alpha_0$, then $\beta$ is least such that $\crit(F)<\hat\lambda(E_\beta^\tree{T})$.
  
Now suppose $\tree{S}$ and $\tree{T}$ be plus trees
 of successor length on a premouse $M$ such that $\tree{S}\restrict\beta+1=\tree{T}\restrict\beta+1$ and if $\beta+1<\lh(\tree{S})$, then $\dom(F)\isneq  M_\beta^\tree{S}|\lh(E_\beta^\tree{T})$.
 We define the quasi-normalization 
\[
\tree{V}=V(\tree{S},\tree{T},F),
\]
a partial extended tree embedding \[\Phi=\Phi^{V(\tree{S}, \tree{T}, F)}\] from $\tree{S}$ into $\tree{V}$, and a nearly elementary map \[\sigma=\sigma^{V(\tree{S},\tree{T}, F)}:\ult(P,F)\to M_\infty^\tree{V},\] where $P$ is the longest initial segment of the last model of $\tree{S}$ to which we can apply $F$. As before,
 we may reach illfounded models in forming $\tree{V}$ and 
stop if we do. Otherwise, we say $\tree{V}$ is wellfounded.

We set \[\tree{V}\restrict\alpha_0+1=\tree{T}\restrict\alpha_0+1\] and
 put 
\[
E_{\alpha_0}^\tree{V}=F.
\]
 For the rest of the definitions, we split into cases, 
as before. Let $Q$ be the initial segment of $M_\beta^\tree{S}$ 
to which $F$ must be applied in a quasi-normal tree.

\paragraph{The dropping case.} Either $\beta+1 = \lh (\itS)$ and 
$Q\isneq M^\tree{S}_\beta$ or
$\beta + 1 < \lh (\itS)$ and
$Q\isneq M^\tree{S}_\beta|\lh (E^\tree{S}_\beta)$.  

In this case we have described all of $\tree{V}$ already:
 \[
 \tree{V}=\tree{T}\restrict\alpha_0+1 \conc \langle F\rangle.
 \]
We also set $\Phi \restriction \beta+1$ to be the identity
 on $\tree{S}\restrict\beta+1$. $\Phi$ is the associated extended tree embedding:
 $u(\beta)=\alpha_0+1$ and $t_\beta = \ihat_{v(\beta),u(\beta)}^\itW \circ s_\beta =
 i_F^Q \circ \mbox{id}$ = $i_F^Q$. In this case $P=Q$ and  
$\ult(P,F) = M_{\alpha +1}^\itV $ is the last model of $\tree{V}$, and so
we set $\sigma=id$.

Note that in the dropping case, $\Phi$ is total exactly when $\beta+1=\lh (\tree{S})$.

\paragraph{The non-dropping case.} Otherwise.

We define the $u$ and $v$ maps of $\Phi$ by
\begin{equation*}
    u(\xi) =
    \begin{cases*}
      \xi & if $\xi<\beta$, \\
      \alpha_0+1+(\xi-\beta) & if $\xi\geq \beta$,
    \end{cases*}
  \end{equation*}
and  
 \begin{equation*}
    v(\xi) =
    \begin{cases*}
      \xi & if $\xi\le \beta$, \\
      \alpha_0+1+(\xi-\beta) & if $\xi > \beta$.
    \end{cases*}
  \end{equation*}
  So $u = v$, except that $v(\beta) = \beta <_V u(\beta)=\alpha_0+1$.
  The remainder of $\itV$ and $\Phi$ are determined by our having set
  $E_{\alpha_0}^\itW = F$ and the rules for tree embeddings. By induction on $\xi\geq \beta$, one shows that there are nearly elementary embeddings $\tilde{\sigma}_\gamma:\ult(M_\gamma^\tree{S}, F)\to M_{u(\xi)}^\tree{V}$ and we let $\sigma$ be the last of these embeddings. 
  
  In the non-dropping case, $\Phi$ is always a total extended tree embedding (assuming that $\tree{V}$ is wellfounded).
  
 This finishes our description of the quasi-normalization. 

The following lemma, due to Siskind, connects the one-step quasi-normalization
 to an analysis of particularly well-behaved tree embeddings.

\begin{lemma} [Factor Lemma]\label{factor lemma}
Let $\Psi: \tree{S}\to \tree{T}$ be an (extended) tree embedding such 
that $\Psi\neq Id$. Let $\beta=\crit (u^\Psi)$ and $\alpha_0+1$ be the successor
 of $v^\Psi(\beta)=\beta$ in $(v^\Psi(\beta), u^\Psi(\beta)]_\tree{T}$. Suppose that $\dom(E^\tree{T}_{\alpha_0})\isneq M_\beta^\tree{S}|\lh(E_\beta^\tree{S})$.
Then $V(\tree{S},\tree{T}\restrict\alpha_0+1, E^\tree{T}_{\alpha_0})$ is defined
 and wellfounded and there is a unique (extended) tree embedding
  $\Gamma: V(\tree{S},\tree{T}\restrict\alpha_0+1, E^\tree{T}_{\alpha_0})\to \tree{T}$
   such that $u^\Gamma\restrict\alpha_0+1=id$ and $\Psi= \Gamma\circ \Phi^{V(\tree{S},\tree{T}\restrict\alpha_0+1,
    E^\tree{T}_{\alpha_0})}$.
\end{lemma}
\begin{proof}[Proof sketch.]
First notice that our hypotheses guarantee that $V(\tree{S},\tree{T}\restrict\alpha_0+1, E^\tree{T}_{\alpha_0})$ is defined. Now let $\tree{V}=V(\tree{S},\tree{T}\restrict\alpha_0+1, E^\tree{T}_{\alpha_0})$ and
 $\Phi= \Phi^{V(\tree{S},\tree{T}\restrict\alpha_0+1, E^\tree{T}_{\alpha_0})}$.
The commutativity condition together with the demand that $u^\Gamma\restrict\alpha_0+1=\id$ totally determine the $u$-map of $\Gamma$:
\begin{equation*}
    u^\Gamma(\xi)=\begin{cases} \xi &\text{ if } \xi<\alpha_0+1\\
    u^\Psi\circ(u^\Phi)^{-1}(\xi)&\text{ if } \xi\geq \alpha_0+1,
    \end{cases}
\end{equation*}
using in the second case that
 $[\alpha_0+1, \lh (\tree{V}))\subseteq \ran (u^\Phi)$.
One must check by induction on $\xi$ that $u^\Gamma\restrict(\xi+1)$
is the $u$-map of a tree embedding from $\tree{V}\restrict (\xi+1)$ into $\tree{T}$. For this, one uses the result on extending tree embeddings,
Proposition 8.2.1 of \cite{nitcis}. \qed
\end{proof}
The Factor Lemma gives us a sense in which the one-step quasi-normalization
$V(\itS,\itT,F)$ is a minimal plus tree that uses $F$ and tree-embeds $\itS$.\footnote{$\itT$
is relevant too, but if we are dealing with trees by a fixed iteration strategy,
then $F$ determines $\itT$.} It is parallel to the fact that
if $j \colon M \to N$ is elementary, and $E$ is an initial segment of the
extender of $j$, then $\ult(M,E)$ embeds into $N$ in a way that makes
the diagram commute. Indeed, one can think of $V(\itS,\itT,F)$ as an
$F$-ultrapower of $\itS$.

A version of the  Factor Lemma also holds for embedding normalization; this appears in in \cite{associativity}.

We can use iterated applications of the Factor Lemma to factor
 appropriate non-identity extended tree embedding $\Psi$ as
  \[\cdots \circ\Phi_{F_{\xi+1}}\circ\Phi_{F_\xi}\circ\cdots\circ \Phi_{F_1}\circ \Phi_{F_0},\]
where the $F_\xi$ are a (non-overlapping) sequence of exit extenders of
 $\tree{T}$ and the $\Phi_{F_\xi}$ are certain one-step quasi-normalizations by $F_\xi$.
 
To do this, we need that the additional hypotheses on the relationship between the domain of $F$ in the Factor Lemma statement obtain at every step. Rather than state this directly, we identify a simpler, stronger property which suffices for our applications.

\begin{definition}
A tree embedding $\Psi:\tree{S}\to \tree{T}$ is \textit{inflationary} iff for any $\xi+1<\lh(\tree{S})$ and $\gamma+1\in (v^\Psi(\xi), u^\Psi(\xi)]_\tree{T}$, letting $\eta=\tree{T}\pred(\gamma+1)$,
\[\dom(E_\gamma^\tree{T})\isneq \lh(s_{\xi, \eta}^\Psi(E_\xi^\tree{S})).\]
\end{definition}
 
It is easy to see that the quasi-normalization tree embeddings $\Phi^{V(\tree{S}, \tree{T}, F)}$ are inflationary, but tree embeddings (even hull embeddings) may fail to be inflationary, in general. Still, one can prove the following useful closure properties of the class of inflationary tree embeddings.

\begin{proposition}
Let $\Phi:\tree{S}\to \tree{T}$ and $\Psi:\tree{T}\to \tree{U}$ be inflationary tree embeddings. Then $\Psi\circ \Phi$ is inflationary.
\end{proposition}

\begin{proposition}
Let $\Phi:\tree{S}\to \tree{T}$ and $\Psi:\tree{T}\to \tree{U}$ be tree embeddings such that $\Phi$ and $\Psi\circ \Phi$ are inflationary and $[\crit(u^\Psi), \lh(\tree{T}))\subseteq \ran(u^\Phi)$. Then $\Psi$ is inflationary.
\end{proposition}

\begin{proposition}
Let $\mathcal{D}=\langle \{\tree{T}_a\}_{a\in A}, \{\Psi_{a,b}\}_{a\preceq b}\rangle$ be a directed system of plus trees such that all tree embeddings $\Psi_{a,b}$ are inflationary. Suppose $\lim \mathcal{D}$ is wellfounded and let $\{\Gamma_a\}_{a\in A}$ be the direct limit tree embeddings from $\tree{T}_a$ into the direct limit tree. Then for all $a\in A$, $\Gamma_a$ is inflationary.
\end{proposition}

We omit the somewhat tedious but straightforward proofs. The reader can find proofs of each proposition for the representative case that all trees are $\lambda$-tight and normal in \cite{associativity}.

Now we return to our iterative factorization result.

\begin{theorem}
Let $\Psi:\tree{S}\to \tree{T}$ be an inflationary extended tree embedding. Then there is a unique sequence of extenders $\langle F_\xi\mid \xi<\lambda\rangle$ such that there is a directed system of plus trees $\mathcal{D}=\langle \{\tree{S}_\xi\}_{\xi\leq\lambda}, \{\Psi_{\eta,\xi}\}_{\eta\leq \xi\leq \lambda}\rangle$ such that \begin{enumerate}
    \item $\tree{S}_0=\tree{S}$, $\tree{S}_\lambda= \tree{T}$, and $\Psi_{0,\lambda}=\Psi$;
    \item for $\xi+1\leq \lambda$, letting $\beta_\xi=\crit(u^{\Psi_{\eta,\xi}})$ and $\alpha_\xi+1$ the successor of $\beta_\xi$ in $(\beta_\xi, u^{\Psi_{\xi, \lambda}}(\beta_\xi)]_\tree{T}$,
    \begin{enumerate}
        \item $F_\xi=E_{\alpha_\xi}^\tree{T}$, 
        \item $\tree{S}_{\xi+1}=V(\tree{S}_\xi, \tree{T}\restrict\alpha_\xi+1, F_\xi)$, and
        \item $\Psi_{\xi, \xi+1}=\Phi^{V(\tree{S}_\xi, \tree{T}\restrict\alpha_\xi+1, F_\xi)}$;\end{enumerate}
        \item for $\gamma\leq \lambda$ a limit ordinal,
        \begin{enumerate}
            \item $\tree{S}_\gamma= \lim \langle \{\tree{S}_\xi\}_{\xi<\gamma}, \{\Psi_{\eta, \xi}\}_{\eta\leq \xi<\gamma}\rangle$, and
            \item for $\xi<\gamma$, $\Psi_{\xi, \gamma}$ is the direct limit tree embedding;
            \item for $\xi\leq \lambda$ and $\eta<\xi$, $u^\Psi_{\xi, \lambda}\restrict \alpha_\eta+1 =id$.
        \end{enumerate}
    \end{enumerate}
\end{theorem}

\begin{proof}[Proof sketch.]
We define $\tree{S}_\xi$, $\Psi_{\eta, \xi}$, and $F_\xi$ by induction, maintaining that we always have additional inflationary extended tree embeddings $\Gamma_\xi: \tree{S}_\xi\to \tree{T}$ which commute appropriately with the rest of our system. Assuming that the tree embeddings produced so far are inflationary, one continues at limits by taking direct limits at limits and at successors by applying the factor lemma to the last $\Gamma_\xi$ (as long as this is not the identity). The preceding propositions guarantees that all the resulting tree embeddings produced in this way are actually inflationary so that this process never breaks down. 
Since we are pulling all the extenders $F_\xi$ from $\tree{T}$, this process must terminate at some stage $\leq lh(\tree{T})$, and so we must end up with a stage $\Gamma_\xi= Id$ (or else we could continue the construction).
A detailed proof in the case that all the trees are normal appears in \cite{associativity}. \qed
\end{proof}

For $\Psi:\tree{S}\to \tree{T}$ an inflationary extended tree embedding, we call the
 sequence $\vec{F}=\langle F_\xi\, |\, \xi<\lambda\rangle$ as in the previous theorem
 the
  \textit{factorization of $\Psi$}.
  
 One might ask whether $\vec{F}$ and $\itS$ determine
  $\Psi$ and $\itT$. It is easy to see that this is not the case;
  for example, consider the one-step normalization
  embeddings  $\Phi \colon \itS \to W(\itS,\itT,F)$ and
  $\Psi \colon \itS \to W(\itS,\itU,F)$, with $\itT$ and $\itU$ being distinct
  normal trees on $M$ that diverge before each of them reaches a model
  with $F$ on its sequence. Then $\Phi \neq \Psi$, but they have the
  common factorization $\langle F \rangle$.  If we are dealing
  with normal trees on $M$ by a fixed iteration strategy, this kind of
  thing can't happen, and in fact it is easy to see that then
  $\itS$ and the factorization $\vec{F}$ of $\Phi \colon \itS \to \itT$
  determine $\Phi$. In general, $\itS$ and $\vec{F}$ determine $\Phi$
  up to ``similarity".

First, some notation from \cite{nitcis}:
\begin{definition}
For $\tree{T}$ a plus tree and $\eta\leq_\tree{T}\xi$, we 
let $e^\tree{T}_{\eta,\xi}$ be the sequence of extenders used along 
$(\eta,\xi)_\tree{T}$, i.e. $e^\tree{T}_{\eta,\xi}=
\langle E_\alpha^\tree{T}\,|\,\alpha+1\in (\eta,\xi]_\tree{T}\rangle$.
\end{definition}

The next proposition says that if
$\tree{T}$ and $\tree{U}$ are plus trees having two common models 
$M,N$ which are tree-related in
 both trees, and have the same branch embeddings between them in both trees, 
 then the extenders used to get from $M$ to $N$ are the same in both trees.
 We omit the well known proof.

\begin{proposition}\label{branch extender prop}
Suppose $\tree{T}$, $\tree{U}$ are plus trees,
 $\eta\leq_\tree{T} \xi$, $\eta^*\leq_\tree{U} \xi^*$,
  $M_{\xi}^\tree{T}=M_{\xi^*}^\tree{U}$, and either
\begin{enumerate}
    \item $\eta$-to-$\xi$ and $\eta^*$-to-$\xi^*$ don't drop,
     $M_{\eta}^\tree{T}=M_{\eta^*}^\tree{U}$, and 
     $i^\tree{T}_{\eta,\xi}= i^\tree{U}_{\eta^*,\xi^*}$, or
    \item $\eta^*$ and $\eta$ are the locations of the last drop along 
    $\eta$-to-$\xi$ and $\eta^*$-to-$\xi^*$.
\end{enumerate} Then $e^\tree{T}_{\eta,\xi}=e^\tree{U}_{\eta^*,\xi^*}$.
\end{proposition}

\begin{definition}\label{similar tree embeddings}
Let $\Phi:\tree{S}\to \tree{T}$ and $\Psi:\tree{S}\to \tree{U}$ be extended
 tree embeddings. We say $\Phi$ and $\Psi$ are \textit{similar}, and write
  $\Phi\equiv \Psi$, iff for all $\xi<\lh (\tree{S})$,
   $e^\tree{T}_{0,u^\Phi(\xi)}=e^\tree{U}_{0,u^\Psi(\xi)}$.
\end{definition}
We can characterize similarity using the Factor Lemma.\footnote{The definitions
and results on factoring tree embeddings are due to Siskind.} 

\begin{theorem}\label{tree embedding uniqueness}
Let $\Phi:\tree{S}\to \tree{T}$ and $\Psi:\tree{S}\to \tree{U}$ are extended inflationary
 tree embeddings. Then $\Phi\equiv \Psi$ iff $\Phi$ and $\Psi$ have the same
  factorization.
    \end{theorem}
We omit the proof; the special but representative case that
 all the trees are $\lambda$-tight and normal appears in \cite{associativity}.

We now briefly describe the quasi-normalization of a length 2, maximal $M$-stack $\langle \tree{T},\tree{U}\rangle$ with $\tree{U}$ a normal, 
$V(\tree{T},\tree{U})$. This is the last tree
 system 
 \[\langle \tree{V}_\xi, \sigma_\xi,F_\zeta,\Phi^{\eta,\xi}\,|\, \eta,\xi,\zeta+1
 <\lh (\tree{U}), \, \eta\leq_\tree{U} \xi\rangle,
 \] 
 defined by induction on $\lh (\tree{U})$ by:
\begin{enumerate}
    \item $\tree{V}_\xi = V(\tree{T},\tree{U}\restrict \xi+1)$;
    \item $\sigma_\xi:M^\tree{U}_\xi\to R_\xi$ is nearly elementary, 
    where $R_\xi$ is the last model of $\tree{V}_\xi$, and 
    if $\xi+1<\lh (\tree{U})$, $F_\xi=\sigma_\xi(E^\tree{U}_\xi)$;
    \item For $\zeta\leq_\tree{U}\eta\leq_\tree{U}\xi$, 
    \begin{enumerate}
        \item $\Phi^{\eta,\xi}:\tree{V}_\eta\to \tree{V}_\xi$ is a partial 
        extended tree embedding, and
        \item $\Phi^{\zeta,\xi}=\Phi^{\eta,\xi}\circ \Phi^{\zeta,\eta}$;
    \end{enumerate}
    \item For $\eta=\tree{U}\pred(\xi+1)$,
    \begin{enumerate}
        \item $\tree{V}_{\xi+1}= V(\tree{V}_\eta,\tree{V}_\xi,F_\xi)$,
        \item $\Phi^{\eta,\xi+1}=\Phi^{V(\tree{V}_\eta,\tree{V}_\xi,F_\xi)}$ and
         $\sigma_{\xi+1}= \sigma^{V(\tree{V}_\eta,\tree{V}_\xi,F_\xi)}$;
    \end{enumerate}
    \item For $\lambda <\lh (\tree{U})$ a limit and $b=[0,\lambda)_\tree{U}$,
    \[\tree{V}_\lambda = \lim \langle\tree{V}_\xi,\Phi^{\eta,\xi}\,|\,
     \eta\leq_\tree{U}\xi\in b\rangle\] and
    $\Phi^{\xi, \lambda},\sigma_\lambda$ the direct limit tree embeddings
    and direct limit map.
\end{enumerate}

\noindent If $\tree{U}$ has successor length, we 
let $\sigma^{V(\tree{T},\tree{U})}$ be the last of the $\sigma_\xi$.

In general, we define $V(s)$ for a stack of plus trees $s$ of length
 $n+1$ and $\sigma_s$ from the last model of $s$ to the last model of
  $V(s)$ by induction on length $n$. For $s=\langle \tree{S}_1,\ldots \tree{S}_n\rangle$
   we put
    \[
    V(s)=V( V(s\restrict n), \sigma_{s\restrict n} \tree{S}_n),
    \]
and, if $\tree{S}_n$ has successor length, we let
 \[
 \sigma_s = \sigma^{V(V(s\restrict n), \sigma_{s\restrict n}\tree{S}_n)}\circ\sigma^*,
 \]
  where $\sigma^*$ is the last copy map from the last model of 
  $\tree{S}_n$ to the last model of $\sigma_{s\restrict n}\tree{S}_n$.
  
  \begin{remark} {\rm The normalizations we are considering here are ``bottom-up",
in the sense of \cite{nitcis}. One could normalize in other orders, for example,
by setting $\itV^*(\langle \itU_0,\itU_1,\itU_2\rangle =
 V(\itU_0,V(\itU_1,\itU_2))\rangle$.
We expect that all such quasi-normalizations are equivalent; in the case of embedding normalization for $\lambda$-tight normal trees, this has been proven by Siskind (see \cite{associativity}). The same proof should pass over to quasi-normalizations, but this has not been checked.}\end{remark}
  
  \subsection{Mouse pairs}
  
        We recall some notions from \cite{nitcis}. 
\begin{definition}
\label{normalizeswell}
        Let $\Sigma$ be a complete $(\lambda,\theta)$-iteration strategy for $M$, where
        $\lambda > 1$. We say that $\Sigma$ {\em quasi-normalizes well} iff whenever
        $s$ is an $M$-stack by $\Sigma$ and $\langle \tree{T}, \tree{U}\rangle$ is a maximal 2-stack by $\Sigma_s$ such that $\tree{U}$ is normal, then
        \begin{itemize}
                \item[(a)] $V(\tree{T}, \tree{U})$ is by $\Sigma_s$, and
                \item[(b)] letting $\tree{V}=V(\tree{T}, \tree{U})$ and $\sigma=\sigma^{V(\tree{T}, \tree{U})}$,
                $\Sigma_{s\conc\langle \tree{T}, \tree{U}\rangle} = (\Sigma_{s\conc\tree{V}})^\sigma$.
                \end{itemize}
                \end{definition}
Here $M$ is a premouse of some kind, coarse, pure extender, or least branch. 
It is easy to see that if $\Sigma$ normalizes well,
then so do all its tails $\Sigma_s$. 


Recall that $\itT$ is a {\em pseudo-hull} of $\itU$ iff there is
a tree embedding of $\itT$ into $\itU$.

\begin{definition}
  \label{stronghull}
Let $\Sigma$ be a complete $(\lambda,\theta)$-iteration strategy for a premouse $M$;
then $\Sigma$ has \textit{strong hull condensation}
        iff whenever $s$ is a stack
of plus trees by $\Sigma$ and $N\is M_\infty(s)$, and
$\itU$ is a plus tree on $N$ by $\Sigma_{s,N}$, then
for any plus $\itT$ on $N$,
\begin{itemize}
\item[(a)] if $\itT$ is a pseudo-hull of $\itU$,
then $\itT$ is by $\Sigma_{s,N}$, and
\item[(b)] if $\Phi \colon \itT \to \itU$ is an extended
tree embedding, with last $t$-map $\pi$, and $Q \unlhd \dom(\pi)$,
then $\Sigma_{s\conc \langle\itT\rangle,Q} = (\Sigma_{s\conc \langle\itU\rangle,\pi(Q)})^\pi$.
        \end{itemize}
\end{definition}

The most important part of \ref{stronghull} is clause (a) in the case
$s=\emptyset$; that is, that every pseudo-hull of a plus tree by
$\Sigma$ is also by $\Sigma$.

\begin{definition}\label{pureextenderpair} Let $\delta$ be regular.
$(M,\Omega)$ is
a {\em pure extender pair} with scope $H_\delta$
 iff
\begin{itemize}
\item[(1)]$M$ is a pure extender premouse, and $M \in H_\delta$,
\item[(2)] $\Omega$ is a complete $(\delta,\delta)$-iteration strategy for $M$,
and
\item[(3)] $\Omega$ quasi-normalizes well, is internally lift consistent\footnote{Roughly,
for all $N \lhd M$, $\Omega_N$ is consistent with $\Omega$ in the
natural sense. See \cite[Def. 5.4.4]{nitcis}.} and has strong
hull condensation.
\end{itemize}
\end{definition}

\begin{definition}\label{lbrhodpair} Let $\delta$ be regular.
$(M,\Omega)$ is a {\em least branch hod pair}
(lbr hod pair) with scope $H_\delta$ iff
\begin{itemize}
\item[(1)] $M$ is a least branch premouse, and $M \in H_\delta$,
\item[(2)] $\Omega$ is a complete $(\delta,\delta)$-iteration strategy
for $M$, 
\item[(3)] $\Omega$ normalizes
well, is internally lift consistent, and has strong hull condensation, and
\item[(4)] $(M,\Omega)$ is {\em pushforward consistent},
in that if $s$ is by $\Omega$ and has last model $N$, then $N$ is an lpm, and
$\hat\Sigma^N \subseteq \Omega_s$.
\end{itemize}
\end{definition}

\begin{definition}\label{mousepair} $(M,\Omega)$ is a {\em mouse pair}
        iff $(M,\Omega)$ either a pure extender pair, or an lbr hod pair.
\end{definition} 

According to these definitions, the strategy in
 a mouse pair with scope $H_\delta$ must be a $(\delta,\delta)$-strategy.
 \cite{nitcis} demanded only an $(\omega,\delta)$-strategy, so as to
 avoid having to talk about normalizing infinite stacks. Here we are
 not avoiding that at all. It is more natural for a strategy with scope
 $H_\delta$ to act on all stacks in $H_\delta$, that is, to be
 a $(\delta,\delta)$-strategy.

The paradigm case here is a mouse pair $(M,\Omega)$ with scope $\hc$, 
considered in a model
of $\adp$. One important way to obtain such pairs, assuming $\adp$,
 is
to start with a coarse $\Gamma$-Woodin pair $(N^*,\Sigma^*)$. Let
$\C$ be a full background construction of of one of the two types
done in  $N^*$, and let $(M,\Omega) = (M_{\nu,k}^\C,\Omega_{\nu,k}^\C)$
be one of its levels. From the point of view of $N^*$,
$(M,\Omega)$ is a mouse pair with scope $H_\delta$, where $\delta$ is
the $\Gamma$-Woodin cardinal of $N^*$. But $\Omega$ can be extended to
act on all stacks $s \in \hc$, because it is induced by $\Sigma^*$, which
acts on all stacks in $\hc$. So from the point of view of $V$,
$(M,\Omega)$ is a mouse pair with scope $\hc$.\footnote{These facts are
proved in \cite{nitcis}; see especially
Theorems 7.6.2 and 10.4.1. Those proofs deal only with $(\omega,\delta)$-iteration
strategies, but they adapt in a routine way to $(\delta,\delta)$-strategies.}
In sum

\begin{theorem}(\cite{nitcis}) 
Assume $\adp$, let $(N^*,\delta,S,T,\lhd,\Sigma^*)$
be a coarse $\Gamma$-Woodin tuple, and
let $\C$ be a $\lhd$-construction ( either pure extender or least branch)
 in
$L[N^*,S,T,\lhd]$ with all $F_\nu^{\C} \in N^*$; then $\C$ does not break
down. Moreover, letting
$M = M_{\nu,k}^\C$, and letting $\Omega$
be the canonical extension of $\Omega_{\nu,k}^{\C}$
to all $M$-stacks in $\hc$,
$(M,\Omega)$ is a mouse pair, with scope $\hc$.
\end{theorem}

The strategy in a mouse pair is actually determined by its action
on countable normal trees:

\begin{theorem}\label{uniqueextension} Let $\delta$ be regular,
and let $(M,\Sigma)$ and $(M,\Omega)$ be mouse pairs with scope
$H_\delta$. Suppose that $\Sigma$ and $\Omega$ agree on all
countable $\lambda$-separated trees; then $\Sigma = \Omega$.
\end{theorem}
This is proved as Lemma 7.6.5 in \cite{nitcis}, for the slightly weaker
notion of mouse pair used in that book. The same proof yields
\ref{uniqueextension}.

In a similar vein, one can show that a strategy for a pure
extender premouse that behaves well
on normal trees can be extended to the strategy component of a pure
extender pair.\footnote{There seems to be no way to show that the extended
strategy is pushforward consistent, so these strategy extension
results do not seem to help one to construct least branch hod pairs.}
One need only assume that every pseudo-hull of a normal tree by $\Sigma$ is by 
$\Sigma$. The following theorem, which combines work of Schlutzenberg
and Siskind, captures this fact.

\begin{theorem} (\cite{farmer}, \cite{associativity})\label{strategyextension}
Let $M$ be a countable premouse, 
and $\Sigma$ an $\omega_1 +1$-iteration
strategy for $M$. Suppose that every pseudo-hull of a countable
tree by $\Sigma$ is by $\Sigma$; then there is a unique $(\omega_1,\omega_1)$-iteration
strategy $\Sigma^*$ for $M$ such that 
\begin{itemize}
\item[(1)] $\Sigma \restriction \hc \subseteq \Sigma^*$, 
\item[(2)] $\Sigma^*$ quasi-normalizes well, and
\item[(3)] $\Sigma^*$ has strong hull condensation.
\end{itemize}
\end{theorem}
That such a strategy $\Sigma$ extends uniquely to a strategy for finite stacks 
of normal trees that normalizes well was proved
independently by Schlutzenberg and Steel. (See \cite{nitcis}, Theorem 7.3.11.)
 The extension to countable stacks
requires significantly more work, and is due to Schlutzenberg. (See \cite{farmer}.)
We shall outline the proof in the next section.
That the resulting strategy $\Sigma^*$ has strong
hull condensation is due to Siskind. (See \cite{associativity}.)

\section{Meta-iteration trees}

The system $\langle \tree{V}_\xi, F_\eta
 \leq_\tree{U}\,|\, \xi,\eta+1<\lh (\tree{U})\rangle$ that arises in the
 construction of $V(\itT,\itU)$ is a kind of 
 tree of iteration trees, with tree-embeddings between tree-order related nodes. 
 This perspective is used in \cite{nitcis} and abstracted in \cite{farmer},
  \cite{jensen} to their notions of \lq\lq factor trees of inflations" and 
  \lq\lq insertion iterations", respectively. Here, we isolate a 
  closely related abstraction, and call the resulting system a
   \textit{meta-iteration tree}, or \textit{meta-tree}.

   \subsection{The definition}

\begin{definition}\label{meta-trees} Let $\itS$ be a plus tree on a
premouse $M$. 
A \textit{meta-iteration tree on $\itS$} is a system 
\[
\mtree{S}=\langle \lh (\mtree{S}), \leq_\mtree{S}, \{\tree{S}_\xi\}_{\xi<\lh (\mtree{S})},
 \{F_\xi\}_{\xi+1<\lh (\mtree{S})}, \{\Phi^{\eta,\xi}\}_{\eta\leq_\mtree{S} \xi}\rangle
 \]
  such that $\itS_0 = \itS$, and
\begin{enumerate}
    \item $\lh(\mtree{S})$ is an ordinal and $\leq_\mtree{S}$ is a 
    tree-order on $\lh(\mtree{S})$;
    \item for all $\zeta\leq_\mtree{S}\eta\leq_\mtree{S}\xi<\lh (\mtree{S}$),
    \begin{enumerate}
        \item $\tree{S}_\xi$ is a plus tree on $M$ of successor length,
        \item if $\xi+1<\lh (\mtree{S})$, $F_\xi^-$ is an extender on 
        the sequence of $M_\infty^{\tree{S}_\xi}$,
        \item $\Phi^{\eta,\xi}$ is a partial extended tree embedding 
        from $\tree{S}_\eta$ into $\tree{S}_\xi$,
        \item $\Phi^{\zeta,\xi}=\Phi^{\eta,\xi}\circ\Phi^{\zeta,\eta}$;
    \end{enumerate}
    \item (Normality) 
    for $\xi+1<\lh(\mtree{S})$, letting $\alpha_\xi=\alpha_0(\tree{S}_\xi, F_\xi)$,
    \begin{enumerate}
        \item for all $\eta<\xi$, $\lh (F_\eta)<\lh (F_\xi)$,
        \item for $\eta=\mtree{S}\pred(\xi+1)$,
        \begin{enumerate}
            \item $\eta$ is least such that $\crit (F_\xi)<\lambda(F_\eta)$,
            \item $V(\tree{S}_\eta,\tree{S}_\xi,F_\xi)\restrict 
\alpha_\xi+2\subseteq\tree{S}_{\xi+1}\subseteq V(\tree{S}_\eta,\tree{S}_\xi,F_\xi)$,
\item $\Phi^{\eta,\xi+1}= \Phi^{V(\tree{S}_\eta,\tree{S}_\xi,F_\xi)}$
 (as a partial extended tree embedding from $\tree{S}_\eta$ into $\tree{S}_{\xi+1}$);
     \end{enumerate}
    \end{enumerate}
   \item for $\lambda<\lh(\mtree{S})$, $b=[0,\lambda)_\mtree{S}$ is a cofinal subset
of $\lambda$ and there is a tail $c\subseteq b$ such that \begin{enumerate}
  \item for all $\eta,\xi\in c$ with $\eta\leq \xi$, $\Phi^{\eta,\xi}$ is total,
        \item $\lim\langle \tree{S}_\xi,\Phi^{\eta,\xi}\,|\,
        \xi\leq_\mtree{S}\xi\in c\rangle$ is wellfounded,
        \item $\lim\langle \tree{S}_\xi,\Phi^{\eta,\xi}\,|\,
        \xi\leq_\mtree{S}\xi\in c\rangle\restrict \sup \{ \alpha_\xi\,|\,
         \xi<\lambda\}\subseteq\tree{S}_\lambda \subseteq \lim\langle \tree{S}_\xi,
         \Phi^{\eta,\xi}\,|\,\xi\leq_\mtree{S}\xi\in c\rangle,$ and
        \item for all $\eta\in c$, $\Phi^{\eta,\lambda}$ is the direct limit
         extended tree embedding and for $\eta\in b\setminus c$, $\Phi^{\eta,\lambda}
         =\Phi^{\xi, \lambda}\circ \Phi^{\eta,\xi},$ where $\xi=\min c$.
    \end{enumerate}
\end{enumerate}
\end{definition}
For a meta-tree $\mtree{S}$ and a branch $b$ of $\mtree{S}\restrict\gamma$,
 we let 
 \[
 \lim_b(\mtree{S}\restrict\gamma) = 
 \lim\langle \tree{S}_\xi,\Phi^{\eta,\xi}\,|\,\eta\leq_\mtree{S}\xi\in c\rangle,
 \]
  where $c$ is any tail of $b$ where the $\Phi^{\eta,\xi}$ are total
   (if there is such a $c$). We also let
\[
\tree{S}_\gamma^+= V(\tree{S}_\eta,\tree{S}_\xi, F_\xi)\text{ if }
\gamma\text{ is a successor }\xi+1,
\]
where $\eta=\mtree{S}\pred(\xi+1)$, and
\[
\tree{S}_\gamma^+ = \lim_b(\mtree{S}\restrict\gamma)\text{ if }\gamma \text{ is a limit},
\]
where $b=[0,\gamma)_\mtree{S}$.

We also put $\alpha_\xi^\mtree{S}= \alpha_0(F_\xi, \tree{S}_\xi)$ 
and $\beta^\mtree{S}_\xi= \beta(F_\xi,\tree{S}_\xi)$.

A meta-tree $\mtree{S}$ \textit{drops along} $(\eta,\xi]_\mtree{S}$
 iff for some $\gamma\in (\eta,\xi]_\mtree{T}$, either $\gamma$ is a
successor, $\lambda+1$, and for $\zeta=\mtree{S}\pred(\lambda+1)$, $V(\tree{S}_\zeta,\tree{S}_{\lambda},F_{\lambda})$
 is in the dropping case, or $\tree{S}_\gamma\subsetneq \tree{S}_\gamma^+$. We
  call the former way of dropping a \textit{necessary drop} and the 
  latter way of dropping a \textit{gratuitous drop}.  When 
  $\eta=\mtree{S}\pred (\xi+1)$ and $V(\tree{S}_\eta, \tree{S}_\xi, F_\xi)$ is
   in the dropping case we'll say that $[\eta,\xi+1)_\mtree{S}$ 
   \textit{is a necessary drop} or that  $\eta$-to-$\xi+1$
    \textit{is a necessary drop in} $\mtree{S}$. We shall use other natural
     variants of this terminology with their obvious meaning
      (e.g. \lq\lq$\eta$-to-$\xi$ drops in $\mtree{S}$").
A meta-tree $\mtree{S}$ is \textit{simple} iff it has no gratuitous drops.

We write $\Phi^{\mtree{S}}_{\eta,\xi}$ for the tree embedding from
$\itS_\eta$ to $\itS_\xi$ given by $\mtree{S}$, when it exists,
that is, when $\mtree{S}$ does not drop along $(\eta,\xi]_\mtree{S}$.
We write $\Phi^{\mtree{S}}$ for $\Phi^{\mtree{S}}_{0,\xi}$, where
$\lh(\mtree{S})=\xi+1$, if this exists.

Here are a few examples of meta-trees. The first is familiar: normal plus trees can be viewed as simple meta-trees with the same tree-order
  and exit extenders.

\begin{example}\label{inseration1}
Let $\tree{T}$ be a normal plus tree. Let $\tree{T}_\xi=
\tree{T}\upharpoonright \xi+1$ and $F_\xi=E^\tree{T}_\xi$ and
 for $\eta=\tree{T}\pred(\xi+1)$; then
\[
\mtree{T}=\langle \tree{T}_\xi, \Phi_{\xi,\eta}, F_\zeta\,|\,
 \xi,\eta,\zeta+1<\lh (\tree{T}), \, \xi\leq_{T}\eta\rangle
 \]
is a meta-tree with underlying tree structure $\tree{T}$
 (i.e. $\lh (\mtree{T})=\lh (\tree{T})$ and $\leq_\mtree{T}=\leq_\tree{T}$).

This is a meta-tree since for $\eta=\tree{T}\pred(\xi+1)$, we have
 \begin{align*}
     \tree{T}_{\xi+1} & =\tree{T}_\xi\conc\langle E^\tree{T}_{\xi}\rangle\\
     & =V(\tree{T}_\xi,\tree{T}_\eta,F_\xi).
\end{align*}

Notice that we haven't explicitly defined the tree embeddings of $\mtree{T}$. 
As remarked above, they are uniquely determined by the extenders 
$F_\zeta$ and the tree order $\leq_{\tree{T}}$.
One can check that for $\eta=\tree{T}\pred (\xi+1)$, the
 (extended) tree embedding $\Phi^{\eta,{\xi+1}}$ is just the unique extended
  tree embedding with associated $u$-map given by
\begin{equation*}
    u(\zeta) =
    \begin{cases}
      \zeta & \text{if\,\,} \zeta<\eta, \\
      \xi+1        & \text{if\,\,} \zeta=\eta.
    \end{cases}
\end{equation*}
\end{example}

\begin{definition}\label{meattreeonM} Let $M$ be a premouse; then a
{\em meta-iteration tree on $M$} is a meta-iteration tree on $\itS$,
where $\itS$ is the trivial iteration tree of length 1 on $M$.
\end{definition}

Because we have allowed gratuitous dropping, we have another familiar example.

\begin{example}
Let $\tree{S}$ and $\tree{T}$ be normal trees of successor length with
 $\tree{T}$ a normal extension of $\tree{S}$; say $\tree{T}=
 \langle M_\xi, E_\xi,\leq_\tree{T}\rangle$ and 
 $\tree{S}=\tree{T}\restrict\gamma+1$. We obtain a gratuitously dropping meta-tree
  $\mtree{S}$ on $\tree{S}$ with trees
$\tree{S}_\xi = \tree{T}\restrict\gamma+\xi$ and extenders 
$F_\xi= E_{\gamma+\xi}$ (in particular
 $\tree{S}_0= \tree{S}$ and $\tree{S}_{\lh (\tree{T})-1}= \tree{T}$).
  Note that at every step here we are letting $\tree{S}_{\xi+1}= \tree{S}_\xi\conc 
  \langle E_\xi\rangle$, which will often be a gratuitous drop.  
\end{example}

The more important example of meta-trees comes from quasi-normalization.

\begin{example}
Using our notation above for the quasi-normalization of a maximal $M$-stack
 $\langle \tree{T},\tree{U}\rangle$, with $\tree{U}$ normal,
  \[
  \mtree{V}(\tree{T},\tree{U}) =
  \langle \tree{V}_\xi, \Phi^{\xi,\eta},
   F_\zeta\,|\,\xi,\eta,\zeta+1<\lh (\tree{U}), \, \xi\leq_\tree{U} \eta\rangle
   \] is a meta-tree with underlying tree structure $\tree{U}$.
\end{example}

We call a meta-tree \textit{countable} if it has countable length
 and all of its component trees are countable.
 
 Clause (3) of Definition \ref{meta-trees} requires that meta-trees be weakly
 normal, in the sense that the exit extenders $F_\xi$ increase in length,
 and are used to inflate the earliest possible tree. That $\mtree{V}(\tree{T},\tree{U})$ is a meta-tree, then, makes use of the normality of $\tree{U}$. We shall also need
 to consider stacks of meta-trees.

\begin{definition}\label{stack of meta-trees} Let $\itT$ be a plus tree
on $M$.
A \textit{meta-stack on $\itT$} is a sequence
 $s = \langle \mtree{S}^\xi\,|\,\xi<\delta\rangle$ of 
 meta-trees $\mtree{S}^\xi$ of successor lengths, with associated normal trees $\itT_\xi$
 on $M$,  such that $\itT_0 = \itT$, and for all $\xi < \delta$
 \begin{itemize}
 \item[(a)] $\mtree{S}^\xi$ is on $\itT_\xi$
 \item[(b)]  if $\xi+1<\delta$, then $\itT_{\xi+1}$ is 
  the last tree of $\mtree{S}^\xi$, and
  \item[(c)] if $\xi$ is a limit ordinal, then $\itT_\xi$ is the direct limit
  of the $\itT_\alpha$ for $\alpha < \xi$ sufficiently large, under
  the tree embeddings $\Phi^{\mtree{S}_\alpha} \colon \itT_\alpha \to
  \itT_{\alpha+1}$.
  \end{itemize}
  If $\itT$ is the trivial tree of length 1, then we say that $s$
  is a meta-stack on $M$.
\end{definition}

Of course, if $s$ is a meta-stack on $\itT$, and $\itT$ is on $M$,
then $\langle \itT \rangle^\frown s$ is a meta-stack on $M$. We emphasize that meta-trees themselves must be weakly normal:
stacks of meta-trees are \textit{not} meta-trees, in general.

\subsection{Meta-iteration strategies}

Let $\itS$ be a normal tree on $M$, and $\theta \in \OR$.
In the {\em meta-iteration game}
$G(\itS,\theta)$, players {\rm I} and {\rm II}
cooperate to produce a meta-iteration tree
$\mtree{S}=\langle \lh (\mtree{S}), \leq_\mtree{S}, \{\tree{S}_\xi\}_{\xi<\lh \mtree{S}},
 \{F_\xi\}_{\xi+1<\lh \mtree{S}}, \{\Phi^{\eta,\xi}\}_{\eta\leq_\mtree{S} \xi}\rangle$
 on $\itT$. Player {\rm I} plays the extenders $F_\xi$, and decides whether to
 drop gratuitously, and if so, how far. Player
 {\rm II} must play cofinal wellfounded branches at limit ordinals, which are then used
 to extend $\mtree{S}$. If {\rm II} fails to do this, or if some
 $\itS_{\xi+1}$ has an illfounded model, then {\rm I} wins. If {\rm I}
 does not win in the first $\theta$ moves, then {\rm II} wins.
 
 Clearly, if $\lh (\itS) =1$, then $G(\itS,\theta)$ is equivalent to the
 usual normal iteration game $G(M,\theta)$.
 
 For $\lambda > 1$, $G(\itS,\lambda,\theta)$ is the corresponding variant
 of $G(M,\lambda,\theta)$. In a play of $G(\itS,\lambda,\theta)$ in
 which {\rm II} has not lost, the output is a stack of meta-trees
 $\langle \mtree{S}_\xi \mid \xi < \lambda \rangle$, with $\mtree{S}_0$
 being on $\itS$, and each meta-tree $\mtree{S}_\xi$ having length $< \theta$.

\begin{definition}
Let $\tree{S}$ be a normal tree on $M$.
  A \textit{$\theta$-iteration strategy for $\tree{S}$} is a 
  winning strategy for {\rm II} in $G(\itS,\theta)$.
  A \textit{$(\lambda,\theta)$-iteration strategy} is a winning strategy
  for {\rm II} in $G(\itS,\lambda,\theta)$.
   for choosing branches
   in meta-trees on $\tree{S}$ of limit length.
\end{definition}
We shall sometimes call such strategies {\em meta-strategies}.
Clearly, if $\itS$ is the trivial tree of length 1 on $M$, then
a meta-strategy for $\itS$ of whatever type is equivalent to an ordinary
iteration strategy for $M$ of the corresponding type.
So it makes sense to speak of meta-strategies for $M$.

The following is the main theorem on the existence of meta-strategies.
The theorem is essentiatlly a generalization, due to Jensen, of
Schlutzenberg's theorem that normal iterability implies
stack iterability. (See \ref{strategyextension} above.) The main ideas
occur in Schlutzenberg's work.

\begin{theorem}\label{meta-iterability}
Suppose that $M$ is a countable premouse and $\Sigma$ is an $\omega_1+1$ strategy 
for $M$ with strong hull condensation. Then there is a unique
 $(\omega_1,\omega_1)$-meta-strategy $\Sigma^*$
 for $M$ such that for any stack
  $\langle \mtree{S}_\xi \mid \xi < \lambda \rangle$,
  \[
  \langle \mtree{S}_\xi \mid \xi < \lambda \rangle
  \text{ is by }\Sigma^* \Leftrightarrow\text{$\forall \xi < \lambda$ (every tree occurring in $\mtree{S}_\xi$ is by }\Sigma).
  \] 
\end{theorem}

We call $\Sigma^*$ the {\em meta-strategy induced by $\Sigma$}.
Note that the restriction of $\Sigma$ to countable normal trees 
is essentially the same as the restriction of $\Sigma^*$ to
meta-trees on $M$.

It follows from Theorem \ref{meta-iterability}
that, assuming $\adp$, whenever $(M,\Sigma)$ is a mouse pair
with scope $\hc$, then $\Sigma$
induces an $(\omega_1,\omega_1)$-meta-strategy $\Sigma^*$
for $M$. Moreover, whenever
$\itS$ is a plus tree by $\Sigma$, the tail strategy
$\Sigma^*_{\itS}$ is an $(\omega_1,\omega_1)$-meta-strategy
for $\itS$.

We shall outline a proof of Theorem \ref{meta-iterability} now. 
A full proof for the normal tree context appears in
\cite{associativity}.


\begin{lemma}\label{meta-strategy lemma}
Let $M$ is a countable premouse and $\Sigma$ an $\omega_1+1$ 
strategy for $M$ with strong hull condensation. Suppose $\tree{T}$ is a countable
plus tree by $\Sigma$ of successor length, 
 and $\mtree{S}=\langle \tree{S}_\xi,F_\xi,\ldots\rangle$ is a meta-tree on
  $\tree{T}$ of limit length 
  $\leq \omega_1$ such that for all $\xi<\lh (\mtree{S})$,
  $\tree{S}_\xi$ is by $\Sigma$. Let $\delta(\mtree{S})= 
  \sup \{\alpha_\xi+1\,|\,\xi<\lh (\mtree{S})\}$.
Then there is a unique branch $b$ of $\mtree{S}$ such that
 $(\lim _b\mtree{S})\restrict\delta(\mtree{S})+1$ is by $\Sigma$.
  Moreover, for this $b$, $\lim_b \mtree{S}$ is wellfounded, and
   is by $\Sigma$ whenever $\lh (\mtree{S})<\omega_1$.
\end{lemma}

\begin{proof}[Proof sketch.] The proof resembles the analysis of branches of $W(\tree{T},\tree{U})$ when $\tree{U}$ has limit length from section 6.6 of \cite{nitcis}.
Let $\itW$ be the plus tree of
length $\delta = \delta(\mtree{S})$ such that 
for all $\xi < \delta$, $\itW \restrict \xi =
\itS_\eta \restrict \xi$ for all sufficiently large $\eta$.
Let $a =
\Sigma(\itW)$.  As in \S6.6 of \cite{nitcis}, we can decode from $a$
a branch $b$ of $\mtree{S}$\footnote{This can also be found in \cite{jensen} for the corresponding notion of meta-tree.}.
The condensation property of $\Sigma$ implies that
if $\mtree{S}$ is countable, then $\itW_b = \lim_b \mtree{S}$ is by $\Sigma$.

  We are left to show that if $\lh (\mtree{S})=\omega_1$, then all models in $\itW_b$
  are wellfounded. This follows by a simple Skolem hull argument, using again
  the condensation property of $\Sigma$.  \qed
\end{proof}

Notice that if we started with $\Sigma$ a strategy for countable trees and
 $\mtree{S}$ of countable length with the properties in the hypothesis of the 
 lemma, the conclusion still holds.

Lemma \ref{meta-strategy lemma} tells us that, under its hypotheses
on $(M,\Sigma)$, for any plus tree $\itT$ by $\Sigma$, the meta-strategy
$\Sigma^*_{\itT}$ induced by $\Sigma$ is well-defined on meta-trees on $\itT$
of length $\le \omega_1$.

\begin{definition}\label{sigma*0}
Suppose that $M$ is a countable premouse, and $\Sigma$ is an $\omega_1+1$ 
strategy for $M$ with strong hull condensation. Let $\itT$ be a countable tree of successor length
that is by $\Sigma$. Then $\Sigma^{*,0}_{\itT}$ is the
$(\omega_1+1)$-meta-strategy for $\itT$ given by: if
$\mtree{S}$ is a meta-tree on
  $\tree{T}$ of limit length $\leq \omega_1$
by $\Sigma^{*,0}_{\itT}$, then
\[
\Sigma^{*,0}_{\itT}(\mtree{S}) = b \text{ iff }
 \lim _b\mtree{S}\restrict\delta(\mtree{S})+1 \text{ is by $\Sigma$.}
 \]
 \end{definition}
 
 Lemma \ref{meta-strategy lemma} shows that the putative
 $(\omega_1,\omega_1)$-meta-strategy $\Sigma^*$ for $M$ defined in the
 statement of Theorem \ref{meta-iterability} does not break down
 at successor rounds. For if $\itT$ is the last tree of $\mtree{S}^\xi$,
 then part of $\Sigma^*$ that is needed to form $\mtree{S}^{\xi+1}$
 in the next round
 is just $\Sigma^{*,0}_{\itT}$. So what is left is to show that $\Sigma^*$
 does not break down at some limit of rounds.
  
One ingredient here is a theorem on normalizing stacks of meta-trees,
due to Schlutzenberg and Siskind.\footnote{The theorem is a variant
 of Schlutzenberg's theorem on
 the commutativity of inflation from \cite{farmer}. As stated,
 it was discovered later but independently by Siskind.}

\begin{theorem} [Schlutzenberg, Siskind] \label{meta-tree full norm} 
Let $M$ be a countable premouse and $\Sigma$ an $\omega_1$-iteration strategy
for $M$ with strong hull condensation.
Let $\tree{S}$ be a plus tree
by $\Sigma$, and $\langle\mtree{S}, \mtree{T}\rangle$
be a stack of (simple) meta-trees on $\tree{S}$ by $\Sigma^*_\tree{S}$
 with last tree $\tree{U}$; then there is a (simple) meta-tree $\mtree{U}$
  on $\tree{S}$ by $\Sigma^*_\tree{S}$ with last tree
   $\tree{U}$ such that $\Phi^\mtree{U}=\Phi^\mtree{T}\circ\Phi^\mtree{S}$.
Moreover, $\tree{S}$-to-$\tree{U}$ drops in $\mtree{U}$ iff
$\tree{S}$-to-$\tree{U}$ drops in $\langle \mtree{S},\mtree{T}\rangle$.
\end{theorem}

For a proof of the theorem in the context of $\lambda$-tight normal trees, see \cite{associativity}.
The theorem says that $\Sigma^*$ fully normalizes well
 for meta-stacks of length 2, but \cite{associativity} shows that in fact
 it does so for arbitrary countable meta-stacks. It is worth noting that
 the proof is entirely combinatorial; no phalanx comparisons like those
 we shall use later to prove that the ordinary strategies in mouse pairs
 fully normalize well are needed.

\begin{lemma}  \label{meta-tree uniqueness}
Suppose $(M,\Sigma)$ is a mouse pair, $\tree{S}$ is a plus tree
 by $\Sigma$, and $\mtree{S}$, $\mtree{T}$ are simple meta-trees
  on $\tree{S}$ by $\Sigma^*$ with the same last tree $\tree{T}$;
then $\mtree{S}=\mtree{T}$.
\end{lemma}
\begin{proof} The proof is really the same as the proof
 that there is a unique normal tree by $\Sigma$ giving
  rise to any normal $\Sigma$-iterate of $M$.
  
Let $\mtree{S} = \langle \tree{S_\xi}, F_\xi, ...\rangle$ and
 $\mtree{T}=\langle \tree{T_\xi}, G_\xi, \ldots\rangle$
 (so $\tree{T}_0=\tree{S}_0=\tree{S}$ and 
 $\tree{T}_{\infty}=\tree{S}_{\infty}=\tree{T}$).
We'll verify by induction that $\tree{T}_\xi=\tree{S}_\xi$ and $F_\xi=G_\xi$ for all $\xi<\lh(\mtree{S})=\lh\mtree{T})$. 

We've assumed the base case, so now suppose $\tree{S}_\xi=\tree{T}_\xi$ and $F_\eta=G_\eta$ for all $\eta<\xi$. Towards a contradiction, suppose $\lh (F_\xi)< \lh (G_\xi)$. Then 
$\lh(F_\xi)< \lh (G_\eta)$ for all $\eta\geq\xi$. It follows that $F_\xi^-$ is
 on the sequence of last model of $\tree{T}_\eta$ for all $\eta\geq \xi$.
  In particular, $F_\xi^-$ is on the sequence of the last model of $\tree{T}$. 
  But $F_\xi$ is used in $\tree{T}$ since it is used in $\tree{S}_{\eta}$ 
  for \textit{every} $\eta>\xi$. A contradiction.
The same argument shows that we can't have $\lh (G_\xi) < \lh (F_\xi)$, either,
 so $F_\xi=G_\xi$. Since $\mtree{S}$ and $\mtree{T}$ are simple,
 we get $\tree{S}_{\xi+1}=\tree{T}_{\xi+1}$.

$\mtree{S}$ and $\mtree{T}$ cannot diverge at a
 limit stage, because both are by $\Sigma^*$,
and are both simple. The argument in the successor case also shows we cannot have $\lh(\mtree{S})<\lh(\mtree{T})$ or vice-versa, so $\mtree{S}=\mtree{T}$, as desired. \qed
\end{proof}

\begin{remark} If you drop the assumption
that the meta-trees are simple, uniqueness can
 fail. However, the above argument shows that the
  extender sequences still must be the same. In particular,
   we get that the partial tree embeddings from $\tree{S}$ into $\tree{T}$ determined by $\mtree{S}$ and $\mtree{T}$ are the same.
\end{remark}

The remaining ingredient
is a comparison theorem for plus trees. It is  
a variant of a theorem of Schlutzenberg (see \cite{farmer}).
 For $\lambda$-tight normal trees, it is proved in \cite{associativity}. We shall
  generalize this theorem in a later section.

\begin{theorem} [Tree comparison; Schlutzenberg, Siskind]
\label{meta-tree comparison v1}
Suppose that $M$ is a countable premouse, and $\Sigma$ is an
 $\omega_1+1$ strategy for $M$ such that every pseudo-hull of a countable
 tree by $\Sigma$ is by $\Sigma$. 
 Suppose $\{\tree{S}_i\,|\, i\in \omega\}$ is a set of countable plus trees of
  successor lengths which are by $\Sigma$; then there is a countable 
  tree $\tree{T}$ by $\Sigma$ and countable simple meta-trees 
  $\mtree{S}^i$ on $\tree{S}_i$ by $\Sigma^*$,
   each with last tree $\tree{T}$, such that for some $i$,
    the branch $\tree{S}_i$-to-$\tree{T}$ of $\mtree{S}^i$ does not drop.
\end{theorem}

The proof is a straightforward modification of comparison by least extender
disagreement, the disagreements in question being those between the 
extender sequences of the last models of the last trees in our current
approximation to the $\mtree{S}^i$'s. It is important here that
the last trees of those approximations cannot diverge at a limit step,
because they are all by $\Sigma$.

Theorem \ref{meta-tree comparison v1} is no longer true if we allow
the $\itS_i$ to be by different strategies, however nice they are.
For example, let $\Sigma$ and $\Omega$ be such that $(M,\Sigma)$
and $(M,\Omega)$ are mouse pairs, and let $\itS^\frown b$ be by $\Sigma$,
$\itS^\frown c$ by $\Omega$, and $b \neq c$. There can be no tree $\itT$
that is by both $\Sigma$ and $\Omega$ such that both $\itS^\frown b$
and $\itS^\frown c$ are pseudo-hulls of $\itT$, by strong hull condensation.

We can now show that the putative meta-strategy $\Sigma^*$ described
in Theorem \ref{meta-iterability} does not break down at limit rounds.

\begin{lemma}\label{limitcase}
Suppose that $M$ is a countable premouse, and $\Sigma$ is an $\omega_1+1$ 
strategy for $M$ such that every pseudo-hull of a countable normal
tree by $\Sigma$ is also by $\Sigma$. Let $\lambda < \omega_1$ be a limit
ordinal, and
$\mtree{S} = \langle \mtree{S}^\xi\,|\,\xi< \lambda\rangle$ be a meta-stack on $M$ such
that for all $\xi < \lambda$,  $\mtree{S}^\xi$ is countable, and all
trees occurring in $\mtree{S}^\xi$ are by $\Sigma$. Let
\[
\itT_\xi = \text{ first tree in } \mtree{S}^\xi.
\]
Then
\begin{itemize}
\item[(a)] for all sufficiently large $\xi<\eta < \lambda$,
$\Phi^{\mtree{S}}_{\xi,\eta}$ exists, and
\item[(b)] letting $\itT$ be the direct limit of the
$\itT_\xi$, for $\xi < \lambda$ sufficiently large,
under the $\Phi^{\mtree{S}}_{\xi,\eta}$, we have
that $\itT$ is by $\Sigma$.
\end{itemize}
\end{lemma}
\begin{proof}
Applying the tree comparison theorem (Theorem 
\ref{meta-tree comparison v1}) to $\{\tree{T}_\xi\,|\,\xi<\lambda\}$, we 
get a countable tree $\itV$ which is by $\Sigma$, and for each
 $\xi<\lambda$ a simple countable meta-tree $\mtree{T}^\xi$ on $\tree{T}_\xi$ by
  $\Sigma^*$ with last tree $\itV$. Moreover, for some 
  $\xi<\lambda$, $\mtree{T}^\xi$ doesn't drop along
   $\tree{T}_\xi$-to-$ \itV$.
   For every $\eta<\lambda$,  
 $\langle \mtree{S}^\eta, \mtree{T}^{\eta+1}\rangle$ is a meta-stack on $M$
 with last tree $\itV$.
 Applying Theorem \ref{meta-tree 
 full norm} to this stack, we get a meta-tree 
 \[
 \mtree{U}^\eta = \mtree{V}(\mtree{S}^\eta, \mtree{T}_{\eta+1})
 \]
 by $\Sigma^*$ on
  $\tree{T}_\eta$ with last tree $\itV$. 
By the remark following Lemma \ref{meta-tree uniqueness}, $\mtree{U}^\eta$ 
and $\mtree{T}^\eta$ 
determine the same partial tree embedding 
from $\tree{T}_\eta\to \itV$.
 (If $\mtree{S}^\eta$ were simple, then we would have $\mtree{U}=\mtree{T}^\eta$.) 
Further, if $\mtree{T}^\eta$ does not drop along $\tree{T}_\eta$-to-$\itV$, 
then $\Phi^{\mtree{T}^\eta}$ is total. 
Since $\Phi^{\mtree{T}^\eta} = \Phi^{\mtree{U}^\eta} = 
\Phi^{\mtree{T}^{\eta+1}}\circ \Phi^{\mtree{S}^\eta}$, 
it follows that $\Phi^{\mtree{S}^\eta}$ and $\Phi^{\mtree{S}^\eta}$ are total as well.

So, we get that $\Phi^{\mtree{T}^\xi}=\Phi^{\eta,\xi}\circ\Phi^{\mtree{T}^\eta}$ 
for all $\eta\leq\xi<\alpha$ and for all sufficiently large $\eta,\xi$, these are 
total extended tree embeddings.

Fixing $\zeta$ above which these are total, we have
 $\lim \langle \tree{S}_\xi, \Phi^{\eta,\xi}\,|\, \zeta\leq\eta\leq\xi<\alpha\rangle$
  exists, and a pseudo-hull of $\itV$.  (See Proposition \ref{direct limit prop}.)
   Since $\itV$ is by $\Sigma$, so is 
    $\lim \langle \tree{S}_\xi, \Phi^{\eta,\xi}\,|\, \zeta\leq\eta\leq\xi<\alpha\rangle$.
    \qed
\end{proof}

Lemmas \ref{meta-strategy lemma} and \ref{limitcase} clearly combine to give our variant of
Schlutzenberg's meta-iterability theorem, Theorem \ref{meta-iterability}.

\begin{remark} We shall show in Theorem \ref{induced strategy theorem} that, assuming $\adp$,
 if $\Lambda$ is a meta-strategy for
a normal tree $\itT$ on $M$, and $\Lambda$ has a certain Dodd-Jensen property,
then $\Lambda$ is the meta-strategy induced by some ordinary strategy for $M$.
So meta-trees and strategies are not something fundamentally new. They are rather
a useful way of organizing the construction of ordinary iteration trees
by an ordinary iteration strategy.\end{remark}

As an immediate corollary to \ref{meta-iterability}, we get Schlutzenberg's
 part of Theorem \ref{strategyextension}.

\begin{theorem}[Schlutzenberg]
\label{stack iterability ii} Let $M$ be a premouse 
and $\Sigma$ an $\omega_1+1$ strategy for $M$ such that every
pseudo-hull of a countable plus tree by $\Sigma$ is by $\Sigma$;
then there is a unique extension of $\Sigma$ 
to a strategy for countable stacks of countable trees which quasi-normalizes well.
\end{theorem}

\begin{proof} Let $\Sigma^*$ be the $(\omega_1,\omega_1)$-meta-strategy
for $M$ induced by $\Sigma$. We define an ordinary $(\omega_1,\omega_1)$-strategy
for $M$, which we call $\Omega$, as follows. We associate inductively to any
stack
\[
t = \langle \itT_\xi \mid \xi < \alpha \rangle
\]
by $\Omega$ a meta-stack
\[
\mtree{S} = \langle \mtree{S}^\xi \mid \xi < \alpha \rangle
\]
that is by $\Sigma^*$. For $\xi < \alpha$, let $P_\xi$ be the base model
of $\itT_\xi$. Let $\mtree{S}^0$ be $\itT_0$, considered as a meta-tree on
$P_0 = M$. For $\xi > 0$, let
\[
\itU_{\xi} = \text{ last tree in $\mtree{S}^\xi$},
\]
and for $\lambda$ a limit,
\[
\itU^0_\lambda = \text{ last tree in $\lim_{\xi<\lambda}\mtree{S}^\xi$},
\]
where the limit is taken under the tree embeddings of $\mtree{S}$. So $\itU_0
= \itT_0$. Let
$Q_0=P_0$, and
\[
Q_{\xi+1} = \text{last model of $\itU_\xi$,}
\]
 and
if $\lambda$ is a limit ordinal,
\[
Q_\lambda = \text{last model of $\itU^0_\lambda$}.
\]
So $Q_1= P_1$.
We maintain inductively that we have nearly elementary maps
\[
\sigma_\xi \colon P_\xi \to Q_\xi
\]
that commute with the maps of $t$ and $\mtree{S}$. (That is, if
$i \colon P_\xi \to P_\eta$ is a map of $t$, then $\Phi^{\mtree{S}}_{\xi,\eta}$
yields a $t$-map $j \colon Q_\xi \to Q_\eta$, and $\sigma_\eta \circ i =
j \circ \sigma_\xi$.) $\sigma_0 = \sigma_1 = \mbox{ id}$.

Let $\xi > 0$, and suppose that the stack $t \restriction \xi$ is by $\Omega$,
 with associated meta-stack
$\mtree{S} \restriction \xi$ has been constructed with the properties above. 

If $\xi$ is a limit ordinal, then because $\mtree{S} \restriction \xi$ is by
$\Sigma^*$, $\itU^0_\xi$ is by $\Sigma$, and $Q_\xi$ exists. Moreover
we have $\sigma_\xi \colon P_\xi \to Q_\xi$ by our commutativity hypothesis.
We define the tail strategy $\Omega_{t \restriction \xi, P_\xi}$ for normal trees
by
\[
 \itV \text{ is by }\Omega_{t \restriction \xi, P_\xi} \Leftrightarrow
\mtree{V}(\itU^0_\xi,
 \sigma_\xi\itV) \text{ is by }\Sigma^*.
 \]
 Equivalently, $\Omega$ plays round $\xi$ by pulling back under $\sigma_\xi$
 the strategy $\Psi_{\itU^0_\xi,Q_\xi}$, where $\Psi$ is the unique extension
 of $\Sigma$ to stacks of length two.
 
 Playing round $\xi$ this way, $\Omega$ does not lose, and we have
 $\itT_\xi$ as the result of our play. Set 
 \[
 \mtree{S}^\xi = \mtree{V}(\itU^0_\xi,\sigma_\xi \itT_\xi).
 \]
 $\mtree{S}^\xi$ is a meta-tree on $\itU^0_\xi$, and $\mtree{S} \restriction \xi+1$
 is by $\Sigma^*$. Moreover, the last tree in $\mtree{S}^\xi$ is
  $\itU_{\xi+1} = V(\itU^0_\xi,\sigma_\xi\itT_\xi)$, and so we 
  have a quasi-normalization map
  $\tau$ from the last model $R$ of $\sigma_\xi\itT_\xi$ to $Q_{\xi+1}$. We also
  have a copy map $\pi$ from $P_{\xi+1}$, which is the last model of $\itT_\xi$,
  to $R$. We then set $\sigma_{\xi+1} = \tau \circ \pi$. Our induction
  hypotheses still hold.
  
  The case $\xi = \gamma+1$ is completely parallel. We have $t(\gamma)$ with
  last model $P_{\gamma+1}$,
  $\mtree{S}^\gamma$ with last tree
  $\itU_\gamma$, whose last model is $Q_{\gamma+1}$, and
  $\sigma_{\gamma+1} \colon P_{\gamma +1} \to Q_{\gamma+1}$ already.
  We put
  \[
 \itV \text{ is by }\Omega_{t \restriction \gamma +1, P_{\gamma+1}} \Leftrightarrow
\mtree{V}(\itU_\gamma,
 \sigma_{\gamma+1}\itV) \text{ is by }\Sigma^*,
 \]
and
\[
\mtree{S}^{\gamma+1} = \mtree{V}(\itU_\gamma,\sigma_{\gamma+1}\itT_{\gamma+1}),
\]
and so on. \qed
\end{proof}

\subsection{Copying meta-trees}

    Some of the basic facts about meta-trees can be understood
    through the analogy:
 
\begin{align*}
    \text{premice} &\leftrightsquigarrow \text{ iteration trees}\\
\text{iteration trees} &\leftrightsquigarrow \text{meta-trees}\\
\text{strategies} & \leftrightsquigarrow \text{meta-strategies}\\
\text{elementary embeddings} &\leftrightsquigarrow \text{tree embeddings}
\end{align*}
Notice that
 one might have found the tree comparison theorem (Theorem \ref{meta-tree comparison v1})
this way: the theorem says any two iteration trees have a common
   meta-iterate, modulo the condition that we started with trees by the same strategy.
   Except for this (necessary) restriction, this is an analogue of the usual 
    comparison theorem for premice.
In this section we prove a copying theorem that conforms to this analogy. 

The usual copying construction 
is about lifting an iteration tree via an elementary embedding. Following 
the analogy, our result will be about lifting a meta-tree via a tree embedding.

The copying construction results from repeated applications of one-step
copying, as codified in the Shift Lemma. The Shift Lemma for meta-iteration
says that, under the right conditions,
 we can complete the 
 following diagram.

\begin{center}
    \begin{tikzcd} \tree{S} \arrow[r, "\Psi"] 
    & \tree{T} \\
    F \arrow[r, "\mapsto", phantom] \arrow[u, "\scriptsize \inup" ,phantom]& G\arrow[u, "\scriptsize\inup" ,phantom]\\
    \itU \arrow[r, "\Pi"] \arrow[d, "F"]
    & \itV \arrow[d,"G"]\\
    V(\tree{S},\itT,F) \arrow[r, "\Gamma", dashed] 
    & V(\itU,\itV,G)
    \end{tikzcd}
\end{center}

We need enough agreement between $\Psi$ and $\Pi$ that we can indeed
complete the diagram. In the case of ordinary premice, the lifting maps
must agree on $\dom(F)$. The requirement here is similar, but more
complicated.
     
Before we state and prove the Shift Lemma, we isolate a preliminary fact which is essential to the proof. This fact captures a way in which quasi-normalization is better behaved than embedding normalization: if one replaces the $\alpha_0$'s by $\alpha$'s in its statement, it becomes false (in general). In \cite{associativity}, a variant of the Shift Lemma is proved where one uses embedding normalization instead of quasi-normalization, but it is considerably more complicated than our present version because of possible failures of the next lemma for $\alpha$'s.

\begin{lemma}\label{key lemma 1}
Let $\Psi:\tree{S}\to \tree{T}$ be an extended tree embeddings and $F$ an extender such that $F^-$ is on the $M_\infty^\tree{S}$-sequence and $M_\infty^\tree{S}|\lh(F)\is\dom( t_\infty^\Psi)$. Let $G=t_\infty^\Psi(F)$, $\alpha_0=\alpha_0(\tree{S}, F)$, $\alpha_0^*=\alpha_0(\tree{T}, G)$, $\beta=\beta(\tree{S}, F)$, and $\beta^*=\beta(\tree{T},G)$. Then \begin{enumerate}
    \item[(a)] $\alpha_0^*\in [v^\Psi(\alpha_0), u^\Psi(\alpha_0)]_\tree{T}$ and
    \item[(b)] $s^\Psi_{\alpha_0, \alpha_0^*}\restrict\lh(F)+1=t_{\alpha_0}^\Psi\restrict\lh(F)+1=t^\Psi_\infty\restrict\lh(F)+1$.
    \item[(c)] $\beta^*\in [v^\Psi(\beta), u^\Psi(\beta)]_\tree{T}$, and
    \item[(d)] $s^\Psi_{\beta, \beta^*}\restrict\dom(F)\cup\{\dom(F)\}=s^\Psi_{\alpha_0, \alpha^*_0}\restrict\dom(F)\cup\{\dom(F)\}$.
\end{enumerate}
\end{lemma}
\begin{proof}
We first work towards verifying (a) and (b). Let $\alpha=\alpha(\tree{S}, F)$ and $\alpha^*=\alpha(\tree{T}, G)$. Since $[\alpha, \alpha_0)$ is contained in delay interval which begins at $\alpha$ (cf. \cite{nitcis} Lemma 6.7.2), $\alpha_0=\alpha+n$ for some $n\in \omega$ and, by the definition of $\alpha_0$, for every $i<n$, $E_{\alpha+i}^\tree{S}$ is not of plus type and $\hat\lambda(E_{\alpha+i}^\tree{S})<\lh(F)$.

The agreement hypotheses between the maps of a tree embedding give that $t_{\alpha+n}^\Psi\restrict\lh(F)+1=t_\infty^\Psi\restrict\lh(F)+1$. In particular, $t_{\alpha+n}^\Psi(F)=G$. So let $j\leq n$ least such that $t_{\alpha+j}^\Psi(F)=G$.

\paragraph{Case 1.} $j<n$. \\

First we'll check that for any $k<n-j$,  \[E_{u(\alpha+j+k)}^\tree{T}=t_{\alpha+j+k}^\Psi(E_{\alpha+j+k}^\tree{S}) \text{ and if $k\geq 1$, then } v^\Psi(\alpha+j+k)=u^\Psi(\alpha+j+k).\] Note that since $E_{\alpha+j+k}^\tree{S}$ is not of plus type, this implies that $E_{u^\Psi(\alpha+j+k)}^\tree{T}$ is not of plus type either. 

For any $i<n$, since $\hat\lambda(E_{\alpha+i}^\tree{S})<\lh(F)$, $u^\Psi(\alpha+i)$ is the least $\xi\in[v^\Psi(\alpha+i), u^\Psi(\alpha+i)]_\tree{T}$ such that $t^\Psi_{\alpha+i}(F)^-$ is on the $M_\xi^\tree{T}$-sequence. 
(This is because all critical points used along $[v^\Psi(\alpha+i), u^\Psi(\alpha+i)]_\tree{T}$ must be below the current image of $\lh(F)$, since they are below the current image of $\lh(E_{\alpha+i}^\tree{S})$, since we are blowing up this extender along this partial branch.) 
For any $i<n$, since $\Psi$ is an extended tree embedding,  $v^\Psi(\alpha+i+1)=u^\Psi(\alpha+i)+1$, and $E_{\alpha+i}^\tree{S}$ is not of plus type, either $E_{u(\alpha+i)}^\tree{T}=t_{u(\alpha+i)}^\Psi(E_{\alpha+i}^\tree{S})$ (and so is \textit{not} of plus type) and $s_{\alpha+i+1}^\Psi\restrict\lh(E_{\alpha+i}^\tree{S})+1=t_{\alpha+i}^\Psi\restrict\lh(E_{\alpha+i}^\tree{S})+1$ or else $E_{u(\alpha+i)}^\tree{T}=t_{u(\alpha+i)}^\Psi(E_{\alpha+i}^\tree{S})^+$ (and so is of plus type), $s_{\alpha+i+1}^\Psi\restrict\hat\lambda(E_{\alpha+i}^\tree{S})=t_{\alpha+i}^\Psi\restrict\hat\lambda(E_{\alpha+i}^\tree{S})$, but $s_{\alpha+i+1}^\Psi(\hat\lambda(E_{\alpha+i}^\tree{S}))>t_{\alpha+i}^\Psi(\hat\lambda(E_{\alpha+i}^\tree{S}))$. 

It follows that for any $i<n$, $\lh(t_{\alpha+i}^\Psi(F))\leq \lh( t_{\alpha+i+1}^\Psi(F))$, with equality obtaining only when $E_{u(\alpha+i)}^\tree{T}=t_{u(\alpha+i)}^\Psi(E_{\alpha+i}^\tree{S})$ (which is \textit{not} of plus type) and if $i<n-1$, $v^\Psi(\alpha+i+1)=u^\Psi(\alpha+i+1)$.\footnote{In the case $i=n-1$, we may still have equality obtaining when $E_{u(\alpha+i)}^\tree{T}=t_{u(\alpha+i)}^\Psi(E_{\alpha+i}^\tree{S})$ but $v^\Psi(\alpha+i+1)<u^\Psi(\alpha+i+1)$.} Since we reach our final image $G$ of $F$ already by $u^\Psi(\alpha+j)$, we get that for all $k$ such that $k<n-j$, $E_{u(\alpha+j+k)}^\tree{T}=t_{\alpha+j+k}^\Psi(E_{\alpha+j+k}^\tree{S})$ and if $k\geq 1$, then $v^\Psi(\alpha+j+k)=u^\Psi(\alpha+j+k)$, as desired.

In particular, $u^\Psi(\alpha+j+k)=u^\Psi(\alpha+j)+k$ for all $k< n-j$. This will let us verify that $\{u^\Psi(\alpha+j)+k\mid k<n-j\}$ is contained in a delay interval of $\tree{T}$ and that for all $k<n-j$, $\hat\lambda(E_{u^\Psi(\alpha+j)+k}^\tree{T})<\lh(G)$. Since $\{\alpha+j+k\mid k<n-j\}$ is contained in a single delay interval in $\tree{S}$ and $\hat\lambda(E_{\alpha+j+k})<\lh(F)$, this is immediate from the agreement of the $t$-maps using that $E_{u(\alpha+j)+k}^\tree{T}=t_{\alpha+j+k}^\Psi(E_{\alpha+j+k}^\tree{S})$.

We have $\alpha^*\leq u^\Psi(\alpha+j)$ by our choice of $j$. The considerations above showed that $u^\Psi(\alpha+j)$ is the least $\xi\in [v^\Psi(\alpha+j), u^\Psi(\alpha+j)]_\tree{T}$ such that $s_{\alpha+j, \xi}^\Psi(F)=G$. It is easy to see that $u^\Psi(\alpha+j)$ is also the least $\xi\in [v^\Psi(\alpha+j), u^\Psi(\alpha+j)]_\tree{T}$ such that $G^-$ is on the $M_\xi^\tree{T}$-sequence and also that $G^-$ is not on the sequence of any $M^\tree{T}_{u^\Psi(\alpha+i)}$ for $i<j$. It follows that $\alpha^*\geq v^\Psi(\alpha+j)$ (as this is trivial if $j=0$). 

\paragraph{Subcase A.} $\alpha^*=u^\Psi(\alpha+j)$.\\

In this case, our preceding observations give that  $\alpha_0^*\geq \sup\{u^\Psi(\alpha+j)+k\mid k<n-j\}=v^\Psi(\alpha+n)$. Since $G$ isn't moved along $[v^\Psi(\alpha+n), u^\Psi(\alpha+n)]_\tree{T}$ we must have that either $v^\Psi(\alpha+n)=u^\Psi(\alpha+n)$ or else $\hat\lambda(E^\tree{T}_{v^\Psi(\alpha+n)}>\lh(G)$ ($\hat\lambda(E^\tree{T}_{v^\Psi(\alpha+n)}<\lh(G)$ would imply $G$ is blown up along $[v^\Psi(\alpha+n), u^\Psi(\alpha+n)]_\tree{T}$, as before, contradicting that $t_{\alpha+n}^\Psi(F)=G$). In the latter case (i.e. that $\hat\lambda(E^\tree{T}_{v^\Psi(\alpha+n)}>\lh(G)$), we get $\alpha^*_0=v^\Psi(\alpha_0)$, giving (a). (b) easily follows so suppose $v^\Psi(\alpha+n)=u^\Psi(\alpha+n)$. Since $\alpha+n=\alpha_0$, we have that $\alpha+n=\lh(\tree{S})$ or else either $\hat\lambda(E_{\alpha+n}^\tree{S})>\lh(F)$ or $E_{\alpha+n}^\tree{S}$ is of plus type. If  $\alpha+n+1=\lh(\tree{S})$ then since $\Psi$ is an extended tree embedding we must have  $v^\Psi(\alpha+n)+1=u^\Psi(\alpha+n)+1=\lh(\tree{T})$, so that $\alpha_0^*=v^\Psi(\alpha+n)=u^\Psi(\alpha+n)$ which immediately gives (a) and easily gives (b). In the remaining cases, $E_{u^\Psi(\alpha+n)}^\tree{T}$ must either be of plus type (either $t_{\alpha+n}^\Psi(E_{\alpha+n}^\tree{S})$ if $E_{\alpha+n}^\tree{S}$ is of plus type or possibly $t_{\alpha+n}^\Psi(E_{\alpha+n}^\tree{S})^+$ even if it isn't) or have $\hat\lambda(E_{u^\Psi(\alpha+n)}^\tree{T})>\lh(G)$. So in these cases, too, we get $\alpha_0^*=v^\Psi(\alpha+n)=u^\Psi(\alpha+n)$, giving (a) and (b).

\paragraph{Subcase B.} $\alpha^*<u^\Psi(\alpha+j)$.\\

We can run the argument of Subcase A once we know that every $\xi\in [\alpha^*, u^\Psi(\alpha+j)]$ is not of plus type and has $\hat\lambda(E_\xi^\tree{T})<\lh(G)$. But we already verified that $E_{u(\alpha+j)}^\tree{T}$ cannot be of plus type (or else we'd violate quasi-normality at the next extender) and has $\hat\lambda( E_{u(\alpha+j)}^\tree{T})<\lh(G)$. Since $\tree{T}$ is quasi-normal, every $\xi\in [\alpha^*, u^\Psi(\alpha+j)]$ must have $\hat\lambda(E_\xi^\tree{T})\leq \hat\lambda( E_{u(\alpha+j)}^\tree{T})<\lh(G)$, too. Moreover, as $G^-$ is on the sequence of $M_{\alpha^*}^\tree{T}$ and $M_\infty^\tree{T}$, all of these extender must have length $>\lh(G)$. But then if any of these $E_\xi^\tree{T}$ is of plus type, the delay interval which begins at $\alpha^*$ would have to end strictly before $E_{u^\Psi(\alpha+j)}^\tree{T}$, forcing $\hat\lambda(E_{u^\Psi(\alpha+j)}^\tree{T})>\lh(G)$, a contradiction.

This finishes Case 1.

\paragraph{Case 2.} $j=n$.\\

Recall that this means that for all $i<n$, $t^\Psi_{\alpha+i}(F)\neq G$. Let $\xi\in [v^\Psi(\alpha_0), u^\Psi(\alpha_0)]_\tree{T}$ be least such that $s_{\alpha_0, \xi}^\Psi(F)= G$. 

First suppose $\xi=v^\Psi(\alpha_0)$. Then we must have $n=0$, i.e. $\alpha_0=\alpha$, or else we'd have $t^\Psi_{\alpha+i}(F)= G$ for some $i<n$ (by considerations similar to those at the start of Case 1). It follows that $\alpha^*=\xi$. But then, since $G$ is in the range of $t_{\alpha_0}^\Psi$, either  $u^\Psi(\alpha_0)=v^\Psi(\alpha_0)=\xi$, in which case the argument at the end of Case 1 Subcase A gives $\alpha_0^*=\xi$, or else $\hat\lambda(E_\xi^\tree{T})>\lh(G)$, in which case, again, $\alpha_0^*=\xi$. So if $\xi=v^\Psi(\alpha_0)$, we're done. So for the remainer of the proof, suppose $\xi>v^\Psi(\alpha_0)$.

\paragraph{Subcase A.} $\xi$ is a limit ordinal.\\

In this case, we can again show that $\xi=\alpha^*=\alpha_0^*$. Towards a contradiction, suppose $\alpha^*<\xi$. Then the delay interval beginning at $\alpha^*$ must end strictly before $\xi$, so that quasi-normality guarantees that for some $\eta+1\in [v^\Psi(\alpha_0), \xi)_\tree{T}$, $\hat\lambda(E_\eta^\tree{T})>\lh(G)$. But then we must use an extender with critical point $>\lh(G)$ along $[v^\Psi(\alpha_0), \xi)_\tree{T}$, contradicting our choice of $\xi$. 

So $\xi=\alpha^*$. If $\xi=u^\Psi(\alpha_0)$, then we get $\xi=\alpha_0^*$ by the argument at the end of Case 1 Subcase A. If $\xi<u^\Psi(\alpha_0)$, we must have $\hat\lambda(E_\xi^\tree{T})>\lh(G)$, so that, again, $\xi=\alpha_0^*$. 

\paragraph{Subcase B.} $\xi$ is a successor ordinal.\\

Let $\xi=\gamma+1$ and $\eta=\tree{T}\pred(\gamma+1)$. Since $\xi>v^\Psi(\alpha_0)$, $\eta\geq v^\Psi(\alpha_0)$. Since we haven't reached the final image $G$ of $F$ along $[v^\Psi(\alpha_0), u^\Psi(\alpha_0)]_\tree{T}$ by $\eta$, we must have $\crit(E_\gamma^\tree{T})<\lh(G)$ which implies $\hat\lambda(E_\gamma^\tree{T})<\lh(G)$. It follows that $\gamma<\alpha_0^*$, i.e. $\gamma+1\leq\alpha_0^*$. If $\gamma+1=u^\Psi(\alpha_0)$, one can show $\gamma+1=\alpha_0^*$ (by the Case 1 Subcase A argument, again). So suppose $\gamma+1=u^\Psi(\alpha_0)$. In this case we must again have $\gamma+1=\alpha_0^*$ because $\gamma+1<\alpha_0^*$ is impossible, since it implies that we must next use an extender with critical point $<\lh(G)$, contradicting that we've already reached $G$ as the final image of $F$ at $\gamma+1$. 

This finishes our verification of (a) and (b).\\

Verifying (c) and (d) is actually much easier and essentially appears in \cite{nitcis}. We just handle the case that $\beta+1<\lh(\tree{S})$ (the remaining case being that $\beta+1=\lh(\tree{S})$, so that actually $\alpha_0=\alpha=\beta$; the argument in this case is basically the same). By definition, $\beta$ is the least $\xi$ such that $\dom(F)\isneq M_\xi^\tree{S}|\hat\lambda(E_\xi^\tree{S})$. 
It is easy to see $\beta\leq \alpha_0$ and that the agreement between the maps of $\Psi$ gives $t_\beta^\Psi\restrict\dom(F)+1=t_{\alpha_0}\restrict\dom(F)+1$. 
It follows that $t_\beta^\Psi(\dom(F))=\dom(G)$. Now let $\gamma\in [v^\Psi(\beta), u^\Psi(\beta)]_\tree{T}$ least such that $s^\Psi_{\beta, \gamma}(\dom(F))=\dom(G)$. We'll show that $\gamma=\beta^*$. 
We must have that $\dom(G)\isneq M_\gamma^\tree{T}|\hat\lambda(E_\gamma^\tree{T})$, since if this were to fail then we must have $\gamma<u^\Psi(\beta)$ and that the extender applied to $M_\gamma^\tree{T}$ along $[v^\Psi(\beta), u^\Psi(\beta)]_\tree{T}$ has critical point $\leq\crit(G)$, 
contradicting that we've reached $\dom(G)$ as the image of $\dom(F)$ at $\gamma$. So $\beta^*\leq \gamma$. It's also easy to verify that $\beta^*\geq v^\Psi(\beta)$ and $\gamma$ is the least $\xi\in[v^\Psi(\beta), u^\Psi(\beta)]_\tree{T}$ such that $\dom(G)\isneq M_\xi^\tree{T}|\hat\lambda(E_\xi^\tree{T})$. So if $\beta^*\neq \gamma$, then $\beta^*<\gamma$ is a successor $\xi+1$ for $\xi\geq\beta^*$, and $\tree{T}\pred(\xi+1)=\eta\geq v^\Psi(\beta)$. Since $\eta<\gamma$, we must have $\crit(E_\xi^\tree{T})<\crit(s_{\beta, \eta}^\Psi(F)$ (i.e. we are still moving up $\dom(F)$, and so actually $\crit(F)$, along $[v^\Psi(\beta), u^\Psi(\beta)]_\tree{T}$. But then that $\hat\lambda(E_\xi^\tree{T})<\crit(s_{\beta, \gamma}(F))=\crit(G)$, contradicting that $\xi\geq \beta^*$, as $\dom(G)\isneq M_\xi^\tree{T}|\hat\lambda(E_\xi^\tree{T})$ is impossible. So $\gamma=\beta^*$, giving (c) and (d).
\qed 
\end{proof}

We now isolate the agreement hypotheses needed for the Shift Lemma in a definition.

\begin{definition}\label{metashiftapplies}
Let $\Psi:\tree{S}\to \tree{T}$ and $\Pi:\itU \to \itV$ be extended tree embeddings, $F$ an extender such that $F^-$ be
 an extender on the $M_\infty^\tree{S}$-sequence, and $G$ an extender such that $G^-$ is on the $M_\infty^\tree{T}$-sequence. We say that 
 \textit{the Shift Lemma applies
  to $(\Psi,\Pi, F, G)$} iff letting $\beta = \beta(\itS,F)$ and $\beta^*=\beta(\itT, G)$, 
  \begin{enumerate}
  \item $M_\infty^\tree{S}|\lh(F)\is\dom( t_\infty^\Psi)$ and $G=t_\infty^\Psi(F)$,
  \item  $\Psi\restrict\beta+1\approx \Pi\restrict\beta+1$,
\item $\tree{T}\restrict \beta^*+1=\tree{V}\restrict\beta^*+1$
\item $\beta^*\in [v^\Pi(\beta), u^\Pi(\beta)]_\tree{V}$ and $t_{\beta}^\Pi\restrict\dom(F)\cup\{\dom(F)\}=s_{\beta, \beta^*}^\Pi\restrict\dom(F)\cup\{\dom(F)\}$,
\item
   if  $\beta+1<\lh(\tree{U})$, then 
  $\dom(F) \isneq M_\beta^\tree{U}|\lh(E^\itU_\beta)$, and
  \item
   if  $\beta^*+1<\lh(\tree{V})$,
  $\dom(G) \isneq M_{\beta^*}^\tree{U}|\lh(E^\itV_{\beta^*})$.
   \end{enumerate}
   \end{definition}
 
 Note that in clause (4) we do not need to say that $\beta^*\in[v^\Psi(\beta), u^\Psi(\beta)]_\tree{T}$ (though this is needed for $s^\Psi_{\beta,\beta^*}$ to makes sense), since this is a consequence of Lemma \ref{key lemma 1}.

   In the ordinary Shift Lemma, the upstairs models must agree up to the
   common image of $\dom(F)$. Clauses (2)-(4) play an analgous role: they say that the upstairs
   trees and embeddings agree up to the place in $\itT$ where $t_\infty^\Psi(\dom(F))$ has been created.  (5) and (6)  additionally ensure that $V(\tree{U}, \tree{S}, F)$ and $V(\tree{V}, \tree{T}, G)$ are defined.

\begin{lemma}[Shift Lemma]\label{Shift Lemma}
Let $\Psi:\tree{S}\to \tree{T}$ and $\Pi:\itU \to \itV$ be extended tree embeddings,
and let $F$ be an extender such that $F^-$ be
 an extender  on the sequence of the last model of $\tree{S}$ and $G$ be an extender such that $G^-$ is on the extender sequence of the last model of $\tree{T}$. Let $\alpha_0=\alpha_0(\tree{S}, F)$ and
 $\alpha^*_0=\alpha_0(\tree{T},G)$.
 
 Suppose that the Shift Lemma applies
  to $(\Psi,\Pi, F, G)$. Then $V(\tree{U},\tree{S},F)$ and $V(\tree{V},\tree{T},G)$ are defined and, letting $\mu$ the greatest ordinal such that 
    $V(\itU,\itS,F)\restrict \mu$ is wellfounded and $\mu^*$ the greatest ordinal
 such that $V(\itV,\tree{T}, G)\restrict\mu^*$ is wellfounded, there is a unique partial tree embedding 
$\Gamma: V(\itU,\itS,F)\restrict\mu \to V(\itV,\tree{T}, G)\restrict\mu^*$ with maximal domain such that
 \begin{enumerate}
        \item $\Gamma\restrict \alpha_0+1\approx \Psi \restrict \alpha_0+1$,
        \item $u^\Gamma(\alpha_0)=\alpha^*_0$, and
        \item $\Gamma\circ \Phi^{V(\itU,\itS,F)} =\Phi^{V(\itV,\tree{T}, G)}\circ \Pi$ (on their common domain).
    \end{enumerate}
 Moreover, if $V(\itV,\tree{T}, G)$ is wellfounded, then
 $V(\itU,\itS, F)$ is wellfounded 
  and $\Gamma$ is a (total) extended tree embedding 
  from $V(\itU, \itS, F)$ into $V(\itV,\itT, G)$. If $V(\itV,\tree{T}, G)$ is wellfounded and also $\Pi$ is non-dropping, then $\Gamma$ is a non-dropping extended tree embedding.
\end{lemma}

\begin{remark}\label{shift remark}
If we assume that $\tree{S},\itT,\itU,$ and $\itV$ are all
 by some strategy $\Sigma$ for $M$ with strong hull condensation that quasi-normalizes well, then
  $V(\itU,\tree{S}, F)$ and $V(\itV,\tree{T}, G)$ are by $\Sigma$, and
$\Gamma$ is a total extended tree embedding.
\end{remark}

\begin{proof} We have $V(\tree{U}, {\tree{S}}, F)$ is defined by hypotheses (2) and (5) in the definition of when the Shift Lemma applies. $V(\tree{V},\tree{T},G)$ is defined by hypotheses (3) and (6). So all of the work is in identifying $\Gamma$ and proving it is as desired, inductively. At bottom, we are able to do this because the $s$-maps of tree embeddings are given by the ordinary Shift Lemma at successors. Let $\tree{W}=V(\tree{U}, {\tree{S}}, F)\restrict\mu+1$, $\tree{W}^*=V(\tree{V},\tree{T},G)\restrict\mu^*+1$, $\Phi=\Phi^{V(\tree{U}, {\tree{S}}, F)}$, and $\Phi^*=\Phi^{V(\tree{V},\tree{T},G)}$.

Notice that $\mu\geq\alpha_0+1$ and $\mu^*\geq \alpha^*_0+1$ since $V({\tree{U}},{\tree{S}}, F)\restrict\alpha_0+1= \tree{S}\restrict\alpha_0+1$ and $V(\tree{V},\tree{T}, F)\restrict\alpha_0^*+1= \tree{T}\restrict\alpha_0^*+1$, so the first models which are possibly illfounded are the new models $M^{V({\tree{U}},{\tree{S}}, F)}_{\alpha_0+1}$ and $M^{V(\tree{V},\tree{T}, F)}_{\alpha_0^*+1}$, which are obtained as ultrapowers by $F$ and $G$, respectively.

Now, the $u$-map of $\Gamma$ is totally determined by what we have demanded in (1)-(3). We must have
\begin{equation*}
u^\Gamma(\zeta)=\begin{cases}
u^{\Psi}(\zeta) & \text{ if } \zeta<\alpha_0 \\
\alpha_0^* & \text{ if } \zeta=\alpha_0\\
u^{\Phi^*\circ \Pi}(\xi)  & \text{ if } \zeta>\alpha_0, \text{ where $\xi$ is such that $u^{\Phi}(\xi)=\zeta$}.
\end{cases}
\end{equation*}
Recall that this third case makes sense since $u^{\Phi}$ maps $[\beta,\lh(\tree{U}))$ onto $[\alpha_0+1, \lh (\tree{W}))$ and $u^{\Phi^*\circ\Pi}(\xi)>\alpha_0^*$ for all $\xi$ such that $u^{\Phi}(\xi)>{\alpha_0}$ (since if $u^{\Phi}(\xi)>{\alpha_0}$, then $\xi\geq \beta$, so $u^{\Phi^*\circ \Pi}(\xi)\geq u^{\Phi^*}(\beta)\geq \alpha_0^*+1$). 

The definition of $u^\Gamma(\zeta)$ just given makes sense for any ordinal $\zeta$, but the actual $u$-map of $\Gamma$ has domain $\{\zeta\mid  \zeta< \mu $ and $u^\Gamma(\zeta)< \mu^*\}$ (since the domain can't include anymore than this and we want the domain of $\Gamma$ to be maximal). In the course of the proof, we'll show that we can drop the condition \lq\lq$\zeta<\mu$" from the description of the domain of $u^\Gamma$.

Since we had to define $u^\Gamma$ as above and tree embeddings are totally determined by their $u$-map, uniqueness of $\Gamma$ is guaranteed. We just need to check that we can actually find a tree embedding with this $u$-map. This amounts to identifying $s$-maps and $t$-maps that make the relevant diagrams commute.

We've stipulated $\Gamma\restrict\alpha_0+1\approx \Psi\restrict\alpha_0 +1$ and $u^\Gamma(\alpha_0)=\alpha_0^*$. By Lemma \ref{key lemma 1} this determines a legitimate extended tree embedding from $\tree{S}\restrict\alpha_0+1$ into $\tree{T}\restrict\alpha_0^*+1$ with last $t$-map $t_{\alpha_0}^\Gamma=s_{\alpha_0, \alpha_0^*}^\Psi$.

Since $u^{\Phi}$ maps $[\beta,\lh ({\tree{U}}))$ onto $[\alpha+1,\lh ({\tree{W}}))$, we just need to find appropriate $s^\Gamma_{u^{\Phi}_\xi}$ and $t^\Gamma_{u^{\Phi}_\xi}$ by induction on $\xi\in[\beta,\lh({\tree{U}}))$. We also show that if $u^{\Phi^*\circ \Pi}(\xi)< \mu^*$, then $u^{\Phi}(\xi)< \mu$. 
We start with the base case.\\

\noindent \textbf{Base case.} $\xi=\beta$.\\ 

We first want to define $s_{\alpha_0+1}^\Gamma$. For $\Gamma$ to be a tree embedding, $s_{\alpha_0+1}^\Gamma$ must be the copy map associated to $(t_{\alpha_0}^\Gamma, s_{\beta, \beta^*}^\Gamma, E_{\alpha_0}^\tree{W}, E_{\alpha_0^*}^{\tree{W}^*})$. 

Note that we have $F=E_{\alpha_0}^\tree{W}$, $G=E_{\alpha_0^*}^{\tree{W}^*}$, and by Lemma \ref{key lemma 1}, $t^\Gamma_{\alpha_0}(F)=G$. Since $s^\Gamma_\beta=s^\Psi_\beta$ and $\beta^*\leq \alpha_0^*$, we have $s^\Gamma_{\beta, \beta^*}=s^\Psi_{\beta,\beta^*}$. Lemma \ref{key lemma 1}, again, gives $t^\Gamma_{\alpha_0}\restrict \restrict\dom(E_{\alpha_0}^\tree{W})\cup\{\dom(E_{\alpha_0}^\tree{W})\}=s_{\beta,\beta^*}^\Gamma\restrict \restrict\dom(E_{\alpha_0}^\tree{W})\cup\{\dom(E_{\alpha_0}^\tree{W})\}$. So the ordinary Shift Lemma does apply to $(t_{\alpha_0}^\Gamma, s_{\beta, \beta^*}^\Gamma, E_{\alpha_0}^\tree{W}, E_{\alpha_0^*}^{\tree{W}^*})$. So we can actually let $s_{\alpha_0+1}^\Gamma$ be the copy map associated to $(t_{\alpha_0}^\Gamma, s_{\beta, \beta^*}^\Gamma, E_{\alpha_0}^\tree{W}, E_{\alpha_0^*}^{\tree{W}^*})$. Here is a diagram for this application of the ordinary Shift Lemma.

\[\begin{tikzcd}[column sep = huge]
M^{{\tree{W}}}_{\alpha_0}=M^{{\tree{S}}}_{\alpha_0} \arrow{r}{ s_{\alpha_0,\alpha_0^*}^{\Psi}} & M^{\tree{T}}_{\alpha_0}=M^{\tree{W}^*}_{\alpha_0^*}\\
 F\arrow[mapsto]{r}{}& G\\
M^{{\tree{W}}}_{\beta}=M^{{\tree{S}}}_{\beta} \arrow{r}{ s_{\beta,\beta^*}^{\Psi}}\arrow{d}{ F} & M^{\tree{T}}_{\beta^*}=M^{\tree{W}^*}_{\beta^*}\arrow{d}{G}\\
M^{{\tree{W}}}_{\alpha_0+1}\arrow[dashed]{r}{s_{\alpha_0+1}^\Gamma} & M^{\tree{W}^*}_{\alpha_0^*+1}
\end{tikzcd}\]

If $\alpha^*_0+1<\mu^*$, then since $s^\Gamma_{\alpha_0}$ is a total elementary embedding from $M_{\alpha_0}^{\tree{W}}$ into $M_{\alpha_0^*}^{\tree{W}^*}$, $\alpha_0+1<\mu$, too.

If $V(\tree{U},\tree{S},F)$ is in the dropping case, we've already defined all of $\Gamma$. So for the remainder of the proof, assume that $V(\tree{U},\tree{S},F)$ is \textit{not} in the dropping case. We need to see that $V(\tree{V},\tree{T},G)$ is not in the dropping case either, or else we may not be able to successfully find $t^\Gamma_{\alpha_0+1}$ (or later maps).

First, suppose that $\beta+1=\lh(\tree{U})$. Then $F$ is total over $M_\beta^\tree{U}=M_\beta^\tree{S}=M_\beta^\tree{W}$, i.e. no proper initial segment of $M_\beta^\tree{U}$ beyond $\dom(F)$ projects to or below $\crit(F)$. It follows that for every $\eta\in [v^\Pi(\beta), u^\Pi(\beta)]_\tree{V}$, no level of $M_\eta^\tree{V}$ beyond $s_{\beta, \eta}^\Pi(\dom(F))$ projects to or below $\crit(s_{\beta, \eta}^\Pi(F))$. By hypothesis (4), $\beta^*$ is such an $\eta$ and $s_{\beta, \beta^*}^\Pi(\dom(F))=\dom(G)$. So $V(\tree{V},\tree{T}, G)$ is not in the dropping case either. 

Now suppose that $\beta+1<\lh(\tree{U})$. Then there is no level $P\isneq M_\beta^\tree{U}|\lh(E_\beta^\tree{U})$ beyond $\dom(F)$ which projects across $\crit(F)$. It follows that there is no level $Q\isneq M_{\beta^*}^{\tree{V}}|\lh(s_{\beta,\beta^*}^\Pi(E_\beta^\tree{U})$ beyond $\dom(G)$ which projects across $\crit(G)$. So the only potential problem would be that some level $Q\isneq M_{\beta^*}^{\tree{V}}| \lh(E_{\beta^*}^\tree{V})$ beyond $\lh(s_{\beta,\beta^*}^\Pi(E_\beta^\tree{U}))$ projects across $\crit(G)$. But in this case we must have $\beta^*<u^\Pi(\beta)$ (or else $E_{\beta^*}^\tree{V}=s_{\beta,\beta^*}^\Pi(E_\beta^\tree{U})$) and so, letting $\eta+1$ the successor of $\beta^*$ in $[v^\Pi(\beta), u^\Pi(\beta)]_\tree{V}$, $\crit(E_\eta^\tree{V})$ cannot be in the interval $[\crit(G), \lh(s_{\beta,\beta^*}^\Pi(E_\beta^\tree{U}))]$ (since nothing here is cardinal of $M_{\beta^*}^{\tree{V}}| \lh(E_{\beta^*}^\tree{V})$). Since we haven't reached the final image of $E_\beta^\tree{U}$ along $[v^\Pi(\beta), u^\Pi(\beta)]_\tree{V}$, we must have $\crit(E_\eta^\tree{V})<\lh(s_{\beta,\beta^*}^\Pi(E_\beta^\tree{U}))$. So actually $\crit(E_\eta^\tree{V})<\crit(G)$, contradicting hypothesis (4) that $t_\beta^\Pi\restrict\dom(F)\cup\{\dom(F)\}=s_{\beta,\beta^*}^\Pi\restrict\dom(F)\cup\{\dom(F)\}$.

So we've shown $V(\tree{V},\tree{T}, G)$ is not in the dropping case either. Now, if $\mu^* \leq u^{\Phi^*\circ\Pi}(\beta)$, we stop, so suppose $\mu^*> u^{\Phi^*\circ\Pi}(\beta)$. We must now put \begin{align*}
    t_{\alpha_0+1}^\Gamma &= \hat\imath^{\tree{W}^*}_{\alpha_0+1, u^{\Gamma}_{\beta}}\circ s_{\alpha_0+1}^\Gamma.
\end{align*}
For this to make sense, we need that $\alpha_0^*+1\leq_{\tree{W}^*}u^\Gamma(\beta)=u^{\Phi^*\circ\Pi}(\beta)$. Hypothesis (4) implies that any extender used along $(\beta^*,u^\Pi(\beta)]_\tree{V}$ has critical point $>\crit(G)$. So since $\Phi^*=\Phi^{V(\tree{V},\tree{T},G)}$, $u^{\Phi^*}$ is tree-order preserving on $[\beta^*,u^\Pi(\beta)]_\tree{V}$. So $\alpha_0^*+1=u^{\Phi^*}(\beta^*)\leq_{\tree{W}^*}u^{\Phi^*\circ\Pi}(\beta)=u^\Gamma(\beta)$, as desired. We have the following picture.

\[\begin{tikzcd}[column sep = huge]
M^{{\tree{W}}}_{\beta}=M^{{\tree{U}}}_{\beta} \arrow{r}{ s_{\beta,\beta^*}^{\Pi}}\arrow{d}{ F} & M^{\tree{V}}_{\beta^*}=M^{\tree{W}^*}_{\beta^*}\arrow{d}{G} \arrow {r}{\tree{V}}& M^{\tree{V}}_{u^\Pi_{\beta}}\arrow{d}{t^{\Phi^*}_{u^\Pi_{\beta}}}\\
M^{{\tree{W}}}_{\alpha_0+1}\arrow[dashed]{r}{s_{\alpha_0+1}^\Gamma} & M^{\tree{W}^*}_{\alpha_0^*+1}\arrow{r}{\tree{W}^*} & M^{\tree{W}^*}_{u^{\Phi^*\circ \Pi}_{\beta}}
\end{tikzcd}\]

As part of conclusion (3), we need to see that this diagram commutes. We already have that the left square commutes, by our choice of $s^\Gamma_{\alpha_0+1}$, so we just need to see that the right one does. We may assume $\beta^*<u^\Pi(\beta)$, as otherwise this is trivial. Now since $\Phi^*=\Phi^{V(\tree{V},\tree{T}, G)}$, for every $\zeta>\beta^*$, $u^{\Phi^*}(\zeta)=v^{\Phi^*}(\zeta)$ and so $t_\zeta^{\Phi^*}=s_\zeta^{\Phi^*}$, where $\zeta+1$ be the least element of $(\beta^*,u^\Pi(\beta)]_{\tree{V}}$, so we can expand the right square as follows.

\[\begin{tikzcd}[column sep = huge]
M^{\tree{V}}_{\beta^*}\arrow{d}{G} \arrow {r}{E^\tree{V}_\zeta}& M^\tree{V}_{\zeta} \arrow{r}{\tree{V}} \arrow{d}{s_\zeta^{\Phi^*}}& M^\tree{V}_{u^\Pi_{\beta}}\arrow{d}{s_{u^\Psi_{\beta}}^{\Phi^*}}\\
 M^{\tree{W}^*}_{\alpha_0^*+1}\arrow{r}{E^{\tree{W}^*}_ {u^{\Phi^*}_{\zeta}}} & M^{\tree{W}^*}_{v^{\Phi^*}_{\zeta}} \arrow{r}{\tree{W}^*}& M^{\tree{W}^*}_{u^{\Phi^*\circ \Psi}_{\beta}}
\end{tikzcd}\]

But these squares commute since $\Phi$ is a tree embedding: the left square commutes since $s^\Phi_\zeta$ is the appropriate copy map and the right square commutes since the $s$-maps of a tree embedding commute with branch embeddings, by definition. This finishes the base case.\\

\noindent \textbf{Successor case.} $\beta<\xi+1$ and $u^{\Phi^*\circ \Pi}(\xi+1)< \mu^*$.\\

Since $v^{\Phi^*\circ \Pi}(\xi)<\mu^*$, we have $M^{\tree{W}^*}_{v^{\Phi^*\circ \Pi}(\xi)}$ is wellfounded and so $M^{{\tree{W}}}_{v^{\Phi}(\xi)}$ is as well (since $s^\Gamma_{u^{\Phi}_\xi}:M^{{\tree{W}}}_{v^{\Phi}(\xi)}\to M^{\tree{W}^*}_{v^{\Phi^*\circ \Pi}(\xi)}$ is a total elementary embedding). If $\xi>\beta$, then $v^{\Phi}(\xi)= u^{\Phi}(\xi)$, so we have $u^{\Phi}(\xi)<\mu$. So $u^{\Phi}(\xi+1)=v^{\Phi}(\xi+1)=u^{\Phi}(\xi)+1\leq \mu$. If $\xi=\beta$, then we have $u^{\Phi}(\xi)=\alpha_0+1$, so we already had $u^{\Phi}(\xi)\leq \mu$.

Let $\eta = {\tree{U}}\pred (\xi+1)$ and  $\eta^*=\tree{V}\pred (u^\Pi(\xi)+1)$. We have $\eta\in[v^\Pi(\eta),u^\Pi(\eta)]_\tree{V}$ since $\Pi$ is a tree embedding. 

We consider two subcases depending on the critical point of $E^{{\tree{U}}}_{\xi}$. The arguments in each case are basically the same.

\paragraph{Subcase A.} $\crit (E^{{\tree{U}}}_{\xi})<\crit(F)$. \\

In this case, $\eta\leq \beta$ and $\eta={\tree{W}}\pred(u^{\Phi}(\xi)+1)$. We also have that $\crit (E^\tree{V}_{u^\Pi(\xi)})<\crit (G)$, so that $\eta^*\leq \beta^*$ and $\eta^*=\tree{W}^*\pred(u^{\Phi^*\circ \Pi}(\xi)+1)$, as well. We have the following picture, in the case that we don't drop.

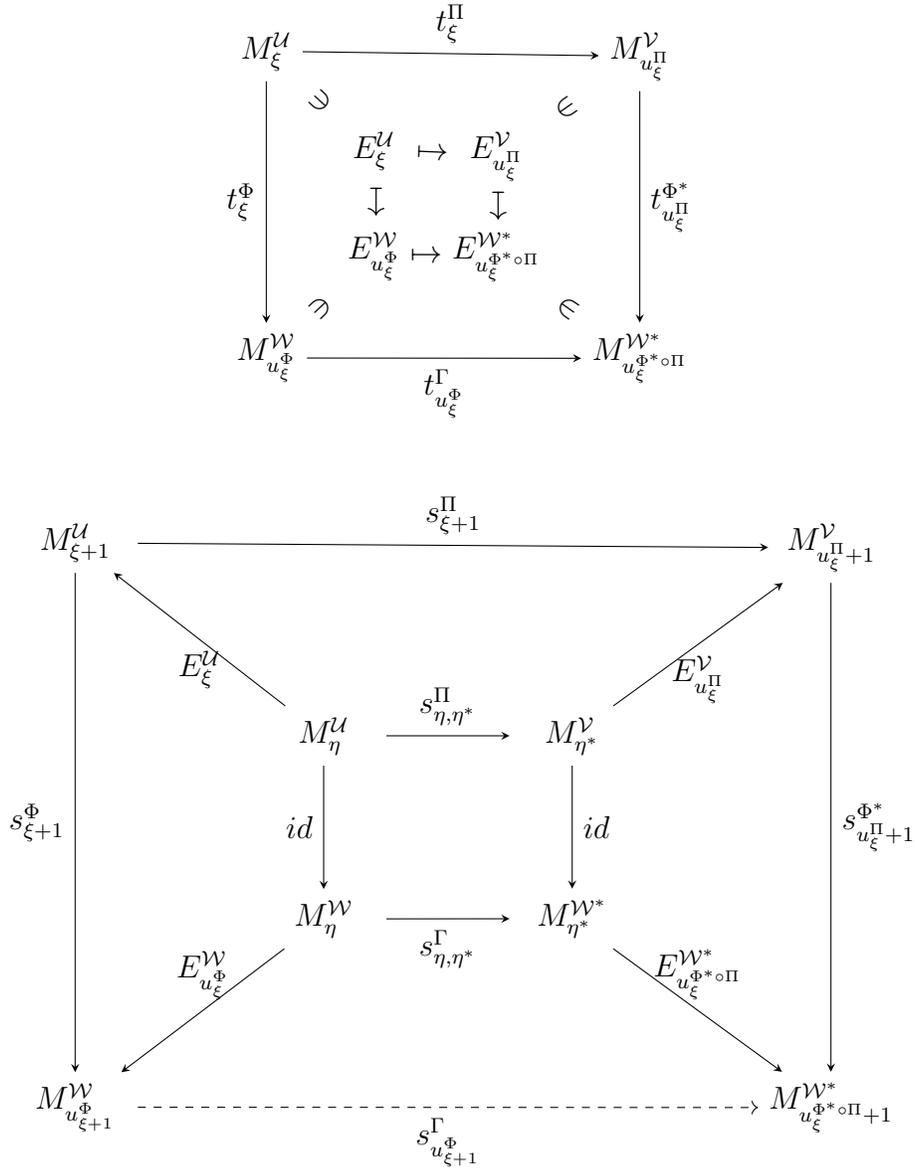
\[\begin{tikzpicture}
  \matrix (m) [matrix of math nodes,row sep=1em,column sep=1em,minimum width=1em]
  {
     M^{{\tree{U}}}_\xi & {} & {} & M^\tree{V}_{u^\Pi_\xi} \\
    {} & E^{{\tree{U}}}_\xi & E^\tree{V}_{u^\Pi_\xi}\\
    {} & E^{{\tree{W}}}_{u^{\Phi}_\xi} & E^{\tree{W}^*}_{u^{\Phi^*\circ\Pi}_\xi} \\
    M^{{\tree{W}}}_{u^{\Phi}_\xi} & {} & {} & M^{\tree{W}^*}_{u^{\Phi^*\circ 
    \Pi}_\xi}\\
   };
  \path[-stealth]
    (m-1-1)
        edge
         node [above] {$t_\xi^\Pi$}
         (m-1-4)
         edge
         node [left] {$t_\xi^{\Phi}$}
         (m-4-1)
    (m-1-4)
        edge
         node [right] {$t_{u^\Pi_\xi}^{\Phi^*}$}
         (m-4-4)
    (m-4-1)
        edge
         node [below] {$t_{u^{\Phi}_\xi}^\Gamma$}
         (m-4-4)     
    
    (m-2-2) edge[draw=none]
             node [sloped, auto=false,
              allow upside down] {$\in$}
        (m-1-1)
    edge[draw=none]
             node [sloped, auto=false,
              allow upside down] {$\mapsto$}
         (m-2-3)
    edge[draw=none]
             node [sloped, auto=false,
              allow upside down] {$\mapsto$}
         (m-3-2)  
    (m-2-3) edge[draw=none]
             node [sloped, auto=false,
              allow upside down] {$\in$}
        (m-1-4)
    edge[draw=none]
             node [sloped, auto=false,
              allow upside down] {$\mapsto$}
         (m-3-3)
    (m-3-2)
    edge[draw=none]
             node [sloped, auto=false,
              allow upside down] {$\mapsto$}
         (m-3-3)
    edge[draw=none]
             node [sloped, auto=false,
              allow upside down] {$\in$}
         (m-4-1)
    (m-3-3)
    edge[draw=none]
             node [sloped, auto=false,
              allow upside down] {$\in$}
         (m-4-4) ;
\end{tikzpicture}\]

\[\begin{tikzpicture}
  \matrix (m) [matrix of math nodes,row sep=4em,column sep=4em,minimum width=4em]
  {
     M^{{\tree{U}}}_{\xi+1} & {} & {} & M^\tree{V}_{u^\Pi_\xi+1} \\
    {} & M^{{\tree{U}}}_{\eta} & M^\tree{V}_{\eta^*}\\
    {} & M^{{\tree{W}}}_{\eta} & M^{\tree{W}^*}_{\eta^*} \\
    M^{{\tree{W}}}_{u^{\Phi}_{\xi+1}} & {} & {} & M^{\tree{W}^*}_{u^{\Phi^*\circ \Pi}_\xi+1}\\
   };
  \path[-stealth]
    (m-1-1)
        edge
         node [above] {$s_{\xi+1}^\Pi$}
         (m-1-4)
         edge
         node [left] {$s_{\xi+1}^{\Phi}$}
         (m-4-1)
    (m-1-4)
        edge
         node [right] {$s_{u^\Pi_\xi+1}^{\Phi^*}$}
         (m-4-4)
    (m-4-1)
        edge[dashed,->]
         node [below] {$s_{u^{\Phi}_{\xi+1}}^\Gamma$}
         (m-4-4)
    (m-2-2) edge
             node
            [below] {$E^{{\tree{U}}}_\xi$}
        (m-1-1)
    edge
             node [above] {$ s_{\eta,\eta^*}^\Pi$}
         (m-2-3)
    edge
             node [left] {$id$}
         (m-3-2)  
    (m-2-3) edge
             node
            [below] {$E^\tree{V}_{u^\Pi_\xi}$}
        (m-1-4)
    edge
             node [right] {$id$}
         (m-3-3)
    (m-3-2)
    edge
             node [below] {$ s_{\eta,\eta^*}^\Gamma$}
         (m-3-3)
    edge
             node
            [above] {$E^{{\tree{W}}}_{u^{\Phi}_\xi}$}
         (m-4-1)
    (m-3-3)
    edge
    node
            [above] {$E^{\tree{W}^*}_{u^{\Phi^*\circ\Pi}_\xi}$}
         (m-4-4) ;
\end{tikzpicture}\]

Each of the maps along the outer square of the bottom diagram are the copy maps associated to the maps along corresponding side of the top square, inner square, and appropriate extenders.
In particular, we let $s_{u^{\Phi}_{\xi+1}}^\Gamma$ be the copy map associated to ($t^\Gamma_{u^{\Phi}_\xi}$, $ s^\Gamma_{\eta,\eta^*}$, $E^{{\tree{W}}}_{u^{\Phi}(\xi)}, E^{{\tree{W}^*}}_{u^{\Phi^*\circ \Pi}(\xi)}$), as we must.  Note that the ordinary Shift Lemma applies in this case because we have assumed that, so far, $\Gamma$ is a tree embedding. In particular, since $\dom(E^{{\tree{W}}}_{u^{\Phi}(\xi)})\is M^{{\tree{W}}}_{\eta}|\hat\lambda(E_{\eta}^{{\tree{W}}})$, the agreement properties of $t$-maps gives that 
\begin{align*} t^\Gamma_{u^{\Phi}(\xi)}\restrict\dom(E^{{\tree{W}}}_{u^{\Phi}(\xi)})\cup \{\dom(E^{{\tree{W}}}_{u^{\Phi}(\xi)})\} & = t^\Gamma_{\eta}\restrict\dom(E^{{\tree{W}}}_{u^{\Phi}(\xi)})\cup \{\dom(E^{{\tree{W}}}_{u^{\Phi}(\xi)})\} \\
& =  s^\Gamma_{\eta,\eta^*}\restrict\dom(E^{{\tree{W}}}_{u^{\Phi}(\xi)})\cup \{\dom(E^{{\tree{W}}}_{u^{\Phi}(\xi)})\},
\end{align*}
using for this second equivalence that $\crit(\hat\imath^{\tree{W}^*}_{\eta,u^{\Gamma}(\eta)})>\crit(E^{\tree{W}^*}_{u^{\Phi^*\circ \Pi}(\xi)})$, as, otherwise, we used an extender $E_\gamma^{\tree{W}^*}$ with $\crit(E_\gamma^{\tree{W}^*})\leq \crit(E^{\tree{W}^*}_{u^{\Phi^*\circ \Pi}(\xi)})<\hat\lambda(E_\gamma^{\tree{W}^*})$, but then $\crit(E^{\tree{W}^*}_{u^{\Phi^*\circ \Pi}(\xi)})$ can't be in $\ran(\hat\imath^{\tree{W}^*}_{\eta^*,u^{\Gamma}(\eta)})\supseteq \ran(t^\Gamma_{\eta})$, a contradiction.

In the lower diagram, the inner square commutes by our induction hypothesis and all the trapezoids commute since the outer maps are copy maps associated to the relevant objects. We now want to see that the full outer square commutes. 
Let's look at the two ways of going around this outer square, $s^{\Phi^*}_{u^{\Pi_\xi +1}}\circ s^\Pi_{\xi+1}$ and $s^\Gamma_{u^{\Phi}_{\xi+1}}\circ s^{\Phi}_{\xi+1}$. 
Since $s^\Pi_{\xi+1}$ and $s^{\Phi^*}_{u^\Pi_\xi +1}$ are copy maps associated to the appropriate objects, Lemma \ref{shift composition} gives that $s^{\Phi^*}_{u^\Pi_\xi +1}\circ s^\Pi_{\xi+1}$ is the copy map associated to ($t^{\Phi^*}_{u^\Pi_\xi}\circ t^\Pi_\xi$, $ s_{\eta,\eta^*}^\Pi$, $E^{{\tree{U}}}_\xi$, $E^{\tree{W}^*}_{u^{\Phi^*\circ \Pi}_\xi}$). 
Similarly, $s^\Gamma_{u^{\Phi}_{\xi+1}}\circ s^{\Phi}_{\xi+1}$ is the copy map associated to ($t^\Gamma_{u^{\Phi}_\xi}\circ t^{\Phi}_\xi$, $ s^\Gamma_{\eta,\eta^*}$, $E^{{\tree{U}}}_\xi$, $E^{\tree{W}^*}_{u^{\Phi^*\circ \Pi}_\xi}$).
But $t^{\Phi^*}_{u^\Pi_\xi}\circ t^\Pi_\xi= t^\Gamma_{u^{\Phi}_\xi}\circ t^{\Phi}_\xi$ and $s_{\eta,\eta^*}^\Pi= s^\Gamma_{\eta,\eta^*}$, so the two ways of going around the outer square are both the copy map associated to the same objects, i.e. $s^{\Phi^*}_{u^\Pi_\xi +1}\circ s^\Pi_{\xi+1}=s^\Gamma_{u^{\Phi}_{\xi+1}}\circ s^{\Phi}_{\xi+1}$.

Note that $u^{\Phi^*\circ \Pi}_\xi+1= v^{\Phi^*\circ\Pi}_{\xi+1}\leq_{\tree{W}^*} u^{\Phi^*\circ\Pi}_{\xi+1}$. We now define \[t_{u^{\Phi}_{\xi+1}}^\Gamma= \hat\imath^{\tree{W}^*}_{u^{\Phi^*\circ \Pi}_\xi+1, u^{\Phi^*\circ \Pi}_{\xi+1}}\circ s_{u^{\Phi}_{\xi+1}}^\Gamma,\] as we must. Finally, we check that this assignment gives us a commuting square of the $t$-maps. We get the following diagram.
\[\begin{tikzpicture}
\matrix (m) [matrix of math nodes,row sep=6em,column sep=6em,minimum width=6em]
  {
     M^{{\tree{U}}}_{\xi+1} & M^\tree{V}_{v^\Pi_{\xi+1}} & M^\tree{V}_{u^\Pi_{\xi+1}} \\
     M^{{\tree{W}}}_{u^{\Phi}_{\xi+1}} & M^{\tree{W}^*}_{v^{\Phi^*}( u^\Pi_\xi+1)} & M^{\tree{W}^*}_{u^{\Phi^*\circ\Pi}_{\xi+1}}\\
 };
  \path[-stealth]
(m-1-1)
    edge
    node
    [above]
    {$s_{\xi+1}^\Pi$}
    (m-1-2)
    edge
    node
    [left]
    {$t_{\xi+1}^{\Phi}$}
    (m-2-1)
(m-1-2)
    edge
    node
    [above]
    {$\tree{V}$}
    (m-1-3)
    edge
    node
    [left]
    {$t_{v^\Pi_{\xi+1}}^{\Phi^*}$}
    (m-2-2)
(m-1-3)
    edge
    node
    [right]
    {$t_{u^\Pi_{\xi+1}}^{\Phi^*}$}
    (m-2-3)
(m-2-1)
    edge
    node
    [below]
    {$s_{u^{\Phi}_{\xi+1}}^{\Gamma}$}
    (m-2-2)
(m-2-2)
    edge
    node
    [below]
    {$\tree{W}^*$}
    (m-2-3);
\end{tikzpicture}\]

We just need to see that this diagram commutes, since $t_{\xi+1}^\Pi$ is just the map going across the top and $t_{u^{\Phi}(\xi+1)}^\Gamma$ is the map going across the bottom (so this really is the relevant square of $t$-maps). The left square is just the outer square of the lower commuting diagram, above (though we used $u^\Pi(\xi)+1= v^\Pi(\xi+1)$ and $u^{\Phi^*}\circ u^\Pi(\xi) +1= v^{\Phi^*}(u^\Pi(\xi)+1)$ to change the labeled indices of the models in the middle column to emphasize how we knew they were tree-related to the appropriate models all the way on the right; we get these equivalences since $\xi+1>\beta$ and $u^\Pi(\xi+1), u^\Pi(\xi)+1>\beta^*$). We've also used that all the vertical $t$-maps are the same as the corresponding $s$-maps (by the equivalence of the indices just mentioned). This last fact (that the vertical $t$-maps are the same as the corresponding $s$-maps) also gives us that the square on the right commutes, since $\Phi$ is a tree embedding.

If we drop when applying any of the $E^{{\tree{U}}}_\xi, E^\tree{V}_{u^\Pi_\xi}$, etc, then we drop applying all of them and the initial segments to which we apply these extenders are all mapped to each other by the relevant maps. In this case, everything remains the same except that we must use the initial segments to which we apply the extenders instead of the models displayed in the above diagrams, e.g. some $P\is M^{{\tree{U}}}_{{\eta}}$ instead of $M^{{\tree{U}}}_{{\eta}}$.

\paragraph{Subcase B.} $\crit (E^{{\tree{U}}}_{\xi})\geq \crit ( F)$. \\

In this case, $\eta\geq \beta$ and $u^{\Phi}(\eta)={\tree{W}}\pred(u^{\Phi}(\xi)+1)$. We also get $\crit (E^\tree{V}_{u^\Pi(\xi)})\geq \crit (G)$, so $\eta^*\geq \beta^*$ and $u^{\Phi^*}(\eta^*)=\tree{W}^*\pred(u^{\Phi^*\circ \Pi}(\xi)+1)$. We now have the model to which $E^{{\tree{U}}}_\xi$ is applied is related to the model to which $E^{{\tree{W}}}_{u^{\Phi}(\xi)}$ is applied by a $t$-map of $\Phi$, whereas they were just the same model in the previous case. Similarly on the $\tree{V}$-$\tree{W}^*$ side. The only thing this changes is that we replace the identity maps in the above previous diagram with these $t$-maps.  This is the diagram for the non-dropping case (as before, dropping makes little difference).

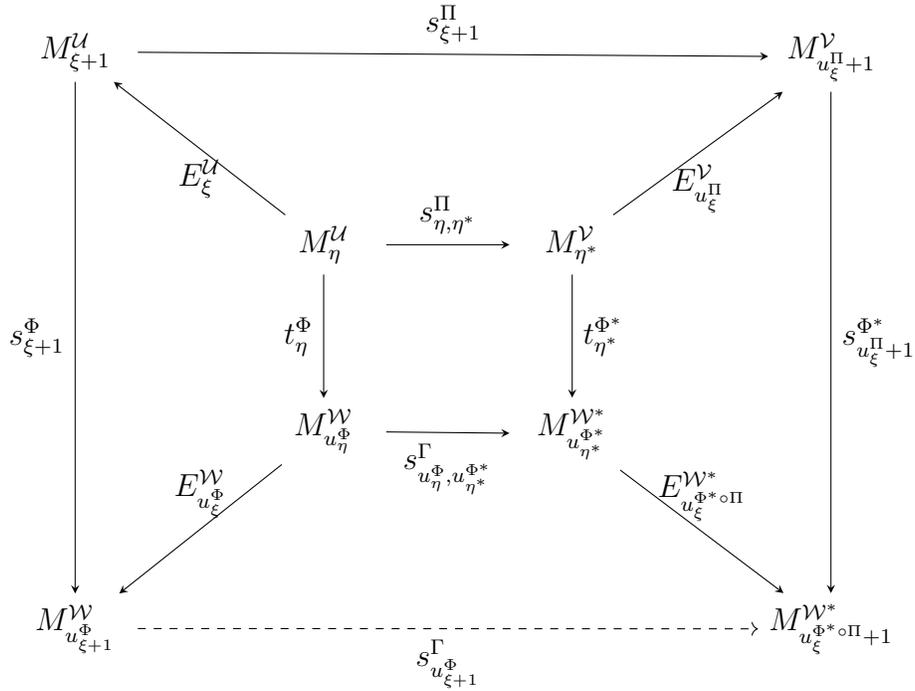
\[\begin{tikzpicture}
  \matrix (m) [matrix of math nodes,row sep=4em,column sep=4em,minimum width=4em]
  {
     M^{{\tree{U}}}_{\xi+1} & {} & {} & M^\tree{V}_{u^\Pi_\xi+1} \\
    {} & M^{{\tree{U}}}_{\eta} & M^\tree{V}_{\eta^*}\\
    {} & M^{{\tree{W}}}_{u^{\Phi}_{\eta}} & M^{\tree{W}^*}_{u^{\Phi^*}_{\eta^*}} \\
    M^{{\tree{W}}}_{u^{\Phi}_{\xi+1}} & {} & {} & M^{\tree{W}^*}_{u^{\Phi^*\circ\Pi}_\xi+1}\\
   };
  \path[-stealth]
    (m-1-1)
        edge
         node [above] {$s_{\xi+1}^\Pi$}
         (m-1-4)
         edge
         node [left] {$s_{\xi+1}^{\Phi}$}
         (m-4-1)
    (m-1-4)
        edge
         node [right] {$s_{u^\Pi_\xi+1}^{\Phi^*}$}
         (m-4-4)
    (m-4-1)
        edge[dashed,->]
         node [below] {$s_{u^{\Phi}_{\xi+1}}^\Gamma$}
         (m-4-4)
    (m-2-2) edge
             node
            [below] {$E^{{\tree{U}}}_\xi$}
        (m-1-1)
    edge
             node [above] {$ s_{\eta,\eta^*}^\Pi$}
         (m-2-3)
    edge
             node [left] {$t_{\eta}^{\Phi}$}
         (m-3-2)  
    (m-2-3) edge
             node
            [below] {$E^\tree{V}_{u^\Pi_\xi}$}
        (m-1-4)
    edge
             node [right] {$t_{\eta^*}^{\Phi^*}$}
         (m-3-3)
    (m-3-2)
    edge
             node [below] {$ s_{u^\Phi_\eta,u^{\Phi^*}_{\eta^*}}^\Gamma$}
         (m-3-3)
    edge
             node
            [above] {$E^{{\tree{W}}}_{u^{\Phi}_\xi}$}
         (m-4-1)
    (m-3-3)
    edge
    node
            [above] {$E^{\tree{W}^*}_{u^{\Phi^*\circ\Pi}_\xi}$}
         (m-4-4) ;
\end{tikzpicture}\]

The rest of the diagrams and arguments are essentially the same as before. This finishes the successor case.

\paragraph{Limit case.} $\lambda>\beta^*$ is a limit and $u^{\Phi^*\circ \Pi}(\lambda)\leq\mu^*$. 

We have $u^{\Phi^*\circ \Pi}(\xi)<\mu^*$ for all $\xi<\lambda$, so that by our induction hypothesis, $u^{\Phi}(\xi)<\mu$.
So, $u^{\Phi}(\lambda)=v^{\Phi}(\lambda)=\sup\{u^{\Phi}(\xi)\mid \xi<\lambda\}\leq \mu$. 

Let $ c=[0,u^{\Phi}(\lambda)_{{\tree{W}}}$ and $c^*=[0,u^{\Phi^*\circ \Pi}(\lambda))_{\tree{W}^*}$. We need to see that $c^*$ is the $\leq_{\tree{W}^*}$-downward closure of $v^\Gamma[c]$. To do this, we just trace $c$, $c^*$ back to the branch $ b=[0,\lambda)_{{\tree{U}}}$. We have \[c=\{\xi\mid  \exists \eta\in b \,(\xi\leq_{{\tree{W}}}v^{\Phi}(\eta))\},\]
and
\[c^*=\{\xi\mid  \exists \eta\in b \,(\xi\leq_{\tree{W}^*} v^{\Phi^*\circ \Pi}(\eta))\}.\]
We also have that $v^\Gamma(v^{\Phi}(\eta))=v^{\Phi^*\circ \Pi}(\eta)$ for every $\eta$, so that $v^\Gamma[c] \subseteq c^*$. This implies $c^*$ is the downward closure of $v^\Gamma[c]$, as desired.

So, we get our map $s^\Gamma_{u^{\Phi}_\lambda}$ commuting with the maps $s^\Gamma_{u^{\Phi}(\xi)}$ since we are taking the direct limits along $c$ and $c^*$ on both sides. From here, we get $t^\Gamma_{u^{\Phi}_\lambda}$ as in the successor case.
This finishes our construction of $\Gamma$. 

For the \lq\lq moreover" clause, we've already shown that if $u^{\Phi^*\circ \Pi}(\xi)\leq\mu^*$, $u^{\Phi}(\xi)\leq\mu$. So, if the full $V(\tree{V},\tree{T}, G)$ is wellfounded, then for all $\xi< \lh({\tree{U}})$, $u^{\Phi}(\xi)<\mu$. But then $\mu+1=\lh (V({\tree{U}},{\tree{S}},F))$, so $\Gamma$ is a total extended tree embedding from $V({\tree{U}},{\tree{S}},F)$ into $V(\tree{V},\tree{T},G)$. 
\qed
\end{proof}

We will carry over our notation for the ordinary Shift Lemma.

\begin{definition} Let $\Psi:\tree{S}\to \tree{T}$ and $\Pi:\itU \to \itV$ be extended tree embeddings, $F$ an extender such that $F^-$ be
 an extender on the $M_\infty^\tree{S}$-sequence, and $G$ an extender such that $G^-$ is on the $M_\infty^\tree{T}$-sequence. Suppose that the Shift Lemma applies
  to $(\Psi,\Pi, F, G)$.
We'll say that an extended tree embedding $\Gamma: V(\tree{U},\tree{S}, F)\to V(\tree{V},\tree{T}, G)$ is \textit{the copy map associated to $(\Psi, \Pi, F, G)$} iff it is the unique tree embedding as in the conclusion of the Shift Lemma.
\end{definition}

We can now carry out the copying construction.

\begin{theorem}[Copying]\label{copying} Let $\Gamma:\tree{S}\to \tree{T}$ be a non-dropping extended tree embedding. Let $\mtree{S}=\langle \tree{S}_\xi, \Phi^{\eta,\xi},F_\zeta\,|\, \xi,\zeta+1<\lh (\mtree{S})\rangle$ be a meta-tree on $\tree{S}$.

Then there is some largest $\mu\leq \lh (\mtree{S})$ such that there is a meta-tree $\Gamma\mtree{S}=\langle \tree{T}_\xi, \Psi^{\eta,\xi}, G_\zeta\,|\,\xi,\zeta+1<\mu\rangle$ on $\tree{T}$ with tree-order $\leq_\mtree{S}\restrict \mu$ and for $\xi<\mu$, non-dropping extended tree embeddings $\Gamma^\xi: \tree{S}_\xi\to \tree{T}_\xi$ with (total) last $t$-map $t_\infty^\xi$ such that 
\begin{enumerate}
    \item $\Gamma=\Gamma^0$, \item$G_\xi=t_\infty^\xi(F_\xi)$, 
    \item  and for all $\eta\leq_\mtree{S}\xi$, $\Gamma^\xi\circ \Phi^{\eta,\xi}=\Psi^{\eta,\xi}\circ \Gamma^\eta$.
\end{enumerate}  

Moreover, \begin{enumerate}
\item if $\mtree{S}$ is simple, so is $\Gamma\mtree{S}$,
    \item if $\mtree{S}=\mtree{V}(\tree{S}, \tree{U})$ for some tree $\tree{U}$, then $\Gamma\mtree{S}= \mtree{V}(\tree{T},t^\Gamma_\infty \tree{U}\restrict\mu)$.\footnote{We can copy $\tree{U}$ by $t^\Gamma_\infty$ since this is a total elementary embedding from the last model of $\tree{S}$ to the last model of $\tree{T}$.}
    \item if $\tree{S},\tree{T}$ are by some strategy $\Sigma$ with strong hull condensation and $\mtree{S}$ is by $\Sigma^*$, then $\mu=\lh (\mtree{S})$ and $\Gamma\mtree{S}$ is by $\Sigma^*$.

\end{enumerate}
\end{theorem}

\begin{proof}
We define $\Gamma \mtree{S}$ by induction, using the Shift Lemma at successors. $\mu$ will just be the least ordinal such that this process breaks down or the full $\lh (\mtree{S})$ if it doesn't break down. We just do the case that $\mtree{S}$ is simple. To deal with gratuitous drops, we just make the corresponding drops on the $\Gamma\mtree{S}$ side as well. That is, if at stage $\xi$, $\tree{S}_\xi= \tree{S}^+_\xi\restrict\eta+1$ with $\eta+1< \lh (\tree{S}^+_\xi)$, we just put $\tree{T}_\xi = \tree{T}^+_\xi\restrict v^{\Gamma^\xi}(\eta)+1$.\footnote{This is somewhat arbitrary. We only need to drop to some level $\tree{T}_\xi^+\restrict \eta^*+1$ such that $v^{\Gamma^\xi}(\eta)\leq_{\tree{T}^+_\xi}\eta^*$ and $(v^{\Gamma^\xi}(\eta),\eta*)_{\tree{T}^+_\xi}$ doesn't drop so that we have an extended tree embedding with a total last $t$-map. (Note that here we technically mean the extension of $\Gamma^\xi$ to an extended tree embedding $\tree{S}^+_\xi\to \tree{T}^+_\xi$.)}

Let $\alpha_\xi=\alpha_0(F_\xi, \tree{S}_\xi)$ and $\beta_\xi=\beta(F_\xi, \tree{S}_\xi)$.
Supposing we've define $\Gamma\mtree{S}\restrict\xi+1$, let $\alpha_\xi^*=\alpha_0(G_\xi, \tree{T}_\xi)$ and $\beta_\xi^*=\beta( G_\xi,\tree{T}_\xi)$. 
We'll maintain the following by induction. For $\eta\leq\xi<\mu$, (1)-(3) hold as well as
\begin{enumerate}
    \item[(4)] $\Gamma^\eta\restrict\alpha_\eta+1\approx\Gamma^\xi\restrict\alpha_\eta+1$,
    \item[(5)] $u^{\Gamma^\xi}(\alpha_\eta)\leq_{\tree{T}_\eta} u^{\Gamma^\eta}(\alpha_\eta)$ and $\tree{T}_\xi \restrict u^{\Gamma^\xi}(\alpha_\eta)+1=\tree{T}_\eta \restrict u^{\Gamma^\xi}(\alpha_\eta)+1$ 
   \item[(6)] $t^\xi_\infty$ is total and $t^\eta_\infty\restrict\lh (F_\eta)+1=t^\xi_\infty\restrict\lh( F_\eta)+1$.
\end{enumerate}

This will allow us to verify that $\Gamma\mtree{S}$ has the same tree order as $\mtree{S}$ and show that we satisfy the hypotheses of the Shift Lemma at successor stages.

Suppose we've defined $\Gamma(\mtree{S}\restrict \xi+1)$, so we have $\Gamma^\xi:\tree{S}_\xi\to \tree{T}_\xi$ an extended tree embedding with total last $t$-map $t^\xi_\infty$, by (6). Let $\eta=\mtree{S}\pred (\xi+1)$. We need to check that $\eta$ is least such that $\crit (G_\xi)<\hat\lambda(G_\eta)$, so that $\eta=\Gamma\mtree{S}\pred(\xi+1)$, according to normality.
$\crit (G_\xi) =t^\xi_\infty(\crit (F_\xi))= t_\infty^\eta(\crit (F_\xi))$, since $\crit (F_\xi)<\hat\lambda(F_\eta)$ and $t^\eta_\infty$ agrees with $t_\infty^\xi$ up to $\lh (F_\eta) +1$ by (6). So $\crit (G_\xi)=t^\eta_\infty(\crit (F_\xi))<t^\eta_\infty(\hat\lambda(F_\eta))=\hat\lambda(G_\eta)$.

Now suppose $\zeta$ is such that $\crit (G_\xi)<\hat\lambda(G_\zeta)$. Then $t^\zeta_\infty(\crit (F_\xi))= t^\xi_\infty(\crit (F_\xi)) = \crit (G_\xi)< \hat\lambda(G_\zeta)=t^\zeta_\infty(\hat\lambda(F_\zeta))$. So $\crit (F_\xi)<\hat\lambda(F_\zeta)$, so $\zeta\geq \eta$, as desired.

We now want to apply our Shift Lemma to, in the notation of that lemma, $\Psi= \Gamma_\xi$, $\Pi=\Gamma_\eta$, $F=F_\xi$ and $G=G_\xi$, and then let $\Gamma_{\xi+1}$ be the resulting copy tree embedding $\Gamma$, assuming $V(\tree{T}_\eta,\tree{T}_\xi, G_\xi)$ is wellfounded.

\begin{claim}\label{copy claim 1} The Shift Lemma applies to $(\Gamma^\xi, \Gamma^\eta, F_\xi, G_\xi)$, i.e. \begin{enumerate}
     \item[(i)] $M_\infty^{\tree{S}_\xi}|\lh(F_\xi)\is\dom( t_\infty^\xi)$ and $G_\xi=t_\infty^\xi(F)$,
  \item[(ii)] $\Gamma_\xi\restrict\beta_\xi+1\approx \Gamma_\eta\restrict\beta_\xi+1$,
\item[(iii)] $\tree{T}_\xi\restrict \beta^*_\xi+1=\tree{T}_\eta\restrict\beta^*_\xi+1$
\item[(iv)] $\beta^*_\xi\in [v^\eta(\beta_\xi), u^\eta(\beta_\xi)]_{\tree{T}_\eta}$ and $t_{\beta_\xi}^\eta\restrict\dom(F_\xi)\cup\{\dom(F_\xi)\}=s_{\beta, \beta^*}^\eta\restrict\dom(F_\xi)\cup\{\dom(F_\xi)\}$,
\item[(v)]
   if  $\beta_\xi+1<\lh(\tree{S}_\eta)$, then 
  $\dom(F_\xi) \isneq M_{\beta_\xi}^{\tree{S}_\eta}|\lh(E^{\tree{S}_\eta}_{\beta_\xi})$, and
  \item[(vi)]
   if  $\beta^*_\xi+1<\lh(\tree{T}_\eta)$,
  $\dom(G_\xi) \isneq M_{\beta^*_\xi}^{\tree{S}_\eta}|\lh(E^{\tree{T}_\eta}_{\beta^*_\xi})$.
\end{enumerate}
\end{claim}
\begin{proof}
First notice that (i) is trivial as $t^\xi_\infty$ is a total elementary embedding. We also have $\beta_\xi\leq\alpha_\eta$, since either $\eta=\xi$ and this is trivial, or else $\eta<\xi$ and $F_\eta=E_{\alpha_\eta}^{\tree{S}_\xi}$, so that this follows by our choice of $\eta$. Similarly, $\beta_\xi^*\leq\alpha^*_\eta$.

These observations and hypothesis (4) imply (ii) and (iii). For (iv) we split into cases depending on whether $\beta_\xi=\alpha_\eta$ or $\beta_\xi<\alpha_\eta$. If $\beta_\xi=\alpha_\eta$, then Lemma \ref{key lemma 1} and hypothesis (5) gives (iv). If $\beta_\xi<\alpha_\eta$, then this just follows from Lemma \ref{key lemma 1} and (4). (v) and (vi) follow by considering similar cases.
\hfill{$\qed$ Claim \ref{copy claim 1}}
\end{proof}

So we may let $\Gamma^{\xi+1}$ be the copy tree embedding associated to $(\Gamma_\xi, \Gamma_\eta, F_\xi, G_\xi)$. We assume that $V(\tree{T}_\eta, \tree{T}_\xi, G_\xi)=\tree{T}_{\xi+1}$ is wellfounded, otherwise we stop. It follows from the Shift Lemma that $\tree{S}_{\xi+1}$ is wellfounded and that $\Gamma^{\xi+1}: \tree{S}_{\xi+1}\to \tree{T}_{\xi+1}$ is a non-dropping extended tree embedding. It is easy to see that, since $\Gamma^{\xi+1}$ is the appropriate copy tree embedding, (1)-(6) still hold at $\xi+1$. This finishes the successor step.

Now let $\lambda<\lh (\mtree{S})$ be a limit ordinal and let $b=[0,\lambda)_\mtree{S}$. We put $\tree{T}_\lambda= \lim\langle \tree{T}_\xi,\Psi^{\eta,\xi}\,|\,\eta\leq_\mtree{S}\xi\in b\rangle$ if this direct limit is wellfounded. Otherwise we put $\mu=\lambda$ and stop. If this direct limit is wellfounded, we let $\Gamma_\lambda$ be the extended tree embedding guaranteed by Proposition \ref{direct limit prop}. We leave it to the reader to check that our induction hypotheses go through.

We know turn to the \lq\lq moreover" clauses. We've already shown (1). For (2), suppose $\mtree{S}=\mtree{V}(\tree{S},\tree{U})$. We verify by induction that $\Gamma\mtree{S}= \mtree{V}(\tree{T},t^\Gamma_\infty\tree{U})$. Let $\tau_\xi: M^\tree{U}_\xi\to M^{t^\Gamma_\infty\tree{U}}_\xi$ be the copy maps (so $\tau_0=t^\Gamma_\infty$) and $\sigma_\xi:M^\tree{U}_\xi\to M^{\tree{S}_\xi}_{\infty}$, $\sigma_\xi^*:M^{\tau_0\tree{U}}_\xi\to M^{\tree{T}_\xi}_\infty$, the quasi-normalization maps. We check by induction that for all $\xi<\mu$, \begin{enumerate}
    \item $\tau_0 \tree{U}\restrict\xi+1$ is a putative plus tree, 
    \item $\Gamma\mtree{S}\restrict \xi+1 = \mtree{V}(\tree{T},\tau_0 \tree{U}\restrict\xi+1)$,
\item $\sigma^*_\xi\circ\tau_\xi = t^\xi_\infty \circ \sigma_\xi$.
\end{enumerate}

(1) follows from (2) by induction, since the last model of $V(\tree{T}, \tau_0 \tree{U}\restrict\xi+1)$ embeds $M^{\tau_0\tree{U}}_\xi$. These hypotheses go through at limits $\lambda<\mu$ since we're taking direct limits everywhere. So we just deal with the successor case. Let $\xi+1<\mu$ and \[\eta=\tree{U}\pred (\xi+1)= \mtree{S}\pred(\xi+1)=\Gamma\mtree{S}\pred(\xi+1).\]
We have  \[E^{\tau_0\tree{U}}_{\xi}=\tau_\xi(E^\tree{U}_\xi)\]
so applying $\sigma^*_\xi$ and using (3), we get
\[\sigma^*_\xi(E^{\tau_0\tree{U}}_{\xi})=\tau_\xi\circ \sigma_\xi(E^\tree{U}_\xi).\]
But $\sigma_\xi(E^\tree{U}_\xi)= F_\xi$, so we get
\[ \sigma^*_\xi(E^{\tau_0\tree{U}}_{\xi}) = G_\xi.\]
So \[\tree{T}_{\xi+1} = V(\tree{T}_\eta, \tree{T}_\xi, G_\xi) = V(\tree{T}, \tau_0\tree{U}\restrict\xi+1),\]
giving us (2) and (1) at $\xi+1$. All that remains is to verify $\sigma^*_{\xi+1}\circ \tau_{\xi+1}=t^{\xi+1}_\infty\circ \sigma_{\xi_1}$. So we just need to see that the outermost square of the following diagram commutes.

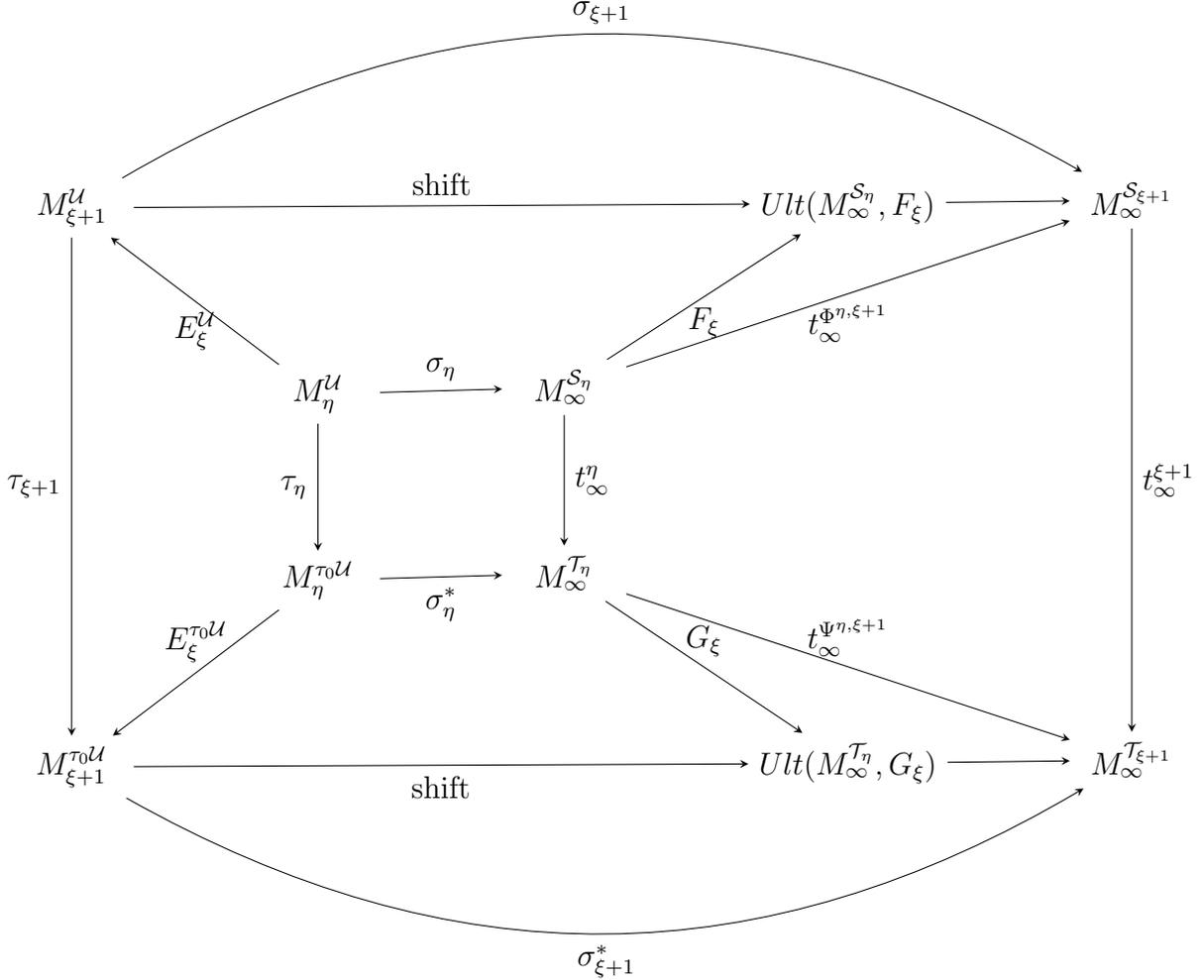
\[\begin{tikzpicture}
  \matrix (m) [matrix of math nodes,row sep=4em,column sep=4em,minimum width=4em]
  {
     M^\tree{U}_{\xi+1} & {} & {} & Ult(M^{\tree{S}_\eta}_\infty, F_\xi) & M^{\tree{S}_{\xi+1}}_\infty \\
    {} & M^\tree{U}_{\eta} & M_\infty^{\tree{S}_\eta}\\
    {} & M^{\tau_0 \tree{U}}_\eta & M_\infty^{\tree{T}_\eta} \\
    M^{\tau_0 \tree{U}}_{\xi+1} & {} & {} & Ult(M_\infty^{\tree{T}_\eta}, G_\xi) & M_\infty^{\tree{T}_{\xi+1}}\\
   };
  \path[-stealth]
    (m-1-1)
        edge
         node [above] {shift}
         (m-1-4)
         edge
         node [left] {$\tau_{\xi+1}$}
         (m-4-1)
         edge [bend left]
         node [above] {$\sigma_{\xi+1}$}
         (m-1-5)
    (m-1-4)
         edge 
         node [right] {}
         (m-1-5)
    (m-1-5) edge node [right] {$t^{\xi+1}_\infty$} (m-4-5)
    (m-4-1)
        edge
         node [below] {shift}
         (m-4-4)
         edge [bend right]
         node [below] {$\sigma^*_{\xi+1}$}
         (m-4-5)
    (m-2-2) edge
             node
            [below] {$E^\tree{U}_\xi$}
        (m-1-1)
    edge
             node [above] {$\sigma_\eta$}
         (m-2-3)
    edge
             node [left] {$\tau_\eta$}
         (m-3-2)  

    (m-2-3) edge
             node
            [below] {$F_\xi$}
        (m-1-4)
    edge
             node [right] {$t^\eta_\infty$}
         (m-3-3)
         edge node [below] {$t^{\Phi^{\eta,\xi+1}}_{\infty}$} (m-1-5)
    (m-3-2)
    edge
             node [below] {$\sigma^*_\eta$}
         (m-3-3)
    edge
             node
            [above] {$E^{\tau_0\tree{U}}_\xi$}
         (m-4-1)
    (m-3-3)
    edge
    node
            [above] {$G_\xi$}
         (m-4-4) 
           edge node [above] {$t^{\Psi^{\eta,\xi+1}}_{\infty}$} (m-4-5)
    (m-4-4)
    edge node {} (m-4-5);
\end{tikzpicture}\]
We have that every region except the rightmost trapezoid commutes by facts about quasi-normalization or the ordinary Shift Lemma. So we just need to see that the rightmost trapezoid commutes, i.e. \[t^{\Psi^{\eta,\xi+1}}_{\infty}\circ t^\eta_\infty =t^{\xi+1}_\infty \circ t^{\Phi^{\eta,\xi+1}}_{\infty}.\] But this follows immediately from the fact that $\Psi^{\eta,\xi+1}\circ\Gamma_{\eta}= \Gamma_{\xi+1}\circ \Phi^{\eta,\xi+1}$. This finishes the successor step of the induction and so establishes part (2) of the \lq\lq moreover" clause.

To finish, we check the \lq\lq moreover" clause part (3), i.e. suppose $\tree{S},\tree{T}$ are by $\Sigma$ and $\mtree{S}$ is by $\Sigma^*$, where $\Sigma$ is some strategy for $M$ with strong hull condensation. We show $\mu=\lh (\mtree{S})$ and $\Gamma\mtree{S}$ is by $\Sigma^*$ simultaneously by induction. 

As long as $\Gamma\mtree{S}\restrict\xi+1$ is by $\Sigma^*$, we know $\xi<\mu$ since the process hasn't broken down. Successor cause no trouble by Remark \ref{shift remark}, so we deal with limits. So we have $\Gamma\mtree{S}\restrict\lambda$ is by $\Sigma^*$ and we need to see that for $b=\Sigma^*(\Gamma\mtree{S}\restrict\lambda)$, $b=[0,\lambda)_\mtree{S}$.  Since we take direct limits of both sides, by Proposition \ref{direct limit prop}, we get a direct limit tree embedding from the last tree of $\mtree{S}\conc b$ to the last tree of $\Gamma( \mtree{S}\restrict\lambda)\conc b$, which is by $\Sigma$. So since $\Sigma$ has strong hull condensation, the last tree of $\mtree{S}\restrict\lambda\conc b$ is by $\Sigma$, hence $b=\Sigma^*(\mtree{S}\restrict\lambda)=[0,\lambda)_\mtree{S}$ by the definition of $\Sigma^*$. Of course, because $b=\Sigma^*(\Gamma\mtree{S}\restrict \lambda)$, $\lambda<\mu$, as well. 

 \qed
\end{proof}

We will also need the analogue of Lemma \ref{shift direct limits}, whose proof we omit. 

\begin{lemma}\label{m-shift direct limits}
Let $\preceq$ be a directed partial order on a set $A$. Suppose we have directed systems of plus trees $\mathcal{C}=\langle \{\tree{S}_\xi\}_{a\in A}, \{\Psi^{a, b}\}_{a\preceq b}\rangle$ and $\mathcal{D}=\langle \{\tree{T}_\xi\}_{a\in A}, \{\Pi^{a, b}\}_{a\preceq b}\rangle$ and extenders $\{F_a\}_{a\in A}$ such that
\begin{enumerate}
    \item for all $a\in A$, $F_a^-$ is on the $M^{\tree{S}_a}_\infty$-sequence and $M^{\tree{S}_a}_\infty|\lh(F_a)\is \dom(t^{\Psi^{a,b}}_\infty)$,
    \item for all $a,b\in A$ such that $a\preceq b$, the Shift Lemma applies to $(\Psi^{a,b}, \Pi^{a,b}, F_a, F_b)$.
\end{enumerate}
For $a,b\in A$ such that $a\preceq b$, let $\Gamma_{a,b}$ be the copy tree embedding associated to $(\Psi^{a,b}, \Pi^{a,b}, F_a, F_b)$. Let $\tree{S}_\infty = \lim\mathcal{C}$, $\tree{T}_\infty=\lim\mathcal{D}$, $\Psi_{a,\infty}:\tree{S}_a\to \tree{S}_\infty$ and $\Pi_{a, \infty}:\tree{T}_a\to \tree{T}_\infty$ the direct limit tree embedding, and $F_\infty$ the common value of $t^{\Psi^{a,\infty}}_\infty(F_a)$.

Let $\tree{V}_a=V(\tree{T}_a, \tree{S}_a, F_a)$, $\Phi^a=\Phi^{V(\tree{T}_a, \tree{S}_a, F_a)}$, $\tree{V}_\infty=\lim \langle \{\tree{V}_a\}_{a\in A}, \{\Gamma^{a,b}\}_{a\preceq b}\rangle$, $\Gamma^{a,\infty}: \tree{V}_a\to \tree{V}_\infty$ the direct limit tree embedding, and $\Phi^\infty: \tree{T}_\infty\to \tree{V}_\infty$ the unique extended tree embedding such that for every $a\in A$, the following diagram commutes.\footnote{such an extended tree embedding is guaranteed by Proposition \ref{direct limit prop}}
\begin{center}
    \begin{tikzcd}
    \tree{T}_a \arrow[r, "\Psi^{a,\infty}"] \arrow[d, "\Phi^a"'] 
    & \tree{T}_\infty \arrow[d,"\Phi^\infty"]\\
    \tree{V}_a \arrow[r, "\Gamma^{a,\infty}"'] 
    & \tree{V}_\infty
    \end{tikzcd}
    \end{center}
for every $\xi<\gamma$ (such an extended tree embedding is guaranteed by Proposition \ref{direct limit prop}).

Then $\tree{V}_\infty=V(\tree{T}_\infty,\tree{S}_\infty,F_\infty)$, $\Phi^\infty= \Phi^{V(\tree{T}_\infty,\tree{S}_\infty, F_\infty)}$, and $\Gamma^{a,\infty}$ is given the copy tree embedding associated to $(\Psi^{a,\infty}, \Pi_{a,\infty}, F_a, F_\infty)$.
\end{lemma}

We end this section with an application of copying which won't be used in the remainder of the paper. It shows that there is some redundancy in the pullback clause (i.e. clause (b)) in the definition of strong hull condensation. Ultimately, nice strategies for plus trees are generated by their action on normal trees and so it seems plausible that this pullback clause is actually redundant.

\begin{proposition} [Siskind]
\label{elementarity prop}
Let $M$ be a premouse and $\Sigma$ an iteration strategy for $M$ which quasi-normalizes well such that every pseudo-hull of a plus tree on $M$ by $\Sigma$ is by $\Sigma$. Let $\tree{S},\tree{T}$ be plus trees on $M$ by $\Sigma$, and $\Phi:\tree{S}\to \tree{T}$ a non-dropping extended tree embedding. Then for any normal tree $\tree{U}$ on $M_\infty^\tree{S}$ by $\Sigma_\tree{S}$, $t^\Phi_\infty \tree{U}$ is by $\Sigma_\tree{T}$.\end{proposition}

\begin{proof}
Suppose $\tree{U}$ is a normal on $M_\infty^\tree{S}$ of limit length $\lambda$ which is by $\Sigma_{\tree{S}}$ such that  $t^\Phi_\infty \tree{U}$ is by $\Sigma_{\tree{T}}$. We need to see that $\Sigma_{\tree{S}}(\tree{U})=\Sigma_\tree{T}(t^\Phi_\infty \tree{U})$. Let $b=\Sigma_\tree{T}(t^\Phi_\infty \tree{U})$.

Then, by our copying result, there is an extended tree embedding $V(\tree{S},\tree{U}\conc b)\to V(\tree{T},\tau\tree{U}\conc b)$ (as these are the last trees of $\mtree{V}(\tree{S},\tree{U}\conc b)$ and $\mtree{V}(\tree{T},t^\Phi_\infty(\tree{U}\conc b))$). Since $V(\tree{T},t^\Phi_\infty\tree{U}\conc b)$ is by $\Sigma$, by quasi-normalizing well, we have that $V(\tree{S},\tree{U}\conc b)$ is by $\Sigma$, since it is a pseudo-hull of a plus tree by $\Sigma$. Hence $b=\Sigma_{\tree{S}}(\tree{U})$.

\qed
\end{proof}
\section{Nice meta-strategies and phalanx comparison}

So far we  have  worked exclusively with meta-strategies
generated by an ordinary strategy, that is, with meta-strategies
of the form $\Sigma^*$, where $\Sigma$ is a strategy acting on  plus trees on some premouse $M$.
In this section we shall start with a plus tree $\tree{S}$ of 
successor length on a
 premouse $M$ and an
 arbitrary meta-strategy $\Sigma$ for
 meta-trees or stacks of meta-trees on $\tree{S}$. 

We shall identify regularity
properties of $\Sigma$ which are the natural analogues of 
strong hull condensation and quasi-normalizing well. We then prove a general comparison theorem 
for pairs of the form $(\tree{S},\Sigma)$ such that $\Sigma$ has these properties.
(See Theorem \ref{main comparison theorem}.) One can think of this as a strategy comparison
theorem for phalanxes of the form $\Phi(\tree{S})$. Not all phalanxes are of this
form, so it is not a truly general strategy comparison theorem for phalanxes.

We then use our tree-phalanx comparison theorem to show characterize those meta-strategies
that are of the form $\Sigma^*$ for some ordinary strategy $\Sigma$. It turns
out every sufficiently nice meta-strategy is of this form. (That is Theorem
\ref{induced strategy theorem}. ) The moral one might draw is that meta-strategies
are not something fundamentally new, but rather a useful way of organizing
constructions and proofs to do with ordinary strategies.

The main step toward Theorem \ref{induced strategy theorem} is Lemma
\ref{induced strategy lemma}, which is a kind of uniqueness theorem
for ordinary iteration strategies $\Sigma$ whose induced meta-strategies
$\Sigma^*$ behave well. 

\subsection{Regularity properties of meta-strategies}

We first isolate the natural notion of tree embedding for meta-trees. 

\begin{definition} Let $\mtree{S}=\langle \tree{S}_\xi, F_\xi, \Phi^{\eta,\xi}\rangle$, $\mtree{T}=\langle \tree{T}_\xi, G_\xi, \Psi^{\eta,\xi}\rangle$, and $\alpha_\xi=\alpha_0(F_\xi,\tree{S}_\xi)$, $\beta_\xi=\beta_0(F_\xi, \tree{S}_\xi)$ and $\alpha^*_\xi=\alpha_0(G_\xi, \tree{T}_\xi)$, $\beta^*_\xi=\beta(G_\xi, \tree{T}_\xi)$.

A \textit{meta-tree embedding} from $\mtree{S}$ to $\mtree{T}$ is a system $\vec{\Delta}=\langle v, u, \{\Gamma_\xi\}_{\xi<\lh\tree{S}}, \{\Delta_\zeta\}_{\zeta+1<\lh(\mtree{S})} \rangle$ such that
\begin{enumerate}
    \item $v:\lh (\mtree{S})\to \lh(\mtree{T})$ is tree-order preserving, $u:\{\eta\mid\eta+1<\lh (\mtree{S})\}\to \lh (\mtree{T})$, $v(\xi)=\sup\{u(\eta)+1\mid \eta<\xi\}$, and for all $\xi+1<\lh(\mtree{S})$, $v(\xi)\leq_{\mtree{T}}u(\xi)$;
    \item For all $\xi$ and $\eta\leq_\mtree{S}\xi$, 
    \begin{enumerate}
        \item $\Gamma_\xi: \tree{S}_\xi\to\tree{T}_{v(\xi)}$ is an extended tree embedding and $\Gamma_0= Id_{\tree{S}_0}$;
        \item $\Psi^{v(\eta),v(\xi)}\circ \Gamma_\eta = \Gamma_\xi\circ \Phi^{\eta,\xi}$,
        \item if $\xi+1<\lh(\mtree{S})$, then $\Delta_\xi= \Psi^{v(\xi),u(\xi)}\circ \Gamma_\xi$ with $M_\infty^{\tree{S}_\xi}|\lh(F_\xi)\is \dom(t_\infty^{\Delta_\xi})$;
    \end{enumerate}
    \item for $\xi+1<\lh (\mtree{S})$, $\eta=\mtree{S}\pred(\xi+1)$, and $\eta^*=\mtree{T}\pred(v(\xi)+1)$,
    \begin{enumerate}
        \item $G_{u(\xi)}=t^{\Delta_\xi}_\infty(F_\xi)$,\footnote{Notice that we haven't built in the option of sending a non-plus-type extender to a plus-type extender extender; this is just because we have no use for such embeddings. Probably one can develop the basics while accommodating this possibility.}
        \item $\eta^*\in[v(\eta),u(\eta)]_\mtree{T}$,
        \item 
        $\Gamma_{\xi+1}\restrict \alpha_\xi +1 \approx \Delta_\xi\restrict \alpha_\xi+1$, and
        \item $u^{\Gamma_{\xi+1}}(\alpha_\xi)=\alpha^*_{u(\xi)}$.
    \end{enumerate}
\end{enumerate}
\end{definition}

One can show, letting $\eta,\xi,\eta^*$ as in (3), that conditions (3) and the commutativity condition (2)(b) imply that the Shift Lemma applies to ($\Delta_\xi$, $\Psi^{v(\eta),\eta^*}\circ \Gamma_\eta$, $F_\xi$, $G_{u(\xi)}$) and that $\Gamma_{\xi+1}$ is the copy tree embedding associated to these objects. At limit ordinals $\lambda$, $\tree{S}_\lambda$ and $\tree{T}_{v(\lambda)}$ are given by direct limits and so (2)(b) implies that $\Gamma_\lambda$ must be the extended tree embedding guaranteed by Proposition \ref{direct limit prop}.

Suppose $\lh (\mtree{S})=\gamma+1$, $\lh (\mtree{T})=\delta+1$, and $v(\gamma)\leq_\mtree{T} \delta$. Then we define the associated \textit{extended meta-tree embedding} by putting $u(\gamma)=\delta$ and $\Delta_\gamma = \Psi^{v(\gamma),\delta}\circ \Gamma_\gamma$.

\begin{proposition}\label{meta-tree agreement}
Let $\vec{\Delta}:\mtree{S}\to \mtree{T}$ be a meta-tree embedding. Then for all $\eta<\xi<\lh (\mtree{S})$,
\begin{enumerate}
    \item $\Delta_\xi\restrict\alpha_\eta+2=\Gamma_\xi\restrict\alpha_\eta+2$ and
    \item $\Delta_\xi\restrict\alpha_\eta+1=\Gamma_\xi\restrict\alpha_\eta+1= \Delta_\eta\restrict\alpha_\eta+1$.
\end{enumerate}

\end{proposition}
\begin{proof}
First, we show for all $\eta<\xi$,
\[\Gamma_\xi\restrict\alpha_\eta+1 = \Delta_\eta\restrict\alpha_\eta+1\]
by induction on $\xi$.
This easily passes through limits, so we check at successors $\xi+1$. But this follows from condition (3) (c) in the definition of meta-tree embedding, since the $\alpha_\eta$'s are increasing.

Now we just need to check that for all $\eta<\xi$, $v^{\Psi^{v(\xi),u(\xi)}}$ is the identity on $v^{\Gamma_\xi}(\alpha_\eta)+1$. 
But this is immediate from the normality condition of meta-trees: letting $\chi+1$ be the successor of $v(\xi)$ along $[v(\xi),u(\xi)]_\mtree{T}$, we have that $\alpha^*_{u(\eta)}<\beta_\chi\leq \alpha^*_\xi$ (using here that $G_{u(\eta)}^*= E^{\tree{T}_{u(\xi)}}_{\alpha^*_{u(\eta)}}$), so $u^{\Psi^{v(\xi),u(\xi)}}$ is the identity on $u^{\Gamma_\xi}(\alpha_\eta)+1$, 
i.e. $v^{\Psi^{v(\xi),u(\xi)}}$ is actually the identity on $v^{\Gamma_\xi}(\alpha_\eta)+2$.

\qed
\end{proof}

\begin{definition}
A meta-strategy $\Sigma$ for $\tree{S}$ has \textit{meta-hull condensation} iff whenever $\mtree{S}$ is by $\Sigma$ and $\vec{\Delta}: \mtree{T}\to \mtree{S}$ is a meta-tree embedding, $\mtree{T}$ is by $\Sigma$.
\end{definition}

\begin{remark}\label{meta-strategy remark}
In the case that $\tree{S}=\{M\}$ (i.e. the trivial tree of length $1$ on a premouse $M$), a meta-strategy is just an iteration strategy for single normal trees. Moreover, meta-hull embeddings are tree embeddings (though not every tree embedding between trees on $M$ is a meta-hull embedding because we only allow mapping non-plus-type extenders to non-plus-type extenders.
\end{remark}

\begin{proposition}\label{nice pulls back}
Let $\Pi:\tree{S}\to \tree{T}$ be a non-dropping extended tree embedding. Let $\Sigma$ be a $(\lambda, \theta)$-strategy for $\tree{T}$. Define $\Sigma^\Pi$ by
\[\mtree{S}\text{ is by }\Sigma^\Pi\text{ iff }\Pi\mtree{S}\text{ is by }\Sigma.\]

Then $\Sigma^\Pi$ is a $(\lambda,\theta)$-strategy for $\tree{S}$. Moreover, if $\Sigma$ has meta-hull condensation, so does $\Sigma^\Pi$.
\end{proposition}

\begin{proof}
It's straightforward to use the copying construction (Theorem \ref{copying}) to get that $\Sigma^\Pi$ is a $(\lambda,\theta)$-strategy for countable meta-trees on $\tree{S}$, so we just handle the \lq\lq moreover" part.

Suppose that $\Sigma$ has meta-hull condensation. Let $\mtree{S}$ be by $\Sigma^\Pi$ and let $\vec{\Delta}=\langle u, v, \Gamma_\xi, \Delta_\xi\rangle:\mtree{T}\to \mtree{S}$ be a meta-tree embedding.

By induction on $\lh (\mtree{T})$ we define a meta-tree embedding $\vec{\Delta}^*:\Pi\mtree{T}\to \Pi\mtree{S}$. This shows that the definition of $\Pi \mtree{T}$ doesn't break down and, since $\Sigma$ has meta-hull condensation, that $\Pi\mtree{T}$ is by $\Sigma$. So $\mtree{T}$ is by $\Sigma^\Pi$, as desired.

Let $\mtree{S}=\langle \tree{S}_\xi, F_\xi, \Phi^{\eta,\xi}\rangle$, $\Pi \mtree{S}= \langle \tree{S}^*_\xi, F^*_\xi, {\Phi^*}^{\eta,\xi}\rangle$, and $\Pi^\mtree{S}_\xi: \tree{S}_\xi\to \tree{S}^*_\xi$ be the copy tree embeddings. Let $\mtree{T}= \langle \tree{T}_\xi,G_\xi, \Psi^{\eta,\xi}\rangle$, $\Pi\mtree{T} = \langle \tree{T}^*_\xi, \tree{T}^*_\xi,G^*_\xi, {\Psi^*}^{\eta,\xi}\rangle$, and $\Pi^\mtree{T}_\xi:\tree{T}_\xi\to \tree{T}^*_\xi$ be the copy tree embeddings.

We define $\vec{\Delta}^*=\langle u,v, \Gamma^*_\xi, \Delta^*_\xi\rangle$ by induction on $\lh (\mtree{T})$ maintaining for all $\xi<\lh (\mtree{T})$,
\begin{enumerate}
    \item $\vec{\Delta}^*\restrict(\Pi\mtree{T}\restrict\xi+1)$ is an extended meta-tree embedding $\Pi\mtree{T}\restrict\xi+1\to \Pi\mtree{S}\restrict u(\xi)+1$ (in particular, $\Pi\mtree{T}\restrict\xi+1$ hasn't broken down) and
    \item $\Pi^\mtree{S}_{v(\xi)}\circ\Gamma_\xi = \Gamma^*_\xi\circ  \Pi^\mtree{T}_\xi$.
    \end{enumerate}
    
Notice that we've demanded $\vec{\Delta}^*$ has the same $u$ and $v$ maps as $\vec{\Delta}$, so we just need to check that we can find the tree embeddings $\Gamma^*_\xi$ (we'll get $\Delta^*_\xi$ for free). Limits cause no trouble, so we'll just handle the successor case.

Suppose $\eta=\mtree{T}\pred (\xi+1)$. Let $\eta^*=\mtree{S}\pred (u(\xi)+1)$ (so $\eta^*\in [v(\eta),u(\eta)]_\mtree{S}$ since $\vec{\Delta}$ is a meta-tree embedding).

We need to see that the process of defining $\tree{T}^*_{\xi+1}$ doesn't break down and that there is a tree embedding $\Gamma^*_{\xi+1}:\tree{T}^*_{\xi+1}\to \tree{S}^*_{u(\xi)+1}$ completing the following diagram.
\[\begin{tikzpicture}
  \matrix (m) [matrix of math nodes,row sep=4em,column sep=4em,minimum width=4em]
  {
    \tree{T}_{\xi+1} & {} & {} & \tree{S}_{u(\xi)+1} \\
    {} & \tree{T}_{\eta} &\tree{S}_{\eta^*}\\
    {} & \tree{T}^*_\eta & \tree{S}^*_{\eta^*} \\
    \tree{T}^*_{\xi+1} & {} & {} & \tree{S}^*_{u(\xi)+1}\\
   };
  \path[-stealth]
    (m-1-1)
        edge
         node [above] {$\Gamma_{\xi+1}$}
         (m-1-4)
         edge
         node [left] {$\Pi^\mtree{T}_{\xi+1}$}
         (m-4-1)
    (m-1-4)
        edge
         node [right] {$\Pi^\mtree{S}_{u(\xi)+1}$}
         (m-4-4)
    (m-4-1)
        edge[dashed,->]
         node [below] {$\Gamma^*_{\xi+1}$}
         (m-4-4)
    (m-2-2) edge
             node
            [below] {$G_\xi$}
        (m-1-1)
    edge
             node [above] {$\Phi^{v(\eta),\eta^*}\circ \Gamma_\eta$}
         (m-2-3)
    edge
             node [left] {$\Pi^\mtree{T}_\eta$}
         (m-3-2)  
    (m-2-3) edge
             node
            [below] {$F_{u(\xi)}$}
        (m-1-4)
    edge
             node [right] {$\Pi^\mtree{S}_{\eta^*}$}
         (m-3-3)
    (m-3-2)
    edge
             node [below] {${\Phi^*}^{v(\eta),\eta^*}\circ\Gamma^*_\eta$}
         (m-3-3)
    edge
             node
            [above] {$G^*_\xi$}
         (m-4-1)
    (m-3-3)
    edge
    node
            [above] {$F^*_{u(\xi)}$}
         (m-4-4) ;
\end{tikzpicture}\]
Of course, $\Gamma^*_{\xi+1}$ will be the copy tree embedding associated to ($\Delta^*_\xi$, ${\Phi^*}^{v(\eta),\eta*}\circ \Gamma^*_\eta$, $G_\xi^*$, $F_{u(\xi)}^*$). Then, since all the trapezoids commute and the outer maps are given by the Shift Lemma, we get that the whole outer square commutes (as in the proof of the Shift Lemma, but now with tree embeddings.) 
We leave it to the reader to check that the Shift Lemma applies. We then set $\Delta_{\xi+1}^*=\Phi^*_{v(\xi+1), u(\xi+1)}\circ \Gamma_\xi^*$, as we must.

At limit $\lambda$, let $b^*= [0,v(\lambda))_\mtree{S}= [0,v(\lambda))_{\Pi\mtree{S}}$. Let $b=[0,\lambda)_\mtree{T}$. We have $v^{-1}[b^*] = b$ since $\vec{\Delta}$ is a meta-tree embedding, so we get a tree embedding $\Gamma^*_\lambda$ from $\tree{T}^*_\lambda=\lim_b (\Pi\mtree{T}\restrict\lambda)$ to $\lim_{b^*} (\Pi\mtree{S}\restrict{\lambda}) = \tree{S}^*_{v(\lambda)}$, as desired. We then continue as in the successor case.
\qed
\end{proof}

We also have the following easy proposition about the existence of meta-strategies with meta-hull condensation.

\begin{proposition}\label{nice meta-strategy existence}
Suppose that $\Sigma$ is a ($\lambda$, $\theta$)-strategy for $M$ with strong hull condensation and $\tree{S}$ is by $\Sigma$. Then the induced meta-strategy $\Sigma^*_\tree{S}$ has meta-hull condensation.
\end{proposition}

\begin{proof}
Let $\vec\Delta= \langle v, u, \Gamma_\xi, \Delta_\xi\rangle :\mtree{T}\to \mtree{S}$ where $\mtree{S}$ is a meta-tree by $\Sigma^*_\tree{S}$. We just want to see that $\mtree{T}$ is by $\Sigma^*_\tree{S}$, so it's enough to show that every tree $\tree{T}_\xi$ is by $\Sigma$ (by the definition of $\Sigma^*$).

For every $\xi<\lh(\mtree{T})$, $\Gamma_\xi$ is a (total) extended tree embedding $\tree{T}_\xi\to\tree{S}_{v(\xi)}$. Since $\tree{S}_{v(\xi)}$ is by $\Sigma$ and $\Sigma$ has strong hull condensation, $\tree{T}_\xi$ is by $\Sigma$. So $\mtree{T}$ is by $\Sigma^*_\tree{S}$, as desired.

\qed
\end{proof}

We'll now see some examples of meta-tree embeddings.

\begin{proposition}\label{meta-tree embedding example prop}

Let $\tree{S}$ be a plus tree of successor length, $\tree{T}$ and $\tree{T}$ normal trees on the last model of $\tree{S}$ and $\Psi: \tree{T}\to \tree{T}^*$ an extended tree embedding. Let $\mu$ greatest such that $\mtree{V}(\tree{S},\tree{T}\restrict \mu+1)$ is wellfounded and $\mu^*$ greatest such that $\mtree{V}(\tree{S},\tree{T}^*\restrict\mu^*+1)$ is wellfounded.

Then $u^\Psi(\mu)\geq \mu^*$ and there is a unique partial meta-tree embedding with maximal domain $\vec \Delta: \mtree{V}(\tree{S},\tree{T}\restrict\mu+1)\to \mtree{V}(\tree{S},\tree{T}^*\restrict\mu^*+1)$ with $u$-map $u^\Psi$.

Moreover, for $\xi\leq\mu$ letting $R_\xi$ be the last model of $V(\tree{S}, \tree{T}\restrict\xi+1)$ and $\sigma_\xi: M^\tree{T}_\xi\to R_\xi$ the quasi-normalization map, for $\xi\leq \mu^*$, letting $R_\xi^*$ be the last model of $V(\tree{S},\tree{T}^*)$ and $\sigma^*_\xi:M^{\tree{T}^*}_\xi\to R_\xi^*$ the  quasi-normalization map, we have the following diagram commutes.
    \begin{center}
    \begin{tikzcd}
    M^{\tree{T}}_\xi \arrow[r, "t^\Psi_\xi"] \arrow[d, "\sigma_\xi"'] 
    & M^{\tree{T}^*}_{u(\xi)} \arrow[d,"\sigma^*_{u(\xi)}"]\\
   R_\xi \arrow[r, "t^{\Delta_\xi}_\infty"'] 
    & R^*_{u(\xi)}
    \end{tikzcd}
    \end{center}
\end{proposition}

\begin{remark}
In particular, if $\Sigma$ is a meta-strategy for $\tree{S}$ with meta-hull condensation and $\mtree{V}(\tree{S},\tree{T}^*)$ is by $\Sigma$, then $\mu^*+1=\lh (\tree{T}^*)$, $\mu+1=\lh(\tree{T})$, $\vec \Delta: \mtree{V}(\tree{S},\tree{T})\to \mtree{V}(\tree{S},\tree{T}^*)$ is a total extended meta-tree embedding, and $\mtree{V}(\tree{S},\tree{T})$ is by $\Sigma$.
\end{remark}

\begin{proof}
Let $\mtree{V}= \langle \tree{V}_\xi, F_\xi, \Phi^{\xi,\eta}\rangle = \mtree{V}(\tree{S},\tree{T})$ and $\mtree{V}^*= \langle \tree{V}^*_\xi, F^*_\xi, {\Phi^*}^{\xi,\eta}\rangle = \mtree{V}(\tree{S},\tree{T}^*)$. We have that $F_\xi= \sigma_\xi(E^\tree{T}_\xi)$, $F^*_\xi= \sigma^*_\xi(E^{\tree{T}^*}_\xi)$.

Our meta-tree embedding $\vec\Delta$ will have $v=v^\Psi$ and $u=u^\Psi$. We just need to see that this works, by induction. Using the notation of the \lq\lq moreover" clause, we have that $R_\xi$ is the last model of $\tree{V}_\xi$ and $R^*_\xi$ is the last model of $\tree{V}^*_\xi$.
Also let $t_{\eta,\xi}$ be the last $t$-map of $\Phi^{\eta,\xi}$ when $\eta\leq_\mtree{V}\xi$ and $t^*_{\eta,\xi}$ the last $t$-map of ${\Phi^*}^{\eta,\xi}$ when $\eta\leq_{\mtree{V}^*}\xi$.

We maintain by induction on $\xi$ that \begin{enumerate}
    \item $\vec\Delta\restrict\xi$ is an extended meta-tree embedding from $\mtree{V}\restrict\xi+1\to \mtree{V}^*\restrict u(\xi)+1$, and
    \item the following diagram commutes. 
    \begin{center}
    \begin{tikzcd}
    M^{\tree{T}}_\xi \arrow[r, "t^\Psi_\xi"] \arrow[d, "\sigma_\xi"'] 
    & M^{\tree{T}^*}_{u(\xi)} \arrow[d,"\sigma^*_{u(\xi)}"]\\
   R_\xi \arrow[r, "t^{\Delta_\xi}_\infty"'] 
    & R^*_{u(\xi)}
    \end{tikzcd}
    \end{center}
\end{enumerate}

Note that the maps in (2) may be partial so we mean that they have the same domain and commute.

We start with the successor case. Let $\eta=\mtree{V}\pred(\xi+1)=\tree{T}\pred(\xi+1)$ and $\eta^*=\mtree{V}^*\pred(u(\xi)+1)=\tree{T}^*\pred(u(\xi)+1)$. Since $\Psi$ is a tree embedding and $\tree{T}^*$ has the same tree order as $\mtree{V}^*$, we get that $\eta^*\in[v(\eta), v(\eta)]_{\mtree{V}^*}$. 
By condition (2) of our induction hypothesis, we get that $t^{\Delta_\xi}_\infty(F_\xi)=F^*{u(\xi)}$. The agreement properties of meta-tree embeddings (Proposition \ref{meta-tree agreement}) imply that the Shift Lemma applies to ($\Delta_\xi$, ${\Phi^*}^{v(\eta),\eta^*}\circ \Gamma_\eta$, $F_\xi$, $G_{u(\xi)}$), so that we may let $\Gamma_{\xi+1}$ be given by the Shift Lemma, as desired. 
We have that $v(\xi+1)=u(\xi)+1\leq_{\mtree{V}^*}u(\xi+1)$, using again that $\tree{T}^*$ has the same tree order as $\mtree{V}^*$ and $\Psi$ is a tree embedding, so we may let $\Delta_{\xi+1}={\Phi^*}^{v(\xi+1),u(\xi+1)}\circ \Gamma_{\xi+1}$, as we must.

This assignment clearly maintains (1), so we just need to see that (2) holds as well. We have the following diagram.

\[\begin{tikzpicture}
  \matrix (m) [matrix of math nodes,row sep=4em,column sep=4em,minimum width=4em]
  {
     M^{\tree{T}}_{\xi+1} & {} & {} & M^{\tree{T^*}}_{v(\xi+1)}  & M^{\tree{T}^*}_{u(\xi+1)} \\
    {} & M^{\tree{T}}_{\eta} & M^{\tree{T}^*}_{\eta^*}\\
    {} & R_\eta & R^*_{\eta^*} \\
    R_{\xi+1}& {} & {} & R^*_{v(\xi+1)} & R^*_{u(\xi+1)}\\
   };
  \path[-stealth]
    (m-1-1)
        edge
         node [above] {$s_{\xi+1}^\Psi$}
         (m-1-4)
         edge
         node [left] {$\sigma_{\xi+1}$}
         (m-4-1)
    (m-1-4)
    
        edge
         node [above] {$\tree{T}^*$}
         (m-1-5)
    
        edge
         node [right] {$\sigma^*_{v(\xi+1)}$}
         (m-4-4)
    (m-1-5)
    
        edge
         node [right] {$\sigma^*_{u(\xi+1)}$}
         (m-4-5)
    (m-4-1)
        edge
         node [below] {$\bar{\tau}_{\xi+1}$}
         (m-4-4)
    (m-4-4)
        edge
         node [below] {$t^*_{v(\xi+1),u(\xi+1)}$}
         (m-4-5)
    (m-2-2) edge
             node
            [below] {$E^{\tree{T}}_\xi$}
        (m-1-1)
    edge
             node [above] {$\hat\imath^{\tree{T}^*}_{v(\eta),\eta^*}\circ s^\Psi_\eta$}
         (m-2-3)
    edge
             node [left] {$\sigma_\eta$}
         (m-3-2)  
    (m-2-3) edge
             node
            [below] {$E^{\tree{T}^*}_{u(\xi)}$}
        (m-1-4)
    edge
             node [right] {$\sigma^*_{\eta^*}$}
         (m-3-3)
    (m-3-2)
    edge
             node [below] {$t^*_{V(\eta),\eta^*}\circ t^{\Gamma_\eta}\infty$}
         (m-3-3)
    edge
             node
            [above] {$t_{\eta,\xi+1}$}
         (m-4-1)
    (m-3-3)
    edge
    node
            [above] {$t^*_{\eta^*,v(\xi+1)}$}
         (m-4-4) ;
\end{tikzpicture}\]

We leave it to the reader to check that everything commutes. This relies quite heavily on properties of the quasi-normalization maps, which can be found in \cite{nitcis}. For example, the leftmost trapezoid commutes for free because of how we define the map $\sigma_{\xi+1}$.
\qed
    \end{proof}

We now turn to another source of meta-tree embeddings: the analogue of embedding normalization for meta-trees. We start with the one-step case.

Given meta-trees $\mtree{S}$ and $\mtree{T}$ of successor lengths, $H$ an extender such that $H^-$ is on the last model of the last tree of $\mtree{T}$, we want to define a meta-tree $\mtree{V}=\mtree{V}(\mtree{S}, \mtree{T}, H)$ and an extended meta-tree embedding $\vec{\Delta}=\vec{\Delta}^{\mtree{V}(\mtree{S},\mtree{T},H)}$ from $\mtree{S}$ into $\mtree{V}$. Moreover, we will have that the last tree of $\mtree{V}$ is just $V(\tree{S}_\infty, \tree{T}_{\infty}, H)$ and $\Delta_\infty$, the last $\Delta$-map of $\vec{\Delta}$, is $\Phi^{V(\tree{S}_\infty, \tree{T}_\infty, H)}$, so that, in particular, this is really producing an analogue of full normalization.

For a meta-tree $\mtree{S}=\langle \tree{S}_\xi, F_\xi\rangle$ of successor length and $H$ on the sequence of the last model of $\mtree{S}$, we define 
\begin{align*}
    a(\mtree{S}, H)&= \text{least $\xi$ such that $H^-$ is on the sequence of the last model of $\tree{S}_\xi$}\\
    b(\mtree{T}, H) &= \text{ least $\xi$ such that $\xi+1=\lh(\mtree{S})$ or $\xi+1<\lh(\mtree{S})$ and $\crit (H)<\hat\lambda(F_\xi).$}
\end{align*}

\noindent Note that $a(\mtree{S},H)$ is also the  least $\xi$ such that $\xi+1=\lh(\mtree{S})$ or $\xi+1<\lh(\mtree{S})$ and $\alpha(\tree{S}_\xi, F_\xi)\leq \alpha(\tree{S}_\infty, G)$.

Let $\mtree{S} = \langle \tree{S}_\xi, F_\xi, \Phi_{\eta,\xi}\rangle$, $\mtree{T}= \langle \tree{T}_\xi, G_\xi,\Psi_{\eta,\xi}\rangle$,  and $H$ on the sequence of the last model of $\mtree{T}$.
Let $a=a(\mtree{T}, H)$ and $b=b(\mtree{T}, H)$. Suppose that $\mtree{S}\restrict b+1= \mtree{T}\restrict b+1$ and if $b+1<\lh(\mtree{S})$, then $\dom(H)\is \lh(F_b)$.

We define $\mtree{V} = \langle \tree{V}_\xi, H_\zeta,\ldots\rangle$, ordinal valued maps $u,v$, extended tree embeddings $\Gamma_\xi: \tree{S}_\xi \to \tree{V}_{v(\xi)}$, and partial extended tree embeddings $\Delta_\xi: \tree{S}_\xi \to \tree{V}_{u(\xi)}$ inductively as follows.

First, we let
$\mtree{V}\restrict a+1= \mtree{T}\restrict a+1$.
We then put
\[u(\xi)=\begin{cases} \xi &\text{if } \xi<b\\
a+1+ (\xi-b) &\text{if }\xi\geq b
\end{cases},\]
and $v(\xi)= \sup\{u(\eta)+1\,|\,\eta<\xi\}$.

We maintain by induction on $\xi\geq b$ that 
\begin{enumerate}
    \item $\vec{\Delta}\restrict\xi+1=\langle u\restrict\xi+1, v\restrict\xi+1, \{\Gamma_\eta\}_{\eta\leq \xi}, \{\Delta_\eta\}_{\eta\leq\xi}\rangle$ is an extended meta-tree embedding from $\mtree{S}\restrict\xi+1$ into $\mtree{V}\restrict u(\xi)+1$,
    \item $\tree{V}_{u(\xi)} = V(\tree{S}_\xi, \tree{T}_a, H)$ and $\Delta_\xi = \Phi^{V(\tree{S}_\xi, \tree{T}_a, H)}$.
\end{enumerate}
There is a lot built into (1); for example, we must have that $\Gamma_\xi = \Delta_\xi$ for $\xi>b$ by our choice of $u$.

We just need to see that this works. It's easy to see that we get the base case $\xi=b$ for free by setting $H_a=H$, so we just need to deal with successors and limits $\xi>b$.

We start with successors.

\paragraph{Successor case.} Suppose our induction hypotheses hold up to $\xi\geq b$. \\

Let $\eta=\mtree{S}\pred (\xi+1)$. There are two subcases depending on the critical point of $F_\xi$.

\paragraph{Subcase 1.}  $\crit (F_\xi)<\crit (H)$.\\

In this case $\eta\leq b$, $\tree{V}_\eta=\tree{S}_\eta$, $\Gamma_\eta=Id$, and $\crit (H_{u(\xi)})=\crit (t_\infty^{\Delta_\xi})$. Note that $t_\infty^{\Delta_\xi}$ has critical point $\crit (H)$, by our induction hypothesis (2).

Since $\mtree{V}\restrict b+1=\mtree{S}\restrict b+1$, we must put $\eta= \mtree{V}\pred (u(\xi)+1)$, as dictated by normality. We let $\tree{V}_{u(\xi)+1}= V(\tree{V}_\eta, \tree{V}_{u(\xi)}, H_{u(\xi)})$ and, following the definition of meta-tree embedding, we let $\Gamma_{\xi+1}$ be the copy tree embedding associated to ($\Delta_\xi$, $\Gamma_\eta = Id$, $F_\xi$, $H_{u(\xi)}$). We leave it to the reader to check that our induction hypothesis (1) guarantees that the Shift Lemma applies. To see (2), we show that $\tree{V}_{u(\xi)+1}= V(\tree{S}_{\xi+1}, \tree{T}_a, H)$ and $\Gamma_{\xi+1}= \Phi^{V(\tree{S}_{\xi+1}, \tree{T}_a, H)}$ simultaneously.

Let $\bar{\tree{V}}= V(\tree{S}_{\xi+1}, \tree{T}_a, H)$ and $\bar \Gamma = \Phi^{V(\tree{S}_{\xi+1}, \tree{T}_a, H)}$. Let $\bar\alpha_\xi = \alpha_0(F_\xi, \tree{S}_\xi)$ and $\alpha_\xi = \alpha_0(H_{\xi}, \tree{V}_\xi)$. Let $\bar\Phi = \Phi^{V(\tree{S}_\eta, \tree{S}_\xi, F_\xi)}$ and $\Phi = \Phi^{V(\tree{V}_\eta,\tree{V}_{u(\xi)}, H_{u(\xi)})}$. We show $\bar\Gamma= \Gamma_{\xi+1}$ by showing it satisfies the conditions which uniquely determine $\Gamma_{\xi+1}$ in the conclusion of the Shift Lemma. First we show that $\bar\Gamma\restrict \bar\alpha +1\approx\Gamma_\xi\restrict\bar\alpha+1$ and $u^{\bar\Gamma}(\bar\alpha)=\alpha$, which guarantees that $\bar\Gamma\restrict\bar\alpha+2 \approx \Gamma_{\xi+1}\restrict\bar\alpha+2$. From here, we'll show simultaneously by induction on $\zeta<\lh (\tree{S}_\eta)$ that $\bar\Gamma \circ \bar\Phi\restrict \zeta+1\approx \Phi\restrict\zeta+1$ and $\bar{\tree{V}}\restrict u^{\Phi}(\zeta)+1 = \tree{V}_{u(\xi)+1}\restrict u^{\Phi}(\zeta)+1$, which establishes $\bar\Gamma =\Gamma_{\xi+1}$. Note that we're using in several places that $\Gamma_\eta= Id$.

We have that $\tree{S}_{\xi+1}\restrict\bar\alpha+1=\tree{S}_\xi\restrict\bar\alpha+1$, so that $\bar {\tree{V}}\restrict \alpha+1= \tree{V}_{u(\xi)}\restrict\alpha+1$ and $\bar\Gamma\restrict\bar\alpha+1\approx\Gamma_\xi\restrict\bar\alpha+1$, as both are just given by one-step quasi-normalization $V(\tree{S}_\xi\restrict\bar\alpha+1, \tree{S}_a, H)$.
Moreover, $u^{\bar\Gamma}(\bar\alpha)= u^{\Gamma_\xi}(\bar\alpha)$, since this is just the $u$-map of embedding normalization by $H$.
We have that $F_\xi=E^{\tree{S}_{\xi+1}}_{\bar\alpha}$ and so $t^{\Gamma_\xi}_{\bar\alpha}=t^{\bar\Gamma}_{\bar\alpha}$ agrees with $t^{\Delta_\xi}_\infty$ on  $F_\xi$. 
It follows that $H_{u(\xi)}= E^{\bar{\tree{V}}}_{u^{\bar\Gamma}(\bar\alpha)}$, so that $u^{\bar\Gamma}(\bar\alpha)=\alpha$. 
This establishes $\bar\Gamma\restrict\bar\alpha+2=\Gamma_{\xi+1}\restrict\bar\alpha+2$.

For the rest, we show by induction on $\zeta<\lh (\tree{S}_\eta)$ that $\bar\Gamma\circ \bar\Phi\restrict\zeta+1 \approx \Phi\restrict\zeta+1$, $u^{\bar\Gamma\circ \bar\Phi}(\zeta)=u^\Phi(\zeta)$, and $\bar{\tree{V}}\restrict u^\Phi(\zeta)+1=\tree{V}_{u(\xi)+1}\restrict u^\Phi(\zeta)+1$.
Let $\bar\beta= \beta(F_\xi,\tree{S}_\xi)$ and $\beta=\beta(H_{u(\xi)}, \tree{V}_{u(\xi)})$. 
We have that $\bar\beta=\beta$ by our case hypothesis and so $v^{\bar\Phi}\restrict \beta +1= v^{\Phi}\restrict\beta+1=id$. 
Also by our case hypothesis we have $v^{\bar\Gamma}\restrict\beta+1=id$, so $\bar\Gamma \circ \bar\Phi\restrict\beta+1 = \Phi\restrict\beta+1= Id_{\tree{S}_\eta\restrict\beta+1}$. 
We have $u^{\bar\Phi}(\beta)=\bar\alpha+1$ so $u^{\bar\Gamma\circ \bar\Phi}(\beta)= u^{\bar\Gamma}(\bar\alpha+1)=\alpha+1$ (since $u^{\bar\Gamma}$ agrees with $v^{\bar\Gamma}$ on $\bar\alpha+1$ and $u^{\bar\Gamma}(\bar\alpha)= \alpha$, as already established). We also clearly have $u^\Phi(\beta) = \alpha+1$ (since $\Phi$ it is just one-step normalization by $H_u(\xi)$), so $u^{\bar\Gamma\circ \bar\Phi}(\beta)= u^\Phi(\gamma)$. Moreover, we already established $\bar{\tree{V}}\restrict\alpha+1 = \tree{V}_{u(\xi)}\restrict\alpha+1 = \tree{V}_{u(\xi)+1}\restrict\alpha+1$, as desired.

Suppose now $\zeta\geq\beta$ and our induction hypothesis holds up to $\zeta$. We have that the exit extenders $E^{\bar{\tree{V}}}_{u^\Phi(\zeta)}$ and $E^{\tree{V}_{\xi+1}}_{u^\Phi(\zeta)}$ are equal since they are both images of $E^{\tree{S}_\eta}_\zeta$ under the same $t$-map (our induction hypothesis implies the $\zeta$th $t$-maps of $\bar\Gamma\circ\bar\Phi$ and $\Phi$ are the same).
It follows that $\bar{\tree{V}}$ and $\tree{V}_{\xi+1}$ agree up to $v^\Phi(\zeta+1)=u^\Phi(\zeta)+1$ and $\bar\Gamma\circ\bar\Phi$ and $\Phi$ agree up to $\zeta+2$. 
Since $\zeta+1\geq \beta$, we get that $u^{\bar\Phi}$ and $u^{\Phi}$ agree with their corresponding $v$-maps on $\zeta+1$. 
Moreover, since $u^{\Phi}(\zeta)>\alpha>\alpha_0(H,\tree{T}_a)$ (using here that $\xi\geq b$ so $u(\xi)>a$), we have that $u^{\bar\Gamma}$ agrees with $v^{\bar\Gamma}$ above $u^{\bar\Phi}(\zeta)$. So $u^{\bar\Gamma\circ \bar\Phi}(\zeta+1)= v^{\bar\Gamma\circ \bar \Phi}(\zeta+1)= v^\Phi(\zeta+1)=u^\Phi(\zeta+1)$ and the trees agree this far, too. 

Both trees must pick the same branches at limit $\alpha$ and at limits $\lambda= v^\Phi(\bar\lambda)= v^{\bar\Gamma\circ\bar\Phi}(\bar\lambda)$, both trees must pick the image of $[0,\bar\lambda)_{\tree{S}_\eta}$ under the same map. So as long as we don't reach illfounded models, this agreement continues through limits. This finishes Subcase 1.

\paragraph{Subcase 2.} $\crit (F_\xi)\geq \crit (H)$.\\

In this case $\eta\geq b$ and $\crit (H_{u(\xi)}) \geq \hat\lambda (G)$. It follows that $a+1\leq u(\eta)= \mtree{V}\pred(u(\xi)+1)$. Where we used $\Gamma_\eta=Id$ above, we must now use $\Delta_\eta$, which is, in particular, \textit{not} the identity. 
We let $\tree{V}_{u(\xi)+1}= V(\tree{V}_{u(\eta)}, \tree{V}_{u(\xi)}, H_{u(\xi)})$ and $\Gamma_{\xi+1}$ the copy tree embedding associated to ($\Delta_\xi$, $\Delta_\eta$, $F_\xi$, $H_{u(\xi)}$). As in the previous case we let $\bar{\tree{V}} = V(\tree{S}_{\xi+1},\tree{T}_a, H)$ and $\bar\Gamma=\Phi^{V(\tree{S}_{\xi+1},\tree{T}_a,H)}$. 
Again we must show that $\bar\Gamma$ satisfies the properties in the conclusion of the Shift Lemma which uniquely determine $\Gamma_{\xi+1}$. 

Getting $\bar\Gamma\restrict\bar\alpha+2=\Gamma_{\xi+1}\restrict\bar\alpha+2$ is the same as before. 
For the rest, we now want to see that $\bar\Gamma\circ \bar\Phi = \Phi\circ \Delta_\eta$, where $\bar\Phi$ is as before and $\Phi=\Phi^{V(\tree{V}_{u(\eta)},\tree{V}_{u(\xi)},H_{u(\xi)})}$. The argument here is pretty much the same as in the previous case: we get agreement up to $\bar\beta+2$ for free and then use that the remainder of our trees and tree embeddings are given by images of $\tree{S}_\eta$ under the same maps. We leave the details to the reader. This finishes Subcase 2 and the successor case.

\paragraph{Limit case.} Suppose our induction holds below $\lambda>b$.\\

We must let $\Gamma_\lambda$ be the unique extemded tree embedding from $\tree{S}_\lambda=\lim_{[0,\lambda)_\mtree{S}} (\mtree{S}\restrict\lambda)$ to $\tree{V}_\lambda= \lim_{[0, v(\lambda))_{\mtree{V}}} (\mtree{V}\restrict v(\lambda))$ which commutes with the rest of our embeddings. By Lemma \ref{m-shift direct limits}, we must have that $\Gamma_\lambda = \Phi^{V(\tree{S}_\lambda, \tree{T}_a, H)}$. This finishes the one-step case.\\

In our main comparison result, we will use analogues of the Factor Lemma and Shift Lemma for meta-tree embeddings, though we don't need to go on to develop analogues of the factorization analysis or copying construction.

\begin{lemma}[Factor Lemma]
Let $\vec{\Psi}: \mtree{S}\to \mtree{T}$ be an (extended) meta-tree embedding such that $\vec{\Psi}\neq Id$. 
Let $b=\crit(u^{\vec{\Psi}})$ and $a+1$ be the successor of $v^{\vec{\Psi}}(b)=b$ in $(v^{\vec{\Psi}}(b), u^{\vec{\Psi}}(b)]_\mtree{T}$. Suppose that $\dom(F_a^\mtree{T})\isneq M_\infty^{\tree{S}_b}|\lh(F_b^\mtree{S})$. Then $\mtree{V}(\mtree{S}, \mtree{T}\restrict a+1, F_a^\mtree{T})$ is defined and wellfounded 
and there is a unique (extended) meta-tree embedding $\vec{\Pi}:\mtree{V}(\mtree{S},\mtree{T}\restrict a+1, F_a^\mtree{T})\to \mtree{T}$ such that $u^{\vec{\Pi}}\restrict a+1=id$ and $\vec{\Psi}=\vec{\Pi}\circ \vec{\Delta}^{\mtree{V}(\mtree{S},\mtree{T}\restrict a+1, F_a^\mtree{T})}$.
\end{lemma}

\begin{proof}
Let $\vec{\Psi} = \langle u^{\vec{\Psi}}, v^{\vec{\Psi}}, \{\Phi_\xi\}, \{\Psi_\xi\}\rangle$. We'll define $\vec{\Pi} = \langle u^{\vec{\Pi}}, v^{\vec{\Pi}}, \{ \Lambda_\xi\}, \{\Pi_\xi\}\rangle$ our desired meta-tree embedding by induction. Note that our hypotheses immediately gives that $\mtree{V}(\mtree{S},\mtree{T}\restrict a+1, F_a^\mtree{T})$ is defined.

Let $\mtree{V}=\mtree{V}(\mtree{S},\mtree{T}\restrict a+1, F_a^\mtree{T})$ and $\vec{\Delta} =  \vec{\Delta}^{\mtree{W}(\tree{S},\tree{T}\restrict a+1, F_a^\mtree{T})}=\langle u^{\vec{\Delta}},v^{\vec{\Delta}}, \{\Gamma_\xi\}, \{\Delta_\xi\}\rangle$.

We can now define our meta-tree embedding $\vec{\Pi}$. To start, we define the $u$-map of $\vec{\Xi}$ like we did in the the tree embedding Factor Lemma: 
\begin{equation*}
    u^{\vec{\Pi}}(\xi)=\begin{cases} \xi &\text{ if } \xi<a+1\\
    u^{\vec{\Psi}}\circ(u^{\vec{\Delta}})^{-1}(\xi)&\text{ if } \xi\geq a+1.
    \end{cases}
\end{equation*}
So we put $\Pi_\xi = \Lambda_\xi = Id$ for $\xi<a+1$ and $\Lambda_{a+1} = Id$. All $\xi\geq a+1$ are in the range of $u^{\vec{\Delta}}$, so we find the remaining $\Lambda_\xi$ and $\Pi_\xi$ as follows.

First
let $\alpha_0=\alpha_0(\tree{T}_a, F_a^\mtree{T})$ and $\beta=\beta(\tree{T}_a, F_a^\mtree{T})$. We have that  $\tree{S}_b=\tree{T}_b$, $\Phi_b= Id$, and $\Psi_{b}=\Phi_{b, u^{\vec{\Psi}}(b)}^\mtree{T}$. So $\Psi_b$ is an inflationary tree  embedding with $\crit(u^{\Psi_b})=\beta$, $u^{\Psi_b}(\beta)=\alpha_0+1$, and first factor $F_a^\mtree{T}$. In particular, as either $\beta+1=\lh(\tree{S}_b)$ or $\beta+1<\lh(\tree{S}_b)=\lh(\tree{T}_b)$, $E_\beta^{\tree{S}_b}=E_\beta^{\tree{T}}$, and $\dom(F_a^\mtree{T})\isneq M_\beta^{\tree{S}_b}|\lh(E_\beta^{\tree{S}_b})$. Now let $\xi>b$. 
The agreement properties of meta-tree embeddings give that $\crit (u^{\Phi_\xi})=\beta$ and $\alpha_0+1$ is the successor of $\beta$ in $(\beta, u^{\Phi_\xi}(\beta)]_{\tree{T}_{v^{\vec{\Psi}}(\xi)}}$. 
Moreover, $E_{\alpha_0}^{\tree{T}_{v^{\vec{\Psi}}(\xi)}}=F_a^\mtree{T}$ and either $E_\beta^{\tree{S}_\xi}=E_\beta^{\tree{S}_b}$ or else $\beta=\alpha_0(\tree{S}_b, F_b^\mtree{S})$ and $E_\beta^{\tree{S}_\xi}=F_b^\mtree{S}$. 
In either case, $\dom(F_a^\mtree{T})\isneq M_\beta^{\tree{S}_\xi}|\lh(E_\beta^{\tree{S}_\xi})$, so that the tree embedding Factor Lemma applies to $\Phi_\xi$, i.e. there is an extended tree embedding $\Lambda_{u^{\vec{\Delta}}(\xi)}:\tree{V}_{u^{\vec{\Delta}}(\xi)}\to \tree{T}_{u^{\vec{\Psi}}(\xi)}$ such that $\Lambda_{u^{\vec{\Delta}}(\xi)}\circ \Delta_\xi = \Phi_\xi$.

Finally, for $\xi\geq b$ we define $\Pi_{u^{\vec{\Delta}}(\xi)} = \Phi^\mtree{T}_{v^{\vec{\Pi}}\circ u^{\vec{\Delta}}(\xi), u^{\vec{\Psi}}(\xi)}\circ \Lambda_{u^{\vec{\Delta}}(\xi)}$, as we must.

To finish, one needs to check that this really is a meta-tree embedding. This is straightforward, by induction. We leave the details to the reader.

\qed 
\end{proof}

Because we demanded that meta-trees are normal, we have used an analogue of embedding normalization instead of the quasi-normalization in the meta-tree normalization procedure just introduced. We mentioned above that one \textit{cannot} prove that $\alpha(\tree{T},G)\in [v(\alpha(\tree{S},F)), u(\alpha(\tree{S},F)]_\tree{T}$ for an arbitrary extended tree embedding $\Phi:\tree{S}\to \tree{T}$ and extenders $F$, $G$ such that $F^-$ is on the $M^\tree{S}_\infty$-sequence and $G=t_\infty(F)$. We needed to move to the $\alpha_0$'s to have this property. Similarly, we cannot prove $a(\mtree{T},G)\in [v(\alpha(\mtree{S},F)), u(\alpha(\mtree{S},F)]_\mtree{T}$ for an arbitrary meta-tree embedding $\vec{\Delta}:\mtree{S}\to \mtree{T}$ and extenders $F$, $G$ such that $F^-$ is on the $M^\mtree{S}_\infty$-sequence and $G=t^{\Delta_\infty}_\infty(F)$. Still, this condition will be met in our application of the meta-tree embedding Shift Lemma, so we just add it as an additional hypothesis, (2), below.

\begin{lemma}[Shift Lemma]
Let  $\vec{\Psi}=\langle u^{\vec{\Psi}}, v^{\vec{\Psi}}, \{\Phi_\xi\}_{\xi<\lh(\mtree{S})},\{\Psi_\xi\}_{\xi<\lh(\mtree{S})}\rangle:\mtree{S}\to \mtree{T}$ and $\vec{\Pi}=\langle u^{\vec{\Pi}}, v^{\vec{\Pi}}, \{\Lambda_\xi\}_{\xi<\lh(\mtree{U})},\{\Pi_\xi\}_{\xi<\lh(\mtree{U})}\rangle:\mtree{U}\to \mtree{V}$ extended meta-tree embeddings and $F,G$ extenders be such that $F^-$ on the $M_\infty^\mtree{S}$-sequence and $G$ such that $G^-$ is on the $M_\infty^\mtree{T}$-sequence. Let $a=a(\mtree{S}, F)$, $b=b(\mtree{S},F)$, $a^*=a(\mtree{T},G)$, and $b^*=b(\mtree{T},G)$.

Suppose
\begin{enumerate}
    \item $M_\infty^\mtree{S}|\lh(F)\is\dom(t^{\Psi_\infty}_\infty)$ and $G=t_\infty^{\Psi_\infty}(F)$,
    \item $a^*\in [v^{\vec{\Psi}}(a), u^{\vec{\Psi}}(a)]_\mtree{T}$ and $G=t_\infty^{{\Phi^\mtree{T}_{v^{\vec{\Psi}}(a),a^*}}\circ \Phi_a}(F)$,
    \item $\vec{\Psi}\restrict b+1\approx \vec{\Pi}\restrict b+1$,
    \item $\mtree{T}\restrict b^*+1=\mtree{V}\restrict b^*+1$,
    \item $b^*\in [v^{\vec{\Pi}}(b), u^{\vec{\Pi}}(b)]_\mtree{V}$ and $t^{\Pi_\infty}_\infty\restrict\dom(F)\cup\{F\}=t^{\Phi^\mtree{V}_{v^{\vec{\Pi}}(b), b^*}\circ \Lambda_b}_\infty \restrict\dom(F)\cup\{F\}$,
    \item if $b+1<\lh(\mtree{U})$, then $\dom(F)\isneq M_\infty^\mtree{U}|\lh(F_b^\mtree{U})$, and
    \item if $b^*+1<\lh(\mtree{V})$, $\dom(G)\isneq M_\infty^\mtree{V}|\lh(F_{b^*}^\mtree{V})$.
\end{enumerate}
Then $\mtree{V}(\mtree{U},\mtree{S},F)$ and $\mtree{V}(\mtree{V},\mtree{T},G)$ are defined and, letting $\mu$ the greatest ordinal such that 
    $\mtree{V}(\mtree{U},\mtree{S},F)\restrict \mu$ is wellfounded and $\mu^*$ the greatest ordinal
 such that $\mtree{V}(\mtree{V},\mtree{T}, G)\restrict\mu^*$ is wellfounded, there is a unique partial meta-tree embedding
$\vec{\Delta}=\langle v^{\vec{\Delta}}, u^{\vec{\Delta}}, \{\Gamma_\xi\}, \{\Delta_\xi\}\rangle: \mtree{V}(\mtree{U},\mtree{S},F)\restrict \mu \to \mtree{V}(\mtree{V},\mtree{T}, G)\restrict\mu^*$ with maximal domain such that
 \begin{enumerate}
        \item $\vec{\Delta}\restrict a+1\approx \vec{\Psi} \restrict a+1$,
        \item $u^{\vec{\Delta}}(a)=a^*$, and
        \item $\vec{\Delta}\circ \vec{\Delta}^{\mtree{V}(\mtree{U},\mtree{S},F)} =\vec{\Delta}^{\mtree{V}(\mtree{V},\mtree{T},G)}\circ \vec{\Pi}$ (on their common domain).
    \end{enumerate}
 Moreover, if $\mtree{V}(\mtree{V},\mtree{T},G)$ is wellfounded, then
 $\mtree{V}(\mtree{U},\mtree{S},F)$ is wellfounded 
  and $\vec{\Delta}$ is a (total) extended meta-tree embedding 
  from $\mtree{V}(\mtree{U},\mtree{S},F)$ into $\mtree{V}(\mtree{V},\mtree{T},G)$. If $\mtree{V}(\mtree{V},\mtree{T},G)$ is wellfounded and also $\vec{\Pi}$ is non-dropping, then $\vec{\Delta}$ is a non-dropping extended meta-tree embedding.

\end{lemma}

We omit the proof because one can essentially copy the proof of the tree embedding Shift Lemma, above, using that lemma where we used the ordinary Shift Lemma everywhere in that proof. Alternatively, one can give a somewhat simpler proof by using that the meta-tree embedding normalization coincides with full normalization, in the sense discussed above.\\

For normalizing a stack of meta-trees $\langle \mtree{S},\mtree{T}\rangle$, we'll need to talk about direct limits of systems of meta-trees under meta-tree embeddings. Our analysis of direct limits of trees under extended tree embeddings from \S1.2 carries over to meta-trees under meta-tree embeddings in the obvious way.

\begin{definition}
A \textit{directed system of meta-trees} is a system $\mathbb{D}=\langle\{\mtree{T}_a\}_{a\in A},\{\vec{\Delta}^{a,b}\}_{a\preceq b}\rangle$, where $\preceq$ is a directed partial order on some set $A$ and \begin{enumerate}
    \item[(a)] for any $a\in A$, $\mtree{T}_a$ is a simple meta-tree of successor length,
    \item[(b)] for any $a,b\in A$ with $a\prec b$, $\vec{\Delta}^{a,b}: \mtree{T}_a\to \mtree{T}_b$ is an extended meta-tree embedding,
    \item[(c)] for any $a,b,c\in A$ such that $a\preceq b\preceq c$, $\vec{\Delta}^{a,c}= \vec{\Delta}^{b,c}\circ \vec{\Delta}^{a,b}$.
\end{enumerate}
\end{definition}

We define $\lim \mathbb{D}$ just as before, except we replace the parts of the tree embeddings with the corresponding parts of our meta-tree embeddings, e.g. we form $u$-threads $x$ using the $u_{a,b}$ and form trees $\tree{T}_x$ by taking direct limits along $\Delta^{a,b}_{x(a)}$ instead of the $t$-maps, provided that enough of these are total. We also define systems $\vec{\Pi}^a$ which, when the direct limit is wellfounded, are extended meta-tree embeddings from $\mtree{T}_a$ into $\lim \mathbb{D}$.

We say $\lim \mathbb{D}$ is \textit{wellfounded} if all the $\tree{T}_x$ are defined and are actually plus trees, the order on $u$-threads is wellfounded, and the direct limit object is a meta-tree. Like in the case of direct limits of trees under extended tree embeddings, the last two conditions follow from the first. 

We get that this construction really identifies the direct limit in the category of meta-trees of successor lengths and extended meta-tree embeddings between them, i.e. we have

\begin{proposition}\label{meta-tree direct limit prop}
Let $\mathbb{D}= \langle \mtree{T}_a, \vec{\Delta}^{a,b}, \preceq\rangle$ be a directed system of meta-trees, where $\preceq$ has field $A$.

Suppose there is a meta-tree $\mtree{S}$ and for all $a\in A$ meta-tree embeddings $\vec{\Psi}^a:\mtree{T}_a\to \mtree{S}$  such that whenever $a\preceq b$, $\vec{\Psi}^b=\vec{\Delta}^{a,b}\circ \vec{\Psi}^a$. 

Then the direct limit $\lim\mathbb{D}$ is wellfounded and there is a unique tree embedding $\vec{\Psi}:\lim \mathbb{D}\to \mtree{S}$ such that $\vec{\Psi}^a= \vec{\Psi}\circ \vec{\Pi}^a$ for all $a\in A$.
\end{proposition}

Now, given a stack of meta-trees $\langle \mtree{S},\mtree{T}\rangle$, with $\mtree{S}= \langle \tree{S}_\xi, F_\xi, \Phi_{\eta,\xi}\rangle$, $\mtree{T}=\langle \tree{T}_\xi, G_\xi, \Psi_{\eta,\xi}\rangle$, we define $\mtree{V}(\mtree{S},\mtree{T})$ as the last meta-tree in a sequence of meta-trees $\mtree{V}^\xi=\langle \tree{V}^\xi_\zeta, F^\xi_\zeta\rangle$ (each of successor length). We also define (partial) extended meta-tree embeddings $\vec{\Delta}^{\eta,\xi}:\mtree{V}^\eta\to \mtree{V}^\xi$ for $\eta\leq_\mtree{T} \xi$. Of course, our construction only makes sense as long as we never reach illfounded models, in which case we'll say that $\mtree{V}(\mtree{S},\mtree{T})$ is  \textit{wellfounded}.

We maintain the following by induction.
 \begin{enumerate}
    \item $\tree{V}^\xi_{\infty}=\tree{T}_\xi$,
    \item for $\eta\leq \xi$, $\mtree{V}^\eta\restrict a(\mtree{V}^\eta, G_\eta)+1 = \mtree{V}^\xi\restrict a(\mtree{V}^\eta, G_\eta)+1$,
    \item for $\eta<\xi$, $G_\eta= F^\xi_{a(\mtree{V}^\eta, G_\eta)}$,
    \item for $\zeta\leq_\mtree{T} \eta\leq_\mtree{T} \xi$, $\vec{\Delta}^{\zeta,\xi}= \vec{\Delta}^{\eta, \xi}\circ \vec{\Delta}^{\zeta,\eta}$.
    \item for $\eta\leq_\mtree{T} \xi$, $\Delta^{\eta,\xi}_{\infty}= \Psi_{\eta,\xi}$.
\end{enumerate}

To start, $\mtree{V}^0=\mtree{S}$. Given everything up to $\mtree{V}^\xi$, let $\eta=\mtree{T}\pred (\xi+1)$. We want to set $\mtree{V}^{\xi+1}= \mtree{V}(\mtree{V}^\eta, \mtree{V}^\xi, G_\xi)$, so we need to see that the agreement hypotheses of the one-step case are met. If $\eta=\xi$, we're good; so assume $\eta<\xi$. By our induction hypothesis (3), we have that $G_\eta = F^\xi_{a(\mtree{V}^\eta,G_\eta)}$. By the normality of $\mtree{T}$, we have that $\crit (G_\xi)<\hat\lambda( G_\eta)$, so $b(\mtree{V}^\xi, G_\xi)\leq a(\mtree{V}^\eta,G_\eta)$. If $b(\mtree{V}^\xi, G_\xi)< a(\mtree{V}^\eta,G_\eta)$ we're done by our induction hypothesis (2). If $b(\mtree{V}^\xi, G_\xi)+1= a(\mtree{V}^\eta,G_\eta)+1= \lh (\mtree{V}^\eta)$, we're also done, since $F^\eta_{b(\mtree{V}^\xi, G_\xi)}$ is undefined. So assume $b(\mtree{V}^\xi, G_\xi)+1= a(\mtree{V}^\eta,G_\eta)+1< \lh (\mtree{V}^\eta)$. Then $F^\eta_{a(\mtree{V}^\eta,G_\eta)}$ is defined but as $G_\eta^-$ is on the sequence of the last models of both $\tree{V}^\eta_{a(\mtree{V}^\eta, G_\eta)}$ and $\tree{T}_\eta$, the last tree of $\mtree{V}^{\eta}$, we must have that $\lh (F^\eta_{a(\mtree{V}^\eta,G_\eta)})\geq \lh(G_\eta)$. So the hypotheses of the one-step case still apply. We also put $\vec{\Delta}^{\eta,\xi+1} = \vec{\Delta}^{\mtree{V}(\mtree{V}^\eta, \mtree{V}^\xi, G_\xi)}$ and $\vec{\Delta}^{\zeta, \xi+1} =\vec{\Delta}^{\eta,\xi+1}\circ \vec{\Delta}^{\zeta,\eta}$ whenever $\zeta\leq_\mtree{T}\eta$. By our work in the one-step case and our induction hypothesis at $\eta$ and $\xi$, it's easy to see all our induction hypotheses still hold at $\xi+1$.

At limit $\lambda$ we take the extended meta-tree embedding direct limit along the branch chosen by $\tree{T}$. That is, letting $\mathbb{D}_\lambda=\langle \mtree{V}^\eta, \vec{\Delta}^{\eta,\xi}\,|\,\eta\leq_\mtree{T}\xi<_\mtree{T}\lambda\rangle$, we let $\mtree{V}^\lambda = \lim \mathbb{D}_\lambda$, if this is wellfounded. The last tree of $\mtree{V}^\lambda$ is the direct limit of the $\tree{T}_\eta$ under $\Psi_{\eta,\xi}$ for $\eta\leq_\mtree{T}\xi<_\mtree{T}\lambda$ by our induction hypotheses (1) and (5), which is just $\tree{T}_\lambda$, since $\mtree{T}$ is a meta-tree. We also let $\vec{\Delta}^{\eta, \lambda}$ be the direct limit meta-tree embeddings. It's easy to see that this maintains the rest of our induction hypotheses.

This finishes the definition of $\mtree{V}(\mtree{S},\mtree{T})$.\\

For a finite stack of meta-trees $\vec{\mtree{S}}=\langle \mtree{S}^0,\ldots, \mtree{S}^n\rangle$, we also define, by induction, $\mtree{V}(\vec{\mtree{S}})=\mtree{V}(\mtree{V}(\vec{\mtree{S}}\restrict n), \mtree{S}_n)$. Notice that this makes sense since, by induction, $\mtree{S}_{n-1}$ and $\mtree{V}(\mtree{S}\restrict n)$ have the same last tree, so $\langle \mtree{V}(\mtree{S}\restrict n), \mtree{S}_n\rangle$ is really a stack of meta-trees.

\begin{definition}
A $(\lambda,\theta)$ strategy $\Sigma$ for $\tree{S}$ \textit{normalizes well} iff for any finite stack of meta-trees $\vec{\mtree{S}}$ by $\Sigma$, $\mtree{V}(\vec{\mtree{S}})$ is by $\Sigma$.
\end{definition}

\begin{remark}
Like with the case of a strategy on a premouse $M$ with strong hull condensation, one can show that a meta-strategy on a tree $\tree{S}$ with meta-hull condensation extends uniquely to a strategy for finite stacks which normalizes well. We have no use for this here and it would take some space to write out, so we'll just assume we are given strategies for finite stacks of meta-trees.
\end{remark}

\begin{proposition}\label{nice meta-strategy existence 2}
Suppose $\Sigma$ is a strategy for $M$ with strong hull condensation and $\tree{S}$ is by $\Sigma$. Then $\Sigma^*_\tree{S}$ normalizes well.
\end{proposition}
\begin{proof}[Proof sketch.]
By induction, we just need to verify this for stacks of length $2$. By our characterization of $\Sigma^*_\tree{S}$, we just need to see that all the trees in all the $\mtree{V}^\xi = \mtree{V}(\mtree{S},\mtree{T}\restrict\xi+1)$ are by $\Sigma$. We do this by induction. At successors all our new trees are all of the form $V(\tree{U}, \tree{V}, G)$ for trees $\tree{U},\tree{V}$ which are by $\Sigma$, so $V(\tree{U},\tree{V}, G)$ is by $\Sigma$, by strong hull condensation and quasi-normalizing well.
At limit $\lambda$, we have all the  $\mtree{V}^\xi$ for $\xi<\lambda$ are by $\Sigma^*$ and we want to see the direct limit along $[0,\lambda)_\mtree{T}$ is by $\Sigma^*$. The trees of the direct limit are either trees of $\mtree{V}^\xi$ or else are (non-trivial) direct limits along one-step embedding normalization tree embeddings by the extenders of $[0,\lambda)_\mtree{T}$. All these non-trivial direct limit trees agree with $\tree{T}_\lambda$ up to $\delta_\lambda +1= \sup_{\xi<\lambda}\{\alpha_0(G_\xi, \tree{T}_\xi)\}$, which is by $\Sigma$. At limit ordinals $\gamma>\delta_\lambda$ in these direct limit trees, the branches are images under the $v$-maps of earlier trees which are by $\Sigma$ and so must be by $\Sigma$ by strong hull condensation.
\qed
\end{proof}

We can now easily prove Theorem \ref{meta-tree full norm}.

\begin{proof}[Proof of Theorem \ref{meta-tree full norm}.]
Let $\tree{S}$ be a countable plus tree by $\Sigma$ and $\langle \mtree{S},\mtree{T}\rangle$ be a stack by $\Sigma^*_\tree{S}$ with last tree $\tree{U}$. Since $\Sigma^*_\tree{S}$ normalizes well, $\mtree{U}=\mtree{V}(\mtree{S},\mtree{T})$ is by $\Sigma^*_\tree{S}$, and is a meta-tree with last tree $\tree{U}$. We get $\Phi^\mtree{U}=\Phi^\mtree{T}\circ \Phi^\mtree{S}$ easily by induction (using the commutativity condition in the definition of meta-tree embeddings). 
\qed 
\end{proof}

We end this section with a few definitions. Let $M$ a premouse and $\tree{S}$ a plus tree on $M$ of successor length. Suppose $\Sigma$ is a strategy for finite stacks of meta-trees on $\tree{S}$ which has meta-hull condensation and normalizes well.

\begin{definition}\label{tail strategies}
 Let $ \mtree{S}$ be a meta-tree by $\Sigma$ of successor length. We define the \textit{tail normal tree strategy} $\Sigma_{\mtree{S}}$ by
\[\tree{U} \text{ is by }\Sigma_{\mtree{S},P} \Leftrightarrow \langle \mtree{S},\mtree{V}(\tree{T}, \tree{U})\rangle\text{ is by }\Sigma.\]
It is easy to see that this is a strategy for single normal trees on $M_\infty^\mtree{S}$.

We can naturally extend this to an internally lift consistent strategy for $M_\infty^\mtree{S}$. By putting, for $Q\is M_\infty^\mtree{S}$, \[\tree{U} \text{ is by }\Sigma_{\mtree{S},Q} \Leftrightarrow \langle \mtree{S}, \mtree{V}(\tree{T},\tree{U}^+)\rangle \text{ is by }\Sigma,\] where $\tree{U}^+$ is the copy of $\tree{U}$ on $Q$ to a normal tree on the full $M_\infty^\mtree{S}$.
\end{definition}
It is easy to see the following.

\begin{proposition}\label{tail strategy prop}
Let $(M,\Sigma)$ be a mouse pair, $\tree{S}$ a plus tree by $\Sigma$ of successor length, and $\mtree{S}$ be a meta-tree by $\Sigma$ of successor length with last tree $\tree{T}$. Then for all $Q\is M_\infty^\mtree{S}$, $(\Sigma^*_\tree{S})_{\mtree{S}, Q}=\Sigma_{\tree{T},Q}$.
\end{proposition}

\begin{definition}
If $M$ is a least branch premouse, $\tree{S}$ is a tree on $M$, and $\Sigma$ is a strategy for finie stacks of meta-trees on $\tree{S}$ which has meta-hull condensation and normalizes well, we say that $\Sigma$ is \textit{pushforward consistent} iff for every $\mtree{S}$ by $\Sigma$ and $Q\is M_\infty^\mtree{S}$, $\dot\Sigma^Q\subseteq \Sigma_{\mtree{S},Q}$.
\end{definition}

By Proposition \ref{tail strategy prop} and pushforward consistency for lbr hod pairs, we immediately get

\begin{proposition}\label{pushforward}
Let $(M,\Sigma)$ be a lbr hod pair, $\tree{S}$ a plus tree by $\Sigma$ of successor length. Then $\Sigma^*_\tree{S}$ is pushforward consistent.
\end{proposition}

A property in the same vein needed for our meta-strategy comparison result is a version of strategy coherence for $\lambda$-separated meta-trees, which follows from normalizing well and meta-hull condensation (essentially by the proof of Corollary 7.6.9 from \cite{nitcis}).

\begin{proposition}\label{strategy coherence prop}
Let $M$ a premouse, $\tree{S}$ a countable 
plus tree on $M$ of successor length, and $\Sigma$ a meta-strategy for finite stacks of countable meta-trees on $\tree{S}$ which normalizes well and has meta-hull condensation. Let $\mtree{S}, \mtree{T}$ be $\lambda$-separated meta-tree on $\tree{S}$ of successor lengths and suppose $Q\is M_\infty^\mtree{S}, M_\infty^\mtree{T}$. Then $\Sigma_{\mtree{S}, Q}=\Sigma_{\mtree{T}, Q}$.
\end{proposition}

We'll consider just one additional property for a meta-strategy, relating it back to an ordinary iteration strategy for the base model.

\begin{definition}\label{doddjensendef}
Let $M$ a premouse, $\Omega$ be an iteration strategy for $M$, $\tree{S}$ a plus tree on $M$ of successor length, and $\Sigma$ a meta-strategy for finite stacks of meta-trees on $\tree{S}$. $\Sigma$ has\textit{ the Dodd-Jensen property relative to $\Omega$} iff for any meta-tree $\mtree{S}$ on $\tree{S}$ by $\Sigma$ of successor length with last tree $\tree{T}$ and last model $P$, 
\begin{enumerate}
    \item if $\tree{T}$ doesn't drop along its main branch, then $(\Sigma_{\mtree{S}})^{i^\tree{T}}\subseteq \Omega$ and
    \item for any $Q\is P$ and nearly elementary map $\pi:M\to Q$ such that  $(\Sigma_{\mtree{S},Q})^\pi\subseteq \Omega$, \begin{enumerate}
    \item $M$-to-$P$ doesn't drop in $\tree{T}$ and $Q=P$,
    \item $i^\tree{T}(\xi)\leq \pi(\xi)$ for any $\xi<o(M)$.
\end{enumerate}
\end{enumerate} 
\end{definition}

\begin{proposition}\label{dodd-jensen}
Let $(M,\Sigma)$ be a lbr hod pair, $\tree{S}$ a plus tree by $\Sigma$ of successor length. Then $\Sigma^*_\tree{S}$ has the Dodd-Jensen property relative to $\Sigma$.
\end{proposition}

\begin{proof}
Let $\mtree{S}$ on $\tree{S}$ by $\Sigma^*_\tree{S}$ of successor length, $\tree{T}$ its last tree, and $P$ its last model. Since $\mtree{S}$ is by $\Sigma^*_\tree{S}$, $\tree{T}$ is by $\Sigma$. If $\tree{T}$ doesn't drop along its main branch, then since $\Sigma$ is pullback consistent, $(\Sigma_\tree{T})^{i^\tree{T}}=\Sigma$. But $(\Sigma^*_\tree{S})_{\mtree{S}}\subseteq \Sigma_\tree{T}$, by Proposition \ref{tail strategy prop}, so $((\Sigma^*_\tree{S})_{\mtree{S}})^{i^\tree{T}}\subseteq \Sigma$, giving (1).

Now suppose we have $Q\is P$ and $\pi:M\to Q$ such that $((\Sigma^*_\tree{S})_{\mtree{S},Q})^\pi\subseteq \Sigma$. By Proposition \ref{tail strategy prop}, $(\Sigma^*_\tree{S})_{\mtree{S},Q}\subseteq\Sigma_{\tree{T},Q}$. It follows that $\Sigma$ and $(\Sigma_{\tree, Q})^\pi$ agree on normal trees, and so by Theorem \ref{uniqueextension}, $\Sigma=(\Sigma_{\tree, Q})^\pi$. That is $\pi$ nearly elementary as a map from $(M, \Sigma)$ into $(Q, \Sigma_{\tree{T},Q})$. So the Dodd-Jensen Lemma for mouse pairs (Theorem 9.3.4 from \cite{nitcis}) immediately gives (2).

\qed
\end{proof}

To review what we've established above, let $(M,\Sigma)$ be a pfs mouse pair or lbr hod pair with 
scope $H_\delta$, $\tree{S}$ a plus tree of successor length by $\Sigma$, and let $\Sigma^*_\tree{S}$ be the induced meta-strategy;
then the tails
 $(\Sigma^*_\tree{S})_{\mtree{S},P}$ are contained in the appropriate 
tails of $\Sigma$. Moreover $\Sigma^*_\tree{S}$ 
\begin{itemize}
\item[(i)] has meta-hull condensation (Proposition \ref{nice meta-strategy existence}), 
\item[(ii)] normalizes well (Proposition \ref{nice meta-strategy existence 2}), 
\item[(iii)] has the Dodd-Jensen property
 relative to $\Sigma$ (Proposition \ref{dodd-jensen}), and 
\item[(iv)] is pushforward consistent, if $(M,\Sigma)$ is an lbr hod pair
(Proposition \ref{pushforward}).
\end{itemize}
 We shall see in the next section that, 
in the right context, these properties uniquely determine $\Sigma^*_\tree{S}$.

\subsection{Comparison}

We need to generalize the tree comparison theorem we proved earlier. In that theorem (Theorem \ref{meta-tree comparison v1}), we showed that, in particular, any two trees by $\Sigma$ had a common meta-iterate (via meta-trees which were by $\Sigma^*$). We now want to compare pairs of the form $(\tree{S}, \Sigma)$, $(\tree{T}, \Lambda)$ where $\Sigma, \Lambda$ are meta-strategies for $\tree{S},\tree{T}$. To do this, we will have to weaken the conclusion that we tree embed both into a common tree. Recall that we reached this conclusion by first arranging that there we reached trees with a common last model. This is what we'll try to arrange in our generalization.

As in \cite{nitcis}, at stages where we reach extender disagreements, we will hit the corresponding plus extender, rather than the disagreement extender itself. A meta-tree $\mtree{S}$ on $\tree{S}$ is \textit{$\lambda$-separated} if for all $\xi+1<\lh(\mtree{S})$, $F_\xi^\mtree{S}$ is of plus type. Notice that if $\tree{T}$ is a plus tree of successor length and $\tree{U}$ is a $\lambda$-separated plus tree on $M_\infty^\tree{T}$, then actually $\tree{U}$ is normal and so $\mtree{V}(\tree{T}, \tree{U})$ is defined and is also $\lambda$-separated meta-tree (assuming it is wellfounded).
\begin{theorem}[Meta-strategy comparison]\label{tree comparison v2}
Assume $\adp$. Let $M, N$ be countable, strongly stable premice of the same type, $\tree{S},\tree{T}$ countable plus trees on $M,N$ of successor lengths, and $\Sigma$, $\Lambda$ meta-strategies for finite stacks of countable meta-trees on $\tree{S}, \tree{T}$ which have meta-hull condensation and normalize well.

Then there are simple $\lambda$-separated meta-trees $\mtree{S}$ by $\Sigma$ and $\mtree{T}$ by $\Lambda$ with last models $P,Q$ such that either
\begin{enumerate}
    \item $(P,\Sigma_{\mtree{S},P}) \is (Q,\Lambda_{\mtree{T},Q})$ and $\mtree{S}$ doesn't drop along its main branch, or
    \item $(Q,\Lambda_{\mtree{T},Q}) \is (P,\Sigma_{\mtree{S},P})$ and $\mtree{T}$ doesn't drop along its main branch.
    \end{enumerate}
\end{theorem}

\begin{remark}
By Remark \ref{meta-strategy remark} and the fact that agreement on normal trees (in fact $\lambda$-separated trees) suffices for full strategy agreement (see \cite{nitcis}), this is a generalization of the main comparison theorem for mouse pairs in \cite{nitcis}.
\end{remark}

This theorem, and a variant we'll need later, follow from a general theorem about comparison with the levels of a background construction. The statement is almost what you would expect, except that we have added that we gratuitously drop however we want, in a sense to follow. To get Theorem \ref{tree comparison v2}, we just need the case $A=\emptyset$, but we will need other $A$ in the variant we need for full normalization.

\begin{definition}\label{parameter set}
Let $\Sigma$ be a meta-strategy for $\tree{S}$ and $A$ any set. A meta-tree $\mtree{S}$ is by $\Sigma, A$ iff it is by $\Sigma$ and for any $\xi<\lh (\mtree{S})$, 
\begin{enumerate}
    \item if there is a unique $\gamma$ such that $\langle \tree{S}_\xi^+,\gamma\rangle\in A$, $\gamma+1<\lh (\tree{S}_\xi^+)$, and $\gamma\geq \sup\{\alpha^\mtree{S}_\eta+1\,|\,\eta<\xi\}$, then $\gamma<\lh(\tree{S}_\xi^+)$ and $\tree{S}_\xi = \tree{S}_\xi^+\restrict\gamma+1$;
\item otherwise, $\tree{S}_\xi = \tree{S}_\xi^+$.

\end{enumerate}
\end{definition}

One can think that, whereas $\Sigma$ is a winning strategy for II in the game where player I must say if and how to gratuitously drop, $(\Sigma, A)$ determines a winning strategy for player II in the game where this information is decided by player II.

We need one more definition.

\begin{definition}
$(P, \Omega)$ a mouse pair and $M$ a premouse of the same type as $P$. Let $\tree{S}$ a plus tree on $M$ of successor length, and $\Sigma$ a $(\lambda, \theta)$-strategy for $\tree{S}$ and $A$ a set. Then $(\tree{S}, \Sigma, A)$ \textit{iterates past} $(P, \Omega)$ iff there is a $\lambda$-separated meta-tree $\mtree{S}$ by $\Sigma$ such that $M_{\eta,j}^\mathbb{C}\is M_\infty^{\mtree{S}}$ and $\Sigma_{\mtree{S}, P}\subseteq \Omega$.\footnote{This inclusion says that $\Sigma_{\mtree{S}, P}$ is the restriction of $\Omega$ to single normal trees of length $<\theta$ on $P$. Recall that $\Omega$ is defined on all plus trees, so includes more information than $\Sigma_{\mtree{S}, P}$, but that it is actually redundant since it is totally determined by its action on normal trees (in fact, $\lambda$-separated trees).} $(\tree{S}, \Sigma, A)$ \textit{iterates strictly past} $(P, \Omega)$ if there is such an $\mtree{S}$ such that, also, either $P\isneq M_\infty^\mtree{S}$ or $\mtree{S}$ has a necessary drop along its main branch.  Finally, $(\tree{S}, \Sigma, A)$ \textit{iterates to} $(P, \Omega)$ if it iterates past $(P,\Omega)$ via an $\mtree{S}$ such that $P=M_\infty^\mtree{S}$ and $\mtree{S}$ doesn't have a necessary drop along its main branch.
\end{definition}

\begin{theorem}\label{main comparison theorem} Suppose $\delta$ is an inaccessible cardinal, $M$ is a strongly stable pfs premouse (or lbr hod mouse), $\tree{S}$ is plus tree on $M$ of successor length, and $\Sigma$ is an $(\omega, \delta)$-strategy for $\tree{S}$ which has meta-hull condensation and normalizes well. 

Let $\mathbb{C}$ be a PFS (or least branch) construction of length $\leq\delta$ such that $\mathcal{F}^\mathbb{C}\subseteq H_\delta$ and for all $E\in \mathcal{F}^\mathbb{C}$, $\crit(E)>o(M), \lh(\tree{S})$, $i_E(\Sigma)\subseteq \Sigma$, and and $i_E(A)\subseteq A$. Let $\langle \nu, k\rangle<\lh(\mathbb{C})$ and suppose that $(\tree{S},\Sigma, A)$ iterates strictly past  $(M_{\eta,j}^\mathbb{C}, \Omega_{\eta,j}^\mathbb{C})$, for all $\langle\eta, j\rangle<_{\text{lex}}\langle \nu, k\rangle$. Then $(\tree{S},\Sigma, A)$ iterates past  $(M_{\nu,k}^\mathbb{C}, \Omega_{\nu,k}^\mathbb{C})$.
\end{theorem}

That Theorem \ref{main comparison theorem} implies Theorem \ref{tree comparison v2} is a simple variant of the analogous argument from \cite{nitcis} (i.e. the proof that Theorem 8.3.5 follows from Theorem 8.3.4). Our proof of Theorem \ref{main comparison theorem} also closely follows the proof of the analogous result, Theorem 8.3.4 from \cite{nitcis}: we compare against levels of the background by least extender disagreement, using plus extenders at every stage, and show by induction on $\langle \nu,k\rangle$ that the background doesn't move and no strategy disagreements show up. The difference is that our comparison process is now producing meta-trees $\mtree{S}$ by $\Sigma, A$ on our base plus tree $\tree{S}$ rather than plus trees by a strategy for a base model $M$.

Let $\mtree{S}_{\nu,k}^*$ be the $\lambda$-separated meta-trees which are by $\Sigma, A$ which are obtained by comparing against the last model of $\tree{S}$ against $M^\mathbb{C}_{\nu,k}$ until we reach an extender disagreement coming from the $M^\mathbb{C}_{\nu,k}$-side \textit{or} until we reach a strategy disagreement. That is, $\mtree{S}_{\nu,k}^*=\langle \tree{S}_\xi, F_\xi, \Phi_{\eta,\xi}\rangle$ is the unique meta-tree by $\Sigma,A$ such that for all $\xi<\lh (\mtree{S}^*_{\nu,k})$, letting $P_\xi$ the last model of $\tree{S}_\xi$,

\begin{itemize}
    \item[(i.)] if $(M_{\nu,k}, \Omega_{\nu,k})\parallel (P_\xi, \Sigma_{\mtree{S}^*_{\nu,k}\restrict\xi+1, P_\xi})$, then  $\xi+1=\lh (\mtree{S}^*_{\nu,k})$, and
    \item[(ii).] If $(M_{\nu,k}, \Omega_{\nu,k})\not\parallel (P_\xi, \Sigma_{\mtree{S}^*_{\nu,k}\restrict\xi+1, P_\xi})$ then for $\langle \eta, l\rangle<\ell(M_{\nu,k})$ least such that \[(M_{\nu,k}, \Omega_{\nu,k})|\langle \eta,l\rangle\not\parallel (P_\xi, \Sigma_{\mtree{S}^*_{\nu,k}\restrict\xi+1, P_\xi})|\langle \eta,l\rangle,\]
    either
    \begin{itemize}
        \item [(a)] $\xi+1<\lh (\mtree{S}^*_{\nu,k})$,
    $l=0$, and $\eta$ is the index of an extender of the $P_\xi$-sequence but not the index of an extender on the $M_{\nu,k}$-seqeunce, \item[(b)] $\xi+1=\lh (\mtree{S}^*_{\nu,k})$ and $l=0$, but $\eta$ is the index of an extender on the $M_{\nu,k}$-sequence which is \textit{not} on the $P_{\nu,k}^*$-sequence, or
    \item[(c)] $\xi+1=\lh(\mtree{S}^*_{\nu,k})$, $M_{\nu,k}|\langle\eta,l\rangle=P_{\nu,k}^*|\langle \eta,l\rangle$ but $\Omega_{\nu,k}|\langle \eta,l\rangle$ disagrees with $\Sigma_{\mtree{S}^*_{\nu,k}, P^*_{\nu,k}}|\langle \eta,l\rangle$ on a normal tree.
    \end{itemize}
\end{itemize}
Here (ii)(a) just says that we are building $\mtree{S}^*_{\nu,k}$ by comparing via least extender disagreements arising from just the $\tree{S}$-side for as long as possible; (ii)(b) says we had to stop because we'd need to hit an extender on the $M_{\nu,k}$-side to continue comparing by least disagreement; (ii)(b) says that we had to stop because we reached a strategy disagreement.

We show by induction on $\langle \nu, k\rangle$ that (i) holds at least until we reach a $\langle \nu, k\rangle$ such that $(\tree{S}, \Sigma, A)$ iterates to $(M_{\nu,k}, \Omega_{\nu,k})$ via $\mtree{S}^*_{\nu,k}$, i.e. a $\langle \nu,k\rangle$ such that $P^*_{\nu,k}=M_{\nu,k}$, $ \Sigma_{\mtree{S}^*_{\nu,k}, P^*_{\nu,k}}\subseteq \Omega_{\nu,k}$, and $\mtree{S}^*_{\nu,k}$ doesn't have a necessary drop along its main branch.

So suppose (i) holds for all $\langle \eta, l\rangle\leq_{\text{lex}}\langle \nu,k\rangle$ and but that we haven't yet reached this situation. We have the following lemma, which is the appropriate generalization of Sublemma 8.3.1.1 of \cite{nitcis}).

\begin{lemma}\label{res map lemma 1}
Suppose $M_{\nu,k}$ is \textit{not} sound. Then $(i)$ holds for $\mtree{S}_{\nu,k+1}^*$ and
\begin{enumerate}
    \item $M_{\nu,k}$ is the last model of $\mtree{S}^*_{\nu,k}$ and $\Sigma_{\mtree{S}^*_{\nu,k}, M_{\nu,k}}\subseteq\Omega_{\nu,k}$,
    \item $\mtree{S}^*_{\nu,k+1}\is \mtree{S}^*_{\nu,k}$,
    \item letting $\eta$ and $\delta_{\nu,k}$ be such that $\mtree{S}^*_{\nu, k+1}=\mtree{S}^*_{\nu,k}\restrict\eta+1$ and $\delta_{\nu,k}+1 = \lh (\mtree{S}^*_{\nu,k})$, we have $\eta\leq_{\mtree{S}^*_{\nu,k}} \delta_{\nu,k}$,
    \item the last $t$-map of the tree embedding $\Phi^{\mtree{S}^*_{\nu,k}}_{\eta,\delta_{\nu,k}}$ is the uncoring map $\pi: M^-_{\nu,k+1}\to M_{\nu,k}$.
\end{enumerate}
\end{lemma}

The proof of this lemma relies on our analysis of drops in meta-trees along with the following easy fact about ordinary normal trees.

\begin{proposition}\label{no projectum overlap}
    Suppose $\tree{T}$ is a normal tree with last model $N$ and $P\is N$ is sound. Then for all $\xi+1<\lh( \tree{T})$,
 \[\lh (E^\tree{T}_\xi)\not\in(\rho(P), o(P)).\]
\end{proposition}
\begin{proof}
We may assume that $\tree{T}$ uses no extenders with lengths $>o(P)$, in which case we just need to show that $\tree{T}$ uses no extenders with lengths $>\rho(P)$.

If $o(P)<o(N)$, then no ordinal $\alpha\in (\rho(P), o(P)]$ is a cardinal of $N$ (since ${|\rho(P)|^+}^Q\geq o(J_1(P))$), so $\lh (E_\xi^\tree{T})\not\in (\rho(P), o(P)]$ as it is a cardinal of $N$.

Now suppose $o(P)=o(N)$, i.e. $P=N|\langle o(N),k(P)\rangle$. Towards a contradiction, assume $\tree{T}$ uses an extender with length $>\rho(P)$. Then there is an extender used along the main branch of $\tree{T}$ with length $>\rho(P)$. 
Let $\delta+1=\lh (\tree{T})$ and let $\eta,\xi$ such that $\eta=\tree{T}\pred(\xi+1)$, $\lh (E_\xi^\tree{T})>\rho(P)$, and $(\xi+1, \delta]_\tree{T}$ doesn't drop.

Then we have that $\pi\defeq\hat\imath^\tree{T}_{\eta, \delta}$ is elementary and cofinal on its domain $Q\is M_\eta^\tree{T}$. Now we claim that $\rho_{k(P)+1}(Q)\leq \crit (E^\tree{T}_\xi)$. To see this, note that if $k(P)<k(N)$, $\pi(\rho_{k(P)+1}(Q))= \rho(P)$ and if $k(P)=k(N)$, $\sup\pi" \rho_{k(P)+1}(Q)=\rho(P)$, so in either case $\rho_{k(P)+1}(Q)>\crit (E^\tree{T}_\xi)$ implies $\rho(P)\geq \lh (E_\xi^\tree{T})$, a contradiction. But then since we apply $E^\tree{T}_\xi$ to $Q$ along the main branch of $\tree{T}$, the last model of $\tree{T}$, $N$, is not $k(P)+1$-sound. Hence $P=N|\langle o(N), k(P)\rangle$ is not sound, a contradiction.

\qed
\end{proof}

\begin{proof}[Proof of Lemma \ref{res map lemma 1}.]
Let $\mtree{S}^*_{\nu,k}= \langle\tree{S}_\xi, F_\xi, \Phi_{\eta,\xi}\rangle$.
Since $M_{\nu,k}$ is an initial segment of the last model of $\mtree{S}^*_{\nu,k}$ but isn't sound, it must be equal to the last model; giving (1). Since we assumed we haven't reached the situation of (1)(c), $\mtree{S}^*_{\nu,k}$ must have a necessary drop along its main branch. 
Let $\eta$-to-$\xi+1$ be the last necessary drop along this branch. By our analysis of droppin in a meta-tree, 
we have that the last $t$-map of $\Phi_{\eta,\delta_{\nu,k}}$ is the uncoring map $\pi:M_{\nu,k+1}^-\to M_{\nu,k}$ and $M_{\nu,k+1}$ is an initial segment of the last model of $\tree{S}_\eta$.

To finish, we just need to see that $\mtree{S}_{\nu,k+1}^*= \mtree{S}^*_{\nu,k}\restrict\eta+1$, since then (2) clearly holds and we've seen (4) and the rest of (3) hold for $\eta$.

Let $\rho=\rho(M_{\nu,k+1})$. We have that $M_{\nu,k+1}|{\rho^+}^{M_{\nu,k+1}}=M_{\nu,k}|{\rho^+}^{M_{\nu,k}}$, so that $\mtree{S}^*_{\nu,k}$ and $\mtree{S}^*_{\nu,k+1}$ use the same extenders with lengths $\leq \rho$. So we just need to see that $\mtree{S}^*_{\nu,k}\restrict\eta+1$ uses no extender with length $>\rho$. 

All of the meta-tree exit extenders $F_\zeta$ for $\zeta< \eta$ are ordinary exit extenders of the normal tree $\tree{S}_\eta$. Since $M_{\nu,k+1}^-$ is a sound initial segment of the last model of $\tree{S}_\eta$, by Proposition \ref{no projectum overlap}, $\tree{S}_\eta$ uses no extenders with lengths in the interval $(\rho, o(M_{\nu,k+1}))$. So, we need to show that $\lh (F_\zeta)<o(M_{\nu,k+1})$ for all $\zeta<\eta$, since then all these $F_\zeta$ must have $\lh (F_\zeta)\leq \rho$.

This just follows by the normality of $\mtree{S}^*_{\nu,k}$ and the quasi-normality of $\tree{S}_{\xi+1}$. Fix $\zeta<\eta$. Since $\eta=\mtree{S}_{\nu,k}^*\pred(\xi+1)$, $\crit (F_\xi)>\hat\lambda (F_\zeta)$. In $\tree{S}_{\xi+1}$, $F_\xi=E^{\tree{S}_{\xi+1}}_{\alpha_\xi}$ and $F_\zeta=E^{\tree{S}_{\xi+1}}_{\alpha_\zeta}$,
 so $\alpha_\zeta<\beta_\xi= \tree{S}_{\xi+1}\pred(\alpha_\xi+1)$. So we have that $\lh (F_\zeta)<\lh (E_{\beta_\xi}^{\tree{S}_{\xi+1}})$. 
 Also, $F_{\xi}$ is applied to $M_{\nu,k+1}^-\is M^{\tree{S}_{\xi+1}}_{\beta_\xi}$, so $M_{\beta_\xi}^{\tree{S}_{\xi+1}}|\lh (E^{\tree{S}_{\xi+1}}_{\beta_\xi})\is M_{\nu,k+1}$.  So $\lh (F_\zeta)<o(M_{\nu,k+1})$, as desired.
 \qed
\end{proof}

Using Lemma \ref{res map lemma 1}, we get that resurrection embeddings are realized as the last $t$-maps of appropriate tree embeddings (which is the appropriate generalization of Lemma 8.3.1 of \cite{nitcis}). We will also use Lemma \ref{res map lemma 1} below to rule out that (ii) (b) ever occurs, i.e. we never stop building $\mtree{S}_{\nu,k}^*$ because we reach an extender disagreement coming from the background.

\begin{lemma}\label{res map lemma 2}
Let $\langle \theta, j\rangle< \langle \nu, k\rangle$ and $P\is M_{\theta, j}$. Let $\langle \theta_0, j_0\rangle = Res_{\theta,j} [P]$ and $\tau=\sigma_{\theta, j}[P]$ (so $\tau:P\to M_{\theta_0, j_0}$).

Let $\xi$ least such that $P\is M^{{\tree{S}^{\theta, j}}^*}_{z_{\theta,j}(\xi)}$. Then
\begin{enumerate}
    \item $\mathbb{S}^*_{\theta_0,j_0}\restrict \xi+1 = \mathbb{S}^*_{\theta,j}\restrict \xi+1$,
    \item $\xi \leq_{\mtree{S}^*_{\theta_0, j_0}} \delta =\lh (\mtree{S}^*_{\theta_0, j_0})-1$,
    \item for $t^{\xi, \delta}$ the last $t$-map along $[\xi, \delta]_{\mtree{S}^*_{\theta_0,j_0}}$, $t^{\xi, \delta}\restrict P= \tau$.
\end{enumerate}
\end{lemma}

This follows from Lemma \ref{res map lemma 1} by, essentially, the argument for Lemma 8.3.1 of \cite{nitcis}. We omit further details.

It will be useful to observe that nice trees on $V$ naturally give rise to tree embeddings in the following way.

\begin{proposition}\label{hull embeddings prop}
Let $\tree{T}$ be a nice, normal tree on $V$ with last model $Q$. Let $i_\tree{T}: V\to Q$ be the iteration map of $\tree{T}$.

Let $M$ be a premouse such that $\tree{T}$ is above $|M|^+$, i.e. all the critical points of extenders used in $\tree{T}$ are $>|M|^+$.

Then there is a (unique) extended tree embedding
$I^\tree{S}_\tree{T}: \tree{S}\to i_\tree{T}(\tree{S})$ with $u$-map $i_\tree{T}\restrict \lh \tree{S}$ and $t$-maps $t_\xi=i_\tree{T}\restrict M^\tree{S}_\xi$.
\end{proposition}

We'll prove this by induction on length $\tree{T}$ using the following lemma at successor stages.
\begin{lemma}\label{hull embeddings lemma}
Let $\tree{T}$ be a nice, normal tree on $V$ and let $Q_\xi = M^\tree{T}_\xi$. Suppose $\nu=\tree{T}\pred(\gamma+1)$ and $\tree{S}\in Q_\nu$ is a normal tree on a premouse $M\in Q_\nu$ such that the extender $F= E^\tree{T}_\gamma$ has critical point $>|M|^+$. Let $\pi=i^{Q_\nu}_F$. 

Then there is a (unique) extended tree embedding
$I^{\tree{T};\nu,\gamma+1}_\tree{S}: \tree{S}\to \pi(\tree{S})$ with $u$-map $\pi\restrict \lh (\tree{S})$ and $t$-maps $t_\xi=\pi\restrict M^\tree{S}_\xi$.
\end{lemma}
\begin{proof}
We'll just check the simplest case $\tree{T}=\langle G\rangle$ and leave the general case to the reader. 

So let $G$ an extender on $V$ and suppose $M$ is a premouse such that $|M|^+<\crit G$ and $\tree{S}$ is a normal tree on $M$. Let $\pi= i^V_G$. We check that there is an extended tree embedding $I=\langle u,v,\{s_\xi\}, \{t_\xi\}\rangle: \tree{S}\to \pi(S)$ with $u$-map $\pi\restrict \lh (\tree{S})$ and $t$-maps $t_\xi = \pi\restrict M_\xi^\tree{S}$.

Since $u\restrict\kappa=id$, we have that $I\restrict(\tree{S}\restrict\kappa)$ is just the identity tree embedding on $\tree{S}\restrict\kappa$. At $\kappa$, we must let $v(\kappa)=\kappa$ and $s_\kappa = id$. To see we can make our desired assignments for $u(\kappa)$ and $t_\kappa$, we check

\begin{claim}
    $\kappa\leq_{\pi(\tree{S})} \pi(\kappa)$ and $\hat \imath^{\pi(\tree{S})}_{\kappa, \pi(\kappa)} = \pi\restrict M^\tree{S}_\kappa$.
\end{claim}

\begin{proof}
This is routine. First, we have $\xi<_{\pi(\tree{S})} \pi(\kappa)$ for every $\xi<_{\tree{S}}\kappa$ since $\pi\restrict\kappa=id$, so $\kappa\leq_{\pi(\tree{S})} \pi(\kappa)$ since branches of iteration trees are closed below their sup. Also, $\tree{S}\restrict\kappa+1 \is \pi(\tree{S})$ since $\pi$ is tree-order preserving (and the identity below $\kappa$), fixes $M_\xi^\tree{S}$ for $\xi<\kappa$ (since $M$ is small), and $M_\kappa^\tree{S}=M_\kappa^{\pi(\tree{S})}$ since they are both given by the same direct limit.

Now let $x\in M_\kappa^\tree{S}$. Let $\xi<_\tree{S}\kappa$ and $\bar x\in M_\xi^\tree{S}$ such that $\hat\imath^\tree{S}_{\xi, \kappa}(\bar x)= x$. Now
\begin{align*}
    \pi(x)&= \pi( \hat\imath^\tree{S}_{\xi, \kappa}(\bar x))\\
    &= \hat\imath^{\pi(\tree{S})}_{\xi, \pi(\kappa)}(\bar x)\\
    &=\hat\imath^{\pi(\tree{S})}_{\kappa, \pi(\kappa)}\circ \hat\imath^\tree{S}_{\xi, \kappa}(\bar x)\\
    &=\hat\imath^{\pi(\tree{S})}_{\kappa, \pi(\kappa)}(x).
\end{align*}
The second equality is given by the elementarity of $\pi$ and using that $\pi(\bar x)= \bar x$ since $\bar x\in M_\xi^\tree{S}$ and $\xi<\kappa$. The third equality uses that $\kappa\in (\xi, \pi(\kappa))_{\pi(\tree{S})}$ and $\hat\imath^\tree{S}_{\xi, \kappa} = \hat\imath^{\pi(\tree{S})}_{\xi, \kappa}$. The final equality is immediate from our choice of $\xi, \bar x$.
\hfill{$\qed$}
\end{proof}
So letting $u(\kappa)= \pi(\kappa)$ and $t_\kappa = \hat\imath^{\pi(\tree{S})}_{\kappa, \pi(\kappa)}$ (as we must), we have that $I\restrict\tree{S}\restrict\kappa+1 \to\pi( \tree{S})\restrict\pi(\kappa)+1$ is an extended tree embedding meeting the desired assignments. We check the rest by induction. 

The general limit case is basically the same as the case at $\kappa$. Suppose that $I\restrict (\tree{S}\restrict\lambda):\tree{S}\restrict\lambda\to \pi(\tree{S})\restrict \pi(\lambda)$ is a tree embedding\footnote{Here we are \textit{not} assuming it is an extended tree embedding, since that doesn't make sense as $\tree{S}\restrict\lambda$ has limit length.} with $u(\xi)=\pi(\xi)$ and $t_\xi = \pi\restrict M^\tree{S}_\xi$ for all $\xi<\lambda$. To extend $I\restrict(\tree{S}\restrict\lambda)$, we must let $v(\lambda)=\sup \pi" \lambda$. We first check
\begin{claim}
$v(\lambda)\leq_{\pi(\tree{S})} \pi(\lambda)$ and $\pi"[0,\lambda)_\tree{T}$ is a cofinal subset of $[0,v(\lambda))_{\pi(\tree{S})}$.
\end{claim}
\begin{proof}
For $\xi<_\tree{S}\lambda$, $\pi(\xi)<_{\pi(\tree{S})} \pi(\lambda)$, so $\pi"[0,\lambda)_\tree{S}\subseteq[0,\pi(\lambda))_{\pi(\tree{S})}$. Since branches of iteration trees are closed below their sup and $[0,\lambda)_\tree{S}$ is cofinal in $\lambda$,  $v(\lambda)=\sup \pi"[0,\lambda)_\tree{S}\leq_{\pi(\tree{S})} \pi(\lambda)$.

\hfill{$\qed$}
\end{proof}

Now $M_\lambda^\tree{S}=\lim_{\xi\in [0,\lambda)_\tree{S}} M_\xi^\tree{S}$ and $M_{v(\lambda)}^{\pi(\tree{S})}=\lim_{\xi\in [0,\lambda)_\tree{S}} M_{\pi(\xi)}^{\pi(\tree{S})}$
so that $s_\lambda$ must be the unique map from $M^\tree{S}_\lambda$ into $M^{\pi(\tree{S})}_{v(\lambda)}$ such that for all $\xi<_\tree{S} \lambda$, $s_\lambda \circ \hat\imath^{\tree{S}}_{\xi, \lambda}=\hat\imath^{\pi(\tree{S})}_{\pi(\xi), v(\lambda)}\circ t_\xi$. Since we've stipulated $u(\lambda)=\pi(\lambda)$, we must put $t_\lambda = \hat\imath^{\pi(\tree{S})}_{v(\lambda),\pi(\lambda)}\circ s_\lambda$. So to finish the limit case we just need to check
\begin{claim}
$\hat\imath^{\pi(\tree{S})}_{v(\lambda),\pi(\lambda)}\circ s_\lambda= \pi\restrict M^\tree{S}_\lambda$.
\end{claim}
\begin{proof}
Let $x\in M_\lambda^\tree{S}$. Let $\xi<_\tree{S}\lambda$ and $\bar x\in M^\tree{S}_\xi$ such that $\hat\imath^\tree{S}_{\xi, \lambda}(\bar x)= x$. Then
\begin{align*}
    \pi(x)&= \pi(\hat\imath^\tree{S}_{\xi, \lambda}(\bar x))\\
    &= \hat\imath^{\pi(\tree{S})}_{\pi(\xi), \pi(\lambda)}(\pi(\bar x))\\
    &=\hat\imath^{\pi(\tree{S})}_{v(\lambda), \pi(\lambda)}\circ\hat\imath^{\pi(\tree{S})}_{\pi(\xi), v(\lambda)}\circ t_\xi(\bar x)\\
    &= \hat\imath^{\pi(\tree{S})}_{v(\lambda), \pi(\lambda)}\circ s_\lambda \circ \hat\imath^{\tree{S}}_{\xi, \lambda}(\bar x)\\
    &= \hat\imath^{\pi(\tree{S})}_{v(\lambda), \pi(\lambda)}\circ s_\lambda(x).
\end{align*}
The second equality is just using the elementarity of $\pi$, the third splits up the branch embedding and uses our induction hypothesis that $\pi(\bar x)=t_\xi(\bar x)$ (since $\bar x\in M^\tree{S}_\xi$), the fourth uses $s_\lambda \circ \hat\imath^{\tree{S}}_{\xi, \lambda}=\hat\imath^{\pi(\tree{S})}_{\pi(\xi), v(\lambda)}\circ t_\xi$, as observed above, and the final equality just uses that $x=\hat\imath^\tree{S}_{\xi, \lambda}(\bar x)$.

\hfill{$\qed$}
\end{proof}

For the successor case, suppose we have $I\restrict (\tree{S}\restrict\xi+1)$ is an extended tree embedding from $\tree{S}\restrict\xi+1$ to $\pi(\tree{S})\restrict \pi(\xi)+1$ and $t_\eta = \pi\restrict M^\tree{S}_\eta$ for all $\eta\leq \xi$. 

Notice that $\pi(\xi+1)=\pi(\xi)+1$, so we have $u(\xi+1)=v(\xi+1)$, and so only need to define $s_{\xi+1}=t_{\xi+1}$. We also have that $E_{\pi(\tree{S})}^{\pi(\tree{S})} = \pi(E^\tree{S}_\xi)=t_\xi(E^\tree{S}_\xi)$ and for $\eta=\tree{S}\pred(\xi+1)$, $\pi(\eta)=\pi(\tree{S})\pred (\pi(\xi)+1)$, so we know that $s_{\xi+1}$ is given by the Shift Lemma applied to $t_\xi, t_\eta, E_\xi$. To finish, we just need to check that
\begin{claim}
$s_{\xi+1}=\pi\restrict M_{\xi+1}^\tree{S}$.
\end{claim}
\begin{proof}
We'll just deal with the case in which we take 0-ultrapowers on both sides. Let $P\is M^\tree{S}_\eta$ the level to which we apply $E^\tree{S}_\xi$. By the elementarity of $\pi$, $t_\eta(P)= \pi(P)$ is the level of $M^{\pi(\tree{S})}_{\pi(\xi)}$ to which we apply $E^{\pi(\tree{S})}_{\pi(\xi)}$. Let $a\in [\hat\lambda(E^\tree{S}_\xi)]^{<\omega}$ and $f:[\crit (E^\tree{S}_\xi)]^{|a|}\to Q$, $f\in  Q$. We have that \begin{align*}
    \pi([a,f]^P_{E^\tree{S}_\xi})&=[\pi(a),\pi(f)]^{\pi(P)}_{E^{\pi(\tree{S})}_{\pi(\xi)}}\\
    &= [t_\xi(a), t_\eta(f)]^{\pi(P)}_{E^{\pi(\tree{S})}_{\pi(\xi)}}\\
    &=s_{\xi+1}([a,f]^P_{E^\tree{S}_\xi}),
\end{align*}
as desired.
\hfill{$\qed$}
\end{proof}

\qed
\end{proof}

This proposition naturally extends to meta-tree embeddings. The proof is straightforward using Proposition \ref{hull embeddings prop}, so we omit it.

\begin{proposition}\label{hull embeddings prop 2}
Let $\tree{T}$ a nice, normal tree on $V$ and let $Q_\xi = M^\tree{T}_\xi$. 
Suppose $\nu\leq_\tree{T}\gamma$ and $\mtree{S}\in M^\tree{T}_\nu$ is a normal
 tree on a premouse $M\in Q_\nu$ such that the extenders in $[\nu, \gamma)_\tree{T}$ 
have critical points $>|M|^+$.

Then there is a (unique) extended meta-tree embedding from $\mtree{S}$ into $i^\tree{T}_{\nu,\gamma}(\mtree{S})$ with $u$-map $i^\tree{T}_{\nu, \gamma}\restrict \lh (\mtree{S})$ and $\Delta$-maps $\Delta_\xi=I^{\tree{T};\nu, \gamma}_{\tree{S}_\nu}$.
\end{proposition}

We now verify that we never stop building $\mtree{S}^*_{\nu,k}$ because we 
reach an extender disagreement coming from the background, i.e. case (ii)(b)
 doesn't occur. So suppose we've built $\mtree{S}^*_{\nu,k}\restrict\xi+1$ with 
last model $P_\xi$ and let $\langle \eta,l\rangle$ least such that 
\[
(M_{\nu,k}, \Omega_{\nu,k})|\langle \eta,l\rangle\not\parallel (P_\xi, \Sigma_{\mtree{S}^*_{\nu,k}\restrict\xi+1, P_\xi})|\langle \eta,l\rangle.
\]
Then we have the following.
\begin{lemma}\label{background doesn't move}
If $l=0$, $E^{M_{\nu,k}}_\eta=\emptyset$.
\end{lemma}
\begin{proof}
Suppose not and let $G= E^{M_{\nu,k}}_\eta$. If $k=m+1$ for some $m$, then by our induction hypothesis and Lemma \ref{res map lemma 1}, $M_{\nu,k+1}$ is an initial segment of the last model of $\mtree{S}^*_{\nu,k+1}$, so $G$ must have been on the $P_\xi$-sequence after all. So $k=0$ and $G$ is the top extender of $M_{\nu,0}$. 

Let $\mtree{S}=\mtree{S}^*_{\nu,0}$. Let $G^*$ be the background for $G$ in $\mathbb{C}$. Let $\kappa=\crit (G)=\crit (G^*)$ and $\pi: V\to Ult(V,G^*)$ the ultrapower embedding. 
Let $\vec{\Psi}=\langle \pi, \Psi_\xi\rangle$ be the extended tree embedding from $\tree{S}$ into $\pi(\tree{S})=(\tree{S}^*_{\pi(\nu), 0})^{\pi(\mathbb{C})}$ given by Proposition \ref{hull embeddings prop 2}. 
By hypothesis,  $M, \tree{S}$ are fixed by $\pi$ and $\pi(\Sigma)\subseteq \Sigma$ and $\pi(A)\subseteq A$. So $\pi(\mtree{S})$ is by $\Sigma, A$. Since $\nu= \lh (G)< \lh (G^*)<\pi(\kappa)$, and $G^*$ is strong to its length, we have that $M_{\nu,0}^{\pi(\mathbb{C})}||\lh (G)=M_{\nu,0}||\lh (G)$. Since background certificates must be Mitchell-order minimal, we get that $M_{\nu,0}^{\pi(\mathbb{C})}$ is passive, i.e. $M_{\nu,0}^{\pi(\mathbb{C})}= M_{\nu,0}||\lh (G)$.
Since $G$ is the trivial completion of $(\kappa, \hat\lambda(G))$ extender of $\pi\restrict M_{\nu,0}$, we get that \[M_{\pi(\nu),0}^{\pi(\mathbb{C})}|\lh (G) = \pi(M_{\nu,0})|\lh (G) = M_{\nu,0}||\lh (G)=M_{\nu,0}^{\pi(\mathbb{C})}.\]

Since $\mtree{S}$ is obtained by comparing $\tree{S}$ and $M_{\nu,0}$ by least extender disagreements and below $\lh (\mtree{S})$ $(ii)(b)(c)$ never occur, by elementarity of $\pi$, $\pi(\mtree{S})$ is obtained by comparing $\tree{S}$ and $M_{\pi(\nu),0}^{\pi(\mathbb{C})}$ by least extender disagreements and below $\lh (\pi(\mtree{S}))$, $(ii)(b)(c)$ never occur. 
By the agreement between $M_{\nu,0}$ and $M_{\pi(\nu),0}^{\pi(\mathbb{C})}$ observed above, and since both trees are by $\Sigma, A$, we get that $\mtree{S}\is \pi(\mtree{S})$. Since $\kappa=\crit (\pi)$, we get $\kappa\leq_{\pi(\mtree{S})}\pi(\kappa)$ and there is no dropping (of any kind) between $\kappa$ and $\pi(\kappa)$. So $\Phi=\Phi_{\kappa,\pi(\kappa)}^{\pi(\mtree{S})}$ is a total extended tree embedding from $\tree{S}_\kappa$ into $\pi(\tree{S}_\kappa)$. Moreover, it is easy to see that $\Phi=\Psi_\kappa$, since both have $u$-map and $t$-maps the appropriate restrictions of $\pi$.\footnote{$\Psi_\kappa$ is the unique tree embedding with this property; to see that $\Phi$ has this property, we just use that $\kappa$ is the critical point of $\pi$ and the directed system of trees whose limit is $\tree{S}_{\pi(\kappa)}$ is the image under $\pi$ of the directed system of trees whose limit is $\tree{S}_\kappa$. We leave the details to the reader. 
}

Now, $\lh (\mtree{S})\leq \nu+1< \pi(\kappa)$ and $M_{\nu,0}^{\pi(\mathbb{C})}$ is an initial segment of the last model of $\mtree{S}\is \pi(\mtree{S})$, so since all meta-tree exit extenders of $\pi(\mtree{S})$ used after $\mtree{S}$ have length strictly greater than $\lh (G)$ (as $G$ is not on the sequence of the last model of $\mtree{S}$), we have that $M_{\nu, 0}^{\pi(\mathbb{C})}$ is an initial segment of the last model of $\tree{S}_\xi^{\pi(\mathbb{C})}$ for any $\xi\geq \lh (\mtree{S})$.

Let $\lambda+1$ be the successor of $\kappa$ in $(\kappa,\pi(\kappa)]_{\pi(\mtree{S})}$. $F_\lambda^{\pi(\mtree{S})}$ is compatible with $G$ since both are initial segments of the extender of $t^\Phi_\kappa$.
\footnote{for $F_\lambda^{\pi(\mtree{S})}$, this is just because $t^{\Phi_{\kappa,\lambda+1}}_\kappa $ is just the ultrapower by $F_\lambda^{\pi(\mtree{S})}$; 
for $G$ this is because $t^\Phi_\kappa$ is just $\pi\restrict M_\kappa^{\tree{S}_\kappa}$ and $G$ is total on $M_{\pi(\kappa)}^{\tree{S}_{\pi(\kappa)}^{\pi(\mathbb{C})}}||\lh (G)$ which is either an initial segment of $M_\kappa^{\tree{S}_\kappa}$ or has the same $P(\kappa)$ as it.}
Now if $\lambda<\lh (\mtree{S})$, then $\lh (F_\lambda^{\pi(\mtree{S})})<\lh (G)$ so that $F_\lambda^{\pi(\mtree{S})}$ is on the $M_{\nu,0}$-sequence by the Jensen initial segment condition. But then $F_\lambda^{\pi(\mtree{S})}$ isn't an extender disagreement, a contradiction. So $\lambda\geq \lh (\mtree{S})$. But then $\lh (F_\lambda^{\pi(\mtree{S})})$ and we get  $G$ is on the sequence of the last model of $\tree{S}_\lambda^{\pi(\mtree{S})}$. But $M_{\pi(\nu),0}^{\pi(\mathbb{C})}$ agrees with the last model of $\tree{S}_\lambda^{\pi(\mtree{S})}$ through $F_\lambda^{\pi(\mtree{S})}$, so $G$ is on the sequence of $M_{\pi(\nu),0}^{\pi(\mathbb{C})}$, a contradiction (as observed above, $M_{\pi(\nu),0}^{\pi(\mathbb{C})}|\langle\lh (G),0\rangle$ is passive).
\qed
\end{proof}

All that remains is to verify that no strategy disagreements show up, i.e. (ii)(c) never occurs. To do this, we generalize the proof of Theorem 8.4.3 of \cite{nitcis}. This generalization is straightforward, but quite involved, so we give a sketch of it here.

Before we start, we can quote that result to prove the present theorem in the case when $\Sigma=\Lambda^*_\tree{S}$ for some $\Lambda$ a strategy for the base model $M$ such that $(M,\Lambda)$ is a mouse pair. 
Suppose we're in this case. Let $\delta+1=\lh (\mtree{S}^*_{\nu,k})$ and $\langle \eta, l\rangle$ such that $P_\delta|\langle \eta, l\rangle = M_{\nu,k}|\langle \eta, l\rangle$. Let $\gamma<\lh (\tree{S}_\delta)$ least such that $M_{\nu,k}|\langle \eta, l\rangle \is M^{\tree{S}_\delta}_\gamma$. 
Then  $\tree{S}_\delta\restrict\gamma+1$ is the unique shortest tree which is by $\Lambda$ and has  $M_{\nu,k}|\langle \eta, l\rangle$ as an initial segment of its last model. So then $\tree{S}_\delta\restrict\gamma+1\is \tree{W}^*_{\nu,k}$, the comparison of $M$ with $M_{\nu,k}$ via $\Lambda$. 
By (the proof of) Theorem 8.4.3 of \cite{nitcis}, we have that \[( P_\delta|\langle \eta, l\rangle, \Lambda_{\tree{S}_\delta\restrict\gamma+1, P_\delta|\langle \eta,l\rangle}) = (M_{\nu,k}, \Omega_{\nu, k})|\langle \eta,l\rangle.\]
But $(M,\Lambda)$ is internally lift consistent, so \[\Lambda_{\tree{S}_\delta, P_\delta}|\langle \eta, l \rangle =\Lambda_{\tree{S}_\delta\restrict\gamma+1, P_\delta|\langle \eta,l\rangle}.\]
By the definition of $\Lambda^*$, $\Sigma_{\mtree{S}_{\nu,k}^*, P_\delta}= \Lambda_{\tree{S}_\delta, P_\delta}$, so
\[( P_\delta, \Sigma_{\mtree{S}_{\nu,k}^*, P_\delta})| \langle\eta, l\rangle = (M_{\nu,k}, \Omega_{\nu, k})|\langle \eta, l\rangle,\]
as desired.

We now turn to the general case of establishing that no strategy disagreements show up. Suppose now that we have built $\mtree{S} = \mtree{S}^*_{\nu,k}\restrict\xi+1$ and lined up it's last model $P$ up to $\langle \eta, l\rangle$ with $M_{\nu,k}$. 
It suffices to show the iteration strategies agree on $\lambda$-separated trees, so let $\tree{U}$ be a $\lambda$-separated plus tree of limit length on $P|\langle \eta, l\rangle= M_{\nu,k} |\langle \eta, l\rangle$ which is by both $\Sigma_{\mtree{S}, P}$ and $\Omega_{\nu,k}$. 
We show that $\Sigma_{\mtree{S}, P}(\tree{U})=\Omega_{\nu,k}(\tree{U})$. Let 
\[c=\langle M_{\nu, k}, id, M_{\nu,k},\mathbb{C}, V\rangle,\]
\[\text{lift}(\tree{U},c)=\langle \tree{U}^*,\langle c_\alpha\mid\alpha<\lh(\tree{U})\rangle\] and 
\[c_\alpha=\langle M_\alpha^\tree{U}, \psi_\alpha, Q_\alpha, \mathbb{C}_\alpha, S_\alpha\rangle\rangle.\] 

A reflection argument, as in \cite{nitcis} Lemma 8.4.27, gives that $\tree{U}^*$ has a cofinal wellfounded branch, so it suffices to show that that for any cofinal wellfounded branch $b$ of $\tree{U}^*$, $b=\Sigma_{\mtree{S}, P}(\tree{U})$. Then it follows that $\Sigma_{\mtree{S}, P}(\tree{U})=\Omega_{\nu,k}(\tree{U})$, since the definition of $\Omega_{\nu,k}$ as a (partial) iteration strategy gives $\Omega_{\nu,k}(\tree{U})=b$ iff $b$ is the unique cofinal branch of $\tree{U}^*$. So this would establishes that there are no strategy disagreements, finishing the proof.

So we just need to prove the following lemma.
\begin{lemma}\label{main comparison lemma}
If $b$ is a cofinal wellfounded branch of $\tree{U}^*$, then $\Sigma_{\mtree{S}, P}(\tree{U})=b$.
\end{lemma}
\begin{proof}[Proof.]
We write $(\mtree{S}^*_{\nu,k})^{S_\gamma}$ for $\langle \eta,j\rangle\leq_{\text{lex}} i_{0,\gamma}^{\tree{U}^*}(\langle \nu, k\rangle)$ to denote $i^{\tree{U}^*}_{0,\gamma}(\langle \zeta, l\rangle \mapsto \mtree{S}^*_{\zeta, l})(\eta, j)$. It is easy to see that our hypotheses give that $M$, $\tree{S}$ are fixed by $i_{0,\gamma}^{\tree{U}^*}$ and that $i^{\tree{U}^*}(\Sigma)\subseteq \Sigma$ and $i^{\tree{U}^*}(A)\subseteq A$, so $(\mtree{S}^*_{\nu,k})^{S_\gamma}$ is by $\Sigma, A$.

As $M_b^{\tree{U}^*}$ is wellfounded, in addition to the $c_\alpha$, we also have a last conversion stage
\[c_b=\langle M_b^\tree{U}, \psi_b, Q_b, \mathbb{C}_b, S_b\rangle\]
in $\text{lift}(\tree{U}\conc b, c)$. For $\gamma<\lh(\tree{U})$ or $\gamma=b$, let 
Let\begin{align*}
		\langle \eta_\gamma, l_\gamma\rangle&=\text{ the unique }\langle \eta, l\rangle\text{ such that }Q_\gamma=M_{\eta, l}^{\mathbb{C}_\gamma},\\
			\mtree{S}^*_\gamma&=(\mtree{S}^*_{\eta_\gamma, l_\gamma})^{S_\gamma},\\
			  N_\gamma&=M_\infty^{\mtree{S}^*_\gamma}.
		\end{align*}
	
		So $Q_\gamma\is N_\gamma$, the last model of $\mtree{S}^*_\gamma$, which is the unique $\lambda$-separated meta-tree on $\tree{S}$ by $\Sigma,A$ which iterates $\tree{S}$ past $Q_\gamma$, and if $\nu<_\tree{U}\gamma$, and $(\nu, \gamma]_\tree{U}$ doesn't drop, then $i^{\tree{U}^*}_{\nu, \gamma}(\mtree{S}^*_\nu)=\mtree{S}^*_\gamma$.
		
		Note that $\mtree{S}^*_0$ is our initial meta-tree $\mtree{S}$. For simplicity, we will assume that $Q_0=N_0$, so that $\tree{U}=\tree{U}^+$ (this makes little difference).

		Let $\tree{T}$ be the last tree of $\mtree{S}=\mtree{S}^*_0$. Let $\mtree{U}=\mtree{V}(\tree{T}, \tree{U}\conc b)$. By the definition of $\Sigma_{\mtree{S}, P}$, to show $b=\Sigma_{\mtree{S}, P}(\tree{U})$, it suffices to show $\langle \mtree{S}, \mtree{U}\rangle$ is by $\Sigma$. To do this, it suffices to show $\mtree{V}(\mtree{S}, \mtree{U})$ is by $\Sigma$, since $\Sigma$ normalizes well. Since $\Sigma$ has meta-hull condensation this follows immediately from the following claim.
		\begin{sublemma}\label{main sublemma}
			There is an extended meta-tree embedding from $\mtree{V}(\mtree{S}, \mtree{U})$ into $\mtree{S}^*_b$.
		\end{sublemma}
		
	\begin{proof}
	
	\end{proof}	 We'll have to consider stages of the meta-tree normalization $\mtree{V}(\mtree{S},\mtree{U})$. Set	\begin{align*}
	\mtree{V}_\gamma 
			&= \mtree{V}(\mtree{S}, \mtree{U}\restrict (\gamma+1))
			\text{ for $\gamma < \lh (\mtree{U})$, and}\\
			\mtree{V}_b
			&= \mtree{V}(\mtree{S}, \mtree{U}).		\end{align*}
		
		For $\gamma<\lh(\tree{U})$ or $\gamma=b$, let $R_\gamma$ the last model of $\mtree{V}_\gamma$. For $\gamma<\lh(\tree{U})$, we have that the last tree of $\mtree{V}_\gamma$ is $V(\tree{T}, \tree{U}\restrict \gamma+1)$ and the last tree of $\mtree{V}_b$ is $V(\tree{T}, \tree{U}\conc b)$. In particular, for $\gamma<\lh(\tree{U})$ or $\gamma=b$, the models $R_\gamma$ are the last models of these quasi-normalization trees so that the associated quasi-normalization maps $\sigma_\gamma$ are nearly elementary maps from $M_\gamma^\tree{U}$ into $R_\gamma$. For $\gamma<\lh (\tree{U})$, let $F_\gamma=\sigma_\gamma(E^\tree{U}_\gamma)$, $a_\gamma = a(F_\gamma, \mtree{V}_\gamma)$, and $b_\gamma=b(F_\gamma, \mtree{V}_\gamma)$.
		
		If $\nu\leq_\tree{U} \gamma$, let $\vec{\Phi}^{\nu, \gamma}$ the partial extended meta-tree embedding from $\mtree{V}_\nu$ into $\mtree{V}_\gamma$ and $\vec{I}^{\nu, \eta}$ the extended meta-tree embedding from $\mtree{S}_\nu^*$ into $\mtree{S}_\gamma^*$ coming from Proposition \ref{hull embeddings prop 2}.
		
		Let $H_\gamma= \psi_\gamma (E^\tree{U}_\gamma)$, $res_\gamma = (\sigma_{\eta_\gamma, l_\gamma}[M_{\eta_\gamma,l_\gamma}| \lh (H_\gamma)])^{S_\gamma}$,  $G_\gamma = res_\gamma(H_\gamma)$ and $G^*_\gamma = E^{\tree{U}^*}_\gamma$. So $G_\gamma$ comes from resurrecting $P= N_\gamma|\lh (H_\gamma)$ inside $S_\gamma$ and $G_\gamma^*$ is the corresponding background extender. Let $\tau_\gamma$ least such that $P$ is an initial segment of the last model of $\tree{S}_{\tau_\gamma}^{\mtree{S}^*_\gamma}$.
		
		By induction on $\gamma<\lh (\tree{U})$ or $\gamma=b$, we build extended meta-tree embeddings $\vec \Delta^\gamma=\langle u^\gamma, v^\gamma, \{\Psi^\gamma_\xi\}, \{\Delta^\gamma_\xi\}\rangle: \mtree{V}_\gamma \to \mtree{S}_\gamma^*$ maintaining that for all $\nu<\eta\leq \gamma$,
\begin{enumerate}
    \item $\vec \Delta^\nu \restrict a_\nu+1 \approx \vec \Delta^\eta \restrict a_\nu +1$
    \item if $\nu\leq_\tree{U} \eta$ and $\nu$-to-$\eta$ doesn't drop, then $\vec{\Delta}^\eta\circ \vec{\Phi}^{\nu, \eta}= \vec{I}^{\nu, \eta}\circ \vec{\Delta}^\nu$ ;
    \item For $s^\eta, t^\eta$ the last $s$-map and $t$-map of $\Delta^\eta_{\infty}$, respectively,
    \begin{itemize}
        \item[(a)]$s^\eta \restrict\lh (F_\nu) +1 = res_\nu \circ t^\nu \restrict \lh (F_\nu) +1$, and
        \item[(b)]$\psi_\eta = t^\eta\circ \sigma^\eta$;
    \end{itemize}
    \item $G_\nu$ is a meta-tree exit extender of $\mtree{S}^*_\eta$ and for $\xi_\eta$ such that $G_\nu = F_{\xi_\eta}^{\mtree{S}^*_\eta}$, 
    \begin{itemize}
        \item [(a)]$\tau_\nu\leq_{\mtree{S}^*_\eta}{\xi_\nu}\leq_{\mtree{S}^*_\eta}v^{\eta}(a_\eta)$,
        \item [(b)] for $t^{\tau_\nu, \xi_\nu}$ the last $t$-map along $(\tau_\nu,{\xi_\nu}]_{\mtree{S}^*_\eta}$, $t^{\tau_\nu, \xi_\nu}\restrict \lh (H_\nu) +1 = res_\nu \restrict \lh (H_\nu) +1$.
    \end{itemize}
\end{enumerate}

The bulk of the work is seeing that these hypotheses pass through the successor case. In our sketch, we'll just handle the non-dropping successor case, omitting discussion about the dropping successor case and limit case.

Let $\nu=\tree{U}\pred(\gamma+1)$ and suppose (1)-(4) hold up to $\gamma$.
Applying Lemma \ref{res map lemma 2} to $P=N_\gamma|\lh (H_\gamma)$ in $S_\gamma$ with $\langle \theta, j\rangle = \langle \eta_\gamma, l_\gamma\rangle$, we get a meta-tree $\mtree{T}_\gamma = (\mtree{S}^*_{\theta_0,j_0})^{Q_\gamma}$ such that $\mtree{T}_\gamma\restrict\tau_\gamma+1= \mtree{S}^*_\gamma\restrict \tau_\gamma+1$\footnote{Here we are relying on the definition $\mtree{S}^*_\gamma= (\mtree{S}^*_{\eta_\gamma, l_\gamma})^{S_\gamma}$.}, $\tau_\gamma \leq_{\mtree{T}_\gamma} \xi_\gamma \defeq \lh (\mtree{T}_\gamma)-1$ 
and for $t^{\tau_\gamma, \xi_\gamma}_\infty$ the last $t$-map of $\Phi^{\mtree{T}_\gamma}_{\tau_\gamma,\xi_\gamma}$, 
$t^{\tau_\gamma, \xi_\gamma}_\infty\restrict P =  res_\gamma$. Note that since $G_\gamma= res_\gamma(H_\gamma)$, $G_\gamma^-$ is on the sequence of the last model of $\mtree{T}_\gamma$. Let $N^*_\gamma =M_{\theta_0,j_0}^{S_\gamma}$.

Now suppose $(\nu, \gamma+1]_\tree{U}$ doesn't drop. As mentioned above, this is the only case we'll consider. It contains most of the ideas needed for the droping case but is somewhat simpler.

\begin{claim} $\mtree{T}_\gamma\is \mtree{S}^*_{\gamma+1}$ and
    $G_\gamma=F^{\mtree{S}^*_{\gamma+1}}_{\xi_\gamma}$ is the first factor of $\vec{I}^{\nu,\gamma+1}$.
\end{claim}
\begin{proof} 

We first show that $N_{\gamma+1} ||\lh (G_\gamma) = N^*_\gamma|| \lh (G_\gamma)$.  Let $\mu =\crit (F_\gamma)$, and $\bar\mu =\crit (E^\tree{U}_\gamma)$. By our case hypothesis, $E^\tree{U}_\gamma$ is total on $M^\tree{U}_\nu$, so no level of $M^\tree{U}_\nu$ beyond $\lh (E^\tree{U}_\nu)$ projects to or below $\bar\mu$.

So, applying $\sigma_\nu$, no level of  $R_\nu$ beyond $\lh (F_\nu)$ projects to or below $\mu = \sigma_\nu(\bar\mu)$ (using here that $\sigma_\nu$ agrees with $\sigma_\gamma$ up to $\lh (F_\nu)+1> \mu$).

Applying $t^\nu$ and using our induction hypothesis (3)(b) at $\nu$, no level of $N_\nu$ beyond $\lh (H_\nu)$ projects to or below $t^\nu(\mu)$ and so $res_\nu$ is the identity on ${t^\nu(\mu)^+}^{N_\nu}$.

Now since $G_\nu = res_\nu(H_\nu)$, $\crit (G_\nu) =res_\nu(t^\nu(\mu))=t^\nu(\mu)$ and so ${t^\nu(\mu)^+}^{N_\nu^*}<\hat\lambda(G_\nu)$.
So we get
\[ N_\nu|{t^\nu(\mu)^+}^{N_\nu}=N^*_\nu|{t^\nu(\mu)^+}^{N_\nu^*}=N_\gamma|{t^\nu(\mu)^+}^{N_\gamma}\]

Let $\lambda$ be this common value of $t^\nu(\mu)^+$. If $\nu=\gamma$, the above shows that $N_\gamma|\lambda= N^*_\gamma|\lambda$. If $\nu<\gamma$, then no proper initial segment of $M_\gamma^\tree{U}$ projects to or below $\lh (E_\nu^\tree{U})$, so no proper initial segment of $N_\gamma$ projects to or below $\lh (H_\nu)$, so $res_\gamma$ is identity 
on all of $\lh (H_\nu)$. So even when $\nu<\gamma$, we get $N_\gamma|\lambda= N^*_\gamma|\lambda$.

Applying $\pi=i^{S_\nu}_{G^*_\gamma}$, we get
\[\pi(N_\gamma|\lambda)= \pi(N^*_\gamma|\lambda)\]
Now by our choice of background extender $G^*_\gamma$, $\pi(N^*_\gamma)|\lh (G_\gamma)+1 = Ult(N^*_\gamma, G_\gamma)|\lh (G_\gamma)+1$. By our case hypothesis, $N_{\gamma+1}= \pi(N^*_\nu)$, so by the agreement between the models identified above gives
\[N_{\gamma+1}||\lh (G_
\gamma) = N^*_\gamma||\lh (G_\gamma),\]
as desired.
It follows that $\mtree{S}^*_{\gamma+1}$ and $\mtree{T}_\gamma$ use the same meta-tree exit extenders of length $<\lh (G_\gamma)$. 
But $G_\gamma^-$ is the top extender of $N^*_\gamma$, so $\mtree{T}_\gamma$ uses no extenders of length $\geq \lh(G_\gamma)$, hence $\mtree{T}_\gamma\is \mtree{S}^*_{\gamma+1}$. As $G_\gamma^-$ is on the sequence of the last model of $\mtree{T}_\gamma$ but not on the $N_{\gamma+1}|\lh (G_\gamma)+1=Ult(N^*_\gamma, G_\gamma)|\lh (G_\gamma)+1$-sequence, the first new exit extender of $\mtree{S}^*_{\gamma+1}$, $F^{\mtree{S}^*_{\gamma+1}}_{\xi_\gamma}$, must have length $\leq \lh (G_\gamma)$. But then we must have $F^{\mtree{S}^*_{\gamma+1}}_{\xi_\gamma} = G_\gamma$, since $\mtree{S}^*_{\gamma+1}$ and $\mtree{T}_\gamma$ use all the same extenders of shorter lengths.

To finish the claim, we just need to see that the meta-tree embedding Factor Lemma applies to $\vec{I}^{\nu, \gamma+1}$ and the $G_\gamma$ is the corresponding factor. Using that the $u$-map of $\vec{I}^{\nu, \gamma+1}$ and $u$-maps and $t$-maps of the $\Delta$-maps of $\vec{I}^{\nu, \gamma+1}$ are restrictions of $\pi$, one can show that the first meta-tree extender used along $(\kappa, \pi(\kappa)]_{\mtree{T}_\gamma}$ is $G_\gamma$ (since it is compatible with $G_\gamma$ and distinct compatible extenders can't be used in a meta-tree, since they can't be used in a plus tree). But then the agreement between $\mtree{T}_\gamma$ and $\mtree{S}_\nu^*$ guarantee that the Factor Lemma applies and so $G_\gamma$ is the desired factor.
\qed
\end{proof}

Since $G_\gamma$ is the first factor of $\vec{I}^{\nu, \gamma+1}$, we let $\vec{\Pi}^{\gamma+1}$ be the extended meta-tree embedding from $\mtree{V}(\tree{S}^*_\nu, \mtree{T}_\gamma, G_\gamma)$ into $\mtree{S}^*_{\gamma+1}$ such that $\vec{I}^{\nu,\gamma+1}= \vec{\Pi}^{\gamma+1}\circ \vec{\Delta}^{\mtree{V}(\tree{S}^*_\nu, \mtree{T}_\gamma, G_\gamma)}$ and $u^{\vec{\Pi}}\restrict \xi_\gamma+1=id$ guaranteed by the meta-tree embedding Factor Lemma (using here that $a(\mtree{S}^*_{\gamma+1}, G_\gamma)=\xi_\gamma$). 

By our induction hypothesis $(3)(b)$ at $\gamma$, $H_\gamma=t^\gamma(F_\gamma)$, so $\vec{\Delta}^\gamma\restrict a_\gamma+1: \mtree{V}_\gamma\restrict a_\gamma+1\to \mtree{S}^*_\gamma\restrict{\tau_\gamma}+1\is \mtree{T}_\gamma$ is a meta-tree embedding. 
Since $\tau_\gamma\leq_{\mtree{T}_\gamma} \xi_\gamma=\lh (\mtree{T}_\gamma)-1$ and $\tau_\gamma\in [v^\gamma(a_\gamma), u^{\gamma}(a_\gamma)]_{\mtree{S}^*_\gamma}$ (using here that these meta-trees are $\lambda$-separated), we can view $\vec{\Delta}^\gamma\restrict a_\gamma+1$ as an extended meta-tree embedding from $\mtree{V}_\gamma\restrict a_\gamma+1$ into $\mtree{T}_\gamma$, which we'll denote $(\vec{\Delta}^\gamma)^*$. That is, we let the last $\Delta$-tree embedding of $(\vec{\Delta}^\gamma)^*$ be $\Phi^{\mtree{T}_\gamma}_{v^\gamma(a_\gamma), \xi_\gamma}\circ \Gamma^\gamma_{a_\gamma}$. 
Now, one can show that the last $t$-map of $\Phi^{\mtree{T}_\gamma}_{V^\gamma(a_\gamma), \tau_\gamma}\circ \Gamma^\gamma_{a_\gamma}$ maps $F_\gamma$ to $H_\gamma$. Since the last $t$-map of $\Phi^{\mtree{T}_\gamma}_{\tau_\gamma, \xi_\gamma}$ agrees with $\text{res}_\gamma$ on $N_\gamma|\lh (H_\gamma)$ (and $\text{res}_\gamma(H_\gamma)=G_\gamma$), we have that $G_\gamma$ is the image of $F_\gamma$ under the last $t$-map of $(\vec{\Delta}^\gamma)^*$.

By our induction hypotheses and observations we've already made about $(\vec{\Delta}^\gamma)^*$, it is straightforward to check that

\begin{claim}
The meta-tree Shift Lemma applies to $((\vec{\Delta}^\gamma)^*$, $\vec{\Delta}^\nu$, $F_\gamma$, $G_\gamma)$.
\end{claim}
So, letting $\vec{\Psi}^{\gamma+1}:\mtree{V}_{\gamma+1}\to\mtree{V}(\mtree{S}^*_\nu, \mtree{T}_\gamma, G_\gamma)$ be the associated copy meta-tree embedding, we finally set $\vec{\Delta}^{\gamma+1}= \vec{\Pi}^{\gamma+1}\circ \vec{\Psi}^{\gamma+1}$. This is clearly an extended meta-tree embedding from $\mtree{V}_{\gamma+1}$ into $\mtree{S}^*_{\gamma+1}$ and it is straightforward to verify our induction hypotheses are maintained. This finishes the non-dropping successor case and our sketch of the proof of Sublemma \ref{main sublemma}, thereby finishing the proof of Lemma \ref{main comparison lemma} and Theorem \ref{main comparison theorem}. \qed\end{proof}

The next theorem is our main application of Theorem \ref{main comparison theorem}: 
a characterization of meta-strategies of the form $\Sigma^*$.

\begin{theorem}\label{induced strategy theorem} Assume $\adp$. Let $(M,\Sigma)$ be 
a strongly stable mouse pair with scope $\hc$. Let $\tree{S}$ be a countable plus
 tree of successor length by $\Sigma$. Let $\Lambda$ be a meta-strategy 
for finite stacks of meta-trees on $\tree{S}$ which has meta-hull condensation,
 normalizes well, has the Dodd-Jensen property relative
 to $\Sigma$, and is pushforward consistent, if $(M,\Sigma)$ is an lbr hod pair.

Then $\Lambda=\Sigma^*_\tree{S}$.
\end{theorem}

This is a consequence of the following lemma.

\begin{lemma}\label{induced strategy lemma}
Assume $\adp$. Let $(M,\Sigma)$ be a strongly stable mouse pair with scope $\hc$.
 Let $\tree{S}$ be a countable plus tree of limit length by $\Sigma$.

Suppose $c$ is a cofinal, wellfounded branch of $\tree{S}$ and $\Lambda$ 
is a meta-strategy for finite stacks of meta-trees on $\tree{S}\conc c$ which has meta-hull condensation,
 normalizes well, has the Dodd-Jensen property relative to $\Sigma$, and is pushforward consistent,
 if $(M,\Sigma)$ is a lbr hod pair.

Then $\Sigma(\tree{S})=c$.\end{lemma}

In the course of the proof, we shall make use of the following fact about meta-trees, which is a variant of Lemma 6.6.14 of \cite{nitcis} with essentially the same proof.

\begin{proposition}\label{incompatibility prop}
Let $\mtree{T}$ be a meta-tree on a plus tree $\tree{T}$ and $\gamma,\delta<\lh(\mtree{T})$ $\leq_\mtree{T}$-incomparable. Let $\eta=\sup([0,\gamma)_\mtree{T}\cap [0,\delta)_\mtree{T})$, $\bar\alpha\geq \crit(u^{\Phi_{\eta,\gamma}^\mtree{T}})$, $\bar\epsilon\geq\crit(u^{\Phi_{\eta, \delta}^\mtree{T}})$, $\alpha=u^{\Phi_{\eta,\gamma}^\mtree{T}}(\bar\alpha)$, and $\epsilon=u^{\Phi_{\eta,\delta}^\mtree{T}}(\bar\epsilon)$. Then $e_{0,\alpha}^{\tree{T}_\gamma} \bot e_{0,\epsilon}^{\tree{T}_\delta}$.
\end{proposition}

\begin{proof}[Proof of Lemma \ref{induced strategy lemma}.]
By the Basis Theorem we may assume $\Sigma$, $\Lambda$ are Suslin-co-Suslin (in the codes) and work in an appropriate coarse $\Gamma$-Woodin tuple such that $\Sigma$ and $\Lambda$ are coded by sets of reals in $\Gamma$.\\

Fix $\tree{S}$, $c$, and $\Lambda$ as above. Let $b=\Sigma(\tree{S})$. We'll apply Theorem \ref{main comparison theorem} twice: once to $\tree{S}\conc b$ with the induced meta-strategy $\Sigma^*_{\tree{S}\conc b}$ and once to $\tree{S}\conc c$ with $\Lambda$, in each case using some sets $A_b$ and $A_c$ which use information coming from the other meta-strategy.

The idea behind the choices of these sets is to capture the following rules for simultaneously comparing $\tree{S}\conc b$ and $\tree{S}\conc c$ against levels of the background.  Fix a level $\langle \nu, k\rangle$. We build meta-trees $\mtree{S}^b$ and $\mtree{S}^c$ by $\Sigma^*$ and $\Lambda$, respectively, in stages. At every stage we'll drop to the shortest tree which witnesses a disagreement between $\Sigma^*$ and $\Lambda$, in the following sense. At successor stages, given $\mtree{S}^b\restrict \xi+1$ and $\mtree{S}^c\restrict \zeta +1$, we extend the meta-trees by hitting the least disagreements $F^b_\xi$ and $F^c_\zeta$ between the current last models and $M_{\nu,k}$ (using here that least disagreements are extender disagreements which don't come from the background). If the new main branches $[0,\xi+1]_{\mtree{S}^b}$ and $[0,\zeta+1]_{\mtree{S}^c}$ have the same exit extenders (in particular, $F^b_\xi=F^c_\zeta$), then we extend $\mtree{S}^b\restrict \xi+1$ and $\mtree{S}^c\restrict \zeta +1$ by these extenders without gratuitously dropping. Otherwise, we drop to  $\alpha_0(\tree{S}^b_\xi, F^b_\xi)+2$ and $\alpha_0(\tree{S}^c_\zeta, F^c_\zeta)+2$ of these extenders (this may be a gratuitous drop or just a necessary drop).
Similarly, if we're at a stage where we've defined $\mtree{S}^b\restrict \xi$ and $\mtree{S}^c\restrict\zeta$ for limit ordinals $\xi$ and $\zeta$, we extend these meta-trees by choosing $\Sigma^*(\mtree{S}^b\restrict \xi)$ and $\Lambda(\mtree{S}^c\restrict \zeta)$ and gratuitously dropping to the supremum of the respective $\alpha_0$'s whenever these branches use different exit extenders.

We can capture this simultaneous comparison by considering two applications of Theorem \ref{main comparison theorem} using appropriately chosen sets $A_b$ and $A_c$.

First, we apply Theorem \ref{main comparison theorem} to $\tree{S}\conc b$, $\Sigma^*$, and $A_b$ the set of pairs of trees and ordinals $\langle \tree{T}, \gamma\rangle$ such that there is a meta-tree $\mtree{T}=\langle \{\tree{T}_\eta\}, \{F_\eta\}\rangle$ by $\Sigma^*$ on $\tree{S}\conc b$ of length $\xi+1$ such that \begin{enumerate}
    \item[(i)] $\tree{T}=\tree{T}_\xi=\tree{T}_\xi^+$, \item[(ii)]$\gamma=\sup\{\alpha_\eta^\mtree{T}+1\,|\,\eta<\xi\}<\lh (\tree{T})$, and \item[(iii)] for $P$ the last model of $\tree{T}$ and for any meta-tree $\mtree{U}$ by $\Lambda$ with some last model $Q$ such that $P|\delta(\mtree{T}) = Q|\delta(\mtree{T})$,
     
        the exit extenders used along the main branch of $\mtree{T}$ are different from the exit extenders used along the main branch of $\mtree{U}$.
  
\end{enumerate}
The closure properties of Suslin sets guarantee that $A_b$ is also Suslin-co-Suslin (in the codes), since $\Sigma$  and $\Lambda$ are. Moreover, we can choose our pointclass $\Gamma$ sufficiently large at the outset so that $A_b$ (and also $A_c$, which is defined similarly) are is coded by sets of reals in $\Gamma$. \\

Now let $\mtree{S}^b_{\nu,k}$ be the resulting meta-trees by $\Sigma^*,A_b$ and, via the universality argument, let $\langle\nu_b,k_b\rangle$ least such that $(\tree{S}\conc b,\Sigma^*,A_b)$ iterates to $(M_{\nu_b,k_b}, \Omega_{\nu_b, k_b})$.
Similarly, applying Theorem \ref{main comparison theorem} to $\mtree{S}^c$, $\Lambda$, and the set $A_c$ which is defined just like $A_b$ but with the roles of $(\tree{S}\conc b, \Sigma^*)$ and $(\tree{S}\conc c, \Lambda)$ switched, we let $\mtree{S}^c_{\nu,k}$ be the resulting meta-trees by $\Lambda,A_c$ and $\langle \nu_c, k_c\rangle$ least such that $(\tree{S}\conc c, \Lambda, A_c)$ iterates to $(M_{\nu_c,k_c}, \Omega_{\nu_c, k_c})$

Let $\langle \nu,k\rangle$ be the lexicographic minumum of $\langle
\nu_b, k_b\rangle$ and $\langle \nu_c, k_c\rangle$. Let $\mtree{S}^b= \mtree{S}^b_{\nu,k}$ and $\mtree{S}^c=\mtree{S}^c_{\nu,k}$.
Without loss of generality we assume that $\langle \nu,k\rangle=\langle \nu_b, k_b\rangle$. We then have that $\mtree{S}^b$ has no necessary drops along its main branch and the last model of $\mtree{S}^b$ is $M_{\nu,k}$. Let $b^*$ and $c^*$ be the main branches of $\mtree{S}^b$ and $\mtree{S}^c$ and $\tree{T}_b$ and $\tree{T}_c$ be their last trees. Let $\Phi_{b}$ and $\Phi_{c}$ be the possibly partial tree embeddings from $\tree{S}\conc b$ into $\tree{T}_b$  and $\tree{S}\conc c$ into $\tree{T}_c$. Finally, let $b^*$ and $c^*$ be the main branches of $\tree{T}_b$ and $\tree{T}_c$.

Ultimately, we'll show the main branches of $\mtree{S}^b$ and $\mtree{S}^c$ use the same meta-tree exit extenders, and use this to show $b=c$. So let $\langle F_\eta\mid \eta<\gamma\rangle$ and $\langle G_\eta\mid \eta<\lambda\rangle$ enumerate the meta-tree exit extenders used along the main branch of $\mtree{S}^b$ and $\mtree{S}^c$, respectively, in increasing order. We show by induction that $F_\eta=G_\eta$ (and ultimately that $\gamma=\lambda$).

Let $\xi_\eta$ such that $F_\eta=F_{\xi_\eta}^{\mtree{S}^b}$ and $\zeta_\eta$ such that $G_\eta=F_{\zeta_\eta}^{\mtree{S}^c}$. Suppose we've shown that $F_\eta=G_\eta$ for all $\eta<\chi$. Let $\xi=\sup\{\xi_\eta+1\mid \eta<\chi\}$ and $\zeta=\sup\{\zeta_\eta\mid \eta<\chi\}$, so that $\xi$ is on the main branch of $\mtree{S}^b$ and $\zeta$ on that of $\mtree{S}^c$. Since at most one side of our comparison has a necessary drop along its main branch, it follows that there can be no necessary drops at all coming from the $F_\eta=G_\eta$ for $\eta<\chi$, i.e. $0$-to-$\xi$ doesn't have a necessary drop in $\mtree{S}^b$ and $0$-to-$\zeta$ doesn't have a necessary drop in $\mtree{S}^c$. By our choice of $A_b$ and $A_c$, we get that there are no gratuitous drops along these branches so far, either. 

Now we want to show $F_\chi=G_\chi$. We consider cases.

\paragraph{Case 1.} $\xi$-to-$\xi_\chi+1$ drops in $\mtree{S}^b$.

\paragraph{Subcase 1A.} $\zeta$-to-$\zeta_\chi+1$ drops in $\mtree{S}^c$.\\

First suppose that $c^*$ doesn't drop. Using the Dodd-Jensen property relative to $\Sigma$, we get that $b^*$ cannot drop either, $M_{\infty}^{\tree{T}_b}=M_\infty^{\tree{T}_c}$ and $i_{b^*}^{\tree{T}_b}=i_{c^*}^{\tree{T}_c}$. Since the exit extender sequences were the same below $\chi$ and both sides dropped (perhaps gratuitously), we actually must have $F_\chi$ and $G_\chi$ are applied to a common model along $b^*$ and $c^*$ from which it is easy to see that $F_\chi, G_\chi$ are compatible, and so equal by the Jensen ISC.

Now suppose that $c^*$ does drop. Then since its last model is not sound, we must have $M_\infty^{\tree{T}_b}=M_\infty^{\tree{T}_c}$ and $b^*$ also drops. 
Without loss of generality, assume that $\mtree{S}^c$ is the side which doesn't have any necessary drop along its main branch. 
It follows that the last drop in $c^*$ is \textit{before} $\beta^{\mtree{S}^c}_{\zeta_\chi}$ in $c^*$ and is, moreover, in the range of $u^{\Phi_{0,\zeta}^{\mtree{S}^c}}$, say $\eta$ the location of the drop in $\tree{S}$ such that $u^{\Phi_{0,\zeta}^{\mtree{S}^c}}(\eta)$ is the last drop in $c^*$. Since the exit extenders are the same, so far, 
Theorem \ref{tree embedding uniqueness} implies that $M_{u^{\Phi_{0,\xi}^{\mtree{S}^b}}(\eta)}^{\tree{S}^b_\xi}=M_{u^{\Phi_{0,\zeta}^{\mtree{S}^c}}(\eta)}^{\tree{S}^c_\zeta}$ and the common core of $M_\infty^{\tree{T}_b}=M_\infty^{\tree{T}_c}$ is an initial segment of this model. 
One can then show that the models and exit extenders used along $b^*$ and $c^*$ are the same up to the minimum of where $F_\chi$ and $G_\chi$ are applied. As before we can show that $F_\chi$ and $G_\chi$ are actually applied to the same model, so compatible, and the Jensen ISC gives $F_\chi=G_\chi$.

\paragraph{Subcase 1B.} $\zeta$-to-$\zeta_\chi+1$ \textit{doesn't} drop in $\mtree{S}^c$.\\

By our definition of $A_c$, it follows that $G_\chi$ must be a meta-tree exit extender of $\mtree{S}^b$, too. By considering cases about drops in $b^*$ and $c^*$, one gets, again, that the resulting branch embeddings, or a tail of them, are the same, using the Dodd-Jensen property relative to $\Sigma$. But then we can reach a contradiction using Proposition \ref{incompatibility prop}. That proposition gives that the first extender used along the main branch of $\tree{S}^b_{\xi_\chi+1}$ which is not already in $\tree{S}^b_{\xi}$ and the first extender used along the main branch of $\tree{S}^c_{\zeta_\chi+1}$ which is not already in $\tree{S}^c_{\zeta}$ must be incompatible. If we drop at the next stage, this incompatibility must persist, so we do not drop, which gives that $G_{\chi+1}$ is used in $\mtree{S}^b$ (otherwise we must gratuitiously drop in $\mtree{S}^c$ at this next stage). But then Proposition \ref{incompatibility prop} implies that the incompatibility persists anyway, letting us get that $G_{\chi+2}$ is used in $\mtree{S}^b$, and so on. In the end, we get the full $\langle G_\eta\mid \eta<\lambda\rangle$ corresponds to a maximal branch of $\mtree{S}^b$ and can show that, because we chose a different branch, the incompatibility must persist till the end, contradicting that the appropriate branch embeddings are the same.

\paragraph{Case 2.} $\xi$-to-$\xi_\chi+1$ \textit{doesn't} drop in $\mtree{S}^b$.\\

Here we just use the argument from Subcase 1B.\\

Now, without loss of generality assume $\gamma\leq \lambda$. So we've shown $F_\eta=G_\eta$ for all $\eta<\gamma$. Suppose $\gamma<\lambda$ and let $\zeta=\sup\{\zeta_\eta\mid \eta<\gamma\}$. Using Theorem \ref{tree embedding uniqueness} and the Dodd-Jensen property, a tail (or all) of the extenders used in $b^*$ already appear in $\tree{S}^c_\zeta$. But applying the remaining $G_\eta$'s must add new extenders to a tail of $c^*$ which don't already appear in $\tree{S}^c_\zeta$. But the Dodd-Jensen property gives these branches $b^*$ and $c^*$ use common extenders on a tail, a contradiction. So we have $\gamma=\lambda$.

By Theorem \ref{tree embedding uniqueness},  $\Phi_b\restrict \tree{S}\equiv \Phi_c\restrict\tree{S}$. It follows that $v^{\Phi_b}[c]$ generates a cofinal wellfounded branch $d$ of $\tree{T}_b\restrict\sup v^{\Phi_b}[c]$ and $\tree{S}\conc c$ tree embeds into $(\tree{T}_b\restrict\sup v^{\Phi_b}[c])\conc d$. One can show, using the Dodd-Jensen property, that actually $d$ is an initial segment of $b$. So that $(\tree{T}_b\restrict\sup v^{\Phi_b}[c])\conc d\is \tree{T}_b$. In particular, since $\tree{T}_b$ is by $\Sigma$, strong hull condensation gives $\tree{S}\conc c$ is by $\Sigma$, i.e. $c=\Sigma(\tree{S})=b$. 

\qed
\end{proof}

\begin{proof}[Proof of Theorem \ref{induced strategy theorem}.]
This is easy from the previous lemma. Let $\tree{S}$, $\Lambda$ be as in the theorem statement. Suppose $\Sigma^*\neq \Lambda$. Let $\mtree{S}$ of limit length on $\tree{S}$ by both $\Sigma^*$ and $\Lambda$. Let $b=\Sigma^*(\mtree{S})$ and $c=\Lambda(\mtree{S})$. Let $\mtree{S}_b$ and $\mtree{S}_c$ be the gratuitously dropping meta-trees where we drop to the common sup of the $\alpha_0$'s of the meta-tree exit extenders of $\mtree{S}$ on both sides. Then $\mtree{S}_b$ and $\mtree{S}_c$ have last trees some $\tree{T}\conc b^*$ and $\tree{T}\conc c^*$, respectively. By the definition of $\Sigma^*$, $b^*=\Sigma(\tree{T})$.

But then applying Lemma \ref{induced strategy lemma} to $\tree{T}\conc c*$ gives that $\Sigma(\tree{T})=c^*$. So $b^*=c^*$. It follows that $b=c$ by Lemma \ref{meta-strategy lemma}.
\qed
\end{proof}

We finish this section with some applications of these comparison results.

\begin{theorem}\label{normalizeswellthm}
Assume $\adp$. Let $(M,\Sigma)$ be a strongly stable mouse pair with scope $\hc$. Then $\Sigma$ normalizes well.
\end{theorem}

Since we've assumed $\Sigma$ quasi-normalizes well and tails of mouse pairs are mouse pairs, this follows immediately from the following theorem.
\begin{theorem}\label{normalizeswellthm2}
Assume $\adp$. Let $(M,\Sigma)$ be a strongly stable mouse pair with scope $\hc$. Let $\tree{T}$ be a plus tree by $\Sigma$ of successor length and $\tree{S}$ be the normal companion of $\tree{T}$. Then 
\begin{enumerate}
    \item $\tree{S}$ is by $\Sigma$ and \item $\Sigma_\tree{S}=\Sigma_\tree{T}$.
\end{enumerate}
\end{theorem}

\begin{proof}[Proof sketch.]
We'll get this by applying Lemma \ref{induced strategy lemma} and Theorem \ref{induced strategy theorem}. So fix $\tree{T}$ a plus tree by $\Sigma$ of successor length and $\tree{S}$ the normal companion of $\tree{T}$. Meta-trees on $\tree{S}$ correspond in a simple way to those on $\tree{T}$; for example one can show that if we take $F$ from the sequence of the common last model of $\tree{S}$ and $\tree{T}$, then $V(\tree{S}, F)=V(\tree{T},F)$. One can use this correspondence to define a meta-strategy $\Lambda$ for meta-trees on $\tree{S}$ determined in a simple way by the meta-strategy $\Sigma^*_\tree{T}$. Using the corresponding properties of $\Sigma^*_\tree{T}$, one can show that $\Lambda$ has meta-hull condensation, normalizes well, has the Dodd-Jensen property relative to $\Sigma$, and,  in the case that $(M,\Sigma)$ is an lbr hod pair, is pushforward consistent. Moreover, one gets that if $\tree{U}$ is a normal tree on $M_\infty^\tree{T}=M_\infty^\tree{S}$, then $\tree{U}$ is by $\Sigma_\tree{T}$ iff $\mtree{V}(\tree{S}, \tree{U})$ is by $\Lambda$.

We can then use Lemma \ref{induced strategy lemma} to conclude that the normal companion of any plus tree by $\Sigma$ is by $\Sigma$, as follows. Suppose $\tree{T}$ is a plus tree of limit length by $\Sigma$ such that its normal companion $\tree{S}$ is also by $\Sigma$. Let $b=\Sigma(\tree{T})$ and $c$ the unique cofinal wellfounded branch of $\tree{S}$ such that $\tree{S}\conc c$ is the normal companion of $\tree{T}\conc b$. Then we can apply Lemma \ref{induced strategy lemma} to $\tree{S}\conc c$ together with the meta-strategy $\Lambda$ to get $\Sigma(\tree{S})=c$. This gives conclusion (1).

Now assuming conclusion (1) we can easily get conclusion (2) by Theorem \ref{induced strategy theorem}. Let $\tree{T}$ be a plus tree of successor length by $\Sigma$ and $\tree{S}$ its normal companion. Since $\tree{S}$ is by $\Sigma$, from (1), Theorem \ref{induced strategy lemma} applies to $\tree{S}$ together with $\Lambda$, so that $\Lambda=\Sigma^*_\tree{S}$.
Now we want to show that $\Sigma_\tree{S}=\Sigma_\tree{T}$. For this, it suffices to show the two strategies agree on normal trees. So let $\tree{U}$ be a normal tree on $M_\infty^\tree{S}=M_\infty^\tree{T}$. Then, by quasi-normalizing well, $\tree{U}$ is by $\Sigma_\tree{S}$ iff $\mtree{V}(\tree{S},\tree{U})$ is by $\Sigma^*_\tree{S}=\Lambda$. But one property of $\Lambda$ we promised to verify is that $\tree{U}$ is by $\Sigma_\tree{T}$ iff $\mtree{V}(\tree{S},\tree{U)}$ is by $\Lambda$, we we get $\tree{U}$ is by $\Sigma_\tree{S}$ iff $\tree{U}$ is by $\Sigma_\tree{T}$, as desired.
\qed
\end{proof}

Next we show that, in the $\adp$ context, iteration strategies for mouse pairs are totally determined by their action on $\lambda$-tight normal trees.

Let $\tree{T}$ be a normal tree on a premouse $M$. One can define 
the \textit{$\lambda$-tight companion} $\tree{S}$ of $\tree{T}$, a 
$\lambda$-tight tree on $M$ with the same last model and branch embedding
 as $\tree{T}$. We will include a formal definition in a
 subsequent draft but here is the basic idea. One takes 
$\tree{T}$ and considers the (possibly non-quasi-normal) iteration
 tree $\hat{\tree{T}}$ on $M$ which splits up each application of a 
plus extender in $\tree{T}$ into two steps: first using the minus 
extender, and then using the order zero measure on the image of 
the critical point (so the difference between $\tree{T}$ and 
$\hat{\tree{T}}$ being that one considers these both as genuine 
exit extenders in $\hat{\tree{T}}$). One gets $\tree{S}$ as 
$W(\langle M\rangle, \hat{\tree{T}})$, the result of a meta-tree-like 
process where we follow the tree-order of $\hat{\tree{T}}$ and use
 as meta-tree exit extenders the extenders of $\hat{\tree{T}}$.
 One can show that the resulting tree is $\lambda$-tight and 
that the steps embedding normalization coincide with full normalization, 
so $\tree{S}$ has the same last model and branch embedding of $\tree{T}$
 (because it has the same last model and branch embedding of $\hat{\tree{T}}$).

Meta-trees on $\tree{T}$ and its normal companion $\tree{S}$ correspond 
in a sufficiently nice way that a meta-strategie for $\tree{S}$ determines
 a meta-strategy for $\tree{T}$. Moreover, if $\tree{S}$ is by $\Sigma$,
 a nice strategy for the base model $M$, if we start with the
 meta-strategy $\Sigma^*_{\tree{S}}$, the resulting meta-strategy for
 $\tree{T}$ has all of the nice properties need to run the argument of Theorem \ref{normalizeswellthm2} to obtain the following.

\begin{theorem}\label{tight trees}
Assume $\adp$. Let $(M,\Sigma)$ be a strongly stable mouse pair with scope $\hc$. Let $\tree{T}$ be a normal tree by $\Sigma$ of successor length and $\tree{S}$ be the $\lambda$-tight companion of $\tree{T}$. Then 
\begin{enumerate}
    \item $\tree{S}$ is by $\Sigma$ and \item $\Sigma_\tree{S}=\Sigma_\tree{T}$.
\end{enumerate}
\end{theorem}

\noindent Details will be added in a later draft.

\section{Full normalization}
In this section we generalize 
the notion of a tree embedding $\Phi:\tree{S}\to \tree{T}$ by relaxing the requirement
that the images of exit extenders of $\tree{S}$ under the $t$-maps of $\Phi$ are exit extenders of $\tree{T}$. We call
the resulting systems {\em weak tree embeddings}.
We also define the full normalization $X(\tree{T},\tree{U})$ of a stack
of normal trees $\langle\tree{T},\tree{U}\rangle$, and show 
that there are weak tree embeddings from $\tree{T}$ to
$X(\tree{T},\tree{U})$ and from $X(\tree{T},\tree{U})$ to $W(\tree{T},\tree{U})$.
Finally we prove our main theorem: the iteration strategy in a
mouse pair condenses to itself under weak tree embeddings. This 
implies that the strategies in mouse pairs fully normalize well, and are therefore 
positional. 

Condensation under  tree embeddings $\Phi:\tree{S}\to \tree{T}$ has to do
with the structure of the iteration process that produced $\tree{S}$ and $\tree{T}$. Condensation
within the hierarchies of the individual models of $\tree{S}$ and $\tree{T}$
is not relevant; indeed, condensation under tree embeddings
 makes sense for iteration strategies for coarse structures.
In contrast, condensation under weak tree embeddings does involve the condensation
properties of the individual models. We shall need the following theorem in this direction.

\begin{theorem}\label{condensationtheorem}[\cite{trang}] Assume $\adp$, and let
	$(M,\Lambda)$ be a mouse pair with scope $\hc$. Let $H$ be a
	sound premouse of type 1, $\pi \colon H \to M$ be nearly elementary,
and suppose that
	\begin{itemize}
		\item[(1)] $\rho(H) \le \crit(\pi)$, and
		\item[(2)] $H \in M$.
	\end{itemize}
	Then either
\begin{itemize}
	\item[(a)] $H \lhd M$, or
	\item[(b)] $H \lhd \ult(M,E_\alpha^M)$, where $\alpha = \crit(\pi)$.
\end{itemize}
\end{theorem}
\cite{trang} proves a slightly more general result, but \ref{condensationtheorem}
will suffice here.\footnote{
One does not need that $H$ is sound, just that it is $\crit(\pi)$-sound
in an appropriate sense. If $\crit(\pi) < \rho^-(H)$, then the external
strategy also condenses properly, in that $(H,\Lambda^\pi) \lhd (M,\Lambda)$
or $(H,\Lambda^\pi) \lhd \ult((M,\Lambda),E_{\crit(\pi)}^M)$.}

\subsection{Dropdown sequences}
       The connection between exit extenders in a weak tree embedding is
mediated by a dropdown sequence, so we shall need some elementary facts
about such sequences.

\begin{definition}\label{dropdownsequencedef}
Let $Q$ be a pfs premouse
        and $N \lhd Q$. The {\em $N$-dropdown sequence of
        $Q$} is given by
        \begin{itemize}
                \item[(1)] $A_0 = N$,
                \item[(2)] $A_{i+1}$ is the least $B \unlhd Q$ such that
                        $A_i \lhd B$ and $\rho^-(B) < \rho^-(A_i)$.
        \end{itemize}
        We write $A_i = A_i(Q,N)$, and let $n(Q,N)$ be the largest $i$ such that $A_i$ is defined.
        Let also $\kappa_i(Q,N) = \rho^-(A_i(Q,N))$. We also let $A_i(Q,\eta) =
	A_i(Q,Q|\eta)$,  $\kappa_i(Q,\eta) = \kappa_i(Q,Q|\eta)$, and $n(Q,\eta) = n(Q,Q|\eta)$.
\end{definition}

One place that dropdown sequences come up is the following. Suppose that
$Q=M_\xi^\tree{T}$ and $N= M|\lh(E_\xi^\tree{T})$ for some plus tree $\tree{T}$;
then the levels of the $N$-dropdown sequence of $Q$
correspond to levels of $Q$ to which we might apply an extender $E_\delta^\tree{T}$ with $\xi=\tree{T}\pred(\delta+1)$.
More precisely, 
\[
M_{\delta+1}^{*,\itT}  = A_i(Q,N)^-,
\]
 where $i$ is least such that
$\crit(E_\delta^\itT) \le \kappa_i(Q,N)$.\footnote{Because of this, it might be more natural
to let the dropdown sequence consist of the $A_i^-$, rather than the $A_i$.}

Maps that are nearly elementary and exact preserve dropdown sequences. Unfortunately,
we must deal with maps that are not exact, and this adds a small mess.

\begin{lemma}\label{preserverhominus} Let $\pi \colon M\to N$ be nearly elementary
and $\nu \le o(M)$; then
        \begin{itemize}
                \item[(a)] if $\nu < \rho^-(M)$, then $\pi(\nu) < \rho^-(N)$,
                \item[(b)] if $P \lhd M$ and $\nu \le \rho^-(P)$, then $\pi(\nu)
                        \le \rho^-(\pi(P))$, and
                \item[(c)] if $Q \lhd R \unlhd M$, and $\forall P( Q \unlhd P \lhd R) \Rightarrow
                        \nu \le \rho^-(P))$, then $\forall P(\pi(Q) \unlhd P \lhd \pi(R) \Rightarrow
                        \pi(\nu) \le \rho^-(P))$.
        \end{itemize}
\end{lemma}
\begin{proof}
The standard conventions that $\pi(M)=N$, $\pi(M|\langle\ohat(M),k\rangle)=N|\langle\ohat(N),k\rangle$,
        and $\pi(o(M)) = o(N)$ are in force here.

	(a) and (b) are immediate from the definition of near elementarity.
	For part (c), if $R \in M$ the statement is clearly preserved. Otherwise
        we have $R=M|\<\ohat(M),n\>$ for some $n <k(M)$. The statement
         $\forall P \in M( Q \unlhd P \lhd R) \Rightarrow
                        \nu \le \rho^-(P))$ is $\Pi_1$, hence preserved. But $\rho_0(M),...,\rho_{n-1}(M)$
        are also preserved because $\pi$ is nearly elementary, and this covers
                        the remaining $P$.
\end{proof}
If one replaces $<$ by $\le$ in \ref{preserverhominus}(a), then 
it becomes false. Similarly, (c) fails if the quantifier is over
$P \unlhd R$, rather than $P \lhd R$.

\begin{lemma}\label{preservedropdown} Let $\pi \colon M \to X$ be nearly elementary.
Let $N \lhd M$ and $n = n(M,N)$; then
\begin{itemize}
\item[(a)] $n(X,\pi(N)) \in \lbrace n-1, n, n+1 \rbrace$,
\item[(b)] for all $i \le n-1$,
$A_i(X,\pi(N)) = \pi(A_i(M,N))$,
\item[(c)] if $n(X,\pi(N))=n-1$, then $A_n(M,N) = M$
	and $\pi(\rho^-(M)) < \rho^-(X)$,
\item[(d)] if $n \le n(X,\pi(N))$, then
$\pi(A_n(M,N)) = A_n(X,\pi(N))$, and
\item[(e)] if $ n+1 = n(X,\pi(N))$, then
 \begin{itemize}
\item[(i)]$A_n(M,N) \lhd M$,
\item[(ii)]$\rho^-(M)=\rho^-(A_n(M,N))$,
\item[(iii)]
$ \sup \pi``\rho^-(M) \le \rho^-(X)  < \pi(\rho^-(M))$, and
\item[(iv)]$A_{n+1}(X,\pi(N)) = X$.
\end{itemize}
\end{itemize}
\end{lemma}
\begin{proof}
  Let $A_i = A_i(M,N)$, $\kappa_i = \rho^-(A_i)$,
                $B_i = A_i(X,\pi(N))$, and $\mu_i = \rho^-(B_i)$.
                 From Lemma \ref{preserverhominus},
                 we get
        \[
                A_i \lhd M \Rightarrow (i \le n(X,\pi(N)) \wedge \pi(A_i)=B_i \wedge \pi(\kappa_i)=\mu_i).
                \]

Now suppose that $A_{n} \lhd M$. The preservation just noted shows that $n \le n(X,\pi(N))$, as well as
        (b), (c) (vacuously), and (d). We are done if $n(X,\pi(N)) =n$, so suppose not.
        We then  have $\rho^-(B_{n+1}) < \mu_n$. Since
$\forall P( A_n \unlhd P \lhd M) \Rightarrow
                        \kappa_n  \le \rho^-(P))$, Lemma \ref{preserverhominus} implies
                        that $\forall P( B_n \unlhd P \lhd X) \Rightarrow
                        \mu_n \le \rho^-(P))$.
         It follows that $B_{n+1}=X$, and $n(X,\pi(N)) = n+1$.
        This gives us the rest of (a). But also, if $\kappa_n < \rho^-(M)$ then
        $ \mu_n < \rho^-(X)$, so we must have $\kappa_n = \rho^-(M)$, and then
        $\rho^-(X) < \mu_n = \pi(\rho^-(M))$. This proves (e).

So we have proved the lemma when $A_n \lhd M$. Suppose now that  $A_{n} = M.$
Since $\pi(A_{i})=B_{i}$ and $\pi(\kappa_{i}) =
\mu_{i}$ for $i \le n-1$, we have (b), and
                that $n-1 \le n(X,\pi(N))$. If $n(X,\pi(N))=n-1$,
		then $\rho^-(X) \ge \pi(\kappa_{n-1}) > \pi(\kappa_n) =
		\pi(\rho^-(X))$, so we have (c), and we are done.
		So assume $n(X,\pi(N)) \ge n$.  Since $A_n =M$,
                $\forall P( A_{n-1} \unlhd P \lhd M) \Rightarrow
                        \kappa_{n-1} \le \rho^-(P))$. By Lemma \ref{preserverhominus},
                        $\forall P( B_{n-1} \unlhd P \lhd X) \Rightarrow
                        \mu_{n-1} \le \rho^-(P))$.
        It follows that $B_n = X$, and $n(X,\pi(N))=n$. So we have (a). We have (d)
        because $A_n = M$ and $B_n =X$. Clause (e) is vacuously true.
\end{proof}

The mess gets smaller if $\pi$ is almost exact, and disappears if $\pi$ is exact.
Recall here that elementary maps, like the branch embeddings of an iteration tree
or the resurrection maps of a PFS construction, are almost exact. Resurrection maps
are actually exact, as are the branch embeddings of an iteration tree whose
	base model is strongly stable.\footnote{See \cite[\S 4.3]{nitcis}. $M$ is strongly
	stable iff $\eta_{k(M)}^M$ is not measurable by the $M$-sequence. In particular,
	if $k(M)=0$, then $M$ is strongly stable.}
	Factor embeddings, such as the lifting maps of a conversion system,
may fail to be almost exact.

\begin{lemma}\label{preservedropdownelem} Let $\pi \colon M \to X$ be nearly elementary.
Let $N \lhd M$ and $n = n(M,N)$; then
        \begin{itemize}
\item[(a)] If $\pi$ is almost exact, then $n(X,\pi(N)) \in \lbrace  n, n+1 \rbrace$,
and  for all $i \le n$,
$A_i(X,\pi(N)) = \pi(A_i(M,N))$.
\item[(b)] If $\pi$ is exact, then $n(X,\pi(N)) = n$.
\end{itemize}
\end{lemma}
       The proof of \ref{preservedropdownelem} is implicit in that of
       \ref{preservedropdown}.

\subsection{Maps that respect drops}

The connection between exit extenders required by a weak
tree embedding is an abstraction from the following situation.
Let $M$ be a premouse, $\eta<o(M)$, and $E^M_\eta\neq \emptyset$.
Let $F$ be  close to $M$, $\crit F=\mu<\eta$, and suppose $N=\ult(M,F)$
makes sense. Let $\lambda=i^M_F(\eta)$. There is a natural factor embedding
\[\sigma:\ult(M|\eta, F)\to i^M_F(M|\eta)=N|\lambda,\]
given by
\[
\sigma([a,g]^{M|\eta}_F)=[a,g]^M_F
\]
in the case $k(M)=0$. One can use $\sigma$ to
to show that $i_F^{M|\eta}(E)$ is on the $N$ sequence. It is the connection 
between $i_F^{M|\eta}(E)$
and $i_F^M(E)$ provided by $\sigma$ that we wish to abstract.

One shows that $i_F^{M|\eta}(E)$ is on the $N$-sequence 
by factoring
$\sigma$, using the $\ult(P,F)$
for $P^+$ in the $(M,M|\eta)$ dropdown
sequence as factors.
We shall apply Theorem \ref{condensationtheorem} to the natural
embedding from $\ult(P,F)$ into $i^{Q}_F(P)$, where $P^+$ and $Q^+$ are
successive elements of the dropdown sequence. The next two lemmas capture the uses of the condensation theorem here.

The lemmas are complicated by the fact that $\ult(P,F)$ might have type 2.
In that case $\ult(P,F)$ cannot be an initial segment of 
$\ult(Q,F)$, and it is $\mfc_k(\ult(P,F))$ that is an initial segment
of $\ult(Q,F)$, where $k = k(P)$.

\begin{lemma}\label{notes lemma 1}
($\adp$) Let $P$ and $Q$ be countable premice such that 
$P\isneq Q$ and $k(P)>0$.
Suppose that $\rho^-(P)$ is a cardinal of $Q$ and
	that $\rho^-(P) \le \rho^-(Q)$. Suppose that $F$ is an extender that is close
	to $P$ with $\crit (F)=\mu<\rho^-(P)$.
Let
	\[
		\sigma:\ult(P,F)\to i^Q_F(P)
\] be the natural map, and suppose there is a
	$\Lambda$ such that $(i^Q_F(P), \Lambda)$ is a mouse pair; then
either
	\begin{itemize}
		\item[(i)]
$\ult(P,F)\isneq i^Q_F(P)$, or
\item[(ii)] $\ult(P,F)$ has type 2, and for $k=k(P)$,
	$\mfc_k(\ult(P,F)) \isneq i^Q_F(P)$.
	\end{itemize}
\end{lemma}

\begin{proof}
	Let $k=k(P)$ and $n=k(Q)$, and
	\begin{align*}
		i \colon P \to \ult_k(P,F) = R\\
		\intertext{ and }
		j\colon Q \to \ult_n(Q,F)
	\end{align*}
	be the canonical embeddings.
We have the factor map
        \[
		\sigma \colon R \to j(P) 
        \]
        given by
                $\sigma([a,f]_F^{P}) = [a,f]_F^{Q}.$
The function $f$ here is $r\Sigma_k^P$,
	hence $r\Sigma_n^Q$, so the formula makes sense.\footnote{
If $o(P)=o(Q)$, so that $P=Q|\langle o(Q), k\rangle$ where $k<n$, then by
	$j(P)$ we mean $\ult(Q,F) \downarrow k$.}

Let 
	\begin{align*}
		\nu& = \rho_k(P),\\
		\mu&=\crit(F).
\end{align*}
By hypothesis,	$\nu$ is a cardinal of $Q$ and $\nu \le \rho_n(Q)$,
so every $r\Sigma_n^Q$ function $f$ with domain $\mu$ and range bounded in $\nu$
belongs to $P$. Thus
        \[
                \rho_k(R) = \sup i``\nu \le \crit(\sigma).
                \]
		We may assume that $j(P)^k$ is not contained in $\ran(\sigma)$,
	 as otherwise $R=j(P)$,
		so conclusion (i) holds.
But then $\sigma$ witnesses that the reduct
	$R^k$ is a proper initial segment of $j(P)^k$, so $R^k \in j(P)^k$,
		so $R \in j(P)$.

Let \begin{align*}
    \gamma&=\sup i"\nu
    =\rho_k(R),\\
    \intertext{and }
\delta&=j(\nu)
=\rho_k(j(P)).
\end{align*}
	That $\delta=\rho_k(j(P))$ comes from the fact that $j\restrict P$ is sufficiently elementary.

\begin{remark}
These things work out if $\nu=o(P)$, but we don't use that case in our
	applications of Lemma \ref{notes lemma 1}, so we might as well assume that $\nu<o(P)$.
\end{remark}

	\begin{claim}\label{nearelementarityclaim}
	$\sigma \colon R\to j(P)$ is nearly elementary.
\end{claim}
\begin{proof}
	The $k$-th reduct of $P$ is $P^{k}=(P||\nu, A)$, where $A$ is (essentially) $Th^P_k(\nu\cup \{w_k(P)\})$.\footnote{
	$w_k(P) = \langle p_k(P), \nu, \eta_k^P \rangle$, where
	$\eta_k^P$ is the $r\Sigma_k^P$ cofinality of $\rho_k(P)$.}
Let \[B=\bigcup_{\alpha<\nu}i(A\cap \alpha)\]
and \[C=j(A),\]
where $j(A)$ is understood properly if
	$o(Q)=o(P)$.\footnote{In that case, $Th^P_k(\nu\cup\{w_k(P)P\})$ is coded
	into $Th^Q_n(\nu\cup \{p_n(Q)\cup \{ w_n(Q)\})$ since $o(Q)=o(P)$, $k<n$, and $\rho_n(Q)=\nu$.
	We take then $C=Th^{\ult(Q,F)}_k\big(j(\nu)\cup\{j(w_k(P))\}\big)$.}
We have for $\alpha <\nu$,
\begin{align*}\sigma\big(i(A\cap \alpha)\big)&=j(A\cap \alpha),\\
B\cap i(\alpha)&=C\cap j(\alpha).\end{align*}
	So $R^{k}\prec_{\Sigma_0} j(P)^k$. Since $k>0$, $\sigma$ is cardinal preserving.
	Thus $\sigma$, which is the $k$-completion of
	$\sigma\restrict R^{k}=id$, is nearly elementary from $R$ into $j(P)$.
\qed
\end{proof}
\paragraph{Case 1.} $\gamma =\delta$.

Then $R^{k}=j(P)^{k}$, so $R=j(P)$ and conclusion (i) of the lemma holds.

\paragraph{Case 2.} $\gamma<\delta$.

$R^{k}=(j(P)^{k}|\gamma, C\cap \gamma)$, 
so $R^{k}\in j(P)|\delta$, so $R\in j(P)|\delta$.
On the other hand, $\sigma(i(\nu))=\delta$. It follows that $\crit(\sigma)\leq i(\nu)$. Set 
\[\alpha=\crit(\sigma);\]
then
\[
\gamma\leq \alpha\leq i(\nu).
\]

We want to apply the Condensation Theorem \ref{condensationtheorem}.
For this, let $m$ largest such that $\alpha<\rho_m(R)$.
$\gamma=\rho_k(R)$, so $m<k$. Let
\begin{align*}
	X&=R \downarrow m,\\
	Y &= j(P) \downarrow m.
\end{align*}
Then $\sigma:X\to Y$ is nearly elementary (in fact, elementary)
and $\gamma = \rho(X)\leq \alpha<\rho^-(X)$. $X$ is a premouse and $X \in Y$.
Let us suppose first that $X$ is sound. (If not, then $R$ has type 2,
and $X=R^-$ is only almost sound.)
Letting $\Omega=\Lambda_{Y}$, \[\sigma:(X, \Omega^\sigma)\to (Y, \Omega)\]
is nearly elementary in the category of mouse pairs. We get from \ref{condensationtheorem}
that either
\begin{itemize}
    \item[(a)] $X \isneq Y$, or
    \item[(b)] $X \isneq \ult_0(Y, E^{Y}_\alpha)$.
\end{itemize}
We will be done if we can rule out (b).

\paragraph{Subcase 2.1} $\gamma<\alpha$.

Then $(\gamma^+)^{X}\leq \alpha$. But
$X \in Y$ and $\rho(X)=\gamma$, so $(\gamma^+)^{X}<(\gamma^+)^{Y}$,
so $\alpha=(\gamma^+)^{X}$. But then $(b)$ cannot hold,
because $\alpha$ is a cardinal of $\ult_0(Y, E^{Y}_\alpha)$ and $X$ defines a
map from $\alpha$ onto $\gamma$.

\paragraph{Subcase 2.2} $\gamma=\alpha$.

Suppose $\nu$ is a limit cardinal of $P$. It follows that $\gamma$ is a limit
cardinal in $X$ and $Y$. But then $E^{Y}_\gamma=\emptyset$, as desired.

Next suppose $\nu=(\kappa^+)^P$, where $\kappa$ is a cardinal of $P$.
It follows that $j$ is continuous at $\nu$, for otherwise we have a
$\mathbf{\Sigma}^Q_n$ function $f:\mu\to \nu$ with cofinal range. We then have a $\mathbf{\Sigma}_n^Q$
function $g$ such that for all $\xi < \mu$,
$g(\xi)$ be a wellorder of $\kappa$ of ordertype $f(\xi)$, and
$\{(\xi, \eta, \lambda)\,|\, \eta<_{g(\xi)}\lambda\}$ witnesses that $\rho_n(P)\leq \kappa$, contradiction.
So $j$ is continuous at $\nu$. On the other hand, our case 2 hypothesis is that
$j$ is discontinuous at $\nu$.
This finishes Case 2 and the Lemma under the assumption that $X$ is sound.
We have (i) of the Lemma in this case.

Suppose now that $X$ is not sound. It follows that $P$ has type 1B, $i$ is discontinuous at
$\nu$, and $R$ has type 2. Thus  
$\rho_k(R) < \rho_{k-1}(R)$, so $m=k-1$.  Let
\[
	Z = \mfc_m(X),
	\]
	so that
	\[
		R^- = \ult(Z,D),
		\]
		where $D$ is the order zero measure of $Z$
on $\hat{\rho}_k(R) = i(\nu)$. We have that $i_D \circ \sigma$ is a
nearly elementary (in fact elementary) map from $Z$ to $Y$,
and $\crit(\sigma \circ i_D) = \crit(\sigma)=\alpha$. By the condensation
theorem \ref{condensationtheorem}, either
\begin{itemize}
	\item[(a)] $Z \isneq Y$ or
	\item[(b)] $Z \isneq \ult_0(Y, E^{Y}_\alpha)$.
\end{itemize}
One can rule out (b) by the same argument that we applied to $X$.
This gives us (a), and hence conclusion (ii) of the Lemma.\footnote{
 A stronger possible conclusion
 for \ref{notes lemma 1} would demand
 condensation for the external part of $\Lambda$. Namely, either
             \begin{itemize}
                \item[(i)]
$(\ult(P,F), \Lambda^\sigma)\isneq(i^Q_F(P), \Lambda)$, or
\item[(ii)] $\ult(P,F)$ has type 2, and for $k=k(P)$ and
        $i_D \colon \mfc_k(\ult(P,F)) \to \ult(P,F)$ the anticore map,
        $(\mfc_k(\ult(P,F)), \Lambda^{\sigma\circ i_D}) \isneq (i^Q_F(P),\Lambda)$.
        \end{itemize}
Using \cite{trang}, we get the stronger conclusion when
$\crit(\sigma) < \rho_k(\ult(P,F)$.}

\hfill{$\qed$ Lemma \ref{notes lemma 1}}
\end{proof}

Lemma \ref{notes lemma 1} concerns the step from $P$ to $Q=R^-$, where $P$ and
$R$ are successive elements of a dropdown sequence. The next lemma concerns the step
from $R^-$ to $R$.

\begin{lemma}\label{notes lemma 2} 
($\adp$) Let $Q$ be a countable premouse with
$k(Q)>0$ and $F$ close to $Q$. Suppose there is
	$\Lambda$ such that $(\ult(Q,F), \Lambda)$ is a mouse pair.
	Let $P=Q^-$,
	and $\sigma:\ult(P,F)\to \ult(Q,F)^-$ be the natural embedding; then 
	either
        \begin{itemize}
                \item[(i)]
$\ult(P,F) \isneq i^Q_F(P)$, or
\item[(ii)] $\ult(P,F)$ has type 2, and for $k=k(P)$ and
	$\mfc_k(\ult(P,F)) \isneq i^Q_F(P)$.
        \end{itemize}
\end{lemma}
\begin{proof}
If $\rho^-(P)\leq \rho^-(Q)$, then this follows from Lemma \ref{notes lemma 1}.
	So let us assume $\rho^-(Q)<\rho^-(P)$. Let $\nu=\rho^-(Q)=\rho_{k+1}(Q)$, where $k+1=k(Q)$.
	(Recall $k(Q)>0$.) So letting $\hat{Q}$ be the bare premouse associated to $Q$,
$\rho_{k+1}(\hat{Q}) < \rho_k(\hat{Q})$ and
	$\sigma$ is the natural embedding from
	$\ult_k(\hat{Q},F)$ to $\ult_{k+1}(\hat{Q},F)$.

Let $R=\ult(P,F)=\ult_k(Q^-,F)$ and $r=p_k(Q)$. Let $S=\ult_{k+1}(Q,F)$ and $s=p_{k+1}(Q)$.

The elements of $R$ are $[a, f_{\tau, r}^{Q^-}]^{Q^-}_F$ for
$\tau$ a $r\Sigma_k$ Skolem term. The elements of $S$
are $[a,f^Q_{\tau,s}]^Q_F$ for $\tau$ a $r\Sigma_{k+1}$ term.
Each $f^{Q^-}_{\tau,r}$ is also of the form $f^Q_{\gamma, s}$ for some $\gamma$. So $\sigma$ is just given by
\[\sigma([a,g]^{Q^-}_F)=[a,g]^Q_F,\]
where the superscripts $Q^-$ and $Q$ indicate the two classes of functions used.

	Letting $i=i_F^{Q^-}$ and $j=i_F^Q$, we have the diagram
	\[
\begin{tikzcd}
Q
\arrow{r}{i} \arrow[swap]{rd}{j}&
R\arrow{d}{\sigma} \\
&
S
\end{tikzcd}
\]

Set \begin{align*}
    \gamma&=\sup i ``\nu\\
    &=\sup j``\nu.
\end{align*}
The equality holds because the 
functions bounded
in $\nu$ used in $\ult(Q^-,F)$ and $\ult(Q,F)$ are the same. We have
\[\gamma=\rho_{k+1}(S),\] and \[j(s)=p_{k+1}(S)\] by the general theory of $\ult_{k+1}$.

\begin{claim} $\gamma=\rho_{k+1}(R)$ and $p_{k+1}(R)=i(s)$.
\end{claim}
	\begin{proof} That $\gamma \le \rho_{k+1}(R)$ and
		$i(s) \le_{\text{lex}} p_{k+1}(R)$ follows from the usual proof that
		solidity witnesses are mapped by $i$ to generalized solidity witnesses.
		The reverse inequalities follow from the fact that
		$R = \mbox{Hull}_{k+1}^R(\gamma \cup i(s))$.
	\end{proof}

\begin{claim} 
$\sigma$ is nearly elementary.
\end{claim}
	\begin{proof} The proof of Claim \ref{nearelementarityclaim} shows that
		$\sigma$ is $\Sigma_0$ elementary as a map from
		$R^k$ to $S^k$. We must see that $\sigma$ is cardinal preserving.
		This follows from the elementarity if $k>0$. If $k=0$,
then for $f \in Q$, 
\begin{align*}
R \models [a,f] \text{ is a cardinal} &\text{ iff for $F_a$ a.e. $u$,
$Q \models f(u)$ is a cardinal}\\
&\text{ iff $S \models [a,f]$ is a cardinal}.
\end{align*}
The first line holds because the function $g(u) =$ first map from $|f(u)|$ onto $f(u)$
belongs to Q. This shows that $\sigma$ is cardinal preserving when $k=0$.
	\end{proof}

\begin{claim}\label{ransigmabounded}
If $\sigma``\rho_k(R)$ is unbounded in $\rho_k(S)$, then $R=S^-$.
\end{claim}
\begin{proof} Let $g$ and $h$ be the $\Sigma_1$ Skolem functions
	of $R^k$ and $S^k$. If $\alpha < \rho_k(S)$, then we have
	$\alpha = h(\beta,s)$ for some $\beta < \rho_{k+1}(S)$.
	By hypothesis, $\alpha = h^{S^k|\sigma(\delta)}(\beta,s)$
	for some $\delta < \rho_k(R)$. (Here the superscript $S^k||\sigma(\delta)$
	indicates a bound on the witness that $h(\beta,s)=\alpha$.)
	But then $\sigma(g(\delta,r)) = \alpha$, so $\alpha \in \ran(\sigma)$.
	So $\rho_k(S) \subseteq \ran(\sigma)$, so $R^k = S^k$,
	so $R = S^-$.
\end{proof}

By \ref{ransigmabounded} we may assume that $\ran(\sigma)$ is bounded in
$\rho_k(S)$. For this it follows easily that $R \in S$.
If $R$ is type 1 
(that is, $k$-sound), then since $R= \text{Hull}_{k+1}^R(\crit(\sigma) \cup
p_{k+1}(R))$, Theorem \ref{condensationtheorem} applies,
and we get that $R \lhd S$. (The case $\crit(\sigma)$ is an index of an extender in $S$
is ruled out as in Lemma \ref{notes lemma 1}.) If $R$ has type 2,
then $R=\ult(\mathfrak{C}_k(R),D)$, where $D$ is the order zero measure
of $\mathfrak{C}_k(R)$ on $i(\rho_k(Q))$. We then get
$\mathfrak{C}_k(R) \lhd S$ by applying Theorem \ref{condensationtheorem}
to $\sigma \circ i_D$. (Once again, we rule out
$\crit(\sigma)$ being an index of an extender in $S$
as in the proof of \ref{notes lemma 1}.) This proves Lemma
\ref{notes lemma 2}.
\end{proof}

We isolate now, in the definition of ``resolution of $\sigma$",
a fairly extensive list
of the properties of sequences of factor embeddings such as those in
Lemmas \ref{notes lemma 1} and \ref{notes lemma 2}. Lemma \ref{main full normalization lemma}
then shows that for certain $M, F,$ and $\xi$, the natural 
$\sigma \colon \ult(M|\xi,F) \to i_F^M(M|\xi)$ admits a resolution.
This example is the entire motivation for the definition.
We say that $\sigma$ ``respects drops" iff it admits a resolution.

\begin{definition}\label{corepair} Let $B$ be a type pfs premouse and $k=k(B)$; then
$C(B) = B$ if $B$ has type 1, and $C(B) = \mfc_k(B)$ if $B$ has type 2.
If $B$ has type 2, $D(B)$ is the order zero measure $D$
on $\hat{\rho}_k(B)$ such that $B= \ult(C(B),D)$. If $B$ has type 1,
$D(B)$ is a principal measure. In both cases, $i_D \colon C(B) \to B$ is
the ultrapower (i.e. anticore) map.
\end{definition}

	\begin{definition}\label{respectsdropsdef}
Let $N$ be a premouse, and $\sigma \colon N|\eta \to N|\lambda$ be nearly elementary,
where $\eta < \lambda \le o(N)$. 
		We say that $\sigma$ \textit{respects drops over $(N,\eta,\lambda)$} iff
		letting $n=n(N,\eta)$, we have $\eta_i$ for $1 \le i \le n+1$ such that
		\[
			\eta = \eta_1 \le \eta_2 < ... \le \eta_{n+1} = \lambda
			\]
		and $n = n(N,\eta_i)$ for all $i$, together pfs premice
$B^i_k$ for $1 \le k \le n$ such that
\[
C(B^i_k) = A_k(N,\eta_i)^-
\]
for all $i,k$, together with nearly elementary maps
\[
\sigma_i \colon B^i_i \to B^{i+1}_i
\]
defined for $i \le n$ such that
\[
\sigma = \sigma_n \circ ... \circ \sigma_1,
\]
and the following hold. 
\begin{itemize}
\item[(a)] $k(B^1_k) = k(B^i_k)$ for all $i$; let $m_k$ be the common value, and set
\[
\gamma^i_k = \rho(B^i_k) = \rho_{m_k+1}(B^i_k).
\]
\item[(b)] $C(B^i_k))$ is sound, that is, $m_k+1$-sound.
\item[(c)] For all $i<k$, $B^i_k = C(B^k_k)$.
\item[(d)] $\sigma_i \restriction \gamma^i_i = \text{ id}$, 
$\sigma_i(p_{m_i+1}(B^i_i)) = p_{m_i+1}(B^{i+1}_i)$, and
$\sigma_i(\eta_i) = \eta_{i+1}$.
\item[(e)] Either
\begin{itemize}
\item[(i)] $\eta_i = \eta_{i+1}$, $B^i_k = B_k^{i+1}$ for all $k$, and
$\sigma_i = \text{ id}$, or
\item[(ii)] $\eta_i < \eta_{i+1}$, $C(B^i_i) \lhd B^{i+1}_i$.
\end{itemize}
\item[(f)] For $i < n$, $\gamma^{i+1}_{i+1} < \gamma^i_i$,
and either
\begin{itemize}
\item[(i)] $\gamma^i_i = \gamma^{i+1}_i$ and $\sigma_i \restriction
\gamma^i_i +1 = \text{ id}$, or
\item[(ii)] $\gamma^i_i$ is a limit cardinal of $B^{i+1}_i$.
\end{itemize}
\end{itemize}
We call $(\vec{\eta}, \vec{B}, \vec{\sigma})$ a {\em resolution} of $\sigma$.
\end{definition}
We are allowing $\eta_i = \eta_{i+1}$ here simply for its bookkeeping value.
Since $C(B^i_i)$ is $m_i+1$-sound, $\sigma_i$ is determined by $B^i_i$ and
$B^{i+1}_i$. The $B^i_k$ are in turn determined by $\eta_i$. Thus
the whole of the resolution is determined by the sequence of $\eta$'s.

\begin{lemma}\label{main full normalization lemma}
 Let $M$ be a premouse, $\bar{\eta} \le o(M)$, and $E= E^M_{\bar{\eta}}\neq \emptyset$.
Let $F$ be  close to $M$, $\crit(F)=\mu<\hat{\lambda}(E)$, and $N=\ult(M,F)$. Let
\[
\sigma:\ult(M|\bar{\eta}, F)\to i^M_F(M|\bar{\eta})
\]
be the natural factor map,  
$\lambda=i^M_F(\bar{\eta})$, and $\eta = o(\ult(M|\bar{\eta}, F))$; then
\begin{itemize}
\item[(a)] $\ult(M|\bar{\eta}, F) = N|\eta$, and
\item[(b)] $\sigma$ respects drops over $N,\eta,\lambda$.
\end{itemize}
\end{lemma}
\begin{proof} Let $n=n(M,\bar{\eta})$ and $A_i = A_i(M,\bar{\eta})$ for $i \le n$. Since $i_F^M$ is
elementary, it is almost exact, so $n(N,\lambda) \in \lbrace n,n+1 \rbrace$. We assume that
$n(N,\lambda) = n$. The other case is quite similar.

We must define a resolution $(\vec{\eta},\vec{B},\vec{\sigma})$ of $\sigma$. We start by
setting $\eta_{n+1} = \lambda$, and then define the $\eta_i$ and associated objects
for $i \le n$ by reverse induction on $i$. In the end we shall have $\eta_1 = \eta$.
As we go we shall verify the properties of a resolution. 

Let for $i \le n$
\begin{align*}
m_i+1 &= k(A_i),\\
\gamma_i &= \rho^-(A_i) = \rho_{m_i+1}(A_i),\\
B^{n+1}_i &= A_i(N,\lambda)^-,\\
\intertext{ and}
\gamma^{n+1}_i &= \rho^-(A_i(N,\lambda)) = \rho(B^{n+1}_i).
\end{align*}
Let $i_{n+1} = i^M_F$. By our preservation lemmas \ref{preservedropdown} and
\ref{preservedropdownelem}, $i_{n+1}(A_i) = A_i(N,\lambda)$ for all $i \le n$, so
$m_i = k(B^{n+1}_i)$ for all $i \le n$. Also,
\[
\gamma^{n+1}_i = \begin{cases} i_{n+1}(\gamma_i) & \text{ if $i < n$}\\
                           \sup i_{n+1}``\gamma_i & \text{ if $i=n$.}
\end{cases}
\]

Now let us define $\eta_n$ and $\sigma_n$. Set
\begin{align*}
B^n_n &= \ult(A_n^-,F) = \ult_{m_n}(A_n^-,F),\\
\intertext{ and }
\eta_n &= i_n(\bar{\eta}),
\end{align*}
where $i_n \colon A_n^- \to B^n_n$ is the canonical embedding.
We have that $k(A_n^-) = k(B^n_n) = m_n$ and $i_n$ is elementary,
that is, $r\Sigma_{m_n+1}$ elementary. It is possible that $B^n_n$ has type 2.\footnote{
This happens iff $A_n^-$ has type 1B and $\eta_{m_n}(A_n^-) = \mu$.} Let
\[
\sigma_n \colon B^n_n \to i_{n+1}(A_n^-) = B^{n+1}_n
\]
be the natural factor map. $\sigma_n$ is nearly elementary as a map
on premice of degree $m_n$, and since
$\sigma_n \circ i_n = i_{n+1} \restrict A_n^-$,
$\eta_{n+1} = \sigma_n(\eta_n)$. 

\bigskip
\noindent
{\em Claim 1.} 
\begin{itemize}
\item[(a)] $\rho_{m_n+1}(B^n_n) = \sup i_n ``\gamma_n = \sup i_{n+1}``\gamma_n \le
\rho_{m_n+1}(B^{n+1}_n)$.
\item[(b)] $\sigma_n \restriction \gamma^n_n = \text{ id}$.
\item[(c)] $\sigma_n(p_{m_n+1}(B^n_n))=p_{m_n+1}(B^{n+1}_n)$.
\end{itemize}

\medskip
\noindent
{\em Proof.} 
For (a): $\rho_{m_n+1}(B^n_n) = \sup i_n ``\rho_{m_n+1}(A_n^-) =
\sup i_n``\gamma_n$ by the general
properties of $\ult_{m_n}$; see the proof of Lemma \ref{notes lemma 2}.
Every $r\Sigma_{k(M)}^M$ function from $\mu$ into $\gamma_n$ with
range bounded in $\gamma_n$ belongs to $M$ because $n=n(M,\bar{\eta})$.
Thus $\sup i_n``\gamma_n = \sup i_{n+1}``\gamma_n$. Finally,
$\rho_{m_n+1}(B^{n+1}_n) = \sup i_{n+1}`` \gamma_n$ if
$A_n = M$, and $\rho_{m_n+1}(B^{n+1}_n) = \sup i_{n+1}(\gamma_n)$ if
$A_n \lhd M$, so $\sup i_{n+1}`` \gamma_n \le \rho_{m_n+1}(B^{n+1}_n)$ in
either case.

(b) follows from the fact that both ultrapowers use the same
functions with range bounded in $\gamma_n$, namely just those
functions belonging to $M$.

For (c) we use the solidity and universality (over $C(B^n_n))$)
of $p_{m_n+1}(B^n_n)$. See the proof of \ref{notes lemma 2}.
\hfill              $\square$

\bigskip
\noindent
{\em Claim 2.} If $A_n$ has type 1B, then $i_n$ is continuous at $\gamma_n$.

\medskip
\noindent
{\em Proof.} If $A_n$ has type 1B, then $\rho^-(A_n) = \rho_{m_n+1}(A_n)
= \gamma_n$
is $\Sigma_0$ regular in $A_n$. If $i_n$ is discontinuous at $\gamma_n$,
then $\gamma_n$ is $r\Sigma_{m_n}$ singular over $A_n$. But
$\gamma_n < \rho_{m_n}(A_n)$ because $A_n$ was the first level of $M$
with $\rho^- \le \gamma_n$. Thus $\gamma_n$ is $\Sigma_0$ singular in
$A_n$, contradiction.
\hfill    $\square$

\bigskip
\noindent
{\em Claim 3.} $C(B^n_n)$ is $m_n+1$-sound; moreover,
$p_{m_n+1}(C(B^n_n)) = i_n(p_{m_n+1}(A_n))$.

\medskip
\noindent
{\em Proof.} If $A_n$ has type 1A, this follows from the usual
properties of $\ult_{m_n}$; see the proof of
\ref{notes lemma 2}. The proof also works if $A_n$ has type 1B,
except for the problem that $C(B^n_n)^+$ may have type 2.
But if $A_n$ has type 1B, then $\rho_{m_n+1}(B^n_n) =
i_n(\rho_{m_n+1}(A_n))$ by Claim 2. This implies that $C(B^n_n)^+$ has type
1B, rather than type 2.
\hfill  $\square$

The next claim is the main step in our proof, the place where we
integrate the condensation arguments in Lemmas \ref{notes lemma 1}
and \ref{notes lemma 2}.

\bigskip
\noindent
{\em Claim 4.} Either $C(B^n_n) = B^n_n = B^{n+1}_n$, $\eta_n = \eta_{n+1}$,
 and $\sigma_n = \text{ id}$, or $C(B^n_n) \lhd B^{n+1}_n$.

\medskip
\noindent
{\em Proof.} Let $D=D(B^n_n)$, $C=C(B^n_n)$, and $m=m_n$. Let
\[
\pi = \sigma_n \circ i_D,
\]
so that $\pi \colon C \to B^{n+1}_n$ is nearly elementary
at degree $m_n$, and $\crit(\pi) \ge \rho_{m+1}(C)$. (Note
$\rho_{m+1}(C) = \rho_{m+1}(B^n_n) = \gamma^n_n$.) If $\pi = \text{ id}$,
then  $C = B^n_n = B^{n+1}_n$, $\eta_n = \eta_{n+1}$,
 and $\sigma_n = \text{ id}$, and we are done. So suppose that
$\pi \neq \text{ id}$; we shall show that $C \lhd B^{n+1}_n$.

We show first that $C \in B^{n+1}_n$. 

For that,
suppose first that $\rho_{m +1}(C) < \rho_{m+1}(B^{n+1}_n)$;
then since $C$ is coded by $\mbox{Th}_{m+1}^C(\rho_{m+1}(C) \cup
\lbrace r \rbrace)$ for some parameter $r$, and
\[
\mbox{Th}_{m+1}^C(\rho_{m+1}(C) \cup
\lbrace r \rbrace) = \mbox{Th}_{m+1}^{B^{n+1}_n}(\rho_{m+1}(C) \cup
\lbrace \sigma_n(r) \rbrace),
\]
we get that
\[
\mbox{Th}_{m+1}^C(\rho_{m+1}(C) \cup
\lbrace r \rbrace) \in B^{n+1}_n,
\]
so $C \in B^{n+1}_n$.  

Suppose next that $\rho_{m+1}(C) = \rho_{m+1}(B^{n+1}_n)$.
Let $X=C^m$ and $Y=(B^{n+1}_n)^m$ be the two reducts. $\pi \restriction X =
\sigma_n \restriction X$ is $\Sigma_0$ elementary from $X$ to $Y$,
$\rho_1(X) \le \rho_1(Y)$, $\pi \restriction \rho_1(X) = \text{ id}$,
and $\pi(p_1(X))=p_1(Y)$. Since $\rho_{m+1}(A_n) < \rho_m(A_n)$,
$\rho_1(X) < o(X)$ and $\rho_1(Y) < o(Y)$. It will be enough to
show that $\ran(\pi)$ is bounded in $o(Y)$, for then
$\mbox{Th}_1^X(\rho_1(X) \cup \lbrace \rho_1(X), p_1(X) \rbrace) \in Y$,
so $X \in Y$, so $C \in B^{n+1}_n$.

So suppose $\ran(\pi)$ is unbounded in $o(Y)$. It follows that $\pi$
is $\Sigma_1$ elementary from $X$ to $Y$. If $A_n$ has type 1A,
then $Y = \mbox{Hull}_1^Y(\rho_1(Y) \cup p_1(Y))$, so $Y$ is $\Sigma_1$-generated
from points in $\ran(\pi)$, so $X=Y$ and $\pi \restriction X = \text{ id}$.
Thus $\pi = \text{ id}$, contradiction. If $A_n$ has type 1B, then
$Y = \mbox{Hull}_1^Y(\rho_1(Y) \cup \lbrace \rho_1(Y),p_1(Y)\rbrace)$,
so we just need to see that $\rho_1(Y) \in \ran(\pi)$ for our contradiction.
But $\rho_1(Y) = \rho_{m+1}(B^{n+1}_n) = i_{n+1}(\rho_{m+1}(A_n))$,
since if $\rho_{m+1}(B^{n+1}_n) < i_{n+1}(\rho_{m+1}(A_n))$
then we must have $A_n = M$ and $\rho_{m+1}(B^{n+1}_n) = \sup
i_{n+1}``\rho_{m+1}(A_n)$. This implies
$\rho_{m+1}(B^{n+1}_n) < \hat{\rho}_{m+1}(B^{n+1}_n) =
i_{n+1}(\rho_{m+1}(A_n))$, so $B^{n+1}_n$ has type 2.
But $B^{n+1}_n \lhd N$, so it has type 1.

Thus $C \in B^{n+1}_n$. 
By the Condensation Theorem \ref{condensationtheorem},
either $C \lhd B^{n+1}_n$ or $C \lhd \ult(B^{n+1}_n,G)$,
where $G$ is on the $B^{n+1}_n$ sequence and $\lh(G) =
\crit(\pi) = \crit(\sigma_n)$. One can rule out the latter
possibility just as in the proofs of
Lemmas \ref{notes lemma 1} and \ref{notes lemma 2}.
Thus $C \lhd B^{n+1}_n$, as desired.
\hfill  $\square$

Before we proceed to the general inductive step, we need some claims
that deal with the possible difference between $B^n_n$ and $C(B^n_n)$.

\bigskip
\noindent
{\em Claim 5.} If $C(B^n_n) \neq B^n_n$, then there is a limit
cardinal $\xi$ of  $B^n_n$ such that $\xi > \rho_{m_n}(B^n_n)$
and $B^n_n|\xi = C(B^n_n)|\xi$.

\medskip
\noindent
{\em Proof.} If $C(B^n_n) \neq B^n_n$, then $B^n_n$ has type 2,
and we can take $\xi = \hat{\rho}_{m_n}(B^n_n)$.
\hfill    $\square$

\bigskip
\noindent
{\em Claim 6.} $\gamma^n_{n-1} \le \rho_{m_n}(B^n_n)$.

\medskip
\noindent
{\em Proof.} Let $m=m_n$.  Note 
\[
\gamma_{n-1} \le \rho_{m}(A_n),
\]
since $m +1 = k(A_n)$ is by definition the least $i$
such that $\rho_i(A_n) < \gamma_{n-1}$. If $\rho_m(A_n) = \gamma_{n-1}$,
then $A_n^- = A_{n-1}$, and
\[
\rho_m(B^n_n) = \sup i_n``\gamma_{n-1} < i_n(\gamma_i) = \gamma_i^n,
\]
for all $i < n-1$, so $\rho_m(B^n_n) = \gamma^n_{n-1}$ and we are done.

Suppose next that 
\[
\gamma_{n-1} < \rho_m(A_n).
\]
 For all $i < m$,
$\gamma_{n-1} \le \rho_i(A_n)$, so for all $i < m$,
$\gamma_{n-1} < \rho_i(A_n)$, It follows that $A_{n-1} \in A_n$,
so
\[
\gamma_{n-1}^n = i_n(\gamma_{n-1}) < \sup i_n``\rho_m(A_n) = \rho_m(B^n_n),
\]
as desired.
\hfill      $\square$

Note also

\bigskip
\noindent
{\em Claim 7.} 
\begin{itemize}
\item[(a)] $n(N,\eta_n)=n$.
\item[(b)] $C(B^n_n)^+ = A_n(N,\eta_n)$.
\item[(c)] For $k < n$, $B^n_k = \sigma_n^{-1}(B^{n+1}_k)$
and $\gamma^n_k = \sigma_n^{-1}(\gamma^{n+1}_k)$.
\item[(d)] $\gamma^n_n \le \gamma^{n+1}_n$.
\end{itemize}

\medskip
\noindent
{\em Proof.} These are all immediate consequences of the claims above.
\hfill    $\square$

We proceed to the inductive step. Suppose $1 \le e < n$, and suppose
that for all $k \ge e+1$, $B^k_k$, $B^{k+1}_k$ and $\sigma_k \colon B^k_k \to B^{k+1}_k$
satisfy Claims 1--7, with $k$ and $k+1$ replacing $n$ and $n+1$. We define
\[
B^e_e = \ult(A_e^-,F) = \ult_{m_e}(A_e^-,F),
\]
and let
\[
i_e \colon A_e^- \to B^e_e
\]
be the canonical embedding. Let
\[
\sigma_e \colon B^e_e \to i_{e+1}(A_e^-) = B^{e+1}_e
\]
be the factor map. $\sigma_e$ is nearly elementary as a map
on premice of degree $m_e$. Let
\begin{align*}
\eta_e &= i_e(\bar{\eta}) = \sigma_e^{-1}(\eta_{e+1}).
\end{align*}

Claims 1-3 hold with $e$ and $e+1$ replacing $n$ and $n+1$,
with the same proofs. We need a stronger version of Claim 4 now.

\bigskip
\noindent
{\em Claim 8.} Either $C(B^e_e) = B^e_e = B^{e+1}_e$ and $\sigma_e = \text{ id}$,
or $C(B^e_e) \lhd B^{e+1}_e \lhd C(B^{e+1}_{e+1})$.

\medskip
\noindent
{\em Proof.} Suppose the first alternative does not hold. The proof of Claim 4
yields that $C(B^e_e) \lhd B^{e+1}_{e}$, and $B^{e+1}_e \lhd B^{e+1}_{e+1}$,
so we may assume that $C(B^{e+1}_{e+1}) \neq B^{e+1}_{e+1}$. By Claims 5 and 6
at $e+1$, we may fix a limit cardinal $\xi$ of $B^{e+1}_{e+1}$ such that
\[
B^{e+1}_{e+1}|\xi = C(B^{e+1}_{e+1})|\xi
\]
and $\gamma^{e+1}_e < \xi$. But $C(B^e_e)$ has cardinality $\gamma^e_e$ in
$B^{e+1}_{e+1}$, and $\gamma^e_e \le \gamma_e^{e+1}$. Thus
\[
C(B^e_e) \lhd B^{e+1}_{e+1}|(\gamma^{e+1}_e)^{+,B^{e+1}_{e+1}} \lhd C(B^{e+1}_{e+1}),
\]
as desired.
\hfill     $\square$

It is easy to check that Claims 5 and 6 hold with $e$ and $e+1$ replacing
$n$ and $n+1$. 

Let us address item (f) in the definition of resolutions. ( Item (f) did not apply
when $e=n$.)

\bigskip
\noindent
{\em Claim 9.} 
$\gamma^{e+1}_{e+1} < \gamma^e_e$,
and either
\begin{itemize}
\item[(i)] $\gamma^e_e = \gamma^{e+1}_e$ and $\sigma_e \restriction
\gamma^e_e +1 = \text{ id}$, or
\item[(ii)] $\gamma^e_e$ is a limit cardinal of $B^{e+1}_e$.
\end{itemize}

\medskip
\noindent
{\em Proof.} This is trivial if $C(B^e_e)=B^e_e = B^{e+1}_e$ and
$\sigma_e = \text{ id}$, so assume otherwise. By Claim 8,
$C(B^e_e) \lhd B^{e+1}_e$, so $B^e_e \in B^{e+1}_e$ and has cardinality
$\gamma^e_e$ in $B^{e+1}_e$.

For the first part, note
\begin{align*}
\gamma^{e+1}_{e+1} & = \sup i_{e+1}``\gamma_{e+1}\\
& = \sup i_{e} ``\gamma_{e+1}
< \sup i_e ``\gamma_e = \gamma^e_e.
\end{align*}
The second equality holds because the two ultrapowers use the same functions
with range bounded in $\gamma_{e+1}$, namely those belonging to $M$.

Suppose that $\gamma_e$ is a limit cardinal of $A_e$.
It follows that $\gamma^e_e = \sup i_e``\gamma_e$ is a  limit
cardinal of $B^e_e$, and since $\gamma^e_e \le \crit(\sigma_e)$,
$\gamma^e_e$ is a limit cardinal of $B^{e+1}_e$. Thus we have (ii).

Suppose that $\gamma_e = \kappa^{+,A_e}$. By the definition
of $A_{e+1}$, letting $m=m_{e+1}$ we have that
$\gamma_e = \kappa^{+,A_{e+1}}$ and $\rho_m(A_{e+1}) \ge \gamma_e$.
If $\gamma_e$ has $r\Sigma_{m}$ cofinality $\mu$ in $A_{e+1}$,
then  one easily gets that $\rho_{m}(A_{e+1}) \le \kappa$,
contradiction. Thus 
\begin{align*}
\gamma^{e+1}_e &= \sup i_{e+1}``\gamma_e = i_{e+1}(\gamma_e)\\
               &= \sup i_e``\gamma_e = i_e(\gamma_e).
\end{align*}
This shows that (ii) holds.
\hfill     $\square$

Finally, the analog of Claim 7 takes more work when $e<n$.

\bigskip
\noindent
{\em Claim 10.} 
\begin{itemize}
\item[(a)] For $k < e$, $A_k(N,\eta_e)^- = \sigma_e^{-1}(A_k(N,\eta_{e+1})^-) =
\sigma_e^{-1}(B^{e+1}_e)$.
\item[(b)] $A_e(N,\eta_e)^- = C(B^e_e)^+$.
\item[(c)] For $k>e$, $A_k(N,\lambda_e)^- = C(B^k_k)^+$.
\item[(d)] $n(N,\eta_e) = n$.
\end{itemize}

\medskip
\noindent
{\em Proof.} Part (a) follows from Lemma \ref{preservedropdown},
and this implies that $\sigma_e(\gamma^e_k) = \gamma^{e+1}_k$
fot $k<e$. 

Part (b) holds because $\gamma^e_e = \rho(C(B^e_e)) =
\rho(B^e_e) = \sup i_e``\gamma_e < i_e(\gamma_{e-1}) = \gamma^e_{e-1}$, 
and $i_e$ preserves the fact that all $P$ such that
$A_e \lhd P \unlhd A_{e-1}$ satisfy $\gamma_{e-1} \le \rho^-(P)$.
(See Lemma \ref{preserverhominus}(c).)

For (c), we show first that $A_{e+1}(N,\eta_e)^- = C(B^{e+1}_{e+1})$. 
We have that $C(B^e_e) \lhd B^{e+1}_e \lhd C(B^{e+1}_{e+1})$ by Claim 8,
and 
\[
\rho(C(B^{e+1}_{e+1})) = \gamma_{e+1}^{e+1} < \gamma^e_e = \rho(C(B^e_e))
\]
by Claim 8. So it will suffice to show that whenever
$C(B^e_e) \lhd Q \lhd C(B^{e+1}_{e+1})$, then $\gamma^e_e \le \rho(Q)$.
If $B^{e+1}_e \unlhd Q \lhd C(B^{e+1}_{e+1})$ then
\[
\gamma^e_e \le \gamma^{e+1}_e < \rho(Q),
\]
as desired. So suppose $C(B^e_e) \lhd Q \lhd B^{e+1}_e$ and $\rho(Q)<\gamma^e_e$.
If $Q \in B^{e+1}_e$, then $\gamma^e_e$ is not a cardinal of $B^{e+1}_e$,
contrary to both alternatives (i) and (ii) in Claim 9. Thus
$Q = B^{e+1}_e \downarrow k$ for some $k < m_e = k(B^{e+1}_e)$.
This means 
\[
\rho_{k+1}(B^{e+1}_e) < \gamma^e_e \le \gamma^{e+1}_e =
\rho_{m_e+1}(B^{e+1}_e,
\]
a contradiction.

Thus $A_{e+1}(N,\eta_e) = C(B^{e+1}_{e+1})^+$. 
It is then easy to see that 
\[
A_{e+2}(N,\eta_e) =
A_{e+1}(N,\eta_{e+1}) = C(B^{e+2}_{e+2})^+,
\]
and so on until we reach $A_n(N,\eta_e) = C(B^n_n)^+$.
Here the value of $\rho^-$ is $\gamma^n_n$, and no higher levels
of $N$ project strictly below that. Thus $n(N,\eta_e) = n$.
This finishes the proofs of (c) and (d).
\hfill      $\square$

In view of Claim 10, we may set $B^e_k = A_k(N,\eta_e)^-$ 
and $\gamma^k_e = \rho(B^e_k)$ for all
$k \le n$, and we have the properties of a resolution that
apply to the $\eta_i$, $\sigma_i$, and $B^i_k$ for $i \ge e$.

Eventually we reach $e=1$. Since $A_0 = M|\bar{\eta}$ is active
and $\gamma_0 = o(N|\bar{\eta})$,
\begin{align*}
\gamma_1 &= \rho_1(N|\bar{\eta}),\\
m_1 &= 0,\\
\intertext{ and }
B^1_1 &= \ult_0(N|\bar{\eta},F).
\end{align*} 
$B^1_0 = C(B^1_1)$ because $m_1 = 0$. Thus $B^1_1 \unlhd N$.
It its clear that $\sigma = \sigma_n \circ ... \circ \sigma_1$.

\end{proof}

\subsection{Weak tree embeddings}

For most purposes, we need only consider weak tree embeddings that act on $\lambda$-tight,
normal plus trees. This is because we have already shown (assuming $\adp$)\footnote{See Theorem \ref{tight trees}.}
that if $(P,\Sigma)$ is a mouse pair with scope HC, and $(N,\Lambda)$ is an iterate
of $(P,\Sigma)$ via the plus tree $\tree{T}$, then $(N,\Lambda)$ is an iterate
of $(P,\Sigma)$ via a $\lambda$-tight, normal plus tree $\tree{U}$.\footnote{$\tree{U}$ is
the normal companion of $\tree{T}$, re-arranged so that each plus extender $E^+$ used in the
normal companion corresponds to two extenders $E, D$ used in $\tree{U}$.}
It is convenient to restrict attention to $\lambda$-tight, normal trees; in particular,
the fact that they are length-increasing means the agreement properties between their
models are simpler.  So
we shall do this.\footnote{As we shall see below, in one case the
natural embedding from the full normalization $X(s)$ of a stack of plus trees
to its embedding normalization $W(s)$ is not quite a weak tree embedding in the
sense of \ref{weaktreeembeddingdef}. We shall describe the slightly more
general notion of weak tree embedding $\Phi \colon \tree{S} \to \tree{T}$
required in this case in a future draft of this paper. It amounts to allowing
the $M_\alpha^{\tree{T}}$ to be of type 2.}

\begin{definition}\label{weaktreeembeddingdef} Let $\tree{S}$ and $\tree{T}$ be plus trees with the same base model.
A \textit{weak tree embedding} $\Phi \colon
\tree{S}\to \tree{T}$ is a system $\langle v,u, \{s_\xi\}_{\xi<\lh\tree{S}},
 \{t_\zeta\}_{\zeta+1<\lh\tree{S}},\{\sigma_\zeta\}_{\zeta+1<\lh \tree{S}}\rangle$ such that
\begin{enumerate}
    \item  $v:\lh \tree{S}\to \lh\tree{T}$ is tree-order preserving, $u:\{\eta\,|\,\eta+1<\lh \tree{S}\}\to \lh \tree{T}$, $v(\xi)=\sup\{u(\eta)+1\,|\, \eta<\xi\}$, and $v(\xi)\leq_\tree{T}u(\xi)$;
    \item For $\eta\leq_\tree{S}\xi$, 
    \begin{enumerate}
        \item $s_\xi: M^\tree{S}_\xi\to M^\tree{T}_{v(\xi)}$ is elementary and $s_0 = id_{M^\tree{S}_0}$;
        \item $\hat\imath^\tree{T}_{v(\eta),v(\xi)}\circ s_\eta=s_\xi\circ \hat\imath^\tree{S}_{\eta,\xi}$,
        \item $t_\xi= \hat\imath^\tree{T}_{v(\xi),u(\xi)}\circ s_\xi$ (so $t_\xi$ is a partial elementary map $M^\tree{S}_\xi\to M^\tree{T}_{u(\xi)}$);
    \end{enumerate}
    \item for $\xi+1<\lh \tree{S}$, $\eta=\tree{S}\pred(\xi+1)$, and $\eta^*=\tree{T}\pred(u(\xi)+1)$,
    \begin{enumerate}
        \item Either \begin{itemize}
            \item [(i)] (X-case) $\sigma_\xi:M^\tree{T}_{u(\xi)}|\lh E^\tree{T}_{u(\xi)}\to
 M^\tree{T}_{u(\xi)}|\lh t_\xi(E^\tree{S}_{\xi})$ respects drops,
$\sigma_\xi(E^\tree{T}_{u(\xi)})=t_\xi(E^\tree{S}_\xi)$, and $\ran(t_\xi)\subseteq \ran(\sigma_\xi)$, or
            \item[(ii)] (W-case) $\sigma_\xi:M^\tree{T}_{u(\xi)}|\lh t_\xi(E^\tree{S}_{\xi})
\to M^\tree{T}_{u(\xi)}|\lh E^\tree{T}_{u(\xi)}$ respects drops
 and $E^\tree{T}_{u(\xi)}=\sigma_\xi\circ t_\xi(E^\tree{S}_\xi)$;
        \end{itemize} 
        \item $\eta^*\in[v(\eta),u(\eta)]_\tree{T}$,
        \item 
    $s_{\xi+1}\restrict \lh E^\tree{S}_\xi = \begin{cases}
 \sigma_\xi\circ t_\xi\restrict \lh E^\tree{S}_\xi+1 & \text{in the W-case}\\
        \sigma_\xi^{-1}\circ t_\xi\restrict \lh E^\tree{S}_\xi+1 & \text{in the X-case.}
        \end{cases}$
    \end{enumerate}
\end{enumerate}
\end{definition}

Note that a tree embedding is just a weak tree embedding in which $\sigma_\xi=id$ 
for all $\xi+1<\lh \tree{S}$. In particular, a tree embedding is in both 
the X-case and the W-case at $\xi$ for every $\xi+1<\lh\tree{S}$.

We now extend our Shift Lemma and Copying Construction for 
meta-trees to the case where we have a weak tree embedding rather than a tree embedding.
We shall just state the relevant results. The calculations involved in their proofs
are quite similar to those we have done in the tree embedding case, so we defer them
to a later draft of this paper.

We first extend our notation to weak tree embeddings.
\begin{definition}\label{metashiftapplies2}
Let $\Psi:\tree{S}\to \tree{T}$ and $\Pi:\itU \to \itV$ be extended weak tree embeddings, $F$ an extender such that $F^-$ be
 an extender on the $M_\infty^\tree{S}$-sequence, and $G$ an extender such that $G^-$ is on the $M_\infty^\tree{T}$-sequence. We say that 
 \textit{the Shift Lemma applies
  to $(\Psi,\Pi, F, G)$} iff letting $\beta = \beta(\itS,F)$ and $\beta^*=\beta(\itT, G)$, 
  \begin{enumerate}
  \item $M_\infty^\tree{S}|\lh(F)\is\dom( t_\infty^\Psi)$ and $G=t_\infty^\Psi(F)$,
  \item  $\Psi\restrict\beta+1\approx \Pi\restrict\beta+1$,
\item $\tree{T}\restrict \beta^*+1=\tree{V}\restrict\beta^*+1$
\item $\beta^*\in [v^\Pi(\beta), u^\Pi(\beta)]_\tree{V}$ and if $\beta^*<u^\Pi(\beta)$, then either
\begin{enumerate}
    \item if $\Pi$ is in the $W$-case at $\beta$, \[\sigma^\Pi_\beta\circ t_\beta^\Pi\restrict \dom(F)\cup\{\dom(F)\}=s_{\beta, \beta^*}^\Pi\restrict\dom(F)\cup\{\dom(F)\},  \text{ or}\]
    \item if $\Pi$ is in the $X$-case at $\beta$, \[(\sigma^\Pi_\beta)^{-1}\circ t_\beta^\Pi\restrict \dom(F)\cup\{\dom(F)\}=s_{\beta, \beta^*}^\Pi\restrict\dom(F)\cup\{\dom(F)\},\]
\end{enumerate} 
\item
   if  $\beta+1<\lh(\tree{U})$, then 
  $\dom(F) \isneq M_\beta^\tree{U}|\lh(E^\itU_\beta)$, and
  \item
   if  $\beta^*+1<\lh(\tree{V})$,
  $\dom(G) \isneq M_{\beta^*}^\tree{U}|\lh(E^\itV_{\beta^*})$.
   \end{enumerate}
   \end{definition}

\begin{lemma}[Shift Lemma for weak tree embeddings]
 $\Psi:\tree{S}\to \tree{T}$ and $\Pi:\itU \to \itV$ be extended weak tree embeddings,
and let $F$ be an extender such that $F^-$ be
 an extender  on the sequence of the last model of $\tree{S}$ and $G$ be an extender such that $G^-$ is on the extender sequence of the last model of $\tree{T}$. Let $\alpha_0=\alpha_0(\tree{S}, F)$ and
 $\alpha^*_0=\alpha_0(\tree{T},G)$.
 
 Suppose that the Shift Lemma applies
  to $(\Psi,\Pi, F, G)$. Then $V(\tree{U},\tree{S},F)$ and $V(\tree{V},\tree{T},G)$ are defined and, letting $\mu$ the greatest ordinal such that 
    $V(\itU,\itS,F)\restrict \mu$ is wellfounded and $\mu^*$ the greatest ordinal
 such that $V(\itV,\tree{T}, G)\restrict\mu^*$ is wellfounded, there is a unique partial weak tree embedding 
$\Gamma: V(\itU,\itS,F)\restrict\mu \to V(\itV,\tree{T}, G)\restrict\mu^*$ with maximal domain such that
 \begin{enumerate}
        \item $\Gamma\restrict \alpha_0+1\approx \Psi \restrict \alpha_0+1$,
        \item $u^\Gamma(\alpha_0)=\alpha^*_0$,  
        \item either \begin{enumerate}
            \item $\sigma^\Gamma_{\alpha_0}=\id$, or 
       \item $\alpha_0+1<\lh(\tree{S})$,  $\alpha_0^*=u(\alpha_0)$, and $\sigma^\Gamma_{\alpha_0}=\sigma^\Psi_{\alpha_0}$; and \end{enumerate}
        \item $\Gamma\circ \Phi^{V(\itU,\itS,F)} =\Phi^{V(\itV,\tree{T}, G)}\circ \Pi$ (on their common domain).
    \end{enumerate}
 Moreover, if $\Psi$ and $\Pi$ are in the $X$-case at everywhere, so is $\Gamma$; if $\Psi$ and $\Pi$ are in the $W$-case everywhere, so is $\Gamma$. If $V(\itV,\tree{T}, G)$ is wellfounded, then
 $V(\itU,\itS, F)$ is wellfounded 
  and $\Gamma$ is a (total) extended weak tree embedding 
  from $V(\itU, \itS, F)$ into $V(\itV,\itT, G)$. If $V(\itV,\tree{T}, G)$ is wellfounded and also $\Pi$ is non-dropping, then $\Gamma$ is a non-dropping extended weak tree embedding. 
  \end{lemma}
  
Note that this Shift Lemma implies the Shift Lemma for tree embeddings, as if $\Psi$, $\Pi$ are tree embeddings, they are in \textit{both} the $X$-case and $W$-case everywhere, so $\Gamma$ is as well, which implies $\Gamma$ is a tree embedding.

\begin{theorem}[Copying]\label{copying2} Let $\Gamma:\tree{S}\to \tree{T}$ be a non-dropping extended weak tree embedding. Let $\mtree{S}=\langle \tree{S}_\xi, \Phi^{\eta,\xi},F_\zeta\,|\, \xi,\zeta+1<\lh (\mtree{S})\rangle$ be a meta-tree on $\tree{S}$.

Then there is some largest $\mu\leq \lh (\mtree{S})$ such that there is a meta-tree $\Gamma\mtree{S}=\langle \tree{T}_\xi, \Psi^{\eta,\xi}, G_\zeta\,|\,\xi,\zeta+1<\mu\rangle$ on $\tree{T}$ with tree-order $\leq_\mtree{S}\restrict \mu$ and for $\xi<\mu$, non-dropping extended weak tree embeddings $\Gamma^\xi: \tree{S}_\xi\to \tree{T}_\xi$ with (total) last $t$-map $t_\infty^\xi$ such that 
\begin{enumerate}
    \item $\Gamma=\Gamma^0$, \item$G_\xi=t_\infty^\xi(F_\xi)$, 
    \item  and for all $\eta\leq_\mtree{S}\xi$, $\Gamma^\xi\circ \Phi^{\eta,\xi}=\Psi^{\eta,\xi}\circ \Gamma^\eta$.
\end{enumerate}  
\end{theorem}
\begin{proof}
Similar to copying meta-trees via ordinary tree embeddings.
\end{proof}

As usual, the copying construction guarantees that we have pullback strategies.

\begin{definition}
For $\tree{S}$, $\tree{T}$ countable plus trees of successor length on 
a premouse $M$, $\Phi:\tree{S}\to \tree{T}$ a non-dropping 
extended weak tree embedding, and $\Sigma$ a strategy for finite 
stacks of countable meta-trees on $\tree{T}$, we define the
 pullback strategy $\Sigma^\Phi$ for finite stacks of countable
meta-trees on $\tree{S}$ by
\[\mtree{S} \text{ is by }\Sigma^\Phi \Leftrightarrow \Phi\mtree{S} \text{ is by } \Sigma.\]
\end{definition}

\begin{lemma}\label{goodmetastrategy} Suppose $(M,\Sigma)$ is a mouse pair, 
$\tree{S}$ and $\tree{T}$ are countable plus trees on $M$, $\tree{T}$ is by $\Sigma$, and $\Phi:\tree{S}\to \tree{T}$
 is a non-dropping extended weak tree embedding. Let $\Lambda= 
\Sigma^*_\tree{T}$ be the induced meta-strategy for $\tree{T}$; then $\Lambda$ 
has meta-hull condensation, normalizes well, has the
 Dodd-Jensen property relative to $\Sigma$, and if $M$ is a lbr premouse, $\Lambda$ is pushforward consistent.
\end{lemma}

As a corollary to Lemma \ref{goodmetastrategy} and our main 
comparison theorem, we get the main theorem of this section.
(More properly, it is Theorem \ref{induced strategy theorem}
and its Lemma \ref{induced strategy lemma} that are relevant.)

\begin{theorem}\label{vshctheorem} Assume $\adp$.
Suppose $(M,\Sigma)$ is a mouse pair such that $\Sigma$ is coded by
a Suslin-co-Suslin set of reals. Let $\tree{S}$ and $\tree{T}$ be
plus tree on $M$. Suppose that
$\tree{T}$ is by $\Sigma$ and there is a weak tree embedding from
$\tree{S}$ into $\tree{T}$; then $\tree{S}$ is by $\Sigma$.
\end{theorem}

\subsection{Full normalization}
 
Finally, we describe the 
full normalization construction. We shall make some minor simplifying
assumptions as we do that, to the effect that certain ultrapowers are
not type 2 premice. These assumptions can be removed in a way that is
conceptually simple but notationally complicated, so we save the general case
for a later draft of this paper.
 
Suppose we are given a maximal stack of $\lambda$-tight, normal trees 
$s=\langle \tree{S}_0, \ldots, \tree{S}_n\rangle$. The first  of our
simplifying assumptions is that each $\tree{S}_i$
has a strongly stable base model. That guarantees that all models
of all $\tree{S}_i$ are of type 1.\footnote{Strong stability for the base model
of $\tree{S}_0$ does not imply strong stability for its last model, the base model of
$\tree{S}_1$, because the main branch of $\tree{S}_0$ might drop.}
We shall define, subject to further such simplifying assumptions to come,
a putative 
$\lambda$-tight, normal tree $X(s)$ and possibly partial weak tree embeddings 
$\Psi:\tree{S}_0\to X(s)$ and $\Gamma: X(s)\to W(s)$. By Theorem \ref{vshctheorem},
if
$s$ is by the strategy of some mouse pair $(P,\Sigma)$, then since $X(s)$ weakly embeds
into $W(s)$, $X(s)$ is
also 
by $\Sigma$. Our construction guarantees that $s$ and $X(s)$ and $s$ have the same last model.\footnote{We
must assume that $s$ is maximal, because normal trees are maximal. For  example, if
$s = \langle E, N, F \rangle$ is a non-maximal stack on $M$ with $\lh(E) < \crit(F)$, then there
is no normal tree on $M$ with the same last model as $s$. Thought of as a single tree instead of a stack,
$s$ is non-overlapping, but has a gratuitous drop, and so it is not normal.}

As with embedding normalization, we first handle 
the one-step case. Let $\tree{S},\tree{T}$ be 
$\lambda$-tight, normal trees of successor length on $M$, where $M$ is strongly stable.
Let $F$ be on the sequence of last model of $\tree{T}$. 
Let $\alpha=\alpha(F,\tree{T})$ and $\beta=\beta(F,\tree{T})$. 
Suppose that $\tree{S}\restrict\beta+1=
\tree{T}\restrict\beta+1$ and 
$\dom F\leq \lambda( E^\tree{S}_\beta)$, if $\beta+1<\lh \tree{S}$. 
Granted our simplifying assumptions,
 we shall define $\tree{X}=X(\tree{S},\tree{T},F)$, and
a 
partial extended weak tree embedding from $\tree{S}$ into $\tree{X}$, 
\begin{align*}
\Psi&=\Psi^{X(\tree{S},\tree{T},F)}\\
     &= \langle v,u, \{s_\xi\}_{\xi<\lh\tree{S}},
 \{t_\zeta\}_{\zeta+1<\lh\tree{S}},\{\sigma_\zeta\}_{\zeta+1<\lh \tree{S}}\rangle.
\end{align*}
We also define a partial extended weak tree embedding
\[
\Gamma=\Gamma^{X(\tree{S},\tree{T},F)}
\] 
 of $X$ into $\tree{W}=W(\tree{S}, \tree{T},F)$ such that
 \[
\Gamma\circ 
\Psi=\Phi=\Phi^{W(\tree{S},\tree{T},F)}.
\]
The component maps of $\Phi$ and $\Gamma$ we shall indicate by adding
superscripts. We shall have $u^\Gamma = \text{ id}$.
In $\Psi$ the $X$-case will occur everywhere, and in $\Gamma$ the $W$-case will
occur everywhere.

As with embedding normalization, we may reach illfounded models in 
forming $\tree{X}$ and stop when we do. We say that $\tree{X}$ is 
wellfounded 
if we never reach illfounded models. When $\tree{S}$ and $\tree{T}$ are
 by a strategy $\Sigma$ which has strong hull condensation,
 $\tree{X}$ will be wellfounded. In this case, $\tree{X}$ will have last 
model $\ult(P, F)$, where $P$ is the longest initial segment of the last model of 
$\tree{S}$ to which $F$ applies.
 and, moreover, the embedding normalization 
map of $W(\tree{S},\tree{T},F)$ will be the last
 $t$-map of $\Gamma$:
\[
\sigma^{W(\tree{S},\tree{T},F)} = t_\infty^\Gamma.
\]

We let $\tree{X}\restrict \alpha+1 =\tree{T}\restrict \alpha+1$ and $E^\tree{X}_\alpha=F$.
For the rest of $\tree{X}$, we consider cases.

\paragraph{The dropping case.}$F$ is applied to a proper initial 
segment $P\isneq M^\tree{S}_\beta|\lh E^\tree{S}_\beta$, if 
$\beta+1<\lh \tree{S}$, or $P\isneq M^\tree{S}_\beta$ if $\beta+1=\lh \tree{S}$. 

In this case we've described all of $\tree{X}$ already:
 \[
\tree{X}=\tree{T}\restrict\alpha+1 \conc \langle F\rangle.
\]
Notice here that $\ult(P,F)$ has type 1. For otherwise, letting
$k=k(P)$, we have that $\crit(F) = \eta_k^P$. Since $P$ is stable,
that implies $\eta_k^P < \rho_{k+1}(P)$, so that $\ult(P^+,F)$ makes sense.
Our case hypothesis is that 
$P\isneq M^\tree{S}_\beta|\lh E^\tree{S}_\beta$ or 
or $P\isneq M^\tree{S}_\beta$, so in either case, $F$ should have been
applied to $P^+$ rather than $P$.
  
Recall that in this case we also let $W(\tree{S},\tree{T},F)=
\tree{T}\restrict\alpha+1 \conc \langle F\rangle$, so $\tree{W}=\tree{X}$.
 
We let $\Psi$ be the identity on $\tree{S}\restrict\beta+1$ except 
we set $u(\beta)=\alpha+1$ and $t_\beta= i^P_F$. We also let $\Gamma$ be the identity
(on $\tree{X}=\tree{W}$). Recalling how we defined $\Phi^{W(\tree{S},\tree{T},F)}$
in the dropping case, we have $\Phi=\Gamma\circ \Psi$.

Note that $\ult(P,F)$ is the last model of $\tree{X}$,
 and $\sigma^{W(\tree{S},\tree{T},F)} = t_{\beta+1}^\Gamma = \id$,
 as desired

\paragraph{The non-dropping case.}

Suppose $F$ is applied to an initial segment $P\is M^\tree{S}_\beta$ 
with $M^\tree{S}_\beta|\lh E^\tree{S}_\beta\is P$,
 if $\beta+1<\lh \tree{S}$, or $F$ is total on $M^\tree{S}_\beta$ 
if $\beta+1=\lh \tree{S}$.  We define $\Psi, \Gamma$ and $\tree{X}$ as follows.

We define $u=u^\Psi$ as we did in embedding normalization:
\begin{equation*}
    u(\xi) =
    \begin{cases*}
      \xi & if $\xi<\beta$, \\
      \alpha+1+(\xi-\beta) & if $\xi\geq \beta$.
    \end{cases*}
  \end{equation*}
The models of $\tree{X}$ are given by setting $M_\theta^{\tree{X}} = M_\theta^{\tree{T}}$
if $\theta \le \alpha$, and 
\[
M_{u(\xi)}^{\tree{X}} = \begin{cases} \ult(P,F) & \text{ if $\xi = \beta$, and }\\
                                    \ult(M_\xi^{\tree{S}},F) & \text{ if $\xi > \beta$.}
\end{cases}
\]
Note that for all $\xi>\beta$,  $F$ is total
 on $M_\xi^\tree{S}$ and $\crit(F) < \rho^-(M_\xi^\tree{S})$. This follows
easily from

\begin{proposition}\label{local hod prop 1}
Let $\tree{U}$ a normal tree, $\xi+1<\lh \tree{U}$, and $\mu=\lh E^\tree{U}_\xi$. Then if $\xi<\theta<\lh\tree{U}$, then $\mu$ is a successor cardinal of $M_\theta^\tree{U}$ and for $k=\deg(M_\theta^\tree{U})$, $\mu<\rho_k(M_\theta^\tree{U})$.
\end{proposition}
\begin{proof} Well known and routine.
\end{proof}

However, it can now happen that some of the $\ult(M_\xi^{\tree{S}},F)$ for $\xi > \beta$
have type 2. Similarly, when $P=M_\beta^{\tree{S}}$ it can happen that $\ult(P,F)$ has type 2.
(The branch to the relevant model must have dropped to a premouse that is
not strongly stable in this case.) We wish to avoid this possibility, because it adds
a notational mess without requiring any important new ideas.

\bigskip
\noindent
{\em Simplifying assumption:} $\ult(P,F)$ and all
                                    $\ult(M_\xi^{\tree{S}},F)$ for
$\xi > \beta$ have type 1.

\bigskip
\noindent
We shall describe how to remove this assumption in a future draft.\footnote{The construction
we are doing produces a system $Y(\tree{S},\tree{T},F)$ that has premice of type 2
on it, together with a (generalized) weak tree embedding $\Phi$ from
$Y(\tree{S},\tree{T},F)$ to $W(\tree{S},\tree{T},F)$ whose $u$-map is the identity.
$\tree{Y}$ is not literally an iteration tree, but we can convert it to one that reaches
the same models, but using additional steps. This involves inserting 
ultrapowers by order zero measures into the full normalization $\tree{X}$, so that
 it can reach premice of
the form $\ult(Q,F)^-$, where $\ult(Q,F)$ has type 2. Note $\ult(Q,F)^-$ has type 1,
so it could be a model in a normal tree on $M$. Adding these order zero ultrapowers
means that $X(\tree{S},\tree{T},F)$ and $W(\tree{S},\tree{T},F)$ may no longer
have exactly the same tree order. Nevertheless, the existence of $\Phi$ implies that
$\tree{X}$ is ``good", granted that $\tree{W}$ is good.}

Under our simplifying assumption, $\tree{X}$ will have the same tree order and length as $\tree{W}$, 
and we set
$u^\Gamma=v^\Gamma=id$.
We can also now specify the remaining $t$-maps of $\Psi$:
\[
t^\Psi_\xi=i^{M_\xi^\tree{S}}_F.
\]
What is left is to find the extenders $E_\xi^{\tree{X}}$ 
that make $\tree{X}$ into
an iteration tree, and to finish defining $\Psi$ and $\Gamma$.

Proposition \ref{local hod prop 1} implies that for $\mu = \lh E_\xi^{\tree{S}}$ and
$\xi+1 \le \eta < \lh S$, $t_{\xi+1} \restrict (\mu +1) = t_\eta \restrict (\mu+1)$.
In general, we do not have that $t_\xi\restrict \mu =t_{\xi+1}\restrict \mu$ . What we have is the following diagram
\[
\begin{tikzcd}
M^\tree{S}_\xi \arrow{r}{t_\xi} & 
M^\tree{X}_{u(\xi)}\arrow[Eq]{r}& \ult(M_\xi^\tree{S}, F)\\
M^\tree{S}_\xi||\mu \arrow[Is]{u}{}
\arrow{r}{t_\xi} \arrow[swap]{rd}{t_{\xi+1}}&
t_\xi(M^\tree{S}_\xi||\mu)\arrow[Is]{u}{} \\
&
t_{\xi+1}(M^\tree{S}_\xi||\mu)\arrow[swap]{u}{\rs_\xi}\arrow[Eq]{r}& \ult(M_\xi^\tree{S}||\mu, F)
\end{tikzcd}
\]
$t_{\xi+1}(M_\xi^\tree{S}||\mu)$ is the ultrapower computed using functions
 in $M_\xi^\tree{S}||\mu$ and $t_{\xi}(M_\xi^\tree{S}||\mu)$ is the 
ultrapower computed using all functions in $M_\xi^\tree{S}$.
 $\rs_\xi$ is the natural factor map. (\lq\lq rs" is meant to
 suggest \lq\lq resurrection".) Having defined the $E_\gamma^{\tree{X}}$ for $\gamma < u(\xi)$,
we get at once from Proposition \ref{local hod prop 1}:

\begin{claim}\label{local hod claim 1}
For any $\gamma<u(\xi)$, $\rs_\xi\restrict\lh E^\tree{X}_{\gamma}+1 = id$. Also, $\rs_\xi 
\restrict\lh F +1 = id$.
\end{claim}

So for any $\theta\geq \xi+1$, $t_\theta\restrict\lh E_\xi^\tree{S}=
t_{\xi+1}\restrict\lh E^\tree{S}_\xi$. We define now  extenders 
$E_\gamma^\tree{X}$ which make $\tree{X}$ into a normal iteration tree. To start, we let
\[
E_\gamma^\tree{X}=\begin{cases} E_\gamma^\tree{T}\text{ if } \gamma<\alpha\\
F \text{ if } \gamma=\alpha\end{cases}.\]

Now let $\gamma>\alpha$, so $\gamma=u(\xi)$ for some $\xi\geq \beta$.
 Assume that $\xi>\beta$; the argument when $\xi=\beta$ is similar, 
but $M_\beta^\tree{S}$ gets replaced by the initial segment $P\is M_\beta^\tree{S}$
indicated above.
Let $\mu=\lh E_\xi^\tree{S}$. We have the following diagram.

\[
\begin{tikzcd}
M^\tree{S}_\xi \arrow{r}{t_\xi} & 
M^\tree{X}_{u(\xi)}\arrow[Eq]{r}& \ult(M_\xi^\tree{S}, F)\\
M^\tree{S}_\xi|\mu \arrow[Is]{u}{}
\arrow{r}{t_\xi} \arrow[swap]{rd}{t_{\xi+1}}&
t_\xi(M^\tree{S}_\xi|\mu)\arrow[Is]{u}{} \\
&
t_{\xi+1}(M^\tree{S}_\xi|\mu)\arrow[swap]{u}{\rs_\xi}\arrow[Eq]{r}& \ult(M_\xi^\tree{S}|\mu, F)
\end{tikzcd}
\]

The difference from the preceding diagram is just that $M_\xi^\tree{S}|\mu$
 has a predicate symbol $\dot F$ for $E^\tree{S}_\xi$, while
 $M_\xi^\tree{S}||$ is passive. The maps remain elementary (i.e. $\Sigma_1$ elementary)
 even with this added predicate.

Applying Lemma \ref{main full normalization lemma}, we have

\begin{claim}\label{local hod claim 2}
$\ult(M_\xi^\tree{S}|\mu, F)\is M_{u(\xi)}^\tree{X}$
 and $rs_\xi$ respects drops over $(M_{u(\xi)}^\tree{X}, t_{\xi+1}(\mu), t_\xi(\mu))$.
\end{claim}

We set
\begin{align*}
    E_{u(\xi)}^\tree{X} &= \dot F ^{\ult(M_\xi^\tree{S}|\mu, F)}\\
    &= \text{the top extender of }\ult(M_\xi^\tree{S}|\mu, F)\\
    &=\bigcup_{\alpha<\mu} t_{\xi+1}\big(E_\xi^\tree{S}\cap M_\xi^\tree{S}|\alpha\big).
\end{align*}

We may sometimes write
\[E_{u(\xi)}^\tree{X} = t_{\xi+1}(E_\xi^\tree{S})\]
though literally $E_\xi^\tree{S}\not\in \dom t_{\xi+1}$. We do
 always literally have $\lh E_\xi^\tree{S}\in \dom t_{\xi+1}$ 
and $t_{\xi+1}(\lh E_\xi^\tree{S})=\lh E^\tree{X}_{u(\eta)}$.
We now let
\begin{align*}
    G&=E_\xi^\tree{S}\\
    H&= t_\xi(G)\\
    \bar H &= t_{\xi+1}(G)=E_{u(\xi)}^\tree{X}.
\end{align*}

\begin{claim}\label{local hod claim 3}
\begin{enumerate}
    \item[(a)] For any $\delta<\xi$, $\lh E^\tree{X}_{u(\delta)}< \lh \bar H$
    \item[(b)] $\lh F<\lambda (\bar H)$
    \item[(c)] For any $\delta<\xi$, $\crit G<\lambda(E_\delta^\tree{S})\Leftrightarrow \crit  H< \lambda(E_{u(\delta)}^\tree{X})\Leftrightarrow \crit  \bar H< \lambda(E_{u(\delta)}^\tree{X})$
    \item[(d)] If $\crit G<\lambda(E_\delta^\tree{X})$, then $\crit H= \crit \bar H$. In fact, $H\restrict\lh E^\tree{X}_{u(\delta)}=\bar H\restrict \lh E^\tree{X}_{u(\delta)}$.
\end{enumerate}
\end{claim}
\begin{proof}
For (a), let $\delta<\xi$. Then, using Claim \ref{local hod claim 1}, $\lh E_\delta^\tree{S}<\lh E^\tree{S}_\xi$, so \begin{align*}
\lh E^\tree{X}_{u(\delta)}&=t_{\delta+1}(\lh E^\tree{S}_\delta)=t_{\xi+1}(\lh E^\tree{S}_\delta)\\&< t_{\xi+1}(\lh E_\xi^\tree{S})=\lh E^\tree{X}_{u(\xi)},
\end{align*}
as desired.

For (b), $\crit F^+<\lambda(E_\xi^\tree{S})$, so \[i_F^{M_\xi^\tree{S}|\mu}(\crit F^+)=\lh F<i_F^{M_\xi^\tree{S}|\mu}(\lambda(E_\xi^\tree{S}))=\lambda(E^\tree{X}_{u(\xi)}).\]

For (c), let $\kappa=\crit G=\crit E^\tree{S}_\xi$. So $t_\xi(\kappa)=\crit H$ and $t_{\xi+1}(\kappa)=\crit \bar H.$ Then for $\delta<\xi$,
\begin{align*}
    \kappa<\lambda(E^\tree{S}_\delta) &\text{\quad iff \quad} t_{\delta+1}(\kappa)<\lambda(E^\tree{X}_{u(\delta)})\\
    &\text{\quad iff \quad}t_\xi(\kappa)<\lambda(E^\tree{X}_{u(\delta)})\\
    &\text{\quad iff \quad}t_{\xi+1}(\kappa)<\lambda(E^\tree{X}_{u(\delta)}),
\end{align*}
using on the second and third lines that $t_{\delta+1}, t_\xi$, and $t_{\xi+1}$ all agree on $\lh E^\tree{S}_\delta+1$.

(d) is clear.
\qed
\end{proof}

By Claim \ref{local hod claim 3}, setting $E_{u(\xi)}^\tree{X}=\bar H$ 
preserves the length-increasing clause in the definition of normality. 
We now need to check that this choice of extender gives rise to the
 appropriate next model when we apply it where we must using the Jensen normality rules.

Let $\delta=\tree{S}\pred(\xi+1)$. We now break into cases.

\paragraph{Case 1.} $\crit G<\crit F$.

In this case, since $\crit F<\lambda(E_\beta^\tree{S})$, 
$\delta\leq \beta$. If $\delta<\beta$ (so that $u(\delta)=\delta$),
 then Claim \ref{local hod claim 3} (b) tells us that $\bar H$ must
 be applied in $\tree{X}$ to the same $Q\is M_\delta^\tree{S}$ that
 $G$ is applied to in $\tree{S}$. In fact, $\crit \bar H=\crit G=\crit H$.


We then have the commutative diagram

\[
\begin{tikzcd}
M^\tree{S}_{\xi+1} \arrow{r}{t_{\xi+1}} & 
M^\tree{X}_{u(\xi+1)}\arrow[Eq]{r}& \ult(M_{\xi+1}^\tree{S}, F)\arrow[Eq]{r}& \ult(Q, \bar H)\\
Q\arrow{u}{G}\arrow[swap]{ur}{\bar H}\arrow[Is]{d}\\
M_\delta^\tree{S}\arrow[Eq]{r} &M_\delta^\tree{X}.
\end{tikzcd}
\]

It is shown in \cite[\S 6.1]{nitcis} that $\ult(M_{\xi+1}^\tree{S}, F)=
\ult(Q, \bar H)$ and that this diagram commutes.
 (See the proof of Claim \ref{local hod claim 5}
 in Case 2, below, for a similar calculation.)
The situation when $\delta=\beta$ is the same:
 $\tree{X}\pred(\xi+1)=\beta$
 and $\bar H$ is applied to the same $Q$ that $G$ was.
 Note that $u(\beta)\neq \beta$, 
so $u$ does not preserve tree order, just as in embedding normalization.
\paragraph{Case 2.} $\crit F\leq \crit G$.

In this case, $\delta\geq \beta$. Also, $\lambda(F)\leq \crit \bar H$,
 so $\bar H$ is applied in $\tree{X}$ to some $Q\is M_\tau^\tree{X}$, 
where $\tau\geq \alpha+1$. Thus $\tau\in \ran(u)$ and, by Claim 
\ref{local hod claim 3}, $\tau=u(\delta)$. 
That is, $\tree{X}\pred(u(\xi+1))=u(\delta)$.

Let $\kappa=\crit G$ and $P\is M_\delta^\tree{S}$ be such that 
$M_{\xi+1}^\tree{S}=\ult(P,G)$.

\begin{claim}\label{local hod claim 4}
$\bar H$ is applied to $\ult(P,F)$ in $\tree{X}$.
\end{claim}

\begin{proof}
We have the following diagram.

\[
\begin{tikzcd}
M^\tree{S}_\delta \arrow{r}{t_\delta} & 
M^\tree{X}_{u(\delta)}\\
P \arrow[Is]{u}{}
\arrow{r}{t_\delta} \arrow[swap]{rd}{i^P_F}&
t_\delta(P)\arrow[Is]{u}{} \\
&
\ult(P,F)\arrow[swap]{u}{l}\\
M_\delta^\tree{S}|\lh E^\tree{S}_\delta \arrow[Is]{uu}{}
\arrow{r}{i^P_F} \arrow[swap]{rd}{t_{\delta+1}}&
i^P_F(M_\delta^\tree{S}|\lh E^\tree{S}_\delta)\arrow[Is]{u}{} \\
&
\ult(M_\delta^\tree{S}|\lh E^\tree{S}_\delta,F)\arrow[swap]{u}{k}\arrow[bend right=60, swap]{uuu}{\rs_\delta}
\end{tikzcd}
\]

$k$ and $l$ are the natural factor maps and
 \[
\rs_\delta=l\circ k.
\]
By Lemma \ref{main full normalization lemma}, 
$\ult(M_\delta^\tree{S}|\lh E^\tree{S}_\delta, F)\is i^P_F(M_\delta^\tree{S}|\lh E^\tree{S}_\delta)$
 and $\rs_\delta$ respects drops over $(M_{u(\delta)}^{\tree{X}}, t_{\delta+1}(\lh E^{\tree{S}}_\delta),
t_\delta(E^{\tree{S}}_\delta)$.
 Note that 
\[
k\restrict t_{\delta+1}((\crit G^+)^{M_\delta^\tree{S}|\lh E^\tree{S}_\delta})= \text{ id}
\]
because $\rho_{deg(P)}(P)>\crit G$. 
So $t_{\delta+1}(\dom G)=i_F^P(\dom G)$ and 
$t_{\delta+1}\restrict(\crit G^+)^P=i^P_F\restrict (\crit G^+)^P$.

We have that $P = A_e(M_\delta^\tree{S}, M_\delta^{\tree{S}}|\lh E^\tree{S})^-$ for some $e \le n$,
where $n = n( M_\delta^\tree{S}, M_\delta^{\tree{S}}|\lh E^\tree{S})$.
By the proof of Lemma \ref{main full normalization lemma},
\begin{align*}
 \ult(P,F) &= B^e_e \is t_\delta(P),\\
          k &= \sigma_{e-1} \circ ...\circ \sigma_1,\\
\intertext{ and } 
   l &= \sigma_n \circ ...\circ \sigma_e,
\end{align*} where $\rs_\delta = \sigma_n \circ ...\circ \sigma_1$ is the factoring
of $\rs_\delta$ induced by the dropdown sequence of 
$(M_\delta^\tree{S}, M_\delta^{\tree{S}}|\lh E^\tree{S})$. Note here that we are using our
simplifying assumption to conclude that $B^e_e \lhd t_\delta(P)$.

Note that 
for $\kappa=\crit G$, $\rho(\ult(P,F))\leq i^P_F(\kappa)$, 
because $\ult(P,F)$ is generated by $i^P_F(\rho(P))\cup i^P_F"\kappa \cup \lh F$,
 and $\lh F<i^P_F(\kappa)$. But for $k=\deg(P)$, $\rho_k(P)\geq 
(\kappa^+)^P$, so $\rho_k(\ult(P,F))\geq i^P_F((\kappa^+)^P)$. It 
follows that $\bar H$, whose domain is $t_{\xi+1}(\dom G)=
t_{\delta+1}(\dom G)=i^P_F(\dom G)$, is applied to $\ult(P,F)$ in $\tree{X}$.
\qed
\end{proof}

\begin{claim}\label{local hod claim 5}
$\ult(\ult(P,F), \bar H)=M^\tree{X}_{u(\xi+1)}= \ult(\ult(P,G), F).$
\end{claim}
\begin{proof}
This is shown in \cite[\S 6.1]{nitcis}, but we repeat the calculations here.

Set $N=\ult(P,G)$ and $Q=\ult(P,F)$. We have the diagram 
\[
\begin{tikzcd}
N \arrow{r}{i^N_F} & 
\ult(N,F)& \ult(Q, \bar H)\\P\arrow{u}{i^P_G}\arrow{r}{i^P_F}& Q\arrow{ur}
\end{tikzcd}
\]

Let $E$ be the extender of $i^N_F\circ i^P_G$; then $\nu(E)\leq \sup i^N_F"\lambda(G)$ 
and for $a\in [\nu(E)]^{<\omega}$, $E_a$ concentrates on $\crit G^{|a|}$. 
Let $K$ be the extender of $i^Q_{\bar H}\circ i^P_F$, so $\nu(K)\leq
 \lh \bar H=\sup i^N_F"\lambda(G)$, and each $K_a$ 
concentrates on $\crit G^{|a|}$. Let $a=[b,g]^N_F$, 
where $g\in N|\lambda(G)=M_\xi^\tree{S}|\lambda(G)$ be a 
typical element of $[\sup i^N_F"\lambda(G)]^{<\omega}$
 and $A\subseteq \crit G^{|a|}$; then
\begin{align*}
    (a,A)\in E &\text{\quad iff \quad} [b,g]^N_F\in i^N_F\circ i^P_G(A)\\
    &\text{\quad iff \quad} \text{for $F_b$-a.e. $u$, $g(u)\in i^P_G(A)$}\\
    &\text{\quad iff \quad} \text{for $F_b$-a.e. $u$, } (g(u),A)\in G\\
    &\text{\quad iff \quad} \big([b,g]^{M_\xi^\tree{S}|\lh G}_F, i^{M_\xi^\tree{S}|\lh G}_F(A)\big)\in \bar H\\
    &\text{\quad iff \quad} \big([b,g]^N_F, i^P_F(A)\big)\in \bar H\\
    &\text{\quad iff \quad} [b,g]^N_F\in i^Q_{\bar H}\circ i^P_F(A)\\
    &\text{\quad iff \quad} (a,A)\in K,
\end{align*}
using in the fifth line that $[b,g]^N_F=[b,g]^{M_\xi^\tree{S}|\lh G}_F$ and $i^{M_\xi^\tree{S}|\lh G}_F(A)=i^P_F(A)$.

So $E=K$ and $\ult(N,F)= \ult(Q,\bar H)$, as desired.
\end{proof}


The proof of the previous proposition shows also that
 $i^N_F\circ i_G=i_{\bar H}\circ i^P_F$. So we have verified that
 as long as all the models $\ult(M_\xi^\tree{S},F)$ are wellfounded, $\tree{X}$ is a normal tree.

Putting together our work above, we have the following.
\begin{proposition}
\begin{enumerate}
    \item $\tree{X}\restrict\alpha+1=\tree{T}\restrict\alpha+1$
    \item $M_{\alpha+1}^\tree{X}=\ult(P,F)$ and $t_\beta=i^P_F$, for $P\is M_\beta^\tree{S}$,
    \item for $\xi>\beta$ $M_{u(\xi)}^\tree{X}=\ult(M_\xi^\tree{S}, F)$ and $t_\xi=i^{M_\xi^\tree{S}}_F$
    \item $t_\xi=\rs_\xi\circ t_{\xi+1}$, $\rs_\xi$ respects drops, and $E^\tree{X}_{u(\xi)}=t_{\xi+1}(E^\tree{S}_\xi)$,
    \item for $\delta\leq_\tree{S}\xi$ either
    \begin{enumerate}
        \item $u(\delta)\leq_\tree{X} u(\xi)$ and $t_\xi\circ \hat\imath^\tree{S}_{\delta, \xi}=\hat\imath^\tree{X}_{u(\delta), u(\xi)}\circ t_\delta$, or
        \item $\delta=\beta$ and for $\zeta+1$ the successor of $\beta$ in $[\beta, \xi]_\tree{S}$,
        
        $\crit E_\zeta^\tree{S}<\crit F$, $\beta\leq_\tree{X} u(\xi)$ and $t_\xi\circ \hat\imath^\tree{S}_{\beta, \xi}=\hat\imath^\tree{X}_{\beta, u(\xi)}$.
    \end{enumerate}
\end{enumerate}
\end{proposition}
For $\xi<\beta$, we let $t_\xi=id$ and $\rs_\xi=id$. We then have that 
$\Psi= \langle u, \{t_\xi\}, \{\rs_\xi\}\rangle$, $\Psi$ is an extended 
weak tree embedding from $\tree{S}$ into $\tree{X}$ which is in the 
$X$-case at every ordinal $\xi+1$.

It just remains to describe the weak tree embedding
 $\Gamma:\tree{X}\to \tree{W}$ and verify that 
$\Gamma\circ \Psi=\Phi$. Recalling the one-step embedding normalization,
 our work above shows that $\tree{X}$ and $\tree{W}$ have 
the same length and tree order and $u=u^\Psi=u^\Phi$.\footnote{Without our
simplifying assumption, this might not be true.}
We let $u^\Gamma=id$. We just need to define the 
$t$-maps of $\Gamma$, which we call
$\gamma_\xi$,  and the $\sigma$-maps $ \rs^*_\xi$.
\[
\Gamma = \langle \text{ id}, \langle \gamma_\xi \mid \xi < \lh \tree{X} \rangle,
\langle \rs^*_\xi \mid \xi+1 < \lh \tree{X} \rangle \rangle.
\]
 $\rs^*_\xi$
will be in the $W$-case at every $\xi$. Let $\pi_\eta=t^\Phi_\eta$.

Since $\tree{W}\restrict\alpha+2=\tree{T}\restrict\alpha+1\conc\langle F\rangle 
=\tree{X}\restrict\alpha+2$, we'll have $\gamma_\xi=id$ and $\rs_\xi=id$ for all
 $\xi\leq \alpha+1$. Now suppose for all $\eta\leq \xi$ we've defined partial
 elementary maps $\gamma_{u(\eta)}:M_{u(\eta)}^\tree{X}\to M^\tree{W}_{u(\eta)}$
 and maps $\rs^*_{u(\eta)}: M^\tree{W}_{u(\eta)}|\gamma_{u(\eta)}
(\lh E^\tree{X}_{u(\eta)})\to M^\tree{W}_{u(\eta)}|\lh E^\tree{W}_{u(\eta)}$
 which respect drops such that \begin{enumerate}
    \item $\pi_\eta = \gamma_{u(\eta)}\circ t_\eta$,
\item $\gamma_{u(\xi)}\restrict\lh F=id$,
\item if $\beta\leq \eta<\zeta$, then $\gamma_{u(\zeta)}\restrict
\lh E^\tree{X}_{u(\eta)}= \gamma_{u(\eta)}\circ \rs_\eta \restrict\lh E^\tree{X}_{u(\eta)}$
\item  $\rs^*_{u(\eta)}\circ\gamma_{u(\eta)}=\gamma_{u(\eta)}\circ \rs_\eta$
\end{enumerate}

Let $G, H, \bar H$ as before. So (1) gives us that $E^\tree{W}_{u(\xi)}=\pi_\xi(E^\tree{S}_\xi)=\gamma_\xi(H)$. Let $\delta=\tree{S}\pred(\xi+1)$. 

\paragraph{Case 1.} $\crit G<\crit F$.

In this case, $\dom \bar H =\dom H=\dom G$, $\delta\leq \beta$, and 
$\delta=\tree{X}\pred(u(\xi+1))=\tree{W}\pred(u(\xi+1))$. 
Suppose $\bar H$ is applied to $P\is M_\delta^\tree{X}=M_\delta^\tree{S}$ in 
$\tree{X}$. Then $H$ would also be applied to $P$ if we had 
set $E_{u(\xi)}=H$. Now, $\bar H$ is a subextender of $H$ under $\rs_\xi$:
\[(a, A)\in \bar H \quad \text{iff} \quad (\rs_\xi(a), A)\in H.\]

So we let $\epsilon$ be the natural map from $M_{u(\xi+1)}^\tree{X}=
\ult(M_{\xi+1}^\tree{S}, \bar H)$ into $\ult(P, H)$, i.e.
\[\epsilon\big([a,f]^P_{\bar H}\big)=[\rs_\xi(a), f]^P_H.\]

So $\epsilon\restrict\lh \bar H=\rs_\xi\restrict\lh \bar H$.
We have the diagram
\[
\begin{tikzcd}
M^\tree{S}_{\xi+1} \arrow{r}{t_{\xi+1}} & 
M^\tree{X}_{u(\xi+1)}\arrow{r}{\epsilon}& \ult(P, H)\\
P\arrow{u}{G}\arrow{ur}{\bar H}\arrow[swap]{urr}{H}\arrow[Is]{d}\\
M_\delta^\tree{S}\arrow[Eq]{r} &M_\delta^\tree{X}.
\end{tikzcd}
\]
Since $\delta\leq \beta$, we have $\gamma_\xi\restrict\dom H=id$, 
so $\gamma_\xi(H)$ is applied to $P$ in $\tree{W}$, and we have the diagram
\[
\begin{tikzcd}
M^\tree{S}_{\xi+1} \arrow{r}{t_{\xi+1}} & 
M^\tree{X}_{u(\xi+1)}\arrow{r}{\epsilon}& \ult(P, H)\arrow{r}{\theta}&M^\tree{W}_{u(\xi+1)}\\
P\arrow{u}{G}\arrow{ur}{\bar H}\arrow[swap]{urr}{H}\arrow[bend right=10, swap]{urrr}{\gamma_{u(\xi)}(H)}\arrow[Is]{d}\\
M_\delta^\tree{S}\arrow[Eq]{r} &M_\delta^\tree{X}\arrow[Eq]{r} &M_\delta^\tree{W}.
\end{tikzcd}
\]

where $\theta$ is given by
 \[
\theta([a,f]^P_H)=[\gamma_{u(\xi)}, f]^P_{\gamma_{u(\xi)}(H)}.
\] 
So $\theta$ agrees with $\gamma_{u(\xi)}$ on $\lh H$. We set
\[
\gamma_{u(\xi+1)}=\theta\circ \epsilon.
\]
So by the agreement of $\theta$ with $\gamma_{u(\xi)}$ and $\epsilon$ with
$\rs_\xi$,
\[
\gamma_{u(\xi+1)}\restrict\lh \bar H=\gamma_{u(\xi)}
\circ \rs_\xi\restrict\lh \bar H,
\] which gives us (3) at $\xi+1$.
 It is easy to see that 
\[
\pi_{u(\xi+1)}=\gamma_{u(\xi+1)}\circ t_{\xi+1},
\]
so we have (1) at $\xi+1$.
We let 
\[
\rs_{u(\xi+1)}^*=\gamma_{u(\xi+1)}(\rs_{\xi+1}).
\]
(It can happen that $\rs_{\xi+1}$ does not belong to
$M_{u(\xi+1)}^{\tree{X}}$, but in any case $\rs_{\xi+1}$ is
definable from parameters over $M_{u(\xi+1)}^{\tree{X}}$
in a way that is simple enough that we can move it
by $\gamma_{u(\xi+1)}$.) This gives us (4).


\paragraph{Case 2.} $\crit F\leq \crit G$.

Suppose first that $\delta<\xi$. This yields that
 $t_\xi\restrict \lh E^\tree{S}_\delta=t_{\xi+1}\restrict\lh E^\tree{S}_{\delta}$,
 so $t_\xi$ and $t_{\xi+1}$ agree on $\dom G$, 
so $\dom \bar H=\dom H$. Moreover, $\rs_\xi\restrict 
\dom \bar H=id$. Let $P\is M^\tree{S}_\delta$ be what
 $G$ is applied to in $\tree{S}$. We have:
\[
i^P_F\restrict\dom G=t_{\delta+1}\restrict\dom G=t_\xi\restrict\dom G=t_{\xi+1}\restrict\dom G.
\]

As in the previous case, $\bar H$ is a subextender of $H$ via $\rs_\xi$, 
so we let $\epsilon$ be the copy
map from $M_{u(\xi+1)}^\tree{X}=\ult(M_{\xi+1}^\tree{S},
 \bar H)$ into $\ult(P, H)$, as in Case 1.
 We again have 
\[
\epsilon\restrict\lh \bar H=\rs_\xi\restrict\lh \bar H.
\]

By (1) at $\xi$, we have that 
\[
E^\tree{W}_{u(\xi)}=
\pi_{\xi}(E^\tree{S}_\xi)=\gamma_{u(\xi)}(t_\xi(E_\xi^\tree{S}0)=
\gamma_{u(\xi)}(H).
\]
 By our case hypothesis, we have that $E^\tree{W}_{u(\xi)}$ is
 applied to $\pi_\delta(P)\is M_{u(\delta)}^\tree{W}$ in
 $\tree{W}$. Let
\begin{align*}
P &= A_e(M_\delta^{\tree{S}}, M_\delta^{\tree{S}}|\lh E_\delta^{\tree{S}})^-,\\
k &=\sigma_{e-1} \circ ...\circ \sigma_1,\\
l &= \sigma_n \circ ...\sigma_e,\\
\intertext{ where}
\rs_\delta &= \sigma_n \circ ...\sigma_1
\end{align*}
is the resolution of the factor map $\rs_\delta$ given by the fact that it respects drops
over $(M_{u(\delta)}^{\tree{X}}, t_{\delta+1}(\lh E_\delta^{\tree{S}},t_\delta
(\lh (E_\delta^{\tree{S}}))$.
We have that
$k:\ult(M_\delta^\tree{S}|\lh E^\tree{S}_\delta, F)\to \ult(P,F)$
 and
$l:\ult(P,F)\to t_\delta(P)$, moreover
\begin{align*}
k\restrict\dom G & = \text{ id},\\
\intertext{ so}
l\restrict\dom G & =\rs_\delta\restrict\dom G.
\end{align*}
So we have
\[
\gamma_{u(\xi)}\restrict\dom G=\gamma_{u(\delta)}\circ \rs_\delta\restrict\dom G=\gamma_{u(\delta)}\circ l\restrict\dom G.
\]
So we may let $\theta$ be given by the Shift Lemma applied to 
$\gamma_{u(\delta)}\circ l$, $\gamma_{u(\xi)}$, and $H$, i.e. 
\[
\theta\big([a, f]^{\ult(P, F)}_H\big)=[\gamma_{u(\xi)}(a), \gamma_{u(\delta)}\circ l(f)]^{\pi_\delta(P)}_{\gamma_{u(\xi)}(H)}.
\]

We have the following commutative diagram.

\[
\begin{tikzcd}
M^\tree{S}_{\xi+1} \arrow{r}{t_{\xi+1}} & 
M^\tree{X}_{u(\xi+1)}\arrow{r}{\epsilon} & \ult(P, H)\arrow{r}{\theta}& M^\tree{W}_{u(\xi+1)}\\
P\arrow[Is]{dd}\arrow{u}{G}\arrow{r}{i^P_F}\arrow{dr}{t_\delta} &\ult(P,F)\arrow{u}{\bar H}\arrow[swap]{ur}{H}\arrow{d}{l}\\
& t_\delta(P) \arrow[Is]{d}\arrow{r}{\gamma_{u(\delta)}}& \pi_\delta(P)\arrow[ swap]{uur}{\gamma_{u(\xi)}(H)}\arrow[Is]{d}\\
M_\delta^\tree{S}\arrow{r}{t_\delta} & M_{u(\delta)}^\tree{X}\arrow{r}{\gamma_{u(\delta)}} & M_{u(\delta)}^\tree{W}
\end{tikzcd}
\]
We now let \[\gamma_{u(\xi+1)}=\theta\circ \epsilon.\]
Since  $\epsilon \restrict\lh \bar H=\rs_\xi\restrict\lh \bar H$
 and $\theta\restrict\lh H=\gamma_{u(\xi)}\restrict\lh H$, we get
\[
\gamma_{u(\xi+1)}\restrict\lh \bar H = \gamma_{u(\xi)}\circ \rs_\xi\restrict\lh \bar H,
\]
giving (3) at $\xi+1$. As in the previous case, we let
 \[
\rs_{u(\xi+1)}^*=\gamma_{u(\xi+1)}(\rs_\xi),
\]
 which guarantees (4) at $\xi+1$. 

It remains to show (1) at $\xi+1$, i.e. 
that $\pi_{\xi+1} = \gamma_{u(\xi+1)}\circ t_{\xi+1}$.
Note first that the two sides agree on $\ran i^P_G$, for 
letting $j=\hat\imath^\tree{W}_{u(\delta), u(\xi+1)} \colon
\pi_\delta(P)\to M^\tree{W}_{u(\xi+1)}$, we have
\begin{align*}
    \theta\circ \epsilon\circ t_{\xi+1}\circ i^P_G&=j\circ \gamma_{u(\delta)}\circ t_\delta\\
    &=j\circ \pi_\delta\\
    &=\pi_{\xi+1}\circ i^P_G,
\end{align*}
using the commutativity properties of the maps in embedding normalization.

But $M^\tree{S}_{\xi+1}$ is generated by $\ran i^P_G \cup \lambda(G)$, so it is enough to see that $\theta\circ \epsilon\circ t_{\xi+1}$ and $\pi_{\xi+1}$ agree on $\lambda(G)$. Since $\pi_{\xi+1}$ agrees with $\pi_\xi$ on $\lambda(G)$, we get
\begin{align*}
    \pi_{\xi+1}\restrict\lambda(G)&=\pi_\xi\restrict\lambda(G)\\
    &=\gamma_{u(\xi)}\circ t_\xi\restrict\lambda(G)\\
    &=\gamma_{u(\xi)}\circ \big(\rs_\xi\circ t_{\xi+1}\big)\restrict\lambda(G)\\
    &= \big(\gamma_{u(\xi)}\circ \rs_\xi\big)\circ t_{\xi+1}\restrict\lambda(G)\\
    &=\gamma_{u(\xi+1)}\circ t_{\xi+1}\restrict\lambda(G),
\end{align*}
as desired.

This finishes Case 2 under the additional hypothesis that $\delta<\xi$. The case $\delta=\xi$ is not different in any important way. In that case, we may have $\crit\bar H<\crit H$, but the relevant diagram is the same. We leave it to the reader to confirm this.

For $\lambda$ a limit ordinal, we define $\gamma_\lambda = t_\lambda^\Gamma$ by setting
$\gamma_\lambda(i_{\eta,\lambda}^{\tree{X}}(a)) = i_{\eta,\lambda}^{\tree{W}}(\gamma_\eta(a))$
for all sufficiently large $\eta$. Again, $\rs^*_\lambda = \gamma_\lambda(\rs_\lambda)$.
We leave it to the reader to check (1)-(4).

If $\tree{W}$ wellfounded (for example if $\tree{S}$ and $\tree{T}$ are
 by a strategy with strong hull condensation), our induction hypotheses imply
 that $\Gamma$ is an extended weak tree embedding from $\tree{X}$ into 
$\tree{W}$ which is in the $W$-case at every $\xi$,
 and $\Phi=\Gamma\circ \Psi$. 
In this case, let $\delta+1=\lh \tree{S}$, so $u(\delta)=\lh\tree{W}=\lh\tree{X}$. 
Recall that $\sigma=\sigma^{W(\tree{S},\tree{T},F)}$ is 
the natural factor map witnessing that $F$ is an initial segment of
 the extender of $\pi_\delta$, which is just to say that $\sigma$ 
is the unique nearly elementary map from $\ult(M_\delta^\tree{S}, F)
=M^\tree{X}_{u(\delta)}$ into $M^\tree{W}_{u(\delta)}$ such
 that $\sigma\restrict\lh F=id$ and $\pi_\delta=\sigma\circ 
i^{M_\delta^\tree{S}}_F=\sigma\circ t_{\delta}$. 
But $\gamma_{u(\delta)}$ has both of these properties, 
so $\gamma_{u(\delta)}=\sigma$, as desired.

This finishes our work in the non-dropping case.

We define now the full normalization
$X(\tree{T}, \tree{U})$ of a maximal stack $\langle \tree{T},\tree{U}\rangle$.

Formally, this will be another variety of meta-tree, where we do one-step 
full normalization at every step instead of one-step embedding normalization. 
We will not formally define this kind of meta-tree, however. 
We define \[\mtree{X}(\tree{T},\tree{U})=\langle \tree{X}_\xi,
 E^\tree{U}_\zeta,\Psi^{\eta,\xi}\,|\, \eta,\xi,\zeta+1<\lh \tree{U}, 
\, \eta\leq_\tree{U} \xi\rangle,\] where
\begin{enumerate}
    \item $\tree{X}_0 = \tree{T}$ and for all $\xi<\lh \tree{U}$, $\tree{X}_\xi$ 
is a normal tree with last model $M^\tree{U}_\xi$,
    \item for $\zeta\leq_\tree{U}\eta\leq_\tree{U}\xi$, 
    \begin{enumerate}
        \item $\Psi^{\eta,\xi}:\tree{X}_\eta\to \tree{X}_\xi$ is
 an $X$-type partial extended weak tree embedding, and
        \item $\Psi^{\zeta,\xi}=\Psi^{\eta,\xi}\circ \Psi^{\zeta,\eta}$;
    \end{enumerate}
    \item For $\eta=\tree{U}\pred(\xi+1)$,
    \begin{enumerate}
        \item $\tree{X}_{\xi+1}= X(\tree{X}_\eta,\tree{X}_\xi,E^\tree{U}_\xi)$,
        \item $\Psi^{\eta,\xi+1}=\Psi^{X(\tree{X}_\eta,\tree{X}_\xi,E^\tree{U}_\xi)}$
    \end{enumerate}
    \item For $\lambda <\lh \tree{U}$ a limit and $b=[0,\lambda)_\tree{U}$,
    \[\tree{X}_\lambda = \lim \langle\tree{X}_\xi,\Psi^{\eta,\xi}\,|\,
 \eta\leq_\tree{U}\xi\in b\rangle\] and
    $\Psi^{\xi, \lambda}$ is the direct limit weak tree embedding.
\end{enumerate}

It is easy to verify by induction that we have the necessary
 agreement hypothesis so that clause (3)(a) makes sense.

Clause (4) relies on taking direct limits of
 systems of $X$-type weak tree embeddings, which we have not yet
 discussed. This exactly parallels direct limits of tree embeddings, up
 to the point at which we choose the exit extenders of the direct limit.
That is, we must define $E_x$ for a $u$-thread $x$. 
Recall that we only do this when $M_x$ is wellfounded 
and there is an $a\in \dom x$ such that for all $b\succeq a$,
 $t^{a,b}_{x(a)}$ (which is defined to be $t_{x(a)}^{\Psi^{a,b}}$) is  total. 
In this case, we may actually assume that $a$ is such that
 for all $b\preceq a$, $t^{a,b}_{x(a)}(E^a_{x(a)})=E^b_{x(b)}$, 
that is, that $\rs^{a,b}_{x(a)}=id$. 
This is because if $\rs^{a,b}_{x(a)} \neq id$, then, since 
we are in the $X$-case at every step, $\lh E^b_{x(b)}< t^{a,b}_{x(a)}(E^a_{x(a)})$,
 so if for all $a$ there is a $b\preceq a$ such that $\rs^{a,b}_{x(a)}\neq id$, 
then the images of the lengths of these exit extenders forms an 
infinite decreasing sequence in the ordinals of $M_x$, contradicting
 that it is wellfounded. 
So we let $E_x$ be the image of the stabilized valued
 of $E^a_{x(a)}$ under the direct limit $t$-map $t^a_x$. 
Finally, for $c\in \dom x$, we also define the $\sigma$-map
 of the direct limit weak tree embedding $\Psi^c$ as the 
common value $\rs^c_x=t^b_x(\rs^{c,b}_{x(c)})$ for
 any $b\preceq a,c$. Note that in dealing with direct limits 
of arbitrary weak tree embeddings abstractly, we must just assume 
that there is an $a\in \dom x$ such that $t^{a,b}_{x(a)}(E^a_{x(a)})=E^b_{x(b)}$
in order to define the direct limit. 
With this additional hypothesis in the definition of a wellfounded direct limit, we get the obvious analogue of Proposition \ref{direct limit prop}.

Let $\mtree{W}(\tree{T},\tree{U})=\langle \tree{W}_\xi, 
F_\zeta,\Phi^{\eta,\xi}, \,|\,\eta,\xi,\zeta+1<\lh \tree{U},
 \eta\leq_\tree{U} \xi\rangle$. Let $\sigma_\xi$ the embedding
 normalization map from $M_\xi^\tree{U}$ into the last model of
 $\tree{W}_\xi$. We define $W$-type tree embeddings $\Gamma^\xi:\tree{X}_\xi\to \tree{W}_\xi=W(\tree{T},\tree{U}\restrict\xi+1)$, by induction, such that \begin{enumerate}
    \item for all $\eta\leq_\tree{U}\xi$, 
$\Phi^{\eta,\xi}=\Gamma^\xi\circ \Psi^{\eta, \xi}$,
\item $\sigma_\xi$ is the last $t$-map of $\Gamma^\xi$.
\end{enumerate}
Note that (2) implies that $F_\xi$ is the image of $E_\xi^\tree{U}$ under the last $t$-map of $\Gamma_\xi$.

To start, we let $\Gamma_0=Id_\tree{T}$. At limits, $\Gamma_\lambda$ is exists 
and is uniquely determined by the commutativity condition
 and the $\Gamma_\xi$ for $\xi\leq_\tree{U}\lambda$. So we just 
handle the successor stages. So suppose  $\eta=\tree{U}\pred(\xi+1)$ 
and we have $\Gamma^\eta:\tree{X}_\eta\to \tree{W}_\eta$ and $\Gamma^\xi:\tree{X}_\xi\to \tree{W}_\xi$.

\begin{claim}
The weak tree embedding Shift Lemma applies to $(\Gamma^\eta, \Gamma^\xi, E^\tree{U}_\xi, F_\xi)$.
\end{claim}

\begin{proof}
This is a routine induction.
\end{proof}

Now let $\Delta:W(\tree{X}_\eta,\tree{X}_\xi, E^\tree{U}_\xi)\to W(\tree{W}_\eta, \tree{W}_\xi, F_\xi)$ be the  copy $W$-type weak tree embedding associated to ($\Gamma^\eta, \Gamma^\xi, E^\tree{U}_\xi, F_\xi$). Also let $\Xi=\Gamma^{X(\tree{X}_\eta,\tree{X}_\xi,E^\tree{U}_\xi)}$.

We have the following commutative diagram.

 \begin{center}
    \begin{tikzcd}
  & \tree{X}_\eta \arrow{r}{\Gamma_\eta}\arrow[swap]{dl}{E_\xi^\tree{U}}\arrow{d}{\Phi^ {W(\tree{X}_\eta, \tree{X}_\xi, E^\tree{U}_\xi)}} & \tree{W}_\eta \arrow{d}{F_\xi}\\
   \tree{X}_{\xi+1}\arrow[swap]{r}{\Xi} &
 W(\tree{X}_\eta, \tree{X}_\xi, E^\tree{U}_\xi) \arrow[swap]{r}{\Delta} &\tree{W}_{\xi+1}
    \end{tikzcd}
    \end{center}

We let $\Gamma^{\xi+1}= \Delta\circ \Xi$, 
which is as desired since the above diagram commutes. This finishes the definition of the $\Gamma_\xi$.

Note that the maps $\Gamma_\xi$ witness that $\mtree{X}(\tree{S},\tree{T})$ is 
a kind of weak meta-hull of $\mtree{W(\tree{T}, \tree{U})}$. 

This finishes our discussion of $X(\tree{T},\tree{U})$. It is a simple matter to extend the definition
so as to define $X(s)$ for finite stacks $s$.

Let us call a stack $s$ of plus trees {\em simple} iff its component
plus trees are all $\lambda$-tight and normal, and the construction
of $X(s)$ defined above never leads to type 2 premice; that is,
our simplifying assumption always applies. What we have shown is:

\begin{theorem}\label{fullnormalizationeheorem} Assume $\adp$, and let
$(M,\Sigma)$ be a strongly stable mouse pair. Let $s$ be a simple stack on $(M,\Sigma)$
with last pair $(N,\Sigma_s)$; then there is a unique normal,
$\lambda$-tight tree $\tree{X}$ on $(M,\Sigma)$ with last pair $(N,\Sigma_s)$.
\end{theorem}
\begin{proof} Let $\tree{W} = W(s)$ and $(R,\Lambda)$ be the last pair of
$\tree{W}$. Let $\sigma \colon N \to R$ be the last $\sigma$-map of the
embedding normalization. Since $(M,\Sigma)$ embedding normalizes well,
\[
\Sigma_s = \Lambda^\sigma.
\]
Let $\tree{X} = X(s)$ be the tree constructed above. $\tree{X}$ has last
model $N$. Moreover, we produced a weak tree embedding $\Gamma$ from
$\tree{X}$ into $\tree{W}$ whose
last $t$-map is
\[
t^\Gamma = \sigma.
\]
Since $\Sigma$ has very strong hull condensation (by Theorem \ref{vshctheorem}),
\[
\Lambda^{t^\Gamma} = \Sigma_{\tree{X},N}.
\]
Putting things together, we get $\Sigma_s = \Sigma_{\tree{X},N}$,
as desired.
\end{proof}

As we indicated in the course of constructing $X(s)$, the hypothesis that
$s$ is simple can be dropped from Theorem \ref{fullnormalizationeheorem},
and the details of that will appear in a future draft of this paper.
From this stronger version of the theorem, we get at once

\begin{corollary}
\label{positionalitytheorem} Assume $\adp$, and let
$(M,\Sigma)$ be a strongly stable mouse pair. Let $s$ and $t$ be finite stacks of
plus trees on
$M$ such that $s$ and $t$ are by $\Sigma$ and have a common last model
$N$; then $\Sigma_{s,N} = \Sigma_{t,N}$.
\end{corollary}

That is, the strategy component of a
 mouse pair is positional.

\end{document}